\documentclass[11pt]{article}
\usepackage{geometry} 
\usepackage{amsfonts}
\usepackage{epsfig}
\usepackage{url}
\usepackage{amssymb,latexsym,amsmath,amsthm,verbatim}
\usepackage{graphicx,epsfig,epstopdf,amssymb,color}
\usepackage{bm}
\newcommand{\N}{{\mathbb N}}

\newcommand{\veps}{\varepsilon}
\newcommand{\R}{\mathbb{R}} 
\newcommand{\C}{\mathbb{C}}
\newcommand{\K}{\mathbb{K}}
\newcommand{\OO}{\mathcal{O}}

\newcommand{\ba}{\begin{array}}
\newcommand{\ea}{\end{array}}
\newcommand{\tPsi}{\widetilde{\Psi}}
\newcommand{\tA}{\widetilde{A}}
\newcommand{\EE}{{\bf E}}
\newcommand{\eps}{\varepsilon}

\renewcommand{\epsilon}{\varepsilon}

\newtheorem{theorem}{Theorem}[section]
\newtheorem{lemma}[theorem]{Lemma}

\newtheorem{remark}[theorem]{Remark}

\DeclareMathOperator*{\supp}{supp}

\def\ualpha{{\alpha}}

\def\eps{\epsilon}
\def\ua{\ualpha}

\begin{document}
\title{Validity of the Nonlinear Schr\"odinger Approximation  
for the Two-Dimensional Water Wave Problem With and Without Surface Tension in the Arc
  Length Formulation}
\author{Wolf-Patrick D\"ull$^{1}$}
\date{\today}

\footnotetext[1]{
       Universit\"at Stuttgart, 
       Institut f\"ur Analysis, Dynamik und Modellierung, 
       Pfaffenwaldring 57, D-70569 Stuttgart, Germany; e-mail:
duell@mathematik.uni-stuttgart.de}
\maketitle

\begin{abstract}
We consider the two-dimensional water wave problem in an infinitely 
long canal of finite depth both with and without surface tension.
In order to describe the evolution of the envelopes of small oscillating wave packet-like solutions to this problem the Nonlinear Schr\"odinger equation can be derived as a formal approximation equation. In recent years, the validity of this approximation
has been proven by several authors for the case without surface tension. In this paper, we rigorously justify the Nonlinear Schr\"odinger approximation for the cases with and without surface tension by proving error estimates over a physically relevant timespan  in the arc length formulation of 
the two-dimensional water wave problem. The error estimates are uniform with respect to the strength of the surface tension, as the height of the wave packet and the surface tension go to zero.
\end{abstract}

\section{Introduction}
\label{sec1}
In this paper, we consider the two-dimensional water wave problem with finite depth of water. The two-dimensional water wave problem consists in finding the flow of an incompressible, inviscid fluid in an infinitely long canal of finite or infinite depth with a free top surface under the influence of gravity and possibly of surface tension. In Eulerian coordinates, 
the two-dimensional water wave problem with finite depth has the following form: The fluid fills a domain $\Omega(t)=\{(x,y) \in \R^2:\, x \in \R, -h < y < \eta(x,t)\}$ in between the bottom  $B = \{(x,y) \in \R^2:\, x \in \R,\, y = -h\}$ and the free top surface $\Gamma(t) = \{(x,y) \in \R^2:\, x \in \R,\, y = \eta(x,t)\}$. The velocity field $V=(v_1, v_2)$ of the fluid is governed by the incompressible Euler's equations 
\begin{align}
V_t + (V\cdot\nabla)V & =  -\nabla p + g \left(\!\begin{array}{c} 0 \\ -1 \end{array} \!\right) \qquad \text{in}\,\, \Omega(t), \label{euler} \\[2mm] 
\nabla \cdot V & =  0 \qquad \text{in}\,\, \Omega(t), \label{incompr}
\end{align}
where $p$ is the pressure and $g$ the constant of gravity. 

Assuming that fluid particles on the top surface remain on the top
surface, that the pressure at the top surface is determined by the
Laplace-Young jump condition and that the bottom is impermeable yields the boundary conditions
\begin{align}
\eta_t & =  V \cdot \left(\!\begin{array}{c} -\eta_x \\ 1 \end{array} \!\right) \qquad \text{at}\,\, \Gamma(t), \label{surf} 
\\[-2mm]  \nonumber
\end{align}
\begin{align}  
p & =  -bgh^2 \kappa \qquad \text{at}\,\, \Gamma(t), \label{press} 
\\[2mm]  
v_2 & =  0   \qquad \text{at}\,\, B \label{bot},
\end{align}
where $b \geq 0$ is the Bond number, which is proportional to the strength of the surface tension, and $\kappa$ is the curvature of $\Gamma(t)$.

If the flow is additionally assumed to be irrotational, the above
system can be reduced to a system defined on $\Gamma(t)$. Due to the
irrotationality of the motion there exists a velocity potential $\phi$
with $V=\nabla \phi$, which is harmonic in $\Omega(t)$ with vanishing
normal derivative at $B$. Moreover, the motion of the vertical component of the velocity is uniquely determined
by the horizontal one, i.e., there exists an operator 
$ \mathcal{K} = \mathcal{K}(\eta) $ such that
\begin{align} \label{Knl}
\phi_y &= \mathcal{K}(\eta) \phi_x.
\end{align}

By using the potential $\phi$, the system (\ref{euler})--(\ref{bot})
can be reduced to  
\begin{align} 
\eta_t & =  V \cdot \left(\!\begin{array}{c} -\eta_x \\ 1 \end{array} \!\right) \qquad \mathrm{at}\,\, \Gamma(t), \label{isurf2} \\[2mm]  
\phi_t &= -\frac{1}{2} ((\phi_x)^{2} + 
(\mathcal{K} \phi_x)^{2}) - g \eta + b g h^2 \left( \frac{\eta_x}{\sqrt{1+\eta_x^2}}\right)_{x} \qquad \text{at}\,\, \Gamma(t) \label{ipot}
\end{align}
or to
\begin{align} 
\eta_t & =  \mathcal{K} v_1 -v_1\eta_x 
 \qquad \text{at}\,\, \Gamma(t), \label{surf2} \\[2mm]  
(v_1)_t &= - g \eta_x -{\frac{1}{2}} (v_1^{2} + 
(\mathcal{K} v_1)^{2})_x + b g h^2 \left( \frac{\eta_x}{\sqrt{1+\eta_x^2}}\right)_{xx} \qquad \text{at}\,\, \Gamma(t). \label{pot}
\end{align}
From now on, let space and time in the above system be rescaled in such a
way that $h=1$ and $g=1$.

Choosing Eulerian coordinates to formulate the equations for the motion of water waves is natural for describing many physical experiments. But there are also alternative coordinate systems which yield appropriate frameworks for formulating the water wave problem. Each of these coordinate systems has its own advantages concerning applicability and mathematical structure of the resulting equations of the water wave problem. Hence, depending on the problem one intends to solve, one has to find out which coordinate system is the most adapted one.

The most known alternative systems are Lagrangian coordinates, see, for example, \cite{Y83},
holomorphic coordinates, see, for example, \cite{HaIT16}, the arc length formulation,
see, for example, \cite{AM05}, and abstract coordinate independent systems which base on the fact that the solutions of \eqref{euler}-\eqref{bot} can be interpreted as the geodesic flow with respect to the potential energy, the kinetic energy and in case of surface tension also the surface energy on the infinite dimensional Riemannian manifold of volume-preserving homeomorphisms of $\Omega(0)$, see, for example, \cite{SZ11}. In this differential geometric variational framework, the boundary conditions \eqref{surf}-\eqref{bot} appear as natural boundary conditions.

For an overview on the local and global well-posedness results for the water wave problem in the various formulations we refer to \cite{D18} and the references therein.

Concerning the qualitative behavior of the solutions, the full
water wave problem is extremely complicated to analyze. A qualitative understanding of the solutions to the full water wave problem being usable for practical applications does not seem within reach for the near future, neither analytically nor numerically. Therefore, it is important to approximate the full model in different parameter regimes by suitable reduced model equations whose solutions have similar but more easily accessible qualitative properties.

The simplest reduced model equation is the linear wave equation. The most famous nonlinear approximation equations are the Korteweg-de Vries (KdV) equation and the Nonlinear Schr\"odinger (NLS) equation. By inserting the ansatz
\[
\begin{pmatrix} \eta \\ v_1 \end{pmatrix} (x,t)=
\varepsilon^2 A \left(\varepsilon (x \pm t), \varepsilon^3t\right) 
\begin{pmatrix} 1 \\ \mp 1  \end{pmatrix} + \mathcal{O} (\varepsilon^3)
\]
with $0 < \varepsilon \ll1$ and $ A: \R^2 \to \R$ into
\eqref{surf2}--\eqref{pot}, expanding the operator $\mathcal{K}$ with respect to $\varepsilon$ and equating the terms with the lowest power
of $\eps$ one obtains that $A$ has to satisfy in lowest order with respect to $\eps$ the KdV equation
\begin{equation} \label{kdv}
A_{\tau} = \pm \Big(\frac16-\frac{b}{2}\Big) A_{\xi\xi\xi} \pm \frac32 A A_{\xi}\,,
\end{equation}
where $\tau = \varepsilon^3 t$ and $\xi = \varepsilon (x\pm t)$, if $b \neq 1/3$.
For further information about the KdV approximation we refer to \cite{D18} and the references therein. 
 
The ansatz for the NLS approximation is 
\[
\begin{pmatrix} \eta \\ v_1 \end{pmatrix} (x,t)= 
\varepsilon A \left(\varepsilon (x -c_gt), \varepsilon^2t\right) e^{i(k_0x-\omega_0 t)}
\varphi(k_0,b) + \mathcal{O} (\varepsilon^2) + \mathrm{c.c.}\,,
\]
where $0 < \varepsilon \ll 1$.
Here $ \omega_0 >0$ is the basic temporal 
wave number associated via the linear dispersion relation of the two-dimensional water wave problem with finite depth, namely
$$\omega(k) = \omega(k,b) = \pm \sqrt{(k+bk^3)\tanh(k)}\,,$$ 
to the basic spatial wave number $ k_0 > 0$ of the underlying carrier wave $ e^{i(k_0 x - \omega_0 t)}$, that means that $\omega_0= \omega(k_0)$, where the branch of solutions
\begin{equation} \label{wwdisp}
\omega(k) = \omega(k,b) = {\rm sgn}(k) \sqrt{(k+bk^3)\tanh(k)}
\end{equation} 
is chosen. Moreover,
$c_g$ is the group velocity, i.e., $c_g= \omega'(k_0)= \partial_k \omega(k_0,b)$, $A$ the complex-valued amplitude, $\varphi(k_0,b) \in \C^2$ and c.c. the complex conjugate.  This ansatz leads to waves moving to the right; to obtain waves moving
to the left, $\omega_0$ and $c_g$ have to be replaced by  $-\omega_0$ and $-c_g$.

\begin{figure}[htbp]
\epsfig{file=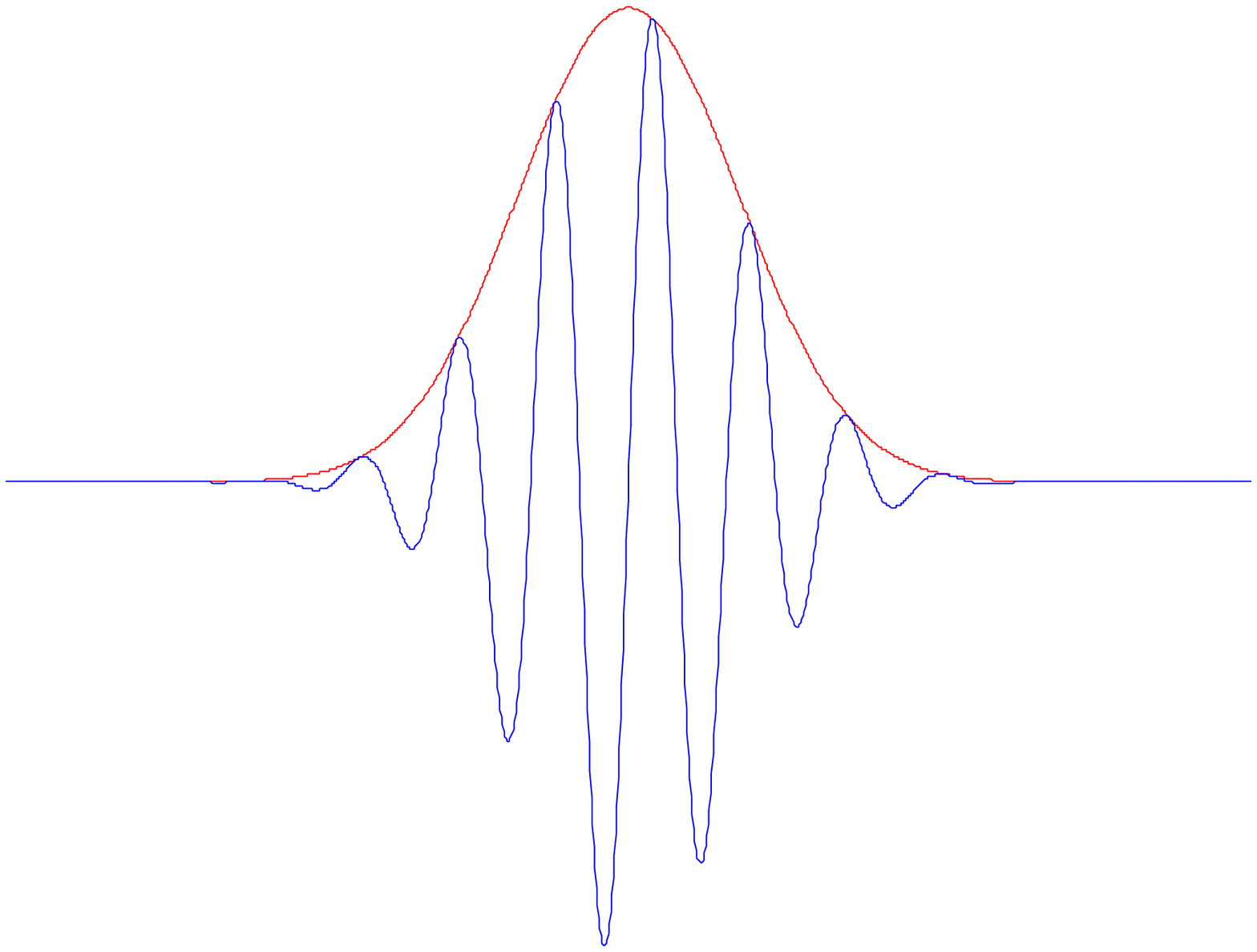,width=15cm,height=6.2cm,angle=0}
\vspace*{-6.3cm}

\hspace{1.9cm}
\epsfig{file=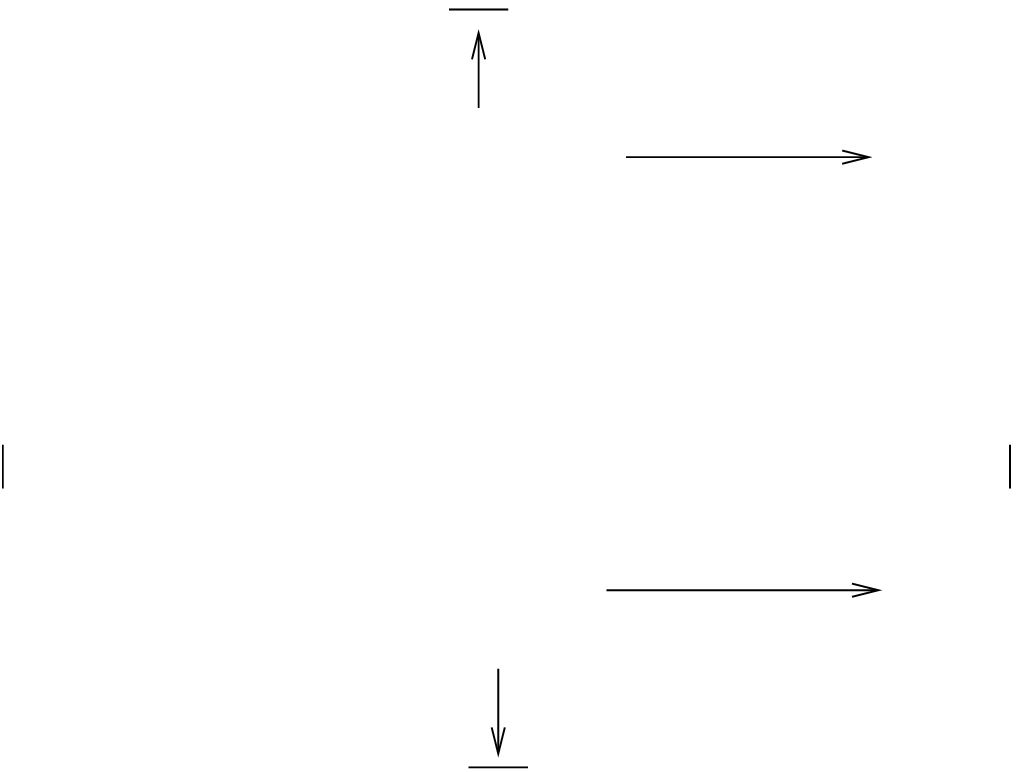,width=10.7cm,height=6.5cm,angle=0}

\vspace*{-5.34cm}
\hspace*{9.7cm}$c_\mathrm{g}$  \vspace*{3.2cm}

\hspace*{9.7cm} $c_\mathrm{p}$

\vspace*{-2.5cm}
\hspace{2.2cm}\hspace{5.0cm}$\varepsilon$
\vspace*{0.3cm}

\hspace{2.65cm}$1/\varepsilon$
\hspace{2.7cm}
\vspace*{2.6cm}

\caption{{\small An envelope (advancing with the group velocity 
$ c_g $) with characteristic length scale of order $\mathcal{O}(1/\varepsilon)$ of an oscillating wave packet $\eta$ of order $\mathcal{O}(\varepsilon)$ 
(advancing with the phase velocity 
$ c_p= \omega_0/k_0  $) is approximately described by the 
amplitude $ A $ which solves  the NLS equation \eqref{NLS}.} \label{fig3} }
\end{figure}

By inserting this ansatz into (\ref{surf2})--(\ref{pot}), one obtains that for an explicitly computable vector $\varphi(k_0,b)$ the amplitude $A$ has to satisfy at leading order in $\eps$ the NLS equation
\begin{equation} \label{NLS}
A_{\tau} = i \frac{\partial_k^{2} \omega(k_0,b)}{2} A_{\xi\xi} +i \nu(k_0,b) A |A|^2 \,,
\end{equation}
where $\tau = \varepsilon^2 t, \xi = \varepsilon (x -c_g t)$ and $\nu(k_0,b) \in \R$.
Hence, the NLS equation \eqref{NLS} approximately describes the dynamics of spatially and temporarily oscillating wave packet-like solutions to the two-dimensional water wave problem; see Figure \ref{fig3}.

In one space dimension, both the KdV equation and the NLS equation are completely integrable Hamiltonian systems which 
can be explicitly solved with the help of inverse scattering schemes; see, for example, \cite{AS81}.  

The first formal derivation of the NLS approximation for the two-dimensional water wave problem was given by Zakharov \cite{Za68} in 1968. 
The NLS approximation is used, for example, in the context of modeling monster waves; see \cite{JO07}.  However, the NLS approximation plays not only an important role for the mathematical description of surface water waves but also in other areas of science and technology, for example, in nonlinear optics to model data transmission via fiber optic cables with the help of light pulses \cite{AS81, Schn11OWbuch}, in biology to model waves in DNA \cite{H97}, in plasma physics \cite{SS99} or in quantum theory \cite{L01}. In numerical simulations, the simulation of the evolution of the envelope with the help of the NLS approximation yields a significant reduction of complexity and consequently an increase of efficiency compared to the simulation of the whole wave packet. 

Although the NLS approximation is very successful in many applications, it should not be taken for granted that the NLS approximation always yields correct predictions of the behavior of the original system. Indeed, there are counterexamples where the NLS approximation fails \cite{Schn05,SSZ15}. Hence, it is important to answer the question of the validity of the NLS approximation for a given system by proving error estimates over a physically relevant timespan. In general, this is a highly nontrivial mathematical problem for the following reasons.

Given the general abstract evolutionary 
problem
\begin{equation} \label{abstrsyst}
  \partial_t W = \mathfrak{L} W + \mathfrak{B}(W,W) + \mathfrak{H}(W)
\end{equation}
with $x,t \in \R$ and $W=(\,W_1(x,t) \;\;\, W_2(x,t))^{T} \in \R^2$.
Here $\mathfrak{L}$ is a linear operator whose symbol is a diagonal matrix of the form 
\begin{equation} \label{linsymb}
\widehat{\mathfrak{L}}(k) = \mathrm{diag}\,(-i\omega(k),\, i\omega(k))\,,
\end{equation}
where $k\in \R$ and $\omega$ is a piecewise smooth real-valued odd function.
Furthermore, $\mathfrak{B}$ is a bilinear operator and $\mathfrak{H}(W)$ consists of terms being at least cubic in $W$ or is equal to $0$. 

The NLS equation \eqref{NLS} can be derived as a formal approximation equation with the help of the ansatz $W=\eps \tilde{\Psi}$, where
\begin{align}\label{approx1}
\eps \tilde{\Psi}(x,t) \,&=\, \eps \Psi_{NLS}(x,t) \left(\begin{array}{c}
1\\
0
\\ \end{array}\right) + \eps^2 \Psi_{h}(x,t)\,,
\end{align}
\begin{align}
\label{ansatz}
\Psi_{NLS}(x,t) 
\,&=\, A(\eps
(x-c_g t),\eps^2t) e^{i( k_0 x - \omega_0 t)} + \mathrm{c.c.} \,,
\\[2mm]
\label{approx1b}
{\Psi}_h(x,t) \,&=\,  \left(\begin{array}{c}
 {\tA}_{01} (\eps (x-c_gt),\eps^2t)\\
{\tA}_{02} (\eps (x-c_gt),\eps^2t)
\\ \end{array}\right) \nonumber \\[2mm] & 
\,\quad \,\, + \left( \left(\begin{array}{c}
 {\tA}_{21} (\eps (x-c_gt),\eps^2t)\\
{\tA}_{22} (\eps (x-c_gt),\eps^2t)
\\ \end{array}\right) e^{2i( k_0 x - \omega_0 t)} + \mathrm{c.c.} \right) , 
\end{align}
$k_0 >0$, $\omega_0 = \omega(k_0)$, $c_g=\omega'(k_0)$, ${\tA}_{01}, {\tA}_{02}$ are real-valued functions and ${\tA}_{21}, {\tA}_{22}$ complex-valued functions. 

Inserting this ansatz into \eqref{abstrsyst} and equating the coefficients in front of the $\veps^m  e^{ji( k_0 x - \omega_0 t)}$ for $m \in \{1,2,3\}$ and  $j \in \{0,1,2\}$ to $0$ yields the NLS equation 
\begin{equation} \label{NLS2}
A_{\tau} = i \frac{\omega''(k_0)}{2} A_{\xi\xi} +i \nu(k_0) A |A|^2 \,,
\end{equation}
where $\tau = \varepsilon^2 t, \xi = \varepsilon (x -c_g t)$ and $\nu(k_0) \in \R$,
if $\omega$ satisfies 
\begin{equation} 
\lim\limits_{k \to 0^{\pm}} \omega(k) \neq 0  
\end{equation}
or
\begin{equation} \label{nrb1}
 \lim\limits_{k \to 0^{\pm}} \omega'(k) \neq \omega'(k_0) 
\end{equation}
as well as
\begin{equation} 
\pm \omega(2k_0)\neq 2 \omega(k_0)
\end{equation}
and
\begin{equation} \label{nrb3}
\omega''(k_0) \neq 0 \,. 
\end{equation}

The above ansatz leads to wave packets moving to the right; to obtain wave packets moving
to the left, $\omega_0$, $c_g$ have to be replaced by  $-\omega_0$, $-c_g$ and $(  \Psi_{NLS}(x,t) \quad 0 )^{T}$ by $( 0 \quad \Psi_{NLS}(x,t))^{T}$.

It is possible to modify  $\eps \tilde{\Psi}$ to make it an even more accurate approximation. Indeed, if there exists an integer $M > 2$ such that
\begin{equation} \label{nrb4}
\pm \omega(mk_0)\neq m \omega(k_0)
\end{equation}
for all integers $m \in [2,M)$, then there 
exists a function $\Psi$ dependent on $\varepsilon$ such that 
\begin{equation}
\lim_{\varepsilon \to 0} \|\Psi(\cdot,t)-\tilde{\Psi}(\cdot,t)\|_{C^{0}} = 0
\end{equation}
and
\begin{equation} 
\mbox{Res}(\veps \Psi) := - \partial_t (\veps \Psi) + \veps \mathfrak{L}\Psi 
+ \mathfrak{B}(\veps \Psi, \veps \Psi) + \mathfrak{H}(\veps \Psi) = \OO(\eps^M)\, . 
\end{equation}

The two-dimensional water wave problem with finite depth can be transformed to an evolutionary system of the form \eqref{abstrsyst} by diagonalizing the linear part of the system\eqref{surf2}--\eqref{pot}.
More precisely, if one makes the linear coordinate transform
\begin{equation}  \label{Www}
W =  \left(
\begin{array}{c} {W}_{1} \\ {W}_{2}
\end{array}
\right) := 
\frac{1}{2}\left(
\begin{array}{cc}
  \sigma &  1   \\
 -\sigma & 1    
\end{array}
\right)
\left(
\begin{array}{c} {y} \\ {v_1}
\end{array}
\right),
 \end{equation}
where $\sigma$ is a linear operator defined by its symbol 
\begin{equation} 
\sigma(k)= \sigma(k,b) = \sqrt{\frac{k+bk^3}{\mathrm{tanh}(k)}} 
\end{equation}
for $k \in \R$, then $W$ satisfies a system of the form \eqref{abstrsyst} with $\omega$ defined by \eqref{wwdisp}. 

For $b=0$, the dispersion relation $\omega$ satisfies \eqref{nrb1}--\eqref{nrb4} for all $k_0 > 0$. For $b >0$, one has to choose a basic wave number $k_0 > 0$ for which the conditions \eqref{nrb1}--\eqref{nrb4} are valid in order to be able to derive a sufficiently accurate NLS approximation. Or, if $k_0 > 0$ is given, then one has to choose those values of $b \geq 0$ for which \eqref{nrb1}--\eqref{nrb4} are valid in order to be able to derive a sufficiently accurate NLS approximation.

To guarantee that qualitative properties of solutions to the NLS equation \eqref{NLS2} like the way pulses interact with each other are also true for solutions to system \eqref{abstrsyst}, it has to be proven that the error
 \begin{equation}
\label{errorequ}
\eps^\beta R := W - \eps {\Psi}
\end{equation}
is of order $\OO(\eps^\beta)$ with $ \beta > 1 $ 
on the characteristic time scale of the NLS equation \eqref{NLS2}, 
this means that there exists a $\tau_0 >0$ such that
$R$ is of order $\OO(1)$ for all $t \in [0,\tau_0/\eps^2]$.
The rescaled error $R$ satisfies for appropriately chosen $\beta$ and $\Psi$ an evolution equation of the form
\begin{equation} \label{errorequation1}
  \partial_tR = \mathfrak{L}R + \eps B({\Psi},R) + \mathcal{O}(\eps^2)\,,
\end{equation}
where 
\begin{equation}  \label{Bsym}
B({\Psi},R) = \mathfrak{B}({\Psi},R) + \mathfrak{B}(R,{\Psi})\,.
\end{equation}
Since the Fourier transform $\veps \widehat{\Psi}_{NLS}$ of $\veps {\Psi}_{NLS}$ is strongly concentrated around the wave numbers $\pm k_0$, the approximation $ \eps \Psi$ can be split into 
\begin{equation}
\eps \Psi = \eps \Psi_c + \eps^{2} \Psi_s = \eps \left(\begin{array}{c}
\psi_{1}\\
0
\\ \end{array}\right)  +  \eps \left(\begin{array}{c}
\psi_{-1}\\
0
\\ \end{array}\right) + \eps^{2} \Psi_s \,,
\end{equation}
where the supports of $\widehat{\psi}_{\pm 1}$ satisfy 
\begin{equation}
\label{supppsic}
\supp \widehat{\psi}_{\pm 1} = \{k \in \R:\, |k \mp k_0| \leq \delta\}
\end{equation}
for a $\delta \in\, (0, k_0)$ sufficiently small, but independent of $\veps$, and $\Psi_s$ is of order $\OO(1)$.
Consequently, we have
\begin{equation}
\label{errorequation}
\partial_tR = \mathfrak{L}R + \eps B({\Psi}_c,R) + \mathcal{O}(\eps^2)\,.
\end{equation}
Hence, the main difficulty is to control the quadratic term $\eps B({\Psi}_c,R)$ over a timespan of order $\mathcal{O}(\varepsilon^{-2})$. 

A classical strategy is to eliminate the quadratic term with the help of a so-called normal-form transform 
\begin{equation} \label{inft} 
\widetilde{R} := R + \eps
N({\Psi}_c,R)\,,
\end{equation}
 where $N$ is an appropriate bilinear mapping, which can be constructed with the help of the Fourier transform. More precisely, let 
\begin{equation}
\label{B}
\widehat{B}_{j_1}(\Psi_c,R)(k) =
\int_{\R} \sum_{\genfrac{}{}{0pt}{}{\ell \in \{\pm 1\},} {j_2 \in \{1,2\}}} \widehat{b}_{j_1 j_2}(k,k-m,m) \widehat{\psi}_{\ell} (k-m) \widehat{R}_{j_2}(m)\, dm 
\end{equation}
and
\begin{equation}
\label{N}
\widehat{N}_{j_1}(\Psi_c,R)(k) =
\int_{\R} \sum_{\genfrac{}{}{0pt}{}{\ell \in \{\pm 1\},} {j_2 \in \{1,2\}}} \widehat{n}_{j_1 j_2 \ell}(k,k-m,m) \widehat{\psi}_{\ell} (k-m) \widehat{R}_{j_2}(m)\, dm \,,
\end{equation}
where $j_1, j_2 \in \{1,2\}$ denote the components of the vectors $\widehat{B}, \widehat{N}$ and $\widehat{R}$. Then, by inserting \eqref{inft} into \eqref{errorequation}, one obtains that 
$\widetilde{R}$ solves an evolution equation of the form
\begin{equation} \label{newerr}
\partial_t \widetilde{R} = \mathfrak{L} \widetilde{R} + \eps^2 g({\Psi}_c,\widetilde{R}) + \mathcal{O}(\eps^2)\,,
\end{equation}
where $g$ is of order $\OO(1)$, if  
\begin{equation} \label{nkernel}
\widehat{n}_{j_1 j_2 \ell}(k,k-m ,m ) = \frac{ \widehat{b}_{j_1 j_2}(k,k-m ,m )}{ 
i(j_1\omega(k) +  \omega(k-m) - j_2\omega(m))}
\end{equation}
and if the normal-form transform $R \mapsto  \widetilde{R}$ is invertible. Furthermore, due to \eqref{supppsic}, it turns out that it is even possible to simplify the kernels $\widehat{n}_{j_1 j_2 \ell}$ to
\begin{equation} \label{nkernel2}
\widehat{{n}}_{j_1 j_2 \ell}(k) = \frac{ \widehat{b}_{j_1 j_2}(k,\ell k_0,k -\ell k_0 )}{ 
i(j_1\omega(k) + \omega(\ell k_0) - j_2\omega(k -\ell k_0))}\,.
\end{equation}
The strategy of using normal-form transforms to eliminate semilinear quadratic terms in hyperbolic systems was introduced in \cite{Sh85}. In the context
of justifying NLS approximations, it was first applied in \cite{Kal88}.

However, there are serious difficulties.
The first one is 
the possible occurrence of resonances.
This means that the denominator of the fraction in \eqref{nkernel2} may have zeros, the so-called resonances or resonant wave numbers (to the wave number $\ell k_0$). Since $\omega$ is odd, any resonance implies further resonances. Namely, if $k$ is resonant to $\ell k_0$, then $-k$ is resonant to $-\ell k_0$. Moreover, if $k$ is resonant to $\ell k_0$ and $j_1=j_2$, then
$\pm (k-\ell k_0)$ is resonant to $ \mp\ell k_0$.

In the case of the two-dimensional water 
wave problem with finite depth, there is a resonance at $k=0$, but the numerator of the fraction in \eqref{nkernel2} also vanishes at $k=0$ and the singularity is removable. Such a resonance is called trivial. Otherwise it is called non-trivial. The resonance at $k=0$
implies resonances at $k =\pm k_0$, which are non-trivial.  
Moreover, for all basic wave numbers $k_0 >0$ there exist some $b \in (0, 1/3)$ such that there are additional non-trivial resonances for $j_1=j_2=-1$.

In the context of the justification of the NLS approximation for an evolutionary system of the form \eqref{abstrsyst} with resonances at $0$ and $\ell k_0$, it is relevant 
if the wave numbers $\ell k_0$ are stable, this means  
that for any wave number  $\ell k_1 \in \R \setminus \{0,\ell k_0\}$ being a non-trivial resonance with respect to $\ell k_0$  for $j_1=j_2=-1$ the NLS subspace in the Three Wave Interaction (TWI) system associated to the wave numbers $k_{0\ell}:=-\ell k_0$, $k_{1\ell}:=\ell k_1$ and $k_{2\ell}:=- \ell (k_1-k_0)$, which then satisfy  
\begin{align}
k_{0\ell} + k_{1\ell} + k_{2\ell} = 0\,, \qquad \omega(k_{0\ell})+ \omega(k_{1\ell}) + \omega(k_{2\ell}) = 0\,,
\end{align}
is stable. More precisely, inserting the  ansatz
\begin{equation}
W_{\ell}(x,t)= 
\sum_{j=0}^2 A_{j\ell}(\veps t) e^{i(k_{j\ell}x-\omega(k_{j\ell})t)} \left(\begin{array}{c}
1\\
0
\\ \end{array}\right) + \mathrm{c.c.} \,,
\end{equation}
where $0 < \veps \ll 1$, in \eqref{abstrsyst} and equating the coefficients of $ \varepsilon^2 
e^{i(k_{j\ell} x - \omega(k_{j\ell}) t)}  $ for $j \in \{0,1,2\}$
to zero yields 
the so-called TWI system 
\begin{align} \label{twi}
\nonumber
\partial_\tau A_{0\ell} &\,=\;\widehat{b}_{1 1}(-\ell k_0,-\ell k_1, \ell k_1 - \ell k_0 )\, \overline{A_{1\ell} A_{2\ell}}\,, \\[1mm]
\partial_\tau A_{1\ell} &\,=\; \widehat{b}_{1 1}(\ell k_1, \ell k_0, \ell k_1 - \ell k_0 )\, \overline{A_{0\ell} A_{2\ell}}\,, \\[1mm] \nonumber
\partial_\tau A_{2\ell}  &\,=\; \widehat{b}_{1 1}(-\ell k_1+\ell k_0,\ell k_0,-\ell k_1)\, \overline{A_{0\ell} A_{1\ell}}\,\,
\end{align}
with $\tau = \veps t$ and $\widehat{b}_{1 1}$ as in \eqref{B}.
This system has three invariant subspaces
consisting of fixed points, namely
$
M_{0\ell} = \{ A_{1\ell} = A_{2\ell} = 0 \} , \, 
M_{1\ell} = \{ A_{0\ell} = A_{2\ell} = 0 \}$  and $M_{2\ell} = \{ A_{0\ell} = A_{1\ell} = 0 \}$.
The so-called NLS subspace $M_{0\ell}$ is stable if and only if 
\begin{equation} \label{stab}
\frac{\widehat{b}_{1 1}(\ell k_1, \ell k_0, \ell k_1 - \ell k_0 )}{\widehat{b}_{1 1}(-\ell k_1+ \ell k_0, \ell k_0,-\ell k_1)} <0 
\end{equation}
and then
\begin{equation} \label{staben}
E:= |A_{1\ell}|^2 - \frac{\widehat{b}_{1 1}(\ell k_1, \ell k_0, \ell k_1 - \ell k_0 ) }{\widehat{b}_{1 1}(-\ell k_1+ \ell k_0, \ell k_0,-\ell k_1)}\, |A_{2\ell}|^2
\end{equation}
is a non negative conserved quantity of the system \eqref{twi}; see \cite{Schn05}. Since \eqref{abstrsyst} is a real-valued system, the wave number $k_0$ is stable if and only if $-k_0$ is stable.
If $k_0$ is unstable, then the corresponding NLS approximation can fail under certain conditions; see \cite{Schn05, SSZ15}.

For the two-dimensional water wave problem with finite depth, the values of the coefficients in the corresponding TWI system \eqref{twi} can be computed explicitly with the help of the method from \cite{SSZ15}. It turns out that $k_0$ is stable if and only if $k_0 < \max\, \{k_{1}, k_0 - k_{1}\}$ for all $k_1 \in \R \setminus \{0,k_0\}$ being a non-trivial resonance with respect to $k_0$ for $j_1=j_2=-1$.

\begin{figure}[tttt]
\epsfig{file=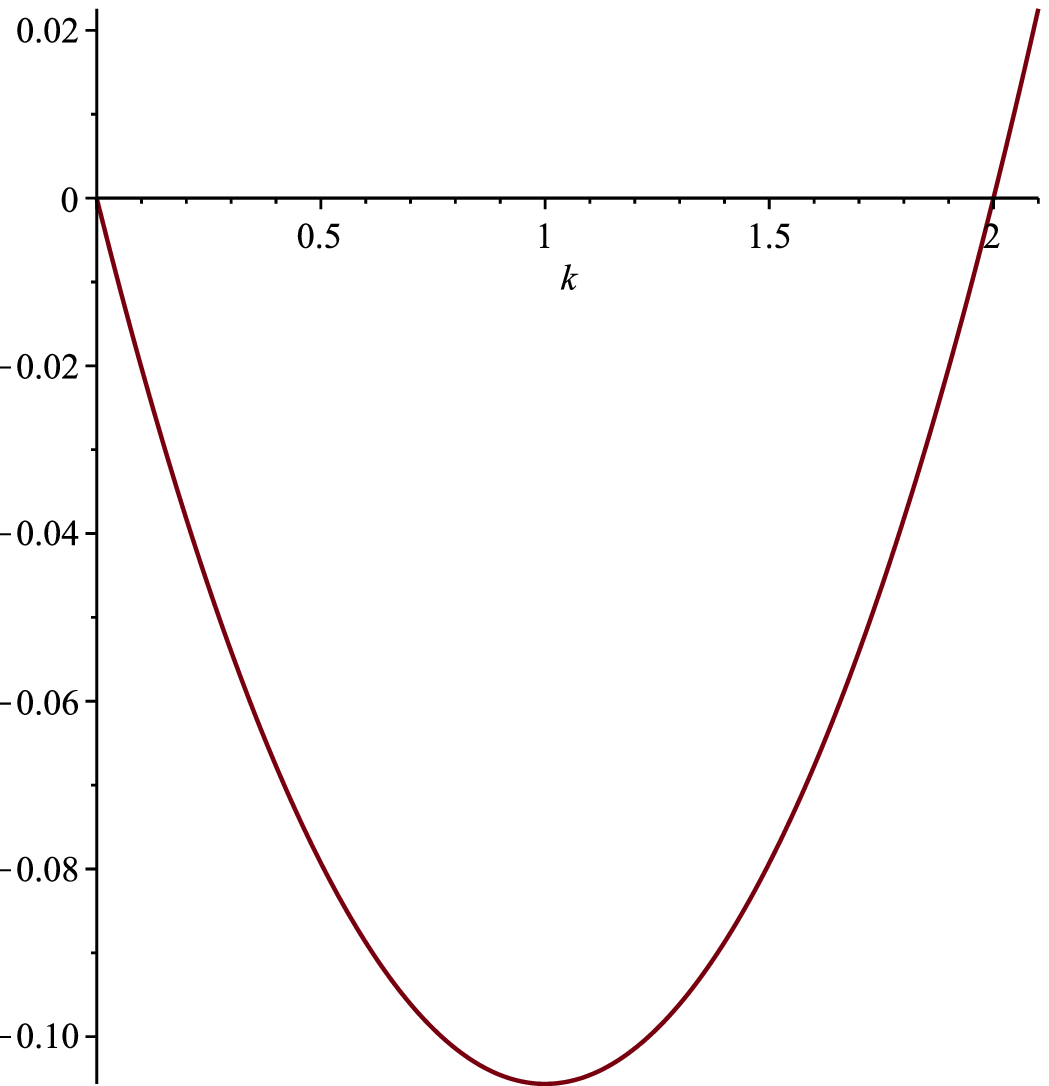,width=4.2cm,height=4.2cm,angle=0}
\qquad
\epsfig{file=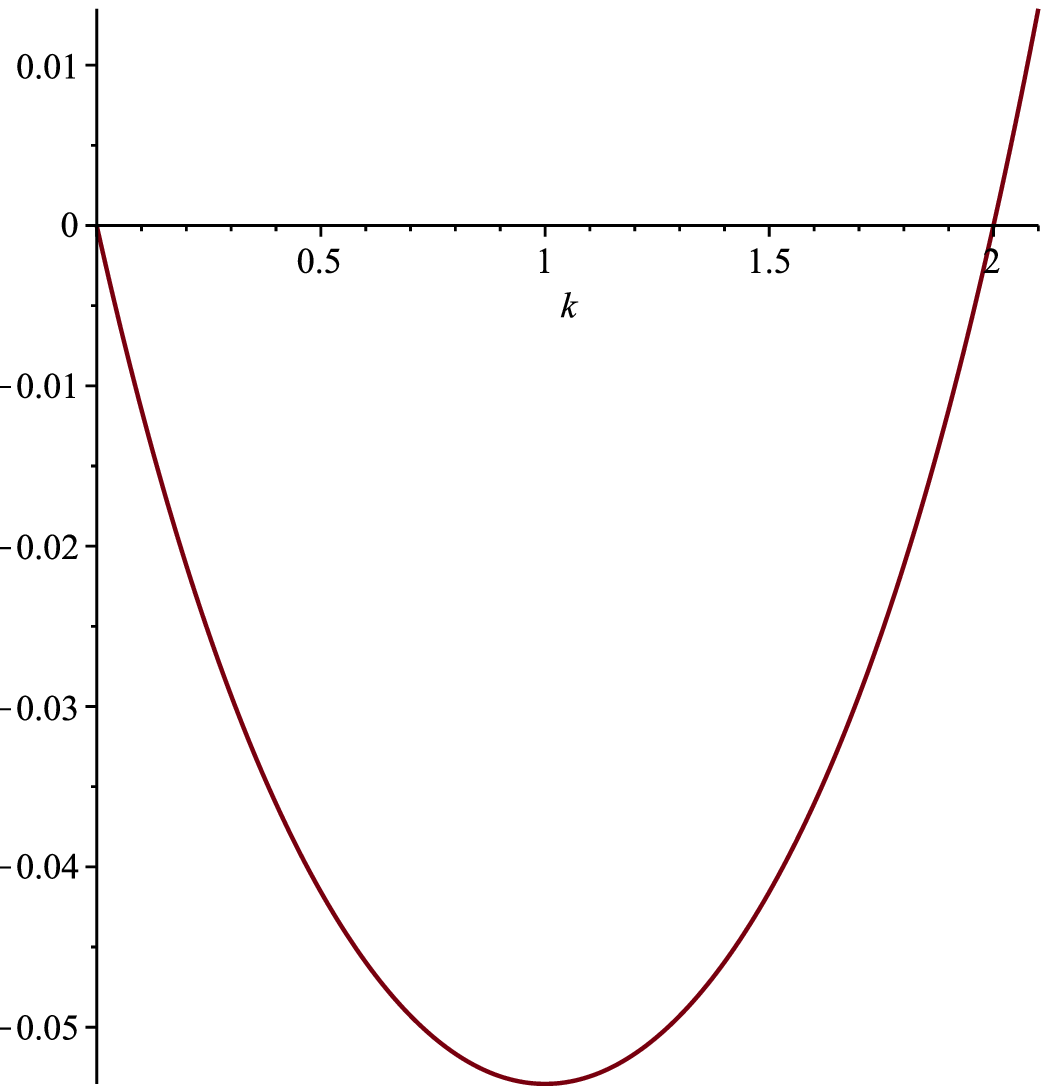,width=4.2cm,height=4.2cm,angle=0}
\qquad
\epsfig{file=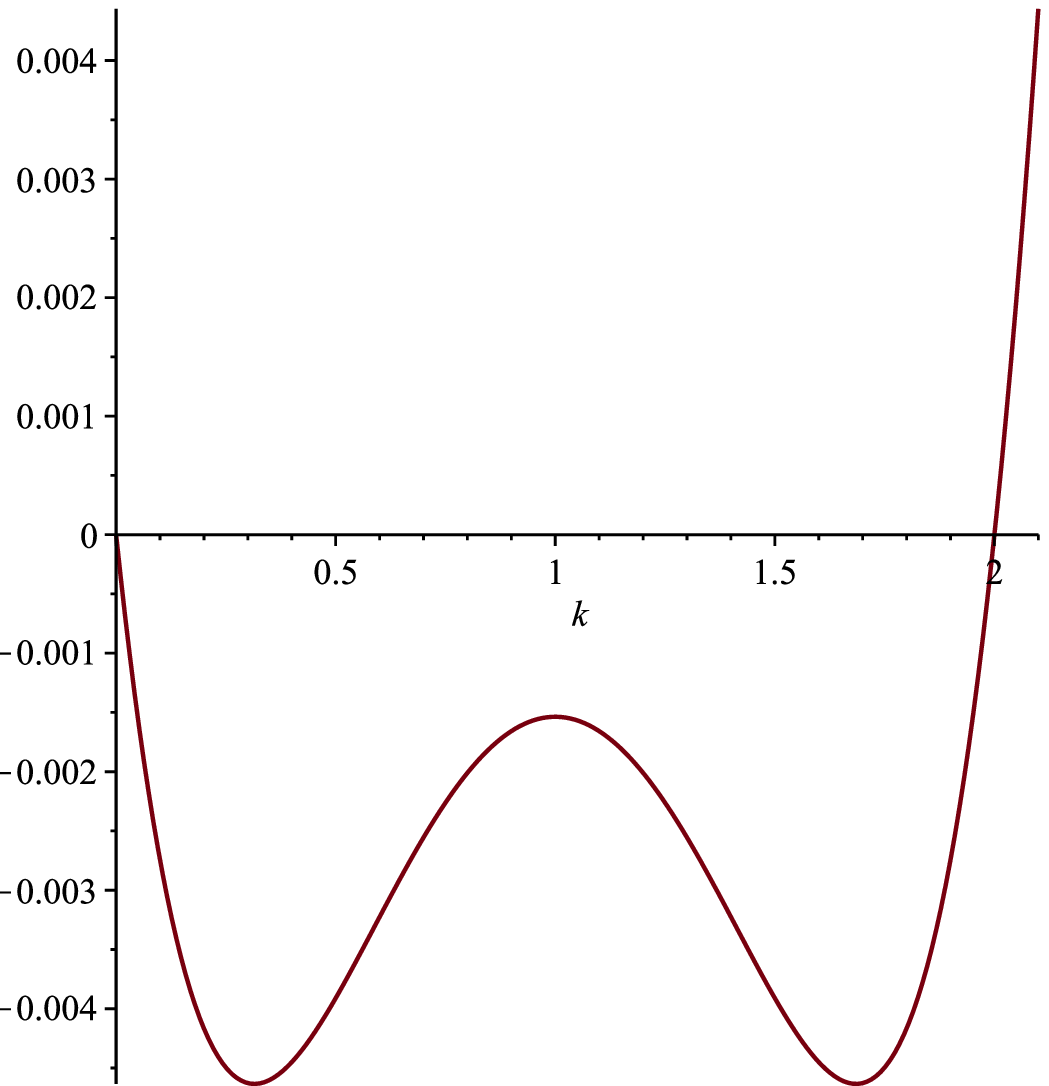,width=4.2cm,height=4.2cm,angle=0}
\\[5mm]
\epsfig{file=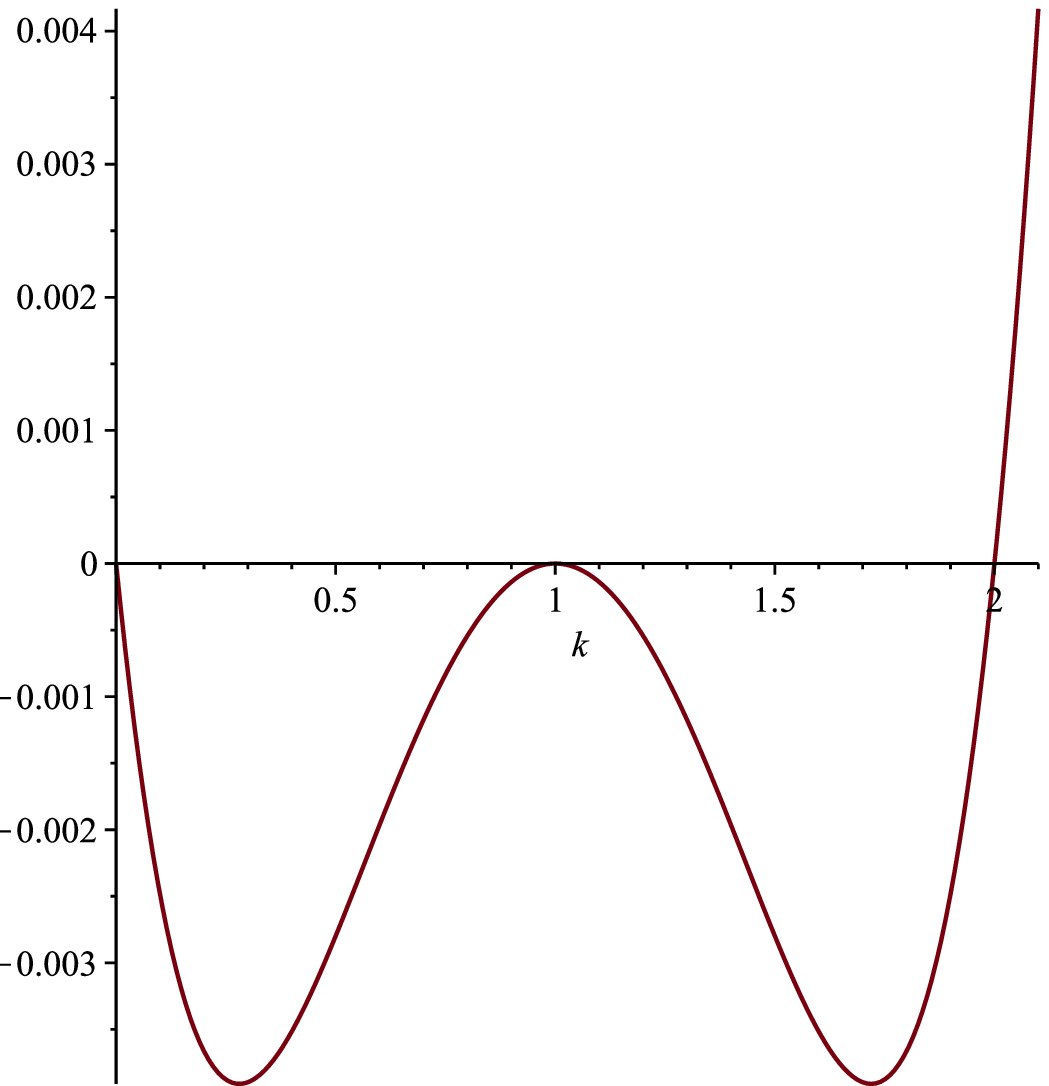,width=4.2cm,height=4.2cm,angle=0}
\qquad
\epsfig{file=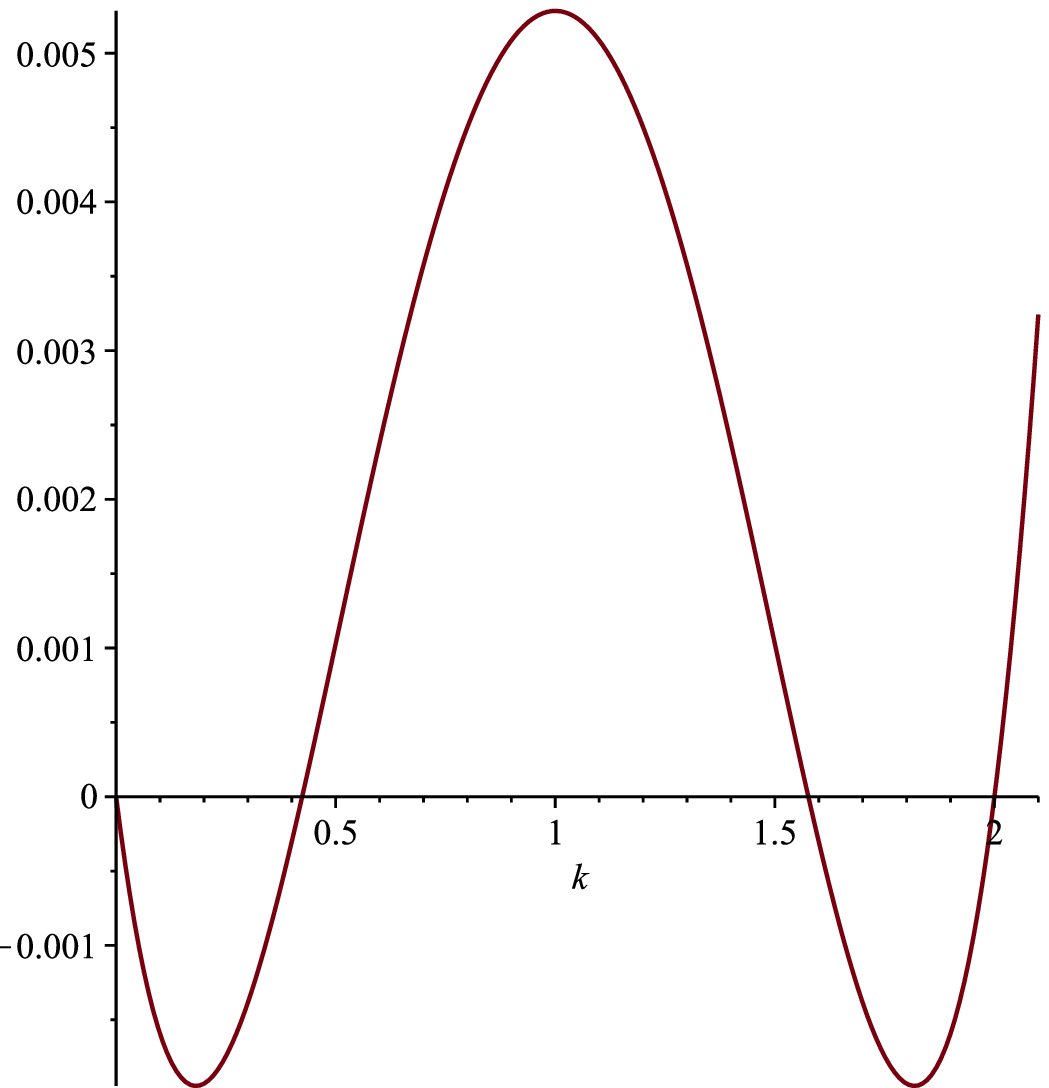,width=4.2cm,height=4.2cm,angle=0}
\qquad
\epsfig{file=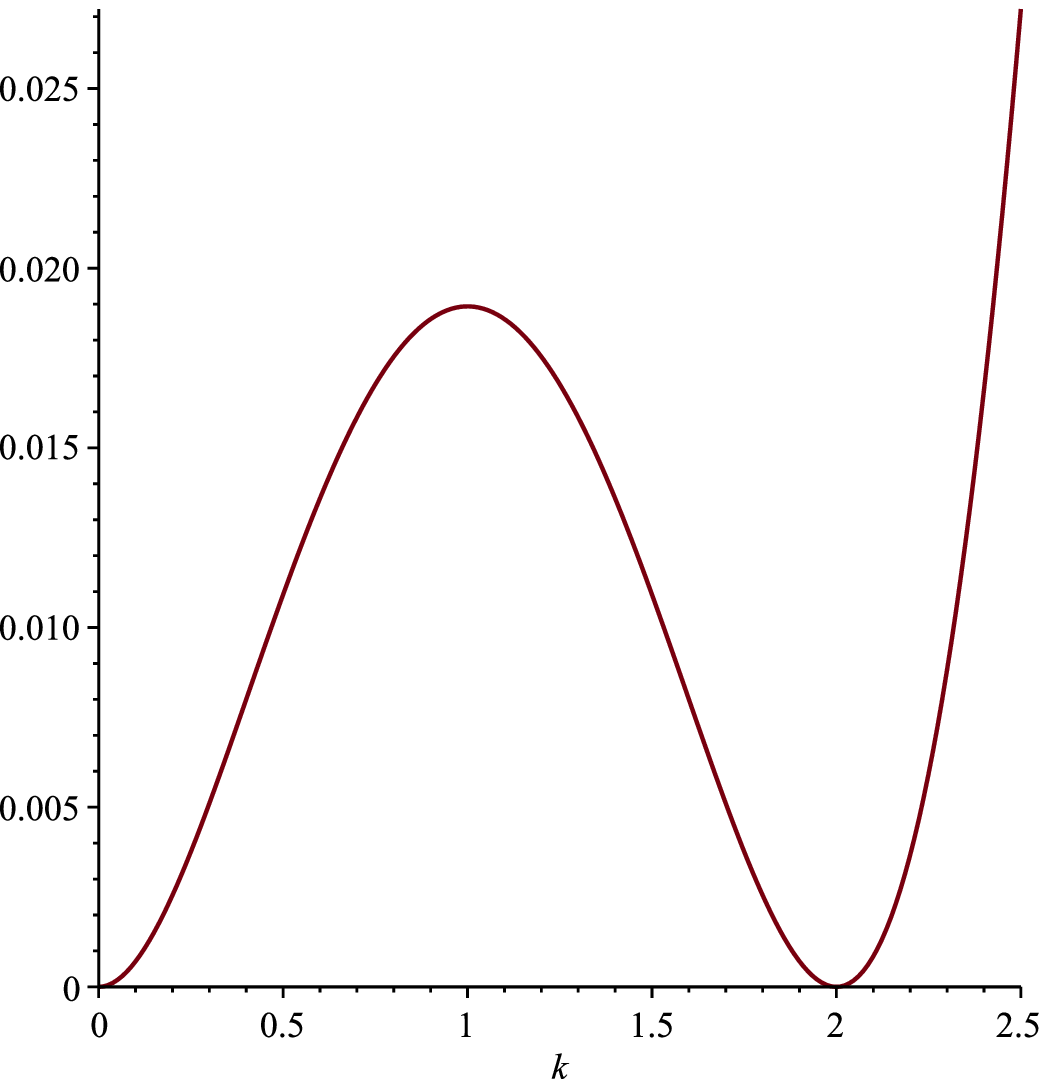,width=4.2cm,height=4.2cm,angle=0}
\\[5mm]
\epsfig{file=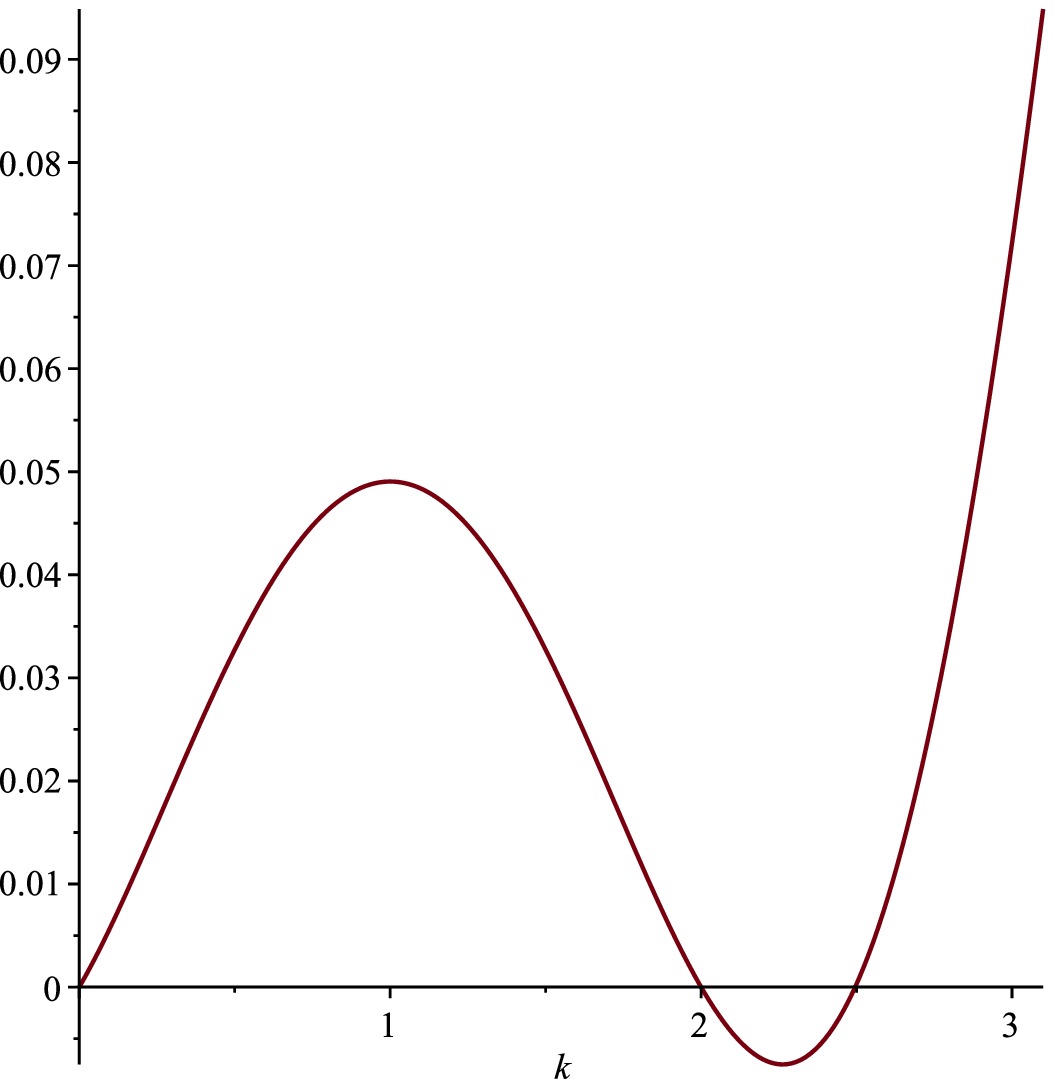,width=4.2cm,height=4.2cm,angle=0}
\qquad
\epsfig{file=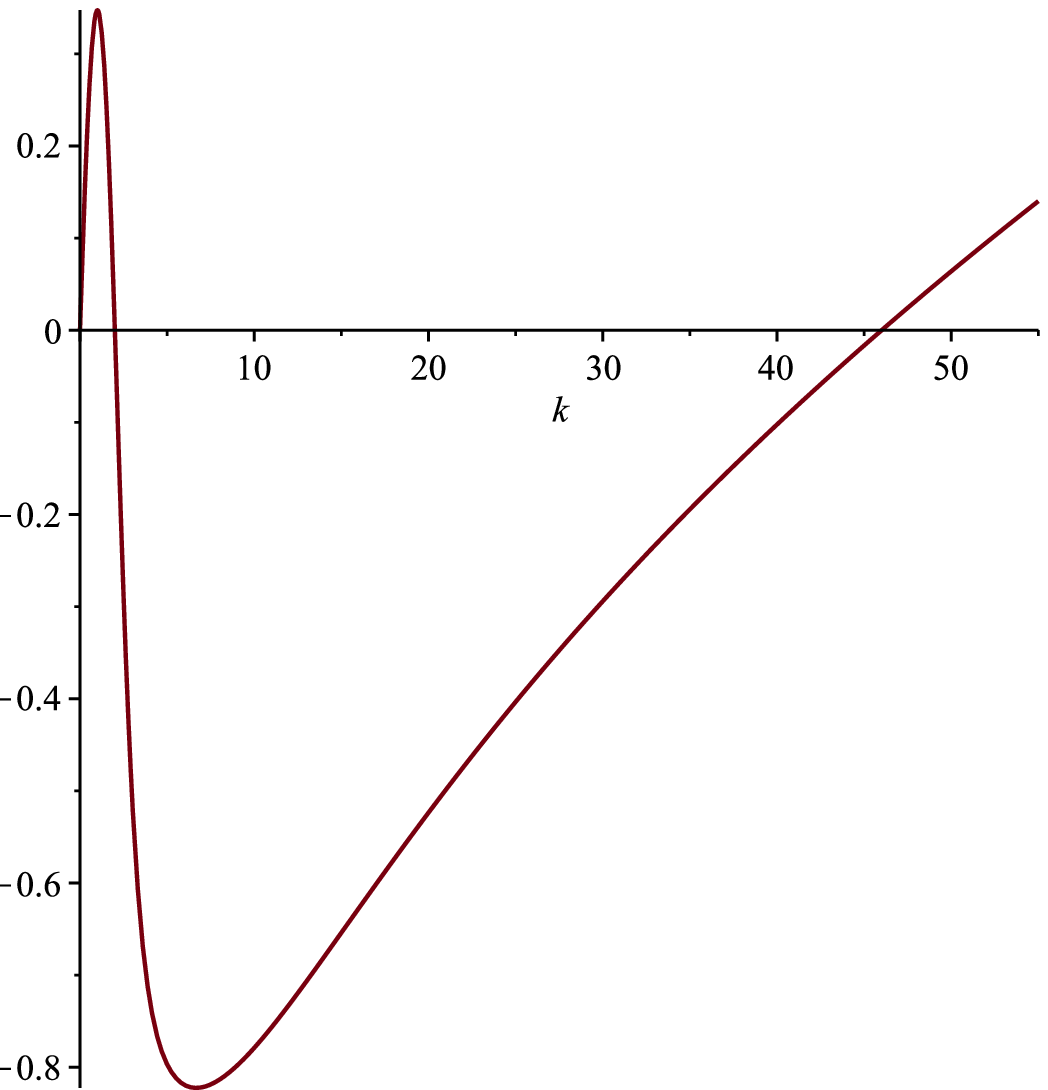,width=4.2cm,height=4.2cm,angle=0}
\qquad
\epsfig{file=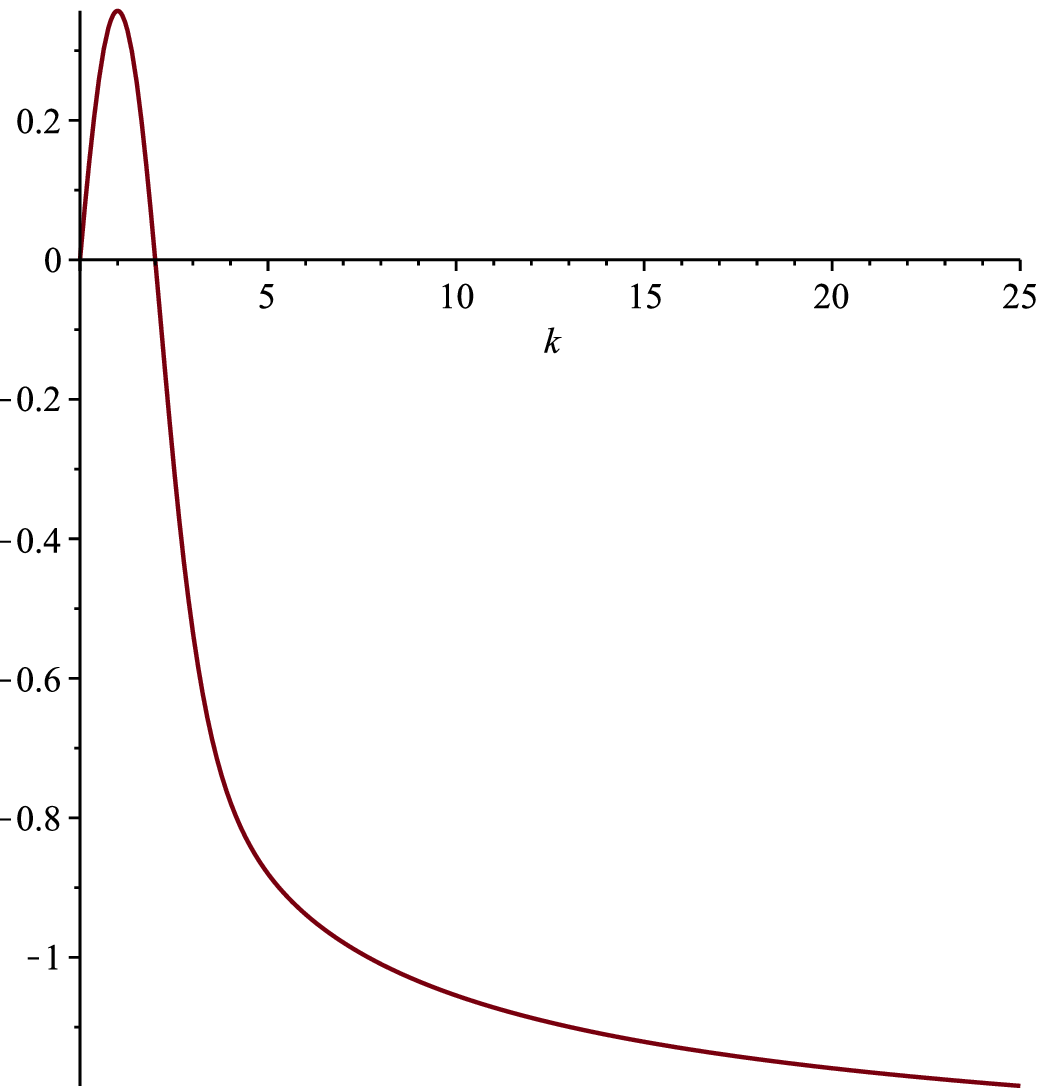,width=4.2cm,height=4.2cm,angle=0}
\caption{ {\small The graph of $k \mapsto  \omega(k,b) - \omega(k-k_0,b) - \omega(k_0,b)$ for $k \geq 0$, $k_0=2$ and different values of $b$. From top left to the bottom right: Panel (i) $b = 1/3\,$, (ii) $b= 1/3.5\,$, (iii) $b= 1/4.15\,$,  (iv) $b= b_{1}(2)=0.2396825654\,$, (v)  $b=1/4.25\,$, (vi) $b= b_{0}(2)= 0.2240838469\,$, (vii) $b=1/5\,$, (viii)  $b=1/200\,$, (ix) $b=0\,$.} \label{fig2} }
\end{figure}

For any $k_0 > 0$ all additional non-trivial resonances and all
values of $b$ for which  $k_0$ is stable can be determined by analyzing the zeros of the function $\widehat{r}$
with $\widehat{r}(k,b) = \omega(k,b) - \omega(k-k_0,b) - \omega(k_0,b)$ 
on $[k_0/2,\infty) \times \R^+_0$ and using the symmetry of $\omega$ as discussed above. It turns out that for all $k_0 > 0$ 
there exist a smallest $b_{1}=b_{1}(k_0) \in (0,1/3)$, a largest $b_{0}=b_{0}(k_0) \in (0,b_1)$ and a strictly monotonically decreasing function  $k_1 \in C^0((0,b_{0}))$ with $k_1(b)> k_0$ for all $b \in (0,b_0)$ and $k_1(b) \to \infty$ for $b \to 0$
such that $\widehat{r}$ has on $[k_0/2,\infty) \times \R^+_0$ no other zeros than $(k_0,b)$ if $b \in \{0\} \cup
(b_1, \infty)$ and exactly two zeros $(k_0,b)$ and $(k_1(b),b)$ if $b \in (0,b_0)$; see Figure \ref{fig2}.

Another difficulty is the fact that the normal-form transform $R \mapsto R + \eps N({\Psi},R)$ may lose regularity, this means that it maps the Sobolev space $H^{n}(\R,\C)$ into $H^{n-j}(\R,\C)$ for a $j \in (0,n]$. A loss of regularity happens, for example, if $B(W,W)$ contains quasilinear terms. A normal-form transform losing regularity may not be invertible. But even if it was invertible, the mapping $\widetilde{R}
\mapsto \eps^2 g({\Psi},\widetilde{R})$ would lose even more regularity such that it would not be possible in general to derive estimates for $\widetilde{R}$ directly from equation \eqref{newerr} - for example, by applying the variation of constants formula and Gronwall's inequality. 

For these reasons the validity of the NLS approximation for 
systems with quasilinear quadratic terms is a highly non-trivial problem, which remained unsolved in general for more than four decades.
The first validity theorems for the NLS approximation in the literature were proven only for systems with special structural properties. 

In \cite{Kal88},
the NLS approximation was justified for quasilinear systems without quasilinear quadratic terms.  Moreover, semilinear quadratic terms were only admitted if they cause no resonances or trivial resonances. In this situation, the method of normal-form transforms discussed above can be successfully used.

In \cite{Schn98Nodea, Schn98JDE, Schn05, DS06}, the method of normal-form transforms was further developed to make it applicable to systems with non-trivial resonances at $k=\pm k_0$, additional non-trivial resonances with the property that the NLS subspace is stable with respect to those resonances or in case of analytic initial data also additional non-trivial resonances with the property that the NLS subspace is unstable with respect to those resonances.  

In \cite{SW10}, the validity of the NLS appro\-ximation was obtained for a quasilinear reduced model of the two-dimensional water wave problem with finite depth and without surface tension. This reduced model shares with the Lagrangian formulation of the two-dimensional water wave problem some of the principal difficulties which have to be overcome for a validity proof for the NLS approximation, for example the fact that the quadratic nonlinearity loses regularity of half a derivative. In this case
the elimination of the quadratic terms is possible with the help of normal-form transforms. The cubic nonlinearity of 
the transformed system then lose one derivative and can be handled by using a Cauchy-Kowalevskaya argument.

For the quasilinear KdV
equation, the NLS approximation was justified by simply applying a
Miura transform \cite{Schn11}.
Another approach to address the problem of the validity of the NLS approximation for a dispersive equation can be found in \cite{MN13}.
In \cite{CDS15}, the NLS approximation of time oscillatory long wave solutions to wave equations with quasilinear quadratic terms was justified. Because of the scaling behavior of the long wave solutions it is not necessary to eliminate the quadratic terms such that a normal-form transform is not needed.
In \cite{DHSZ15}, it was proven that analytic solutions of a two-dimensional wave equation with a quadratic nonlinearity can be approximated with the help of a two-dimensional NLS equation if the set of resonances is separated from the set of integer multiples of the basic wave vector $k_0 \in \R^{2}$ of the underlying carrier wave.  

In \cite{D17}, the NLS approximation was justified for a nonlinear Klein-Gordon equation with a quasilinear quadratic term, which is the first validity proof of the NLS approximation 
of a nonlinear hyperbolic equation with a quasilinear quadratic term losing regularity of more than half a derivative by error estimates in Sobolev spaces. The linear dispersion relation of the Klein-Gordon equation causes no resonances. The loss of regularity is overcome by using the so-called modified energy method. The main idea of this method is as follows. Instead of performing the normal-form transform \eqref{inft} explicitly and estimating the transformed error $\widetilde{R}$ the normal-form transform is only used to construct an energy $\mathcal{E}_s$ which is an appropriate adaptation of 
\begin{equation} \label{ien}
\mathcal{E}_s = \sum_{l=0}^s  \sum_{j=1}^2   \Big( \frac{1}{2}   \|\partial_x^{l} R_j\|_{L^2}^2  +  \veps \int_{\R}\partial_x^{l} R_j \, \partial_x^{l} N_j(\Psi_c,R)\,dx \Big)
\end{equation}
for a sufficiently large $s >0$. Since $\mathcal{E}_s$ differs from $\|\widetilde{R}\|_{H^s}^2$ only by terms of order $\mathcal{O}(\veps^2)$, the evolution equations of $\mathcal{E}_{s}$ and $\|\widetilde{R}\|_{H^s}^2$ share the property that their right-hand sides are of order $\mathcal{O}(\veps^2)$. The energy $\mathcal{E}_s$ has the advantage that in the case of a normal-form transform which loses regularity the right-hand side of
the evolution equation of $\mathcal{E}_{s}$ has better regularity properties than the right-hand side of the evolution equation of $\|\widetilde{R}\|_{H^s}^2$.

An early version of a modified energy can be found in \cite{D12} as an ingredient to simplify and generalize the proof of the error estimates for the KdV approximation of the water wave problem compared with the alternative proofs in \cite{SW00, SW02}.

The first modified energy which was used to overcome regularity problems in quasilinear equations was constructed in \cite{HITW13}. The modified energy from this article was further developed in \cite{HIT14, IT14a, IT14b, HaIT16}  to apply it to prove large time and global existence results for the water wave problem in holomorphic coordinates. 

A similar modified energy as in \cite{D17} was constructed in \cite{DH18} to justify the NLS approximation for a quasilinear equation whose linear dispersion relation causes resonances. In \cite{CW17}, another modified energy was introduced to improve the NLS approximation result from \cite{SW10}.

In \cite{H19}, the modified energies from \cite{D17, DH18} were combined and extended to prove the validity of the NLS appro\-ximation for two further quasilinear quadratic dispersive systems. One system is a reduced model of the two-dimensional water wave problem with finite depth and $b \geq 1/3$,  which shares with the arc length formulation of the two-dimensional water wave problem some of the principal difficulties which have to be overcome for a validity proof for the NLS approximation, for example the fact that the nonlinearity loses regularity of one derivative. The other system is the first dispersive system containing a quasilinear quadratic nonlinearity that loses regularity of $m$ derivatives with an arbitrary $m > 0$ for which the NLS approximation was justified.

For the water wave problem, all justification results for the NLS approximation in the previous literature are restricted to the case without surface tension.
For the two-dimensional water wave problem with infinite depth and without surface tension in Lagrangian coordinates, the NLS approximation was justified in \cite{TW11} by finding an alternative kind of a transform adapted to the special structure of that problem. 
For the three-dimensional water wave problem with infinite depth and without surface tension in Lagrangian coordinates, the two-dimensional NLS approximation was justified in \cite{T14} in an analogous way.

In \cite{DSW12}, the validity of the NLS approximation was proven for the two-dimensional
water wave problem with finite depth and without surface tension in Lagrangian coordinates. In these coordinates, the evolutionary system has a quasilinear quadratic nonlinearity losing regularity of only half a derivative in the case without surface tension. The occurring resonances are handled with the help of the same strategy as in \cite{DS06}. Despite the loss of regularity the normal-form transform can be inverted, which is proven by interpreting the normal-form transform as a system of diffe\-rential equations whose solvability is obtained with the help of appropriate a priori estimates in Sobolev spaces. The loss of one derivative in the evolutionary system for the transformed error can be handled with the help of the same Cauchy-Kowalevskaya argument as in \cite{SW10}.

In \cite{IT18}, the NLS approximation for the two-dimensional water wave problem with infinite depth and without surface tension in holomorphic coordinates was justified by using the modified energy method.

In the present paper, we solve the open problem of justifying the NLS approximation for the full two-dimensional water wave problem with finite depth and with surface tension.  Our approximation result is valid both for the case without surface tension and for the case with surface tension if there are no other non-trivial resonances than $\pm k_0$ or $k_0$ 
is stable. Our error estimates are uniform with respect to the strength of the surface tension as the height of the wave packet and the surface tension go to zero.
We prove the following 

\begin{theorem}
\label{mainresult}
Let $\omega$ be the dispersion relation \eqref{wwdisp} of the two-dimensional water wave problem (\ref{surf2})--(\ref{pot}). Moreover, let $k_0 >0$ and $s \geq 10$.  Then there exist 
$b_0,\, b_1 \in \R$ with $0 < b_0 < b_1 < 1/3$ such that the following holds.
For all ${\tau}_0, C_0 > 0$ there exist an $\epsilon_0 > 0$ and a function $C \in C^0(\mathcal{B}, \R^+)$, where $\mathcal{B}$ is the set of all $b \in  \R^+_0 \setminus [b_0,b_1]$ for which $k \mapsto \omega(k,b)$  satisfies \eqref{nrb1}, \eqref{nrb3} and \eqref{nrb4} with $M=6$, such that for all $b \in \mathcal{B}$, all solutions $A \in
C^0([0,{\tau}_0],H^{s}(\R,\C))$ of the NLS equation (\ref{NLS})
with 
$$ 
\sup_{\tau \in [0,{\tau}_0]} \| A(\cdot,\tau) \|_{H^{s}(\R,\C)} \leq C_0
$$  
and all $\epsilon \in (0,\epsilon_0)$ there exists a solution 
$$
(\eta,v_1)\in
C^0([0,{\tau}_0 \epsilon^{-2}],(H^{s}(\R,\R))^2)
$$
of 
(\ref{surf2})--(\ref{pot}) which satisfies
\begin{align} \label{errest}
  \sup_{t\in[0,{ {\tau}_0} {\epsilon^{-2}}]} \left\| \left( \ba{c}
\eta  \\ v_1 \ea \right)\!(\cdot,t) -
    {\epsilon}\, \Psi_{NLS}(\cdot,t)\,\varphi(k_0,b)\right\|_{(H^{s}(\R,\R))^2} 
\leq C(b)\, {\epsilon^{3/2}}\,,
\end{align}
where
$$
\Psi_{NLS}(x,t) 
\,=\, A(\eps
(x-\partial_k \omega(k_0,b) t),\eps^2t) e^{i( k_0 x - \omega(k_0,b) t)} + \mathrm{c.c.} 
$$
 and
$\varphi(k_0,b) \in \R^2$ is an explicitly computable vector. In particular, the error estimate \eqref{errest} is uniform with respect to $b$ as $b$ and $\eps$ go to zero.
\end{theorem}

The error of order $\OO(\eps^{3/2})$ is small compared with the solution 
$(\eta,v_1)$ and the approximation $\eps \Psi_{NLS}$, which are both of order
$\OO(\eps)$ in $ L^{\infty} $ such that the dynamics 
of the NLS equation can be found in the two-dimensional water wave problem, too. Our theorem guarantees that, for instance, parts of 
the dynamics of time-periodic solitary wave solutions  present  in 
the NLS equation for $ \partial_k^{2} \omega(k_0,b) $ and $ \nu(k_0,b) $ 
having the same sign can be found approximately 
in the water wave problem.
For a discussion of the values of $ \nu(k_0,b) $ in \eqref{NLS},
see also \cite[Figure 4.15, p.~321]{AS81}. 

It should be noted that the smoothness in our error bound is equal to the assumed smoothness of the amplitude. This can be achieved by using a modified approximation which has compact support in Fourier space but differs only slightly from $\eps \Psi_{NLS}$. Such an approximation can be constructed because the Fourier transform of $\eps \Psi_{NLS}$ is sufficiently strongly concentrated around the wave numbers $\pm k_0$.

The constants $b_0$ and $b_1$ from Theorem \ref{mainresult} can be chosen in such a way that
the following holds. $b_1$ is the smallest number such that for all $b \in (b_1,\infty)$ there are
no other non-trivial resonances than $\pm k_0$ and $b_0$ is the largest number such that $k_0$ is stable for all $b \in (0,b_0)$. For $k_0=2$ this choice of $b_0$ and $b_1$ is presented in Figure \ref{fig2}. One can see that the length of the interval $[b_0, b_1]$, which contains all values of $b$ for which $k_0=2$ is unstable and therefore the validity of the NLS approximation can not be expected for all sufficiently small initial data in the Sobolev space $H^{s}(\R,\C)$, is very small. The same is true for the corresponding interval for any other $k_0>0$. Moreover, for all $k_0 >0$ the number of values of $b$ for which the corresponding dispersion relation $k \mapsto \omega(k,b)$ does not satisfy \eqref{nrb1}, \eqref{nrb3} and \eqref{nrb4} with $M=6$
 is finite.

Now, we explain the main ideas of the proof of Theorem \ref{mainresult} and the plan of the paper. Like in many other proofs of related estimates in the literature we will assume in our proof that $s$ is an integer in order to simplify the analysis by using Leibniz's rule,  but our proofs can be generalized to be valid for all $s \geq 10$.

We perform our proof in the arc length formulation of the two-dimensional water wave problem.
The main advantage of this formulation is that in the corresponding evolutionary system
the surface tension dependent term with the most derivatives is linear,
which allows us to prove the desired uniform error estimates. Transferring the estimates into Eulerian coordinates, we do not lose powers of $\eps$ since in the scaling regime of the NLS equation, the coordinates of the free surface in arc length parametrization are very close to Eulerian coordinates. The same advantages have already been used in the proof of the validity of the KdV approximation for the two-dimensional water wave problem in the arc length formulation in \cite{D12}.

In Section 2 we review the arc length formulation and identify the linear terms, the quadratic terms and the terms losing regularity in the corresponding evolutionary system. Then we diagonalize the linear part of the system to obtain a system which has the structure of \eqref{abstrsyst}.
In Section 3 we present the formal derivation of the NLS approximation for this system. Section 4 is devoted to the error estimates.

In order to perform the error estimates we use the modified energy method. The modified energy we construct is a subtle generalization of the energies in \cite{D17, DH18}.
The normal-form transform behind our energy is an extension of a normal-form transform of the form \eqref{inft} and \eqref{nkernel2} in order to handle the non-trivial resonances. 

The problems with the resonances at $\pm k_0$ are circumvented by rescaling the error in Fourier space as in \cite{DS06, SW10, DSW12} with the help of the weight function
\begin{align*}
\widehat{\vartheta}(k) = \left\{ \begin{array}{ll} 1\,, \qquad & |k| >\delta_{0} \,,\\[2mm] \eps + (1-\eps) |k|/\delta_{0} \,,  \qquad 
& |k| \leq \delta_{0}  \end{array} \right.
\end{align*}
with a $\delta_0 = \delta_{0}(b) \in\, (0, k_0/20)$ sufficiently small, but independent of $\veps$. The choice of the weight function makes sense because the quadratic terms 
in the evolutionary system of the two-dimensional water wave problem in the arc length formulation vanish at $k=0$ such that the Fourier transform of the error can grow only slowly for $|k| \ll 1$. But since $\widehat{\vartheta}^{-1}(k) = \mathcal{O}(\eps^{-1})$ for $|k| \leq \delta_{0}$, the normal-form transform has to be extended by an appropriately chosen trilinear mapping.

To control the additional non-trivial resonances the terms of the form \eqref{ien} in our energy are slightly modified by weight functions and correction functions similar to those in 
the final energy in \cite{DS06}, which are motivated by the conserved quantity \eqref{staben}.  

Due to the structure of the evolutionary system of the two-dimensional water wave problem in the arc length formulation all terms
on the right-hand side of the evolution equation of our energy can either directly be estimated by the energy or be identified as time derivatives of time dependent integrals. By adding these integrals to the energy we obtain our final energy, which can be bounded with the help of Gronwall's inequality over the desired timespan of order $\mathcal{O}(\eps^{-2})$. Since this energy controls a Sobolev norm of the error, we finally obtain our approximation result.

The methods of proof developed in the present paper can also be used to prove the validity of the NLS approximation for other dispersive systems with quasilinear quadratic terms. 
\medskip

{\bf Notation}. 
We denote the 
Fourier transform of a function $u \in L^2(\R,\K)$ with $\K=\R$ or $\K=\C$ by
$$\mathcal{F}(u)(k) = \widehat{u}(k) = \frac{1}{2\pi} \int_{\R} u(x) e^{-ikx} dx. $$
Let $H^{s}(\R,\K)$ be
the space of functions mapping from $\R$ into $\K$
for which
the norm
$$ \| u \|_{H^{s}(\R,\K)} = \left(\int_{\R} |\widehat{u}(k)|^2 (1+|k|^2)^{s} 
dk \right)^{1/2} $$ 
is finite. We also write $L^2$ and $H^{s}$ instead of $L^2(\R,\R)$ and $H^{s}(\R,\R)$.
Moreover, we use the space
$ L^p(m)(\R,\K) $ defined by $ u \in L^p(m)(\R,\K) \Leftrightarrow u \sigma^m \in L^p(\R,\K)$, where
$ \sigma(x) = (1+x^2)^{1/2}$.

Furthermore, we write $A
\lesssim B$ if $A \leq C B$ for a constant $C>0$ which does not depend on $A$ and $B$, as well as $A = \mathcal{O}(B)$ if 
$|A| \lesssim B$. 

\section{The water wave problem in the arc length formulation}
\label{sec2}
In the following we review the essential points of the arc length formulation of the two-dimensional water wave problem with finite depth. Let $P(t): \R \rightarrow \Gamma(t),\,\alpha \mapsto P(\alpha,t) = (x(\alpha,t), y(\alpha,t))$ be a parametrization of the free top surface $\Gamma(t)$ by arc length, that means, we have 
\begin{equation} \label{arcxy}
(x_{\alpha}^2+y_{\alpha}^2)^{1/2}=1.
\end{equation} 
Let $U$ and $T$ be the normal and the tangential velocity on the free top surface measured in the coordinates of the arc length parametrization, that means that
\begin{equation} \label{xy}
(x,y)_{t}(\alpha,t) = U(\alpha,t) \hat{n}(\alpha,t) + T(\alpha,t) \hat{t}(\alpha,t),
\end{equation}
where $\hat{n} =(-\sin \theta, \cos \theta)$ and $\hat{t}=(\cos \theta, \sin \theta)$ are the upward unit normal vectors and the unit tangential vectors to the free top surface and $\theta = \arctan (y_{\alpha}/x_{\alpha})$ are the tangent angles on the free top surface. 
Because of \eqref{arcxy}, 
$T$ satisfies
\begin{eqnarray}
T_{\alpha} - \theta_{\alpha} U = 0.
\end{eqnarray}
Integrating this relation determines $T$ depending on $\theta$ and $U$ up to a constant. Since arc length parametrizations are invariant under translations, this constant can be set to $0$ without loss of generality. This implies
\begin{equation} \label{T} 
T(\alpha,t) = \int_{-\infty}^{\alpha} \theta_{\alpha}(\beta,t)U(\beta,t) \,d\beta\,.
\end{equation}
The normal velocity $U$ is governed by the incompressible Euler's
equations (\ref{euler})--(\ref{incompr}), the boundary conditions
(\ref{surf})--(\ref{bot}) and the form of the free top surface. 

From now on, we consider irrotational flows. Then the normal velocity
$U$ can be expressed in terms of the free top surface and the physical tangential velocity $v$ of the fluid particles on the free top surface, where the evolution of $v$ is determined by (\ref{euler})--(\ref{bot}) and the form of the free top surface. Moreover, as long as $y(\cdot,t)$, $\theta(\cdot,t)$ and $v(\cdot,t)$ are sufficiently regular and localized, for example, $y(\cdot,t), v(\cdot,t) \in L^2$ and 
$\theta(\cdot,t) \in H^2$, then, due to \eqref{arcxy}, the evolution of $x$ 
is completely determined by the evolution of $\theta$ and therefore
$U(\cdot,t)$ can be  
represented as a function of
$y(\cdot,t)$, $\theta(\cdot,t)$ and $v(\cdot,t)$. 

Finally, using all the above information, one can derive the following 
evolutionary system:
\begin{eqnarray}
y_{t} &\!\!=\!\!& U \cos \theta + {T} {y}_{\alpha}\,, \label{ytilde}\\[2mm]
{v}_{t} & \!\!=\!\!& - {y}_{\alpha} + b {\theta}_{{\alpha}\, {\alpha}} -  {\delta} {\delta}_{{\alpha}} + U {\theta}_{t}\,, \label{vtilde}
\\[2mm]
\theta_{t} &\!\!=\!\!& U_{\alpha}  + {T} {\theta}_{{\alpha}}\,, \label{thetatilde}
\\[2mm]
{\delta}_{{\alpha}\, t} &\!\!=\!\!& -  {c}{\theta}_{{\alpha}} + b {\theta}_{{\alpha}\,{\alpha}\,{\alpha}} -  ({\delta}{\delta}_{{\alpha}})_{{\alpha}} + (U_{\alpha}+ v {\theta}_{{\alpha}})^{2}\,, \label{deltaalphatilde}
\\[2mm]
y_{\alpha} &\!\!=\!\!& \sin \theta\,, \label{reltheta} 
\\[2mm]
{\delta} & \!\!=\!\! & {v} - {T}\,, \label{deltatildevtilde}
\end{eqnarray}
where
\begin{eqnarray}
\label{ctilde}
{c} &\!\!=\!\! & {U}_{t} + v {\theta}_{{t}}+ {\delta} {U}_{{\alpha}}+ {\delta} v {\theta}_{{\alpha}} + \cos{\theta}\,. 
\end{eqnarray}

For further details of the derivation of this system, an explicit formula for $U$ and the local well-posedness of the system in Sobolev spaces, we refer to \cite{AM05,D12}.

The evolution equations \eqref{thetatilde} and \eqref{deltaalphatilde} are included because they have better regularity properties than the evolution equations for the spatial derivatives of $y$ and $v$. The evolutionary system
\eqref{ytilde}--\eqref{ctilde} could be posed entirely in terms of the variables $\theta$ and $\delta_{\alpha}$ as it is done in \cite{AM05}. But since $y$ and $v$ are the physically relevant variables for which we would like to derive and justify the NLS approximation,
we keep these variables. Consequently, we need the equations \eqref{reltheta}--\eqref{deltatildevtilde} as consistency conditions between the systems \eqref{ytilde}--\eqref{vtilde} and \eqref{thetatilde}--\eqref{deltaalphatilde}.

The main advantage of system \eqref{ytilde}--\eqref{ctilde} is that in the case of surface tension, i.e., for $b>0$,
the term with the most derivatives in \eqref{ytilde}--\eqref{ctilde} is linear.

In order to derive the NLS approximation and to prove the error estimates we need to extract the linear and the quadratic components of system \eqref{ytilde}--\eqref{ctilde}. In this context, the linear operator $K_0$ defined by its symbol
\begin{equation} \label{K0}
\widehat{K}_0(k)\, =\, -i\, \mathrm{tanh}(k)
\end{equation}
for all $k \in \R$ plays an important role. The operator $K_0$ is the linearization
of the operator $\mathcal{K}$ from \eqref{Knl} around the trivial solution $(\eta, \phi_x) =(0,0)$. We present some properties of $K_0$ which we will need below. We have the following 
\begin{lemma}
Let $s \geq 0$ and $q > \frac12$. Then we have
 \begin{align}
\|K_0 f\|_{H^s} & \,\lesssim \;\|f\|_{H^{s}} \label{K01} \,,\\[2mm]
\|[K_0,g]f\|_{H^s} & \,\lesssim \;\|g\|_{H^{s+q}} \|f\|_{H^{0}}\,,\\[2mm]
\|[K_0,g]f\|_{H^s}  & \, \lesssim \; \|g\|_{H^{s}} \|f\|_{H^{q}}\,,\\[1mm]
\|(1+K_0^2)f\|_{H^s} & \, \lesssim \; \|f\|_{H^{0}} \label{K04}\,.
\end{align}
\end{lemma}

\textbf{Proof.}
The lemma is a special case of Lemma 3.7 and Lemma 3.8  in \cite{D12}.
\qed
  
\bigskip
With the help of $K_0$ one obtains the following expansion of the system
\eqref{ytilde}--\eqref{ctilde}. 
\begin{lemma} \label{expsyst}
\begin{eqnarray} \label{y1}
\quad {y}_t &\!\!=\!\!& K_0{v} + \big(K_0[K_0,{y}]{v} - (1+K_0^2)({y}{v})
 \big)_{\ua} + m_{I}\,,\\[2mm]
\qquad {v}_t &\!\!=\!\!& -{y}_{\ualpha} + b {y}_{\ualpha\ualpha \ualpha} - 
{\frac{1}{2}} ({v}^2)_{\ualpha} + 
{\frac{1}{2}} ((K_0{v})^2 )_{\ualpha}
+ m_{I\!I}\,, \label{v} \\[2mm]
\qquad {\kappa}_t &\!\!=\!\!& K_0{\delta}_{\ua\ua} - 
({\delta}{\kappa})_{\ualpha} + \big(K_0[K_0,{y}]{\delta}_{\ualpha} -
(1+K_0^2)({y}{\delta}_{\ualpha})\big)_{\ua\ua} 
\nonumber\\[2mm] &&
+ \big(K_0[K_0,{\theta}]{\delta}_{\ualpha} -
(1+K_0^2)({\theta}{\delta}_{\ualpha}) \big)_{\ualpha}  +  (m_{I\!I\!I})_{\alpha}\,, \label{kappa1} \\[2mm]
\qquad {\delta}_{\ua\ualpha t} &\!\! = \!\!& -\kappa_{\ua} + b {\kappa}_{\ualpha\ualpha\ualpha} 
+ \big( (K_0{\theta}-b K_0{\kappa}_{\ualpha} + {c}_0) {\kappa} \big)_{\ualpha} \nonumber \\[2mm]&&
- ({\delta}{\delta}_{\ualpha\ua})_{\ualpha}  
- (({\delta}_{\ualpha})^2)_{\ualpha} +  ((K_0{\delta}_{\ualpha})^2)_{\ualpha} + (m_{I\!V})_{\alpha}\,, \label{deltaalphaalpha1} \\[2mm]
\qquad y_{\ua} &\!\! = \!\!& \theta + m_{V}, \label{yalphatheta1} \\[2mm]
\qquad {\theta}(\ualpha,t) & \!\!=\!\!&  \int_{-\infty}^{\ualpha} 
{\kappa}({\beta},t)\, d{\beta}\,,
\label{thetakappa1} \\[2mm]
\qquad {\delta}(\ualpha,t) &\!\! =\!\!& {v}(\ualpha,t)  - 
\int_{-\infty}^{\ualpha} ((K_0{v})
{\kappa})({\beta},t)\, d{\beta} + m_{V\!I}(\ualpha,t)\,, \label{deltav1}
\end{eqnarray}
where 
\begin{equation} \label{estmI}
\|m_{I}\|_{H^s} \lesssim \, (\|{y}\|_{L^{2}}^2 + \|{\theta}\|_{H^{s}}^2) (\|{v}\|_{L^{2}} + \|{\delta}_{\ualpha}\|_{L^{2}}) 
\end{equation}
for $s\geq 1$, as long as $\|{y}\|_{L^{2}},  \|{\theta}\|_{H^{s}} \lesssim 1$,
\begin{equation}
\|{m}_{I\!I}\|_{H^s} \lesssim \, (\|{y}\|_{L^{2}} + \|{\theta}\|_{H^{s+1}}) (\|{v}\|_{L^{2}}^2 + \|{\delta}_{\ualpha}\|_{H_{\ualpha}^{s}}^2) 
\end{equation}
for $s\geq 2$, as long as $\|{y}\|_{L^{2}}, \|{\theta}\|_{H^{s+1}} \lesssim 1$, 
\begin{equation} \label{estmIII}
\|(m_{I\!I\!I})_{\alpha}\|_{H^s} \lesssim \, (\|{y}\|_{L^{2}}^2 + \|{\theta}\|_{H^{s+1}}^2) (\|{v}\|_{L^{2}} + \|{\delta}_{\ualpha}\|_{H^{2}}) 
\end{equation}
for $s\geq 2$, as long as $\|{y}\|_{L^{2}}, \|{\theta}\|_{H^{s+1}} \lesssim 1$,
\begin{equation}
\|(m_{I\!V})_{\alpha} \|_{H^{s}} \lesssim \, (\|{y}\|_{L^{2}} + \|{\theta}\|_{H^{s+1}}) (\|{v}\|_{L^{2}}^2 + \|{\delta}_{\ualpha}\|_{H^{s+1}}^2) 
\end{equation}
for $s \geq 2$, as long as $\|{y}\|_{L^{2}}, \|{\theta}\|_{H^{s+1}} \lesssim 1$, 
\begin{equation} \label{estmV}
\|m_{V}\|_{H^{s}} \lesssim \, \|{\theta}\|_{H^{s}}^3 
\end{equation}
for $s \geq 1$, as long as $\|{\theta}\|_{H^{s}} \lesssim 1$,
\begin{equation} \label{estmVI}
\|m_{V\!I}\|_{C_{b}^0} + \|(m_{V\!I})_{\ua}\|_{H^{s-1}}\lesssim \, (\|{y}\|_{L^{2}}^2 + \|{\theta}\|_{H^{s}}^2) (\|{v}\|_{L^{2}} + \|{\delta}_{\ualpha}\|_{H^{s-2}}) 
\end{equation}
for $s \geq 2$, as long as $\|{y}\|_{L^{2}},  \|{\theta}\|_{H^{s}} \lesssim 1$, and
\begin{equation} \label{estc1tilde}
\|{c}_0 \|_{H^{s}} \lesssim \, \|{y}\|_{L^{2}}^{2} + \|{\theta}\|_{H^{s}}^{2} +
 b \|{\theta}\|_{H^{s+1}}^{2} + 
\|{v}\|_{L^{2}}^{2} + \|{\delta}_{\ualpha}\|_{H^{s}}^{2} 
\end{equation}
for $s \geq 6$, as long as $\|{y}\|_{L^{2}}, \|{\theta}\|_{H^{s}},
 \sqrt{b} \|{\theta}\|_{H^{s+1}}, 
\|{v}\|_{L^{2}}, \|{\delta}_{\ualpha}\|_{H^{s}} \lesssim 1$. 
Moreover, we have
\begin{equation}
\|\theta\|_{H^s} \lesssim \, \|y\|_{L^{2}} + \|\kappa\|_{H^{s-1}} \label{estthetaykappa}
\end{equation}
for $s\geq 1$, as long as $\|y\|_{L^{2}} + \|\kappa\|_{L^{2}} \ll 1$, and
\begin{equation}
\|\delta_{\alpha}\|_{H^s} \lesssim \,
(1+\|y\|_{L^{2}} + \|\kappa\|_{H^{1}})(\|v\|_{L^{2}} + \|\delta_{\alpha\alpha}\|_{H^{s-1}})
 \label{estdeltaa}
\end{equation}
for $s\geq 1$, as long as $\|y\|_{L^{2}} + \|\kappa\|_{H^{1}} \ll 1$.
\\[2mm]
All bounds are uniform with respect to $b \lesssim 1$.
\end{lemma}

\textbf{Proof.}
The expansions \eqref{y1}--\eqref{deltav1} follow directly from Lemma 3.9 in \cite{D12}. The bounds \eqref{estmI}--\eqref{estdeltaa} follow directly from the
bounds in the Lemmas 3.1--3.9 in \cite{D12}, which are also uniform with respect to $b \lesssim 1$, and the well-known interpolation inequality
 \begin{align}
\|f_{\alpha}\|_{L^2} & \,\lesssim \; \mu \|f\|_{L^{2}} + \mu^{-1} \|f_{\alpha\alpha}\|_{L^2} \label{intH1}
\end{align}
for all $f \in H^2(\R)$ and all $\mu >0$.
\qed  

\medskip
In Lemma \ref{expsyst} we replaced the evolution equations for $\theta$ and $\delta_{\alpha}$ by the evolution equations for the respective spatial derivative, where we used that the spatial derivative of the tangent angle $\theta$ is the curvature $\kappa$. Due to this replacement, the Fourier transform of all quadratic terms in the resulting evolutionary system \eqref{y1}--\eqref{deltav1} vanishes at $k=0$ such that we can handle the resonances at $k=0$ and $k=\pm k_0$ as explained in the introduction. 

To emphasize the geometric meaning of the terms in system \eqref{y1}--\eqref{deltav1} we included both $\theta$ and $\kappa$ in the system. Therefore, we need the additional equation \eqref{thetakappa1} as a consistency condition. The expansion of \eqref{ctilde} is incorporated directly into \eqref{y1}--\eqref{deltav1}.

In contrast to the justification of the NLS approximation, the occurring resonances need not to be taken into account in the context of proving local 
well-posedness of the water wave problem. Hence, local well-posedness can be shown by analyzing the evolution equations for $\theta$ and $\delta_{\alpha}$, whereas for proving the validity of the NLS approximation we will need the evolution equations for $\kappa$ and $\delta_{\alpha\alpha}$.

\medskip
We diagonalize system \eqref{y1}--\eqref{deltav1} by 
\begin{align}  
 \left(
\begin{array}{c} {y} \\ {v}
\end{array}
\right) & = \,
\left(
\begin{array}{cc}
  \sigma^{-1} &  -\sigma^{-1}   \\
 1 & 1    
\end{array}
\right)
\left(
\begin{array}{c} {u}_{-1} \\ {u}_{1} 
\end{array}
\right),\\[2mm]
 \left(
\begin{array}{c} {\kappa} \\ {\delta_{\ua\ua}}
\end{array}
\right) & =\, 
\left(
\begin{array}{cc}
  \sigma^{-1} &  -\sigma^{-1}   \\
 1 & 1    
\end{array}
\right)
\left(
\begin{array}{c} {u}_{-2} \\ {u}_{2} 
\end{array}
\right)\,,
 \end{align}
where $\sigma^{-1}$ is the inverse of the linear operator $\sigma$ with the symbol
\begin{equation} \label{rhosym}
\sigma(k)= \sigma(k,b) = \sqrt{\frac{k+bk^3}{\mathrm{tanh}(k)}} \,.
\end{equation}
Then we have
\begin{align}  
 \left(
\begin{array}{c} {u}_{-1} \\ {u}_{1}
\end{array}
\right) & = \, 
\frac{1}{2}\left(
\begin{array}{cc}
  \sigma &  1   \\
 -\sigma & 1    
\end{array}
\right)
\left(
\begin{array}{c} {y} \\ {v}
\end{array}
\right), \\[2mm]
 \left(
\begin{array}{c} {u}_{-2} \\ {u}_{2}
\end{array}
\right) & =\, 
\frac{1}{2}\left(
\begin{array}{cc}
  \sigma &  1   \\
 -\sigma & 1    
\end{array}
\right)
\left(
\begin{array}{c} {\kappa} \\ {\delta_{\ua\ua}}
\end{array}
\right),
 \end{align}
and Lemma \ref{expsyst} yields
\begin{eqnarray}
{(u_{-1})}_t &\!\!=\!\!&  - i \omega u_{-1} -  \frac{1}{4} (( u_{-1} + u_{1})^2)_{\ua} +    
 \frac{1}{4} ((K_0( u_{-1} + u_{1}))^2)_{\ua} 
\nonumber\\[2mm] &&
+  \frac{1}{2} \big(\sigma K_0[K_0, \sigma^{-1} ( u_{-1} - u_{1})]( u_{-1} + u_{1})\big)_{\ua} 
\nonumber\\[2mm] &&
-  \frac{1}{2} \big(\sigma (1+K_0^2)(\sigma^{-1} ( u_{-1} - u_{1})( u_{-1} + u_{1}))\big)_{\ua} 
\nonumber\\[2mm] &&
+ m_{-1}\,, \label{w-1} \\[3mm]
{(u_{1})}_t &\!\!=\!\!&  i \omega u_{1} - \frac{1}{4} (( u_{-1} + u_{1})^2)_{\ua} +    
\frac{1}{4} ((K_0( u_{-1} + u_{1}))^2)_{\ua} 
\nonumber\\[2mm] &&
-  \frac{1}{2} \big(\sigma K_0[K_0, \sigma^{-1} ( u_{-1} - u_{1})]( u_{-1} + u_{1})\big)_{\ua} 
\nonumber\\[2mm] &&
+  \frac{1}{2} \big(\sigma (1+K_0^2)(\sigma^{-1} ( u_{-1} - u_{1})( u_{-1} + u_{1}))\big)_{\ua} 
\nonumber\\[2mm] &&
+ m_1\,, \label{w1} 
\\[4mm]
{(u_{-2})}_t &\!\!=\!\!&  - i \omega u_{-2} 
-(\partial_{\ua}^{-2}(u_{-2}+u_{2})u_{-2})_{\ua}
\nonumber\\[2mm] &&
-  \frac{1}{2} ([\sigma,\partial_{\ua}^{-2}(u_{-2}+u_{2})]  \sigma^{-1} ( u_{-2} - u_{2}))_{\ua}
\nonumber\\[2mm] &&
+  \frac{1}{2} (K_0 \partial_{\ua}^{-1} \sigma^{-1}(u_{-2}-u_{2})  \sigma^{-1} ( u_{-2} - u_{2}))_{\ua}
\nonumber\\[2mm] &&
-  \frac{1}{2}b  (\sigma^{-1}(u_{-2}-u_{2})  K_0 \sigma^{-1} ( u_{-2} - u_{2})_{\ua})_{\ua}
\nonumber\\[2mm] &&
-  \frac{1}{2} ((\partial_{\ua}^{-1} ( u_{-2} + u_{2}))^2)_{\ua} +    
 \frac{1}{2} ((K_0  \partial_{\ua}^{-1} ( u_{-2} + u_{2}))^2)_{\ua} 
\nonumber\\[2mm] &&
+  \frac{1}{2} \big(\sigma K_0[K_0, \sigma^{-1} ( u_{-1} - u_{1})] \partial_{\ua}^{-1} ( u_{-2} + u_{2})\big)_{\ua\ua} 
\nonumber\\[2mm] &&
-  \frac{1}{2} \big(\sigma (1+K_0^2)(\sigma^{-1} ( u_{-1} - u_{1})\partial_{\ua}^{-1}( u_{-2} + u_{2}))\big)_{\ua\ua} 
\nonumber\\[2mm] &&
+  \frac{1}{2} \big(\sigma K_0[K_0, \partial_{\ua}^{-1} \sigma^{-1} ( u_{-2} - u_{2})] \partial_{\ua}^{-1} ( u_{-2} + u_{2})\big)_{\ua} 
\nonumber\\[2mm] &&
-  \frac{1}{2} \big(\sigma (1+K_0^2)(\partial_{\ua}^{-1} \sigma^{-1} ( u_{-2} - u_{2})\partial_{\ua}^{-1}( u_{-2} + u_{2}))\big)_{\ua} 
\nonumber\\[2mm] &&
+  \frac{1}{2} (c_1 \sigma^{-1} ( u_{-2} - u_{2}))_{\ua}
+(m_{-2})_{\ua}\,, \label{u-1}
\end{eqnarray}
\begin{eqnarray}
{(u_{2})}_t &\!\!=\!\!&  i \omega u_{2} 
-(\partial_{\ua}^{-2}(u_{-2}+u_{2})u_{2})_{\ua}
\nonumber\\[2mm] &&
+  \frac{1}{2} ([\sigma,\partial_{\ua}^{-2}(u_{-2}+u_{2})]  \sigma^{-1} ( u_{-2} - u_{2}))_{\ua}
\nonumber\\[2mm] &&
+  \frac{1}{2} (K_0 \partial_{\ua}^{-1} \sigma^{-1}(u_{-2}-u_{2})  \sigma^{-1} ( u_{-2} - u_{2}))_{\ua}
\nonumber\\[2mm] &&
-  \frac{1}{2}b  (\sigma^{-1}(u_{-2}-u_{2})  K_0 \sigma^{-1} ( u_{-2} - u_{2})_{\ua})_{\ua}
\nonumber\\[2mm] &&
-  \frac{1}{2} ((\partial_{\ua}^{-1} ( u_{-2} + u_{2}))^2)_{\ua} +    
 \frac{1}{2} ((K_0  \partial_{\ua}^{-1} ( u_{-2} + u_{2}))^2)_{\ua} 
\nonumber\\[2mm] &&
-  \frac{1}{2} \big(\sigma K_0[K_0, \sigma^{-1} ( u_{-1} - u_{1})] \partial_{\ua}^{-1} ( u_{-2} + u_{2})\big)_{\ua\ua} 
\nonumber\\[2mm] &&
+  \frac{1}{2} \big(\sigma (1+K_0^2)(\sigma^{-1} ( u_{-1} - u_{1})\partial_{\ua}^{-1}( u_{-2} + u_{2}))\big)_{\ua\ua} 
\nonumber\\[2mm] &&
-  \frac{1}{2} \big(\sigma K_0[K_0, \partial_{\ua}^{-1} \sigma^{-1} ( u_{-2} - u_{2})] \partial_{\ua}^{-1} ( u_{-2} + u_{2})\big)_{\ua} 
\nonumber\\[2mm] &&
+  \frac{1}{2} \big(\sigma (1+K_0^2)(\partial_{\ua}^{-1} \sigma^{-1}( u_{-2} - u_{2})\partial_{\ua}^{-1}( u_{-2} + u_{2}))\big)_{\ua} 
\nonumber\\[2mm] &&
+  \frac{1}{2} (c_1 \sigma^{-1} ( u_{-2} - u_{2}))_{\ua}
+(m_2)_{\ua} \label{u1}
\end{eqnarray}
as well as
\begin{eqnarray}
\label{thetaneu} 
\partial_{\ua}^{-1} \sigma^{-1}( u_{-2} - u_{2}) &\!\!=\!\!& \sigma^{-1}( u_{-1} - u_{1})_{\ua} + m_{-3}\,, 
\\[2mm]
\partial_{\ua}^{-2} ( u_{-2} + u_{2}) &\!\!=\!\!& ( u_{-1} + u_{1}) 
- \partial_{\ua}^{-1} \big(K_0 ( u_{-1} + u_{1}) \sigma^{-1} ( u_{-2} - u_{2})\big)
+ m_3\,, \label{deltaneu} 
\end{eqnarray}
where $\omega$ is the linear operator with the symbol
\begin{equation} \label{omegasym}
\omega(k) = \omega(k,b) = \mathrm{sgn}(k) \sqrt{(k+bk^3) \mathrm{tanh}(k)}\,,
\end{equation}
$\partial_{\ua}^{-1}$ defined by the symbol $-ik^{-1}$  and
\begin{equation} \label{estm1}
\|m_{-1}\|_{H^s} \lesssim \, \|{u_{-1}}\|_{L^{2}}^3 + \|u_{1}\|_{L^{2}}^3+\|{u}_{-2}\|_{H^{s-1/2}}^3 + \|{u}_{2}\|_{H^{s-1/2}}^3 
\end{equation}
for $s\geq 2$, as long as $\|{u}_{\pm1}\|_{L^{2}},  \|{u}_{\pm2}\|_{H^{s-1/2}} \ll 1$,
\begin{equation} \label{estm3}
\|m_{1}\|_{H^s} \lesssim \, \|{u_{-1}}\|_{L^{2}}^3 + \|u_{1}\|_{L^{2}}^3+\|{u}_{-2}\|_{H^{s-1/2}}^3 + \|{u}_{2}\|_{H^{s-1/2}}^3 
\end{equation}
for $s\geq 2$, as long as $\|{u}_{\pm1}\|_{L^{2}},  \|{u}_{\pm2}\|_{H^{s-1/2}} \ll 1$,
\begin{equation} \label{estm4a}
\|(m_{-2})_{\ua}\|_{H^s} \lesssim \, \|{u_{-1}}\|_{L^{2}}^3 + \|u_{1}\|_{L^{2}}^3+\|{u}_{-2}\|_{H^{s}}^3 + \|{u}_{2}\|_{H^{s}}^3 
\end{equation}
for $s\geq 2$, as long as $\|{u}_{\pm1}\|_{L^{2}},  \|{u}_{\pm2}\|_{H^{s}} \ll 1$,
\begin{equation} \label{estm5a}
\|(m_{2})_{\ua}\|_{H^s} \lesssim \, \|{u_{-1}}\|_{L^{2}}^3 + \|u_{1}\|_{L^{2}}^3+\|{u}_{-2}\|_{H^{s}}^3 + \|{u}_{2}\|_{H^{s}}^3 
\end{equation}
for $s\geq 2$, as long as $\|{u}_{\pm1}\|_{L^{2}},  \|{u}_{\pm2}\|_{H^{s}} \ll 1$, 
\begin{equation} \label{estm6}
\|m_{-3}\|_{H^s} \lesssim \, \|{u_{-1}}\|_{L^{2}}^3 + \|u_{1}\|_{L^{2}}^3+\|\sigma^{-1}({u}_{-2}-u_2)\|_{H^{s-1}}^3 
\end{equation}
for $s\geq 2$, as long as $\|{u}_{\pm1}\|_{L^{2}},  \|\sigma^{-1}({u}_{-2}-u_2) \|_{H^{s-1}} \ll 1$, 
\begin{equation} \label{estm7}
\|m_{3}\|_{C^0} + \|(m_{3})_{\ua}\|_{H^{s-1}} \lesssim \, \|{u_{-1}}\|_{L^{2}}^3 + \|u_{1}\|_{L^{2}}^3 + \|\sigma^{-1}({u}_{-2}-u_2) \|_{H^{{s-1}}}^3 + \|{u}_{-2}+{u}_{2}\|_{H^{{s-3}}}^3 
\end{equation}
for $s\geq 3$, as long as $\|{u}_{\pm1}\|_{L^{2}}, \|\sigma^{-1}({u}_{-2}-u_2)\|_{H^{s-1}}, \|{u}_{-2}+u_2\|_{H^{s-3}}\ll 1$, 
as well as
\begin{equation} \label{estc0}
\|{c}_1 \|_{H^{s}} \lesssim \, \|{u_{-1}}\|_{L^{2}}^{2} + \|u_1\|_{L^{2}}^{2} +
 \|{u}_{-2}\|_{H^{s-1}}^{2} + 
\|{u}_{2}\|_{H^{s-1}}^{2}  
\end{equation}
for $s \geq 6$, as long as $\|{u}_{\pm1}\|_{L^{2}}, \|{u}_{\pm2}\|_{H^{s-1}} \ll 1$. 
All bounds are uniform with respect to $b \lesssim 1$.

In the diagonalized system 
\eqref{w-1}--\eqref{deltaneu} both $\theta$ and $\kappa$ are expressed in terms of the
variables $u_{-2}$ and $u_2$ such that the consistency condition \eqref{thetakappa1} is not needed anymore.

\medskip
We close this section by collecting some properties of the operator $\sigma$, 
which will be useful for our further argumentation. 
We have the identities
\begin{align}
& 
\sigma K_0[K_0,\sigma^{-1}f]g-\sigma(1+K_0^2)(\sigma^{-1}fg) \nonumber \\
& \qquad 
= \; -gf -K_0g K_0f - [\sigma,g]\sigma^{-1}f - [K_0\sigma,K_0 g]\sigma^{-1}f \label{K0-id}
\end{align}
as well as
\begin{align}
[\sigma,f]\sigma^{-1}g \;
= \;& \sigma^{-1}g \sigma f - gf + [\sigma,\sigma^{-1}g]f  \label{sigma-id} \,.
\end{align}

Moreover, a direct computation using 
\begin{align} \label{rootsigma}
\sigma(k) - \sigma(l) = (\sigma(k)+ \sigma(l))^{-1} (\sigma^{2}(k)- \sigma^{2}(l))
\end{align}
 for all $k,l \in \R$ and the mean value theorem
 yields 
\begin{align}
\|[\sigma,g]f\|_{H^s} & \,\lesssim \;  \|\sigma g_{\alpha}\|_{H^{s-1}} \| f\|_{H^{s-1}} + \|g_{\alpha}\|_{H^{s-1}} \|\sigma f\|_{H^{s-1}} 
\label{sigmak1} 
\end{align}
for $s > 3/2$.  

\section{The derivation of the  NLS approximation}
\label{sec3}

In order to derive the NLS approximation for system \eqref{w-1}--\eqref{deltaneu}, we introduce the vector-valued function
\[
\mathcal{U}: = \left(
\begin{array}{c} \mathcal{U}_{1} \\ \mathcal{U}_{2}
\end{array}
\right)
\; \text{with} \;\,\, \mathcal{U}_j: = \left(
\begin{array}{c} {u}_{-j} \\ {u}_{j}
\end{array}
\right)
\,\, \text{for} \;\, j=1,2,
\]
and make the ansatz
\begin{equation} \label{ans1c}
\mathcal{U} = \veps \tPsi = \veps \left( \ba{c}
\tPsi_{1}\\ \tPsi_2 \ea \right)
\end{equation}
with
\begin{align}
\nonumber
\veps \tPsi_j &\, =\,
\eps {\tPsi}^0_{j1} + \eps {\tPsi}^0_{j-1} +  \eps^2 {\tPsi}^0_{j0} + \eps^2{\tPsi}^0_{j2} 
+ \eps^2{\tPsi}^0_{j-2}\,,
\\[3mm] \nonumber
\eps {\tPsi}^0_{j\pm 1}(\alpha,t) 
& =  \eps {\tA}_{-j\pm 1}^0 (\eps
(\alpha-c_g t),\eps^2 t)
\,\EE^{\pm 1} \!\left( \ba{c} 
1 \\ 0  \ea \right),\\[2mm] \nonumber
\eps^2 \tPsi^0_{j0}(\alpha,t)  
& \,=\,  \left(\begin{array}{c}
\eps^2 {\tA}_{-j0}^0 (\eps (\alpha-c_gt),\eps^2t)\\[1mm]
\eps^2 {\tA}_{j0}^0 (\eps (\alpha-c_gt),\eps^2t)
\\ \end{array}\right),
\\[2mm] \nonumber
\eps^2 {\tPsi}^0_{j\pm2}(\alpha,t)  
&\, =\,  \left(\begin{array}{c} 
\eps^2 {\tA}_{-j(\pm2)}^0 (\eps
(\alpha-c_gt),\eps^2t)\,\EE^{\pm2}\\[1mm]
\eps^2 {\tA}_{j(\pm2)}^0( \eps
(\alpha-c_gt),\eps^2t)\,\EE^{\pm2}
 \end{array} \right),
\end{align}
where $0 < \eps \ll 1$, $j \in \{1,2\}$, $\EE = e^{i(k_0 \alpha -  \omega_0 t)} $, $k_0 >0$, $ \omega_0 = \omega(k_0,b) $, $c_g = \partial_k \omega(k_0,b)$ and 
${\tA}_{m-\ell}^0 = \overline{{\tA}_{m\ell}^0}$.   

Our ansatz leads to waves moving to the right.
For waves moving to the left one has to replace in the above ansatz
the vector $
(1 , 0 )^T $ by $
(0 , 1 )^T $ as well as $ -\omega_0$ by $ \omega_0$
and $c_g$ by $-c_g$.

First, we insert the ansatz \eqref{ans1c} for $\mathcal{U}_1$ into \eqref{w-1}--\eqref{w1}. Then we
replace the dispersion relation $k \mapsto \omega(k,b) $ in all terms of the form $\omega {\tA}_{m\ell}^0 \EE^{\ell}$ by their Taylor expansions 
 around $k=\ell k_0$. (Details of these expansions are contained in Lemma 25 of \cite{SW10}, for example.) After that, we
equate the coefficients of the $\eps^p\EE^{\ell}$ to zero. 

We find that
the coefficients of $\eps \EE^1$ and $\eps^2 \EE^1$ vanish identically
due to the definition of $ \omega$ and 
$ c_g$.
For $\eps^3\EE^1 $ we obtain
$$
\partial_{\tau} {\tA}_{-11}^0 = \frac{1}{2}i\, \partial_k^2 \omega(k_0,b)\, \partial^2_{\underline{\alpha}} {\tA}_{-11}^0 +
 g_1\,,
$$
where $\tau = \eps^2t$, $\underline{\alpha}= \eps
(\alpha-c_gt)$ and $g_1$ is a sum of multiples of ${\tA}_{-11}^0 |{\tA}_{-11}^0|^2 $, 
$ {\tA}_{-11}^0 {\tA}_{m0}^0 $ and $ {\tA}_{-1-1}^0 {\tA}_{m2}^0 $ with $m \in \{\pm 1\}$.
In the next steps we obtain algebraic relations
such that  the $ {\tA}_{m2}^0 $ can be expressed in terms of $({\tA}_{-11}^0)^2$ and 
the $ {\tA}_{m0}^0 $ in terms of $ |{\tA}_{-11}^0|^2 $, respectively.

For $\eps^2 \EE^2$ we obtain
\begin{align*}
( -2 \omega_0 + \omega(2 k_0,b) ) {\tA}_{-12}^0 & =  \gamma_{-12} ({\tA}_{-11}^0)^2\,,
 \\
( -2 \omega_0 -  \omega(2 k_0,b) ) {\tA}_{12}^0 & =  \gamma_{12} ({\tA}_{-11}^0)^2\,
 \nonumber
\end{align*}
with coefficients $ \gamma_{m2} \in \C $.
For all $b \geq 0$ with   
\begin{equation}
\label{2om}
-2 \omega_0 \pm \omega(2 k_0,b) \neq 0 \,,
\end{equation}
the $ {\tA}_{m2}^0 $ are well-defined in terms of $({\tA}_{-11}^0)^2$.
All terms vanish identically for $\epsilon^2 \EE^0.$
This is obvious for the linear terms. For the quadratic  terms 
the calculations  are analogous to those of Appendix A of \cite{SW10} (see
specifically equation (94)).
The nonlinear terms in $\epsilon^3\EE^0$
must be perfect derivatives with respect to $\underline{\alpha}$ since no other combination of terms
in the approximation \eqref{ans1c} leads to terms proportional to 
$ \epsilon^3 \EE^0 $.
So we find 
\begin{align}  \nonumber
-c_g  \partial_{\underline{\alpha}} {\tA}_{-10}^0& =  - \partial_k \omega(0,b) \partial_{\underline{\alpha}} {\tA}_{-10}^0 
+ \gamma_{-10} \partial_{\underline{\alpha}} ({\tA}_{-11}^0 {\tA}_{-1-1}^0),\\ \nonumber
-c_g \partial_{\underline{\alpha}} {\tA}_{10}^0& =  \partial_k \omega(0,b) \partial_{\underline{\alpha}} {\tA}_{10}^0 
+ \gamma_{10} \partial_{\underline{\alpha}} ({\tA}_{-11}^0 {\tA}_{-1-1}^0 ),
\end{align}
where now $ \gamma_{m0} \in \R $ according to the fact that 
we consider a real-valued problem.
For all $b \geq 0$ with
\begin{equation}
c_g \neq \pm \partial_k \omega(0,b)\,,
\end{equation} 
 we can divide the equations for $\epsilon^3 \EE^0 $ by $\partial_{\underline{\alpha}}$
and can express the $ {\tA}_{m0}^0 $ in terms of $  |{\tA}_{-11}^0|^2 $.
 
As mentioned above the nonlinear terms in the equation for $\epsilon^3 \EE^1 $ 
include  ${\tA}_{-11}^0 |{\tA}_{-11}^0|^2$ as well as terms  consisting of combinations of ${\tA}_{-11}^0$ with 
the ${\tA}_{m0}^0$ and of ${\tA}_{-1-1}^0$ with the ${\tA}_{m2}^0$. Eliminating
${\tA}_{m0}^0$ and ${\tA}_{m2}^0$ by the algebraic relations obtained 
for $\epsilon^3 \EE^0 $
and $\epsilon^2 \EE^2 $ gives finally the NLS equation
 \begin{equation}\label{nlsderive}
 \partial_{\tau} {\tA}_{-11}^0 = i \frac{\partial_k ^{2}\omega(k_0,b)}{2} \partial^2_{\underline{\alpha}} {\tA}_{-11}^0 
+ i\nu_2 (k_0,b) {\tA}_{-11}^0 |{\tA}_{-11}^0|^2\,
\end{equation}
with a $\nu_2(k_0,b) \in \R$. 

An explicit formula for $\nu_2$ can be found in \cite[p.~504]{CSS92}. It can be seen with the help of that 
formula if the NLS equation \eqref{nlsderive} is defocusing or focusing for a given basic wave number $k_0>0$. Since we will consider solutions of \eqref{nlsderive} on time intervals $[0,{\tau_0}]$ with $\tau_0 \sim 
{\cal{O}}(1)$, this will not affect our analysis.

The approximation function $\veps\tPsi_2$ is determined by inserting \eqref{ans1c} into \eqref{thetaneu}--\eqref{deltaneu}, using the formulas for $\veps\tPsi_1$
derived above and equating the coefficients of the $\eps^p\EE^{\ell}$ to zero. It turns out that $\veps\tPsi_2$ can be expressed in terms of the components of $\veps \tPsi_1$ and its derivatives. In particular, we have
\begin{equation*}
{\tPsi}^0_{2\pm1} = \partial_{\alpha}^2 {\tPsi}^0_{1\pm1}\,.
\end{equation*}

To prove the approximation property of the NLS equation \eqref{nlsderive} it will be helpful to extend the approximation $\veps \widetilde{\Psi}$ by higher order correction terms in order  
to make the residual of the resulting approximation of the equations \eqref{w-1}--\eqref{deltaneu}
smaller in Sobolev norms. 
The residual of an approximation $a$ of an algebraic or differential equation 
\begin{equation} \label{reseqdef}
A(f) = B(f) \,,
\end{equation}
where $A,B$ are functions depending on the function $f$ and in the case of a differential equation also on derivatives of $f$,
is defined by
\begin{equation*}
{\rm Res}(a) := B(a)-A(a) \,.
\end{equation*}
Hence, ${\rm Res}(a)$ contains all terms that do not cancel after inserting the ansatz $f=a$ in \eqref{reseqdef} and quantifies how much $a$ fails to be a solution of \eqref{reseqdef}. There holds ${\rm Res}(a)=0$ if and only if $a$ is an exact solution of \eqref{reseqdef}.

We introduce the notation
$${\rm Res}(\veps \tPsi) =  \left( \ba{c}
{\rm Res}_{1}(\veps \tPsi)\\ {\rm Res}_{2}(\veps \tPsi)  \\ {\rm Res}_{3}(\veps \tPsi)      \ea \right)  \; \text{with} \;\,\, {\rm Res}_{j}(\veps \tPsi) = \left(
\begin{array}{c}  {\rm res}_{-j}(\veps \tPsi)  \\ {\rm res}_{j}(\veps \tPsi)
\end{array}
\right)
\,\, \text{for} \;\, j=1,2,3,$$ 
where ${\rm res}_{-1}(\veps \tPsi)$ is the residual of $\veps \tPsi$ of equation \eqref{w-1}, ${\rm res}_{1}(\veps \tPsi)$ the residual of $\veps \tPsi$ of equation \eqref{w1}, ${\rm res}_{-2}(\veps \tPsi)$ the residual of $\veps \tPsi$ of equation \eqref{u-1}, ${\rm res}_{2}(\veps \tPsi)$ the residual of $\veps \tPsi$ of equation \eqref{u1}, ${\rm res}_{-3}(\veps \tPsi)$ the residual of $\veps \tPsi$ of equation \eqref{thetaneu} and ${\rm res}_{3}(\veps \tPsi)$ the residual of $\veps \tPsi$ of equation \eqref{deltaneu}.

In order to replace $\veps \tPsi$ in ${\rm Res}(\veps \tPsi)$ by 
a better approximation $\veps \Psi$ we proceed analogously as in Section 2 of \cite{DSW12}. In a first step we construct an extended approximation 
\begin{equation} 
\mathcal{U} = \veps \widetilde{\Psi}^{ext} = \veps \left( \ba{c}
\widetilde{\Psi}^{ext}_{1}\\ \widetilde{\Psi}^{ext}_2 \ea \right)
\end{equation}
with
\begin{align*}
\veps \widetilde{\Psi}^{ext}_j &\, =\,
\eps \widetilde{\Psi}_{j} +
\eps^2 \widetilde{\Psi}^{add}_{j}\,,
\end{align*}
where $\eps \widetilde{\Psi}_{j}$ is as above and $\eps^2 \widetilde{\Psi}^{add}_{j}$ is of the form
\begin{eqnarray*}
\eps^2 \widetilde{\Psi}^{add}_{j} \!\!& = &\!\!
\sum_{\ell \in  \{\pm 1 \}} \, \sum_{n=1}^{4} \left(
\begin{array}{c} \epsilon^{1+n} \widetilde{A}^n_{-j\ell} (\epsilon(\alpha -c_gt),\epsilon^2t)\, \EE^{\ell} \\[1mm]
\epsilon^{1+n} \widetilde{A}^{n}_{j\ell} (\epsilon(\alpha -c_gt),\epsilon^2t)\, \EE^{\ell}
\end{array}\right)\\[1mm]
&&\!\! + \sum_{n=1}^{3}  
\left(\begin{array}{c}\epsilon^{2+n}\widetilde{A}^n_{-j0} (\epsilon(\alpha -c_gt),\epsilon^2t)\\[1mm]
\epsilon^{2+n} \widetilde{A}^n_{j0} (\epsilon(\alpha -c_gt),\epsilon^2t)
\end{array}
\right)\\[1mm]
&&\!\! + \sum_{\ell \in  \{\pm 2 \}} \, \sum_{n=1}^{3}
\left(\begin{array}{c}\epsilon^{2+n}{\tA}^n_{-j\ell} (\epsilon(\alpha -c_gt),\epsilon^2t)\, \EE^{\ell}\\[1mm]
\epsilon^{2+n} \widetilde{A}^n_{j\ell} (\epsilon(\alpha -c_gt),\epsilon^2t)\, \EE^{\ell}
\end{array}\right)\\[1mm]
&&\!\! + \sum_{\ell \in  \{\pm 3 \}} \, \sum_{n=0}^{2}
\left(\begin{array}{c}\epsilon^{3+n} \widetilde{A}^n_{-j\ell} (\epsilon(\alpha -c_gt),\epsilon^2t)\, \EE^{\ell}\\[1mm]
\epsilon^{3+n} \widetilde{A}^n_{j\ell} (\epsilon(\alpha -c_gt),\epsilon^2t)\, \EE^{\ell}
\end{array}\right)\\[1mm]
&&\!\! + \sum_{\ell \in  \{\pm 4 \}} \, \sum_{n=0}^{1}
\left(\begin{array}{c}\epsilon^{4+n} \widetilde{A}^{{n}}_{-j\ell} (\epsilon(\alpha -c_gt),\epsilon^2t)\, \EE^{\ell}\\[1mm]
\epsilon^{4 +{n}}  \widetilde{A}^{{n}}_{j{\ell}} (\epsilon(\alpha -c_gt),\epsilon^2t)\, \EE^{\ell}
\end{array}\right)\\[1mm]
&&\!\!  + \sum_{\ell \in  \{\pm 5 \}} 
\left(\begin{array}{c}\epsilon^5 \widetilde{A}^{0}_{-j\ell} (\epsilon(\alpha -c_gt),\epsilon^2t)\, \EE^{\ell}\\[1mm]
\epsilon^5 \widetilde{A}^{0}_{j\ell} (\epsilon(\alpha -c_gt),\epsilon^2t)\, \EE^{\ell}
\end{array}\right)
\end{eqnarray*}
with 
${\tA}_{\mp j -\ell}^{n} = \overline{{\tA}_{\mp j \ell}^{n}}$. Then we have
\begin{align*}
\veps \widetilde{\Psi}_j^{ext} &
 \,=\,\sum_{|\ell| \le 5}\, \sum_{\beta(\ell,n) \le 5} \varepsilon^{\beta(\ell,n)} \widetilde{\Psi}_{j\ell}^{n}   \;,
\end{align*}
where $j \in \{1,2\}$, $\ell \in \mathbb{Z}$, $n \in \N_0$ and $\beta(\ell,n)=1+ \vert \vert \ell \vert -1 \vert +n$ as well as
\begin{align*}
\eps \widetilde{\Psi}^0_{1\pm 1} 
& =  \eps \widetilde{\psi}_{-1\pm 1}^{0} \!\left( \ba{c} 
1 \\ 0  \ea \right),\\[2mm] 
\eps \widetilde{\Psi}^0_{2\pm 1} 
& =  \eps \partial_{{\alpha}}^2 \widetilde{\psi}_{-1 \pm 1}^{0} \!\left( \ba{c} 
1 \\ 0  \ea \right),
\\[2mm] 
\varepsilon^{\beta(\ell,n)} \widetilde{\Psi}^n_{j\ell} 
&\, =\, \varepsilon^{\beta(\ell,n)} \left(\begin{array}{c} 
\widetilde{\psi}^n_{-j\ell} \\[1mm]
\widetilde{\psi}^n_{j\ell}
 \end{array} \right) \quad \text{for} \; (\ell,n) \neq (\pm 1,0)\,,
\\[2mm]
\widetilde{\psi}_{\mp j \ell}^{n}(\alpha,t) &\,=\, \widetilde{A}_{\mp j \ell}^{n}(\varepsilon(\alpha - c_g t), \varepsilon^2 t)\, \textbf{E}^{\ell}\,.
\end{align*}
The functions
$\widetilde{A}_{\mp j \ell}^{n}$ in $\eps^2 \widetilde{\Psi}^{add}_{j}$ are computed by a similar procedure as $\widetilde{A}_{-1\pm 1}^{0}$, $\widetilde{A}_{\mp j 0}^{0}$ and $\widetilde{A}_{\mp j 2}^{0}$ in $\eps \widetilde{\Psi}_{j}$. More precisely, inserting $\veps \widetilde{\Psi}^{ext}_1$ into \eqref{w-1}--\eqref{w1} and equating the coefficients in front of the $\varepsilon^{\beta(\ell,n)} \textbf{E}^{\ell}$ to zero yields 
a system of algebraic equations and inhomogeneous linear Schr\"odinger equations that can be solved recursively. For all $b \geq 0$ with $-\ell \omega_0 \pm \omega(\ell k_0,b) \neq 0 $ for $\ell \in \{2,3,4,5\}$ and $ c_g \neq \pm \partial_k \omega(0,b)$ the functions $\widetilde{A}_{p \ell}^{n}$ with $|p| =1$ 
and $(p,|\ell|) \neq (-1,1)$
 are uniquely determined by the algebraic equations. The functions $\widetilde{A}_{-1\pm 1}^{n}$  satisfy the inhomogeneous linear Schr\"odinger equations. Moreover, since the functions $\widetilde{A}_{-1 \pm 1}^4$  do not appear in the equations for any other $\widetilde{A}_{-1 \pm 1}^{n}$  , we can set  $\widetilde{A}_{-1 \pm 1}^4=0$. 

The approximation function $\veps \widetilde{\Psi}^{ext}_2$ is determined by inserting $\veps \widetilde{\Psi}^{ext}_1$ into \eqref{thetaneu}--\eqref{deltaneu}, using the formulas for $\veps \widetilde{\Psi}^{ext}_1$
derived above and equating the coefficients of the $\eps^{\beta(\ell,n)} \EE^{\ell}$ to zero. It turns out that $\veps \widetilde{\Psi}^{ext}_2$ can be expressed in terms of the components of $\veps \widetilde{\Psi}^{ext}_1$ and their first three spatial derivatives.


In a second step we apply a Fourier truncation procedure to the approximation $\veps \widetilde{\Psi}^{ext}$. More precisely, we define
\begin{align}
{\psi}_{\mp j \ell}^{n} & := \, \mathcal{F}^{-1}( \chi_{[-\delta_0,\delta_0]}\, \mathcal{F}( \widetilde{\psi}_{\mp j \ell}^{n} \,\textbf{E}^{-\ell}))\, \textbf{E}^{\ell} \,,
\end{align}
where $\chi_{[-\delta_0,\delta_0]}$ is the characteristic function on $[-\delta_0,\delta_0]$ and
$\delta_0 = \delta_0(b) \in (0, k_0/20)$ will be determined in Section \ref{sec4}. Then the functions ${\psi}_{\mp j \ell}^{n}$ 
have the compact support 
\begin{align}
&\{k \in \R : |k- \ell k_0| \le \delta_0 < k_0/20\}
\end{align}
in Fourier space for all $0 < \eps \ll 1$ and they are 
of the form
\begin{align}
\psi_{\mp j \ell}^{n}(\alpha,t) &\,=\, A_{\mp j \ell}^{n}(\varepsilon(\alpha - c_g t), \varepsilon^2 t)\, \textbf{E}^{\ell}\,
\label{pm1-ans}
\end{align}
with 
${A}_{\mp j -\ell}^{n} = \overline{{A}_{\mp j \ell}^{n}}$.
By using these functions, we construct our final approximation 
\begin{equation} 
\mathcal{U} = \veps \Psi = \veps \left( \ba{c}
\Psi_{1}\\ \Psi_2 \ea \right)\,
\end{equation}
with
\begin{align}
\label{ans-high}
\veps \Psi_j &\, =\,
\sum_{|\ell| \le 5}\, \sum_{\beta(\ell,n) \le 5} \varepsilon^{\beta(\ell,n)} \Psi_{j\ell}^{n}   \;,
\end{align}
where $j \in \{1,2\}$, $\ell \in \mathbb{Z}$, $n \in \N_0$ and $\beta(\ell,n)=1+ \vert \vert \ell \vert -1 \vert +n$ as well as
\begin{align}
\eps {\Psi}^0_{1\pm 1} 
& =  \eps \psi_{-1\pm 1}^{0} \!\left( \ba{c} 
1 \\ 0  \ea \right),\\[2mm] 
\eps {\Psi}^0_{2\pm 1} 
& =  \eps \partial_{{\alpha}}^2 \psi_{-1\pm 1}^{0} \!\left( \ba{c} 
1 \\ 0  \ea \right),
\\[2mm] 
\varepsilon^{\beta(\ell,n)} {\Psi}^n_{j\ell} 
&\, =\, \varepsilon^{\beta(\ell,n)} \left(\begin{array}{c} 
{\psi}^n_{-j\ell} \\[1mm]
{\psi}^n_{j\ell}
 \end{array} \right) \quad \text{for} \; (\ell,n) \neq (\pm 1,0)\,.
\end{align}

For later purposes we set
\begin{align}
\label{pm1-ans2}
\psi_{\pm 1} & \,:=\, \psi_{-1 \pm 1}^{0}\,,
\\[2mm]
\eps^2 \Psi^h_{j} 
& \,=\, \eps^2 \left(\begin{array}{c}
\psi^h_{-j}\\[1mm]
\psi^h_{j}
\\ \end{array}\right) := \, \veps \Psi_j -
(\eps {\Psi}^0_{j1} + \eps {\Psi}^0_{j-1}) \,.
\end{align}

Since the Fourier transform of the functions $\widetilde{\psi}_{\mp j \ell}^{n}$ in $\eps \widetilde{\Psi}^{ext}$ are strongly concentrated around the wave numbers $\ell k_0$ if the functions $\widetilde{A}_{\mp j \ell}^{n}$ are sufficiently regular, the approximation $\eps \widetilde{\Psi}^{ext}$ is only changed slightly by the Fourier truncation procedure. This fact is a consequence of the estimate
\begin{equation}  \label{cut} 
 \| (\chi_{[-\delta_0,\delta_0]}-1)\, \eps^{-1} 
\widehat{f} ( \eps^{-1} \cdot) \|_{L^2(m)}  \leq C(\delta_0)\,
\eps^{m+{M}-1/2} \| f \|_{H^{m+{M}}} 
\end{equation}

for all $M,m \geq 0$. The Fourier truncation procedure will give us a simpler control
of the error and makes our final approximation $\eps \Psi$ an analytic function. For related strategies compare \cite{Schn98JDE, SW10, DSW12, GMWZ14}.

As in Section 2 of \cite{DSW12}, 
the following estimates for the modified residual hold.
\begin{lemma} \label{lem2}
Let $ s_A \geq {10} $ and
$\widetilde{A}_{-1 1}^{0} \in C^0([0,\tau_0], H^{s_A}(\R,\C)) $ be a solution
of the NLS equation \eqref{nlsderive}  with 
$$ \sup_{\tau \in [0,\tau_0]} \| 
\widetilde{A}_{-1 1}^{0}\|_{H^{s_A}} \leq C_A . $$ 
Then for all $s \geq 0$ there exist $ C_{Res}, C_{\Psi}, \eps_0>0  $ depending on  $C_A$ such that for all
$ \eps \in (0,\eps_0) $ the approximation $\eps \Psi$ satisfies
\begin{align} \label{RES1}
& \sup_{\tau \in [0,\tau_0/\eps^2]} \|  {\rm Res}(\eps \Psi)
\|_{(H^s)^6}
\leq  C_{Res}\, \eps^{{11/2}}, \\ \label{RES2}
& \sup_{\tau \in [0,\tau_0/\eps^2]} \|\eps \Psi_1 - (\eps {\tPsi_{11}^0} + \eps {\tPsi_{1-1}^0})
\|_{(H^{{s_A}})^2}
\leq  {C_{\Psi}}\, \eps^{3/2}, 
\\
\label{RES3}
& \sup_{\tau \in [0,\tau_0/\eps^2]} (\|\widehat{\Psi}_{j \pm 1}^0 \|_{(L^1({s+1})(\R,\C))^2}
+ \|\widehat{\Psi}_j^h \|_{(L^1({s+1})(\R,\C))^2})
\leq {C_{\Psi}}\,
\end{align}
for $j \in \{1,2\}$.
\end{lemma}

\textbf{Proof.}
The first extended approximation $\veps \widetilde{\Psi}^{ext}$ is constructed in a way that formally we have ${\rm Res} (\veps \widetilde{\Psi}^{ext} ) = \OO(\epsilon^{6})$ and $\eps \widetilde{\Psi}^{ext}_1 - (\eps {\tPsi_{11}^0} + \eps {\tPsi_{1-1}^0}) = \OO(\epsilon^{2})$ on the time interval $[0,\tau_0/\epsilon^2]$ if $\widetilde{A}_{-11}^0$
is a solution of the NLS equation \eqref{nlsderive} for $T \in [0,\tau_0]$. 

It can be shown analogously as in the proof of Theorem 2.5 in \cite{DSW12} that the regularity condition 
 $ \widetilde{A}_{-11}^0 \! \in C^0([0, \tau_0],H^{s_A}(\R,\C)) $ with $s_A \geq 8$ implies
$ \widetilde{A}_{- j \pm 1}^n \! \in C^0([0,\tau_0],H^{s_A-n-2-3(|j|-1)}(\R,\C))$ for $n \in \{1,2,3\}$
and $ \widetilde{A}_{p \ell}^n \in C^0([0,\tau_0],H^{s_A-{{n}}-3(|p|-1)}(\R,\C)) $ if $(p,|\ell|) \neq (-j,1) $, where the
respective Sobolev norms are uniformly bound\-ed by the $H^{s_A}$-norm of  $ \widetilde{A}_{-11}^0$.

Therefore, by taking into account that
$\| f(\epsilon\, \cdot)\|_{L^2} = \epsilon^{-1/2} \| f \|_{L^2}$, we conclude that there exist constants $C_1, C_2 >0$ depending on $C_A$ such that
\begin{align} \label{REStilde1}
& \sup_{ \tau \in [0, \tau_0/\eps^2]} \|  {\rm Res}(\eps \widetilde{\Psi}^{ext})
\|_{(H^{s_A -10})^6}
\leq  C_{1}\, \eps^{{11/2}}, \\ \label{REStilde2}
& \sup_{ \tau \in [0, \tau_0/\eps^2]} \|\eps \widetilde{\Psi}^{ext}_1 - (\eps {\tPsi_{11}^0} + \eps {\tPsi_{1-1}^0})
\|_{(H^{{s_A-8}})^2}
\leq  {C_{2}}\, \eps^{3/2} 
\end{align}
if we have $ \widetilde{A}_{-11}^0 \in C^0([0, \tau_0],H^{s_A}(\R,\C)) $ with $s_A \geq 10$ (because two additional 
spatial derivatives of $\widetilde{A}_{-11}^0$ are needed to bound ${\rm Res}(\veps \widetilde{\Psi}^{ext})$).

Since the Fourier transform of the final approximation $\eps \Psi$ has a compact support whose size depends on $k_0$, there exists a constant $C=C(k_0) >0$ such that $\|\Psi\|_{H^s} \leq C \|\Psi\|_{L^2}$ and 
$\|\widehat{\Psi}\|_{L^1({s})} \leq C \|\widehat{\Psi} \|_{L^1}$ for all $s \geq 0$.
Hence, by using the above $L^2$-estimates for $\veps \widetilde{\Psi}_{ext}$ as well as the estimate \eqref{cut}
for $f=\widetilde{A}_{\mp j \ell}^n$ for each $\widetilde{A}_{\mp j \ell}^n$ with $m=0$, $M=M(\ell,n)$ determined by the maximal Sobolev regularity of the respective $\widetilde{A}_{\mp j \ell}^n$ and $\delta_0$ as above, we obtain \eqref{RES1} and
\begin{equation}  \label{RES4}
\sup_{\tau \in [0,\tau_0/\eps^2]} \|\eps \Psi_1 - (\eps \Psi_{11}^0 + \eps \Psi_{1-1}^0)
\|_{(H^{{s_A}})^2}
 \leq  {C}_3\, \eps^{3/2}
\end{equation}
for a constant $C_3=C_3(C_A) >0$ if we have $ s_A  \geq 10$, which yields $\beta(\ell,n)+ M(\ell,n) \geq 6$.  
By combining \eqref{RES4} and \eqref{cut} for $f=\eps {\widetilde{\Psi}_{11}^0} + \eps {\widetilde{\Psi}_{1-1}^0}$, $m=s_A$, $M=0$ and $\delta_0$ as above, we obtain \eqref{RES2}.

Finally, since $\| \epsilon^{-1} \widehat{f}(\epsilon^{-1}  \cdot)\|_{L^1} = \| \widehat{f} \|_{L^1}$,
estimate \eqref{RES3} follows by construction of ${\Psi}_{j\pm 1}^0$ and ${\Psi}_j^h$.
\qed

\begin{remark} 
{\rm 
Due to the estimate \eqref{RES1} of the residual obtained with the Fourier truncation procedure, the Sobolev index of the error estimate \eqref{errest} in Theorem \ref{mainresult} is equal to the Sobolev index of the solution to the NLS equation. If we did not apply the Fourier truncation procedure, we could only use estimate \eqref{REStilde1} in the proof of the error estimates in Section \ref{sec4}. Then, for a given solution of the NLS equation in $H^{s_A}(\R,\C)$, the resulting error could only be estimated in $H^{s_A-10}(\R,\R)$.}  
\end{remark}

\begin{remark} \label{remneucam2}
{\rm 
The bound \eqref{RES3} will be  used for instance to estimate 
$$ \| \psi f \|_{H^{s}} \leq C \| \psi  \|_{C^{s}_b}\| f 
\|_{H^{s}} \leq  C \| \widehat{\psi}  \|_{L^1(s)(\R,\C)}\| f
\|_{H^{s}} $$
without loss of powers in $ \eps $ as it would be the case with $\| \widehat{\psi}  \|_{L^2(s)(\R,\C)}$.}
\end{remark}

Moreover, by an analogous argumentation as in the proof of Lemma 3.3 in \cite{DSW12} we obtain the fact that $\partial_t \psi_{\pm 1}$ can be approximated by
$-i \omega \psi_{\pm 1}$. More precisely, we obtain 
\begin{lemma}\label{deriv_psi} For all $s \geq 0$
there exists a constant $C_{\psi} > 0$ such that
\begin{equation} \label{dtpsi}
\| \partial_t \widehat{\psi}_{\pm 1} + i \omega \widehat{\psi}_{\pm 1} \|_{L^{1}(s)} \le C_{\psi}\, \eps^2\,.
\end{equation}
\end{lemma}

\section{The error estimates} 
\label{sec4}
In this section, we justify the NLS approximation for system \eqref{w-1}--\eqref{deltaneu}. 

\subsection{The structure of the evolutionary system for the error and the approximation result in the arc length formulation}

Our first step in justifying the NLS approximation is writing the exact solution $\mathcal{U}$ of \eqref{w-1}--\eqref{deltaneu} as the sum of the NLS approximation and the error. To avoid problems arising from the non-trivial resonances at $k=\pm k_0$, we rescale the error
with the help of the weight function
\begin{align*}
\widehat{\vartheta}(k) = \left\{ \begin{array}{ll} 1\,, \qquad & |k| >\delta_{0}\,,\\[2mm] \eps + (1-\eps) |k|/\delta_{0}\,,  \qquad 
& |k| \leq \delta_{0}\,, \end{array} \right.
\end{align*}
where $0 < \eps \ll 1$ and $\delta_0  =\delta_0(b)  \in (0, k_0/20)$ will be defined below. 
That means, we write
\begin{equation} \label{ans1cc}
\mathcal{U} = \veps \Psi + \eps^{\beta} \vartheta \mathcal{R} \,,
\end{equation}
where $\beta = 5/2$ and
\[
\vartheta \mathcal{R} := \left(
\begin{array}{c} \vartheta \mathcal{R}_{1} \\ \vartheta \mathcal{R}_{2}
\end{array}
\right)
\; \text{with} \;\,\, \vartheta \mathcal{R}_j := \left(
\begin{array}{c} \vartheta {R}_{-j} \\ \vartheta {R}_{j}
\end{array}
\right)
\,\, \text{for} \;\, j=1,2\,,
\]
and $\vartheta R_{\mp j}$ is defined by $\widehat{\vartheta R}_{\mp j} = \widehat{\vartheta} \widehat{R}_{\mp j}$. 
 By this choice $\widehat{\vartheta R}_{\mp j}(k) 
$ is small at the wave numbers close to zero reflecting 
the fact that the quadratic terms of the evolutionary system of $\mathcal{U}$ vanish at $k = 0$. Hence, we have
\begin{align*}
\qquad u_{-1} \,&=\,   \eps \psi_c + \eps^2 \psi_{-1}^h + \eps^{5/2} \vartheta R_{-1}\,,  \\[2mm]
\qquad u_{1} \,&=\,  \eps^2 \psi_{1}^h + \eps^{5/2} \vartheta R_{1} \,, \\[2mm]
\qquad u_{-2} \,&=\,   \eps \partial_{\ua}^2 \psi_c + \eps^2 \psi_{-2}^h + \eps^{5/2} \vartheta R_{-2} \,, \\[2mm]
\qquad u_{2} \,&=\,  \eps^2 \psi_{2}^h + \eps^{5/2} \vartheta R_{2}\,,  
\end{align*}
where  $\psi_c = \psi_{-1} + \psi_1$. 
\\
The definition of $\vartheta$ directly implies
\begin{align}
\label{theta-1-eps-1}
& \sup_{k \in \R} |\widehat{\vartheta}^{-1}(k)| = \eps^{-1}\,, \\[2mm]
\label{theta-1-eps-0}
& \sup_{k \in \R}\, |(1-\chi_{[-\delta_{0},\delta_{0}]})(k)\,\widehat{\vartheta}^{-1}(k)| = 1\,,
\end{align} 
where the operator $\vartheta^{-1}$ is defined by its symbol $\widehat{\vartheta^{-1}}(k)= \widehat{\vartheta}^{-1}(k)               =(\widehat{\vartheta}(k))^{-1}$.
Moreover, we have
\begin{align*}
|k\, \widehat{\vartheta}^{-1}(k)|= 
\begin{cases}
|k| & \quad {\rm for }\; |k| > \delta_{0} \,,\\[1mm]
\dfrac{|k|}{\varepsilon +(1- \varepsilon) \frac{|k|}{\delta_{0}}} & \quad{\rm for}\; |k| \le \delta_{0} \,.
\end{cases}
\end{align*}
Since 
\begin{align*}
\frac{|k|}{\varepsilon +(1- \varepsilon) \frac{|k|}{\delta_{0}}} \;=\; 
\frac{1}{\frac{\varepsilon}{|k|} + \frac{(1- \varepsilon)}{\delta_{0}}} \;\le\; 
\frac{1}{\frac{\varepsilon}{\delta_{0}} + \frac{(1- \varepsilon)}{\delta_{0}}}\;=\; \delta_{0}
\end{align*}
for $0 \neq |k| \le \delta_{0}$, we obtain
\begin{align}
\label{theta-1}
\sup_{k \in \R}\, |k\widehat{\vartheta}^{-1}(k)| \,=\, \max\{\delta_{0}, |k|\}\,.
\end{align}
Furthermore, we have
\begin{align}
\label{theta_m-k}
\widehat{\vartheta}^{-1}(k)\, \widehat{\vartheta}(m)\, \chi_c(k-m) \, =\, \mathcal{O}(1) 
\end{align}
for $|k| \to \infty$ uniformly with respect to $m \in \R$, where
$\chi_c$ is the characteristic function on $\mathrm{supp} \,\widehat{\psi}_c$.
Using \eqref{K01}--\eqref{K04}, \eqref{w-1}--\eqref{estc0}, \eqref{K0-id} for $f= \vartheta (R_{-1}-R_1)$ and $g=\psi_c$, \eqref{sigmak1}, \eqref{RES3} and \eqref{theta-1-eps-1}--\eqref{theta_m-k}, we obtain
\begin{eqnarray}
\partial_t R_{\mp 1} &\!\!=\!\!& \mp i \omega R_{\mp 1}  \nonumber \\[2mm]
&& -\eps \vartheta^{-1} \partial_{\ua} (\psi_c \vartheta  R_{\mp 1})  
+ \eps \vartheta^{-1} \partial_{\ua} (K_0\psi_c \vartheta K_0 R_{\pm 1}) \nonumber 
\\[2mm]
&& +\eps \vartheta^{-1} \partial_{\ua} B_{\mp1 -1}  (\psi_c ,\vartheta  R_{-1}) +\eps \vartheta^{-1} \partial_{\ua} B_{\mp1 1}  (\psi_c ,\vartheta  R_{1})\nonumber 
\\[2mm]
&& +\eps^2 \vartheta^{-1} \mathcal{C}_{\mp1 -1 -1 -1}  (\psi_{-1} ,\psi_{-1} ,\vartheta  R_{ -1}) +\eps^2 \vartheta^{-1} \mathcal{C}_{\mp1 -1 1 -1}  (\psi_{-1} ,\psi_{1} ,\vartheta  R_{ -1}) \nonumber 
\\[2mm]
&& +\eps^2 \vartheta^{-1} \mathcal{C}_{\mp1 1 1 -1}  (\psi_{1} ,\psi_{1} ,\vartheta  R_{ -1})
+\eps^2 \vartheta^{-1} \mathcal{C}_{\mp1 -1 -1 1}  (\psi_{-1} ,\psi_{-1} ,\vartheta  R_{ 1}) 
\nonumber 
\\[2mm]
&& +\eps^2 \vartheta^{-1} \mathcal{C}_{\mp1 -1 1 1}  (\psi_{-1} ,\psi_{1} ,\vartheta  R_{ 1}) +\eps^2 \vartheta^{-1} \mathcal{C}_{\mp1 1 1 1}  (\psi_{1} ,\psi_{1} ,\vartheta  R_{ 1})
\nonumber 
\\[2mm] \label{Rmp1}
&& +\eps^2 \mathcal{M}_{\mp1}(\Psi, \mathcal{R}) 
+ \eps^{-5/2} \vartheta^{-1} \mathrm{res}_{\mp1}(\eps \Psi)\,,
\end{eqnarray}
\begin{eqnarray}
\partial_t R_{\mp 2} &\!\!=\!\!& \mp i \omega R_{\mp 2}  \nonumber \\[2mm]
&& -\eps \vartheta^{-1} \partial_{\ua} (\psi_c\, \vartheta  R_{ \mp 2})  
\nonumber \\[2mm]
&&+ \frac{\eps}{2}  \vartheta^{-1} \partial_{\ua}((K_0 \sigma^{-1} \partial_{\ua} \psi_c) \, \sigma^{-1} \vartheta (R_{-2}-R_{2})) \nonumber 
\\[2mm]
&&-  \frac{\veps}{2} b \vartheta^{-1} \partial_{\ua} ((\sigma^{-1} \partial_{\ua}^2 \psi_c)\, K_0 \sigma^{-1} \partial_{\ua} \vartheta (R_{-2}-R_{2})) \nonumber 
\\[2mm]
&&\mp \frac{\veps}{2} \vartheta^{-1} \partial_{\ua} ([\sigma, \partial_{\ua}^{-2}\vartheta (R_{-2}+R_{2})]\sigma^{-1} \partial_{\ua}^2 \psi_c) \nonumber 
\\[2mm]
&& +\eps \vartheta^{-1} \partial_{\ua} B_{\mp2 -2}  (\psi_c ,\vartheta  R_{-2}) +\eps \vartheta^{-1} \partial_{\ua} B_{\mp2 2}  (\psi_c ,\vartheta  R_{2})\nonumber 
\\[3mm]
&& -\eps^2 \vartheta^{-1} \partial_{\ua}(\partial_{\ua}^{-2} g_{+}(\Psi_{2}^{h},  \mathcal{R}_{2})\, \vartheta R_{\mp2})
\nonumber 
\\[2mm]
&& + \frac{\eps^2}{2} \vartheta^{-1} \partial_{\ua}((K_0 \sigma^{-1} \partial_{\ua}^{-1} g_{-}(\Psi_{2}^{h},\mathcal{R}_{2})+c(\Psi, \mathcal{R}))\, \sigma^{-1}\vartheta (R_{-2}-R_{2}))
\nonumber 
\\[2mm]
&& -\frac{\eps^2}{2} b \vartheta^{-1} \partial_{\ua}( \sigma^{-1} g_{-}(\Psi_{2}^{h}, \mathcal{R}_{2})\, K_0 \sigma^{-1} \partial_{\ua} \vartheta (R_{-2}-R_{2}))
\nonumber 
\\[2mm]
&& \mp \frac{\eps^2}{2} \vartheta^{-1} \partial_{\ua}( [\sigma, \partial_{\ua}^{-2}\vartheta (R_{-2}+R_{2})] \sigma^{-1} g_{-}(\Psi_{2}^{h}, \mathcal{R}_{2}))
\nonumber 
\\[2mm] \label{Rmp2}
&& +\eps^2 \mathcal{M}_{\mp2}(\Psi, \mathcal{R})  + \eps^{-5/2} \vartheta^{-1} \mathrm{res}_{\mp2}(\eps \Psi) 
\end{eqnarray}
as well as
\begin{eqnarray}
\partial_{\ua}^{-1} \sigma^{-1}( R_{-2} - R_{2}) &\!\!=\!\!& \sigma^{-1}( R_{-1} - R_{1})_{\ua}
\nonumber
\\[2mm]
&&  \label{R2R1m}
+ \eps \mathcal{M}_{-3}(\Psi, \mathcal{R})  + \eps^{-5/2} \vartheta^{-1} \mathrm{res}_{-3}(\eps \Psi)\,, 
\\[3mm]
\partial_{\ua}^{-2} ( R_{-2} + R_{2}) &\!\!=\!\!&  R_{-1} + R_{1} 
\nonumber
\\[2mm]
&& 
- \eps \vartheta^{-1} \partial_{\ua}^{-1} ((\sigma^{-1} \partial_{\ua}^{2} \psi_c)\,            
K_0 \vartheta ( R_{-1} + R_{1}))
\nonumber 
\\[2mm]
&& 
- \eps \vartheta^{-1} \partial_{\ua}^{-1}  ((K_0 \psi_c) \sigma^{-1} \vartheta ( R_{-2} - R_{2}))
\nonumber 
\\[2mm]
&& \label{R2R1p}
+ \eps \mathcal{M}_{3}(\Psi, \mathcal{R})  + \eps^{-5/2} \vartheta^{-1} \mathrm{res}_{3}(\eps \Psi)\,, 
\end{eqnarray}
where 
\begin{equation}
g_{\pm}(\Psi_{2}^{h}, \mathcal{R}_{2}) = \psi_{-2}^{h} \pm \psi_{2}^{h}
+ \eps^{1/2} \vartheta (R_{-2} \pm R_{2})\,,
\end{equation}
$B_{j_1 j_2}$ with $j_1 \in \{\pm 1, \pm 2\}$ and $j_2 \in \{\pm j_1\}$ are bilinear real-valued mappings, $\mathcal{C}_{j_1 m n j_2}$ with $j_1,m,n,j_2 \in \{\pm 1\}$  trilinear real-valued mappings as well as $\mathcal{M}_{j}$ with $j \in  \{\pm 1, \pm 2, \pm 3 \}$ and $c$ nonlinear real-valued functions which satisfy   
\begin{equation}
 \| B_{\mp1 j_2}  (\psi_c ,\vartheta  R_{j_2})\|_{H^1}   \lesssim  \| R_{j_2}\|_{L^2}\,, 
\end{equation}
as long as $\eps^{5/2} \| R_{j_2}\|_{L^2}  \lesssim 1$,
\begin{equation}
\|\partial_{\ua} B_{\mp2 j_2}  (\psi_c ,\vartheta  R_{j_2})\|_{H^{s}} \lesssim \| R_{-1}\|_{L^2}
+ \| R_{1}\|_{L^2} + \| R_{-2}\|_{H^{s}} + \| R_{2}\|_{H^{s}} 
\end{equation}
for $s \geq 2$, as long as $\eps^{5/2} \| R_{\mp 1}\|_{L^2},\,  \eps^{5/2} \| R_{\mp 2}\|_{H^{s}} \lesssim 1$,
\begin{equation}
\label{c1}
 \| \mathcal{C}_{j_1 m n j_2}  (\psi_m ,\psi_n ,\vartheta  R_{ j_2})\|_{H^1}  
 \lesssim  \| R_{j_2}\|_{L^2} + \| R_{2j_2}\|_{H^{3/2}}\,,
\end{equation}
as long as $\eps^{5/2} \| R_{j_2}\|_{L^2},\, \eps^{5/2} \| R_{2j_2}\|_{H^{3/2}} \lesssim 1$,
\begin{equation}
\label{calM1}
 \| \mathcal{M}_{\mp1}(\Psi, \mathcal{R}) \|_{H^1}  \lesssim  \| R_{-1}\|_{L^2}  + \| R_{1}\|_{L^2}  + \| R_{-2 }\|_{H^{3/2}} + \| R_{2}\|_{H^{3/2}}  \,, 
\end{equation}
as long as $\eps^{5/2} \| R_{\mp1}\|_{L^2} ,\, \eps^{5/2} \| R_{\mp 2}\|_{H^{3/2}} \lesssim 1$,
\begin{equation}
\label{estcalMmp2}
\| \mathcal{M}_{\mp2}(\Psi, \mathcal{R}) \|_{H^s} \lesssim \| R_{-1}\|_{L^2}  + \| R_{1}\|_{L^2}  + \| R_{-2}\|_{H^s} + \| R_{2}\|_{H^s}  
\end{equation}
for $s \geq 2$, as long as $\eps^{5/2} \| R_{\mp 1}\|_{L^2},\, \eps^{5/2} \| R_{\mp 2}\|_{H^{s}} \lesssim 1$,
\begin{equation} \label{estcalMm3}
\| \mathcal{M}_{-3}(\Psi, \mathcal{R}) \|_{H^s} \lesssim \, \|{R_{-1}}\|_{L^{2}} + \|R_{1}\|_{L^{2}} +\|\sigma^{-1}({R}_{-2}-R_2)\|_{H^{s-1}} 
\end{equation}
for $s\geq 2$, as long as $\eps^{5/2} \|{R}_{\pm1}\|_{L^{2}},\,  \eps^{5/2} \|\sigma^{-1}({R}_{-2}-R_2) \|_{H^{s-1}} \ll 1$, 
\begin{align} \label{estcalMp3}
\|\vartheta \mathcal{M}_{3}(\Psi, \mathcal{R}) \|_{C^0} + \|\partial_{\ua} \mathcal{M}_{3}(\Psi, \mathcal{R})   \|_{H^{s}} \lesssim &\; \|{R_{-1}}\|_{L^{2}} + \|R_{1}\|_{L^{2}} 
 \nonumber\\[2mm]  &\; 
+ \|\sigma^{-1}({R}_{-2}-R_2) \|_{H^{{s}}} + \|{R}_{-2}+{R}_{2}\|_{H^{{s-2}}}
\end{align}
for $s\geq 2$, as long as $\eps^{5/2} \|{R}_{\pm1}\|_{L^{2}},\, \eps^{5/2} \|\sigma^{-1}({R}_{-2}-R_2)\|_{H^{s-1}},\, \eps^{5/2}\|{R}_{-2}+R_2\|_{H^{s-2}}\ll 1$ and
\begin{equation}
\| c(\Psi, \mathcal{R})\|_{H^s} \leq C( \| R_{-1}\|_{L^2}, \| R_{1}\|_{L^2}, \| R_{-2}\|_{H^{s-1}}, \| R_{2}\|_{H^{s-1}})  
\end{equation}
for $s \geq 6$. All bounds are uniform with respect to $b \lesssim 1$ and $\eps \ll 1$.

Moreover, \eqref{R2R1m}--\eqref{R2R1p}, \eqref{estcalMm3}--\eqref{estcalMp3} and Lemma \ref{lem2} imply  
\begin{align} \label{dxR1R2}
\| R_{-1}\|_{H^s} + \|R_{1} \|_{H^s}  \lesssim  &\; 
\| R_{-1}\|_{L^2}  + \| R_{1}\|_{L^2}  + \| R_{-2}\|_{H^{s-2}} + \| R_{2}\|_{H^{s-2}}  \nonumber\\[2mm]  &\;   + \eps \|\partial_{\alpha}^{s-1} R_{-2}\|_{L^{2}} + \eps \|\partial_{\alpha}^{s-1} R_{2}\|_{L^{2}}+ \eps^{2} 
\end{align}
for $s \geq 3$, as long as $\eps^{5/2} \| R_{\mp 1}\|_{L^2},\, \eps^{5/2} \| R_{\mp 2}\|_{H^{s-1}} \ll 1$,
\begin{equation} \label{dx-1R2}
\|\partial_{\ua}^{-1} (R_{-2} - R_2) \|_{H^{s}} \lesssim \| R_{-1}\|_{L^2}  + \| R_{1}\|_{L^2}  + \| R_{-2}\|_{H^{s-1}} + \| R_{ 2}\|_{H^{s-1}} + \veps^{2}
\end{equation}
for $s \geq 3$,
as long as $\eps^{5/2} \| R_{\mp 1}\|_{L^2}, \,\eps^{5/2} \| R_{\mp 2}\|_{H^{s-1}} \ll 1$, and
\begin{align} \label{dx-2R2}
\|\partial_{\ua}^{-2} \vartheta (R_{-2} + R_2) \|_{C^0} +  \|\partial_{\ua}^{-1} (R_{-2} + R_2) \|_{H^{s}}  \lesssim &\; \| R_{-1}\|_{L^2}  + \| R_{1}\|_{L^2} \nonumber \\[2mm] & \;+ \| R_{-2}\|_{H^{s-1}} + \| R_{2}\|_{H^{s-1}} + \veps^{2}
\end{align}
for $s \geq 3$,
as long as $\eps^{5/2} \| R_{\mp 1}\|_{L^2},\,  \eps^{5/2} \| R_{\mp 2}\|_{H^{s-1}} \ll 1$.
These bounds are also uniform with respect to $b \lesssim 1$ and $\eps \ll 1$.

Local existence and uniqueness of solutions $\mathcal{R}$ to \eqref{Rmp1}--\eqref{R2R1p} in $(L^{2}(\R,\R))^2 \times (H^{s}(\R,\R))^2$ with $s+2 = s_A \geq 10$  follows directly from the local existence and uniqueness results in Sobolev spaces for the arc length formulation of the two-dimensional water wave problem \eqref{ytilde}--\eqref{ctilde} and the NLS equation.

Now, we discuss the structure of the above 
evolution equations for the error $\mathcal{R}$. These equations are of the form 
\begin{equation}
\label{Rmps}
\partial_t \mathcal{R}_j \,=\, \mathcal{L} \mathcal{R}_j + \eps Q_j(\psi_c) \mathcal{R}_j + \eps^2 W_j(\Psi, \mathcal{R}) + \eps^{-\beta} \vartheta^{-1} \mathrm{Res}_{j}(\eps \Psi)
\end{equation}
for $j \in \{1,2\}$, with linear operators $\mathcal{L}$ and $Q_j(\psi_c)$ and nonlinear functions $W_j$ having the following properties. $\mathcal{L}$ can be represented by the diagonal matrix
\begin{equation}
\mathcal{L} \,=\, \mathrm{diag}(-i\omega, i\omega)\,. \\[2mm]
\end{equation}
The operators $Q_j(\psi_c)$ are of the form
\begin{align}
Q_1(\psi_c) \mathcal{R}_1 \;=& \; \left(
\begin{array}{cc} Q_{-1-1}(\psi_c) & Q_{-11}(\psi_c)  \\ Q_{1-1}(\psi_c) & Q_{11}(\psi_c)
\end{array}
\right) \left(
\begin{array}{c} R_{-1} \\ R_{1}
\end{array}
\right) \nonumber \\[2mm]
& \; + \left(
\begin{array}{cc} \mathcal{C}_{-1-1}(\psi_c, \psi_c) & \mathcal{C}_{-11}(\psi_c, \psi_c)  \\ \mathcal{C}_{1-1}(\psi_c, \psi_c) & \mathcal{C}_{11}(\psi_c, \psi_c)
\end{array}
\right) \left(
\begin{array}{c} R_{-1} \\ R_{1}
\end{array}
\right)
\end{align}
and
\begin{align}
Q_2(\psi_c) \mathcal{R}_2 \;=& \; \left(
\begin{array}{cc} Q_{-2-2}(\psi_c) & Q_{-22}(\psi_c)  \\ Q_{2-2}(\psi_c) & Q_{22}(\psi_c)
\end{array}
\right) \left(
\begin{array}{c} R_{-2} \\ R_{2}
\end{array}
\right)
\end{align}
respectively,
with
\begin{align} \label{q}
& Q_{j_1j_2}(g) f \;=\; \sum\limits_{\mu=1}^{2|j_1|+1} Q^{\mu}_{j_1j_2}(g) f\,,\\[2mm]
& (\widehat{Q}^{\mu}_{j_1 j_2}(g)f)(k) \;=\;  \int_{\R} \widehat{\vartheta}^{-1}(k)\, \widehat{q}^{|j_1|,\mu}_{j_1 j_2}(k,k-m,m)\, \widehat{g}(k-m)\, \widehat{\vartheta}(m) \widehat{f}(m)\,dm\,, \\[4mm]
\label{q11}
& \widehat{q}_{j_1 j_2}^{1,1}(k,k-m,m) \; = \;  -\delta_{j_1j_2}\, ik\,,\\[4mm] 
\label{q12}
& \widehat{q}_{j_1 j_2}^{1,2}(k,k-m,m)\; = \;  \,\delta_{j_1-j_2}\, ik \,\widehat{K}_0(k-m)  
\widehat{K}_0(m)\,, \\[4mm]
\label{q13}
& \widehat{q}_{j_1 j_2}^{1,3}(k,l,m)\, \chi_c(l) \, \chi_c(k-m) \; = \;  \left\{ \begin{array}{ll} \mathcal{O}(|k|) \quad & \text{for}\,\, |k| \to 0\,,\\[2mm]  \mathcal{O}(1) \quad 
& \text{for}\,\, |k| \to \infty\,, \end{array} \right.
\\[4mm]
\label{q1}
& \widehat{q}_{j_1 j_2}^{2,1}(k,k-m,m) \; = \;  -\delta_{j_1j_2}\, ik\,,\\[4mm] 
\label{q2}
& \widehat{q}_{j_1 j_2}^{2,2}(k,k-m,m)\; = \;  -\frac12\, \mathrm{sgn}(j_2)\, ik \,\widehat{K}_0(k-m)  \sigma^{-1}(k-m) i(k-m)\, \sigma^{-1}(m)  \,, 
\\[4mm]
\label{q3}
& \widehat{q}_{j_1 j_2}^{2,3}(k,k-m,m)\; = \; 
 -\frac{b}{2}\, \mathrm{sgn}(j_2)\, ik\, \sigma^{-1}(k-m) (k-m)^2\, \widehat{K}_0(m)\sigma^{-1}(m) im\,, 
\\[4mm]
\label{q4}
& \widehat{q}_{j_1 j_2}^{2,4}(k,k-m,m)\; = \; 
 \,\frac12\,  \mathrm{sgn}(j_1)\, ik\, i\, \frac{\sigma(k)-\sigma(k-m)}{k-(k-m)} \, \sigma^{-1}(k-m) (k-m)^2\, (im)^{-1}\,, \\[4mm]
\label{q5}
& \widehat{q}_{j_1 j_2}^{2,5}(k,l,m)\, \chi_c(l) \, \chi_c(k-m) \; = \;  \left\{ \begin{array}{ll} \mathcal{O}(|k|) \quad & \text{for}\,\, |k| \to 0\,,\\[2mm]  \mathcal{O}(1) \quad 
& \text{for}\,\, |k| \to \infty\,, \end{array} \right.
\\[4mm]
\label{frakm}
& \partial_n\, \widehat{q}_{j_1 j_2}^{2,5}(k,l,m)\, \chi_c(l) \, \chi_c(k-m)\; = \;   \mathcal{O}(|k|^{-1}) \quad \text{for}\,\, |k|  \to \infty   \;\, \text{and}\;\,   n \in \{1,3\}\,,
\end{align}
where the bounds are uniform with respect to $m \in \R$ and $b \lesssim 1$,
and
\begin{equation} 
\label{c2}
\mathcal{C}_{j_1 j_2}(g,h)f \;=\; \sum_{m,n\in \{\mp1\}} \eps {\vartheta}^{-1} \mathcal{C}_{j_1 m n j_2}(g_{m},h_{n}, \vartheta f)
\,.
\end{equation}
Here, $\partial_n$ denotes the partial derivative with respect to the $n$th variable and the functions $p_{\ell}$ with $p \in \{g,h\}$ and $\ell \in \{m,n\}$, are defined by $\widehat{p}_{\ell} = \widehat{p}\, \chi_{\R^+_0}$ if $\ell=1$ and $\widehat{p}_{\ell} = \widehat{p}\, \chi_{\R^-_0}$ if $\ell=-1$.
By using \eqref{w-1}--\eqref{deltaneu} and the Taylor expansion of $\sigma$ as function of $m$ around $m=k$, the symbols $\widehat{q}_{j_1 j_2}^{1,3},  \widehat{q}_{j_1 j_2}^{2,5}: \R^{3} \to i\R$ can be computed explicitly. But for simplicity we only present those properties of these symbols that we need for the proof of the error estimates.       

For later purposes we set
\begin{align}
& \widehat{q}_{j_1 j_2}(k,k-m,m) \,:=\, \sum_{\mu=1}^{2|j_1|+1} \widehat{q}_{j_1 j_2}^{|j_1|,\mu} (k,k-m,m)  \,.
\end{align} 

The symbols $\widehat{q}_{j_1 j_2}^{2,\mu}$ with $\mu \in \{1,2,3,4\}$ have the following symmetry properties, which will be essential for the proof of our error estimates. There holds
\begin{align} \label{oddq}
\widehat{q}_{j_1j_2}^{2,\mu} (-k,k-m,-m) \,&=\, - \widehat{q}_{j_1j_2}^{2,\mu} (k,k-m,m) \,,\\[2mm] \label{odd-jq}
 \widehat{q}_{-jj}^{2,\mu} (k,k-m,m) \,&=\, - \widehat{q}_{j-j}^{2,\mu} (k,k-m,m) 
\end{align}
for all $\,j,j_1,j_2 \in \{\pm 2\}$, $\mu \in \{1,2,3,4\}$, $k \in \R$ and $m \in \R \setminus \{0\}$.

The functions $W_j(\Psi, \mathcal{R})$ are of the form 
\begin{equation}
 W_1(\Psi, \mathcal{R}) \,=   \left(
\begin{array}{c} \mathcal{M}_{-1} (\Psi, \mathcal{R}) \\ \mathcal{M}_{1} (\Psi, \mathcal{R})
\end{array}
\right) 
\end{equation}
and
\begin{equation}
 W_2(\Psi, \mathcal{R}) \,= \left(
\begin{array}{cc} W_{-2-2}(\Psi, \mathcal{R}) & W_{-22}(\Psi, \mathcal{R})  \\ W_{2-2}(\Psi, \mathcal{R}) & W_{22}(\Psi, \mathcal{R})
\end{array}
\right) \left(
\begin{array}{c} R_{-2} \\ R_{2}
\end{array}
\right) +  \left(
\begin{array}{c} \mathcal{M}_{-2}(\Psi, \mathcal{R}) \\ \mathcal{M}_{2}(\Psi, \mathcal{R})
\end{array}
\right)   
\end{equation}
with
\begin{equation} 
W_{j_1j_2}(\Psi, \mathcal{R}) R_{j_2}\, =\, \sum\limits_{\mu=1}^4 W^{\mu}_{j_1j_2}(\Psi, \mathcal{R}) R_{j_2}\,,\\
\end{equation}
 \begin{align}
W^{1}_{j_1j_2}(\Psi, \mathcal{R}) R_{j_2}
\; = &\;\, 
Q^{1}_{j_1j_2}(\partial_{\ua}^{-2} g_{+}(\Psi_{2}^{h}, \mathcal{R}_{2})        ) R_{j_2}
\,,\\[2mm] 
W^{2}_{j_1j_2}(\Psi, \mathcal{R}) R_{j_2}
\; = &\;\,  
Q^{2}_{j_1j_2}(\partial_{\ua}^{-2} g_{-}(\Psi_{2}^{h}, \mathcal{R}_{2}) + (K_0 \sigma^{-1} \partial_{\ua})^{-1} c(\Psi,\mathcal{R})) R_{j_2}
\,, \\[2mm]
W^{3}_{j_1j_2}(\Psi, \mathcal{R}) R_{j_2}
\; = &\;\,  
Q^{3}_{j_1j_2}(\partial_{\ua}^{-2} g_{-}(\Psi_{2}^{h}, \mathcal{R}_{2})) R_{j_2}
\,, \\[2mm]
W^{4}_{j_1j_2}(\Psi, \mathcal{R}) R_{j_2} \; = &\;\,  Q^{4}_{j_1j_2} (\partial_{\ua}^{-2} g_{-}(\Psi_{2}^{h}, \mathcal{R}_{2})) R_{j_2}\,.
\end{align}

The size of the Fourier transform of the terms in the above evolutionary system depends on whether $k$ is close 
to zero or not. To separate the behavior in these two
regions more clearly, we define projection operators $P_{0,\alpha}$ and $P_{\alpha,\infty}$ for $\alpha >0$ by the Fourier
multipliers 
\begin{align}
&\widehat{P}_{0,\alpha}(k) = \chi_{[-\alpha, \alpha]}(k)\,,\\[2mm]
&\widehat{P}_{\alpha,\infty}(k)= (1- \chi_{[-\alpha, \alpha]})(k)\,,
\end{align}
where $\chi_{[-\alpha, \alpha]}$ is the characteristic
function on $[-\alpha, \alpha]$.

In order to control the evolution of the error we will use a suitable energy. For the construction of this energy we have to take into account the resonances generated by $\omega$. For any $b \geq 0$  and $j_1,j_2 \in \{\pm 1\}$ let the functions $\widehat{r}_{j_1j_2}$  be defined by 
\begin{align} \label{nkdef0}
& \widehat{r}_{j_1 j_2}(k,l,m) \;=\; \widehat{r}_{j_1 j_2}(k,l,m,b) \;=\;
i\, (\, j_1 \omega(k,b) + \omega(l,b) - j_2 \omega(m,b)) 
\end{align}
for all $k,l,m \in \R$. We analyze the zeros of $(k,m) \mapsto 
\widehat{r}_{j_1 j_2}(k,k-m,m,b)$. We will see later that because of $\mathrm{supp}\, \widehat{\psi}_{\ell} = [\ell k_0-\delta_{0}, \ell k_0+\delta_{0}]$ and $
\eqref{theta-1-eps-0}$ we can additionally prescribe $|k-m \mp k_0|  \leq {\delta}_0$ and if $|k| > \delta_{0}$, we confine ourselves to considering the zeros of 
$k \mapsto \widehat{r}_{j_1 j_2}(k,\pm k_0,k \mp k_0,b)$.  

By the mean value theorem we have
\begin{align}
\label{rmean02c}
& \widehat{r}_{j_1 j_2}(k,k-m,m,b) 
\nonumber \\[1mm]
& \qquad  
\,=\, i \big( j_1 \partial_k \omega(\theta_0(k,b)k,b) + \partial_k \omega(k-m-\theta_1(k,m,b)k,b) \big) k - i (1+j_2)\omega(m,b) \,,
\end{align}
with $\theta_0(k,b),\theta_1(k,m,b) \in [0,1]$, for all $k,m \in \R$.
Since $k \mapsto \partial_k \omega(k,b)$ is a continuous even function which satisfies \eqref{nrb1} for all $b \in \R^{+}_0 \setminus \{{b}_{\ast}\}$, where ${b}_{\ast} \in (0,1/3)$,
there exist 
a function $\widetilde{\delta}_0 \in C^{0} (\R^{+}_0 \setminus \{ {b}_{\ast} \} , (0,k_0/20))$ such that $(k,m) \mapsto 
\widehat{r}_{j_1 j_2}(k,k-m,m,b)$ has zeros
satisfying $|k| \leq \widetilde{\delta}_0(b)$ and $|k-m \mp k_0|  \leq \widetilde{\delta}_0(b)$ 
if and only if $j_2 = -1$ and then the zeros 
are $(0,\mp k_0)$. Moreover, there exist a
function $\gamma_0 \in C^{0} (\R^{+}_0 \setminus \{ {b}_{\ast} \}, \R^{+})$ such that
\begin{align} \label{P0resbound}
& |\widehat{r}_{\pm1 -1}(k,k-m,m,b)| \geq \gamma_0(b)|k| 
\end{align}
for all $(k,m) \in \R^{2}$ with $|k| \leq \widetilde{\delta}_0(b)$ and $|k-m \mp k_0|  \leq \widetilde{\delta}_0(b)$.

Next, we analyze the zeros of $k \mapsto \widehat{r}_{j_1 j_2}(k,\pm k_0,k \mp k_0,b)$. 
We have
\begin{align} \label{0res}
& \widehat{r}_{\pm1 -1}(0,\pm k_0,\mp k_0,b) = 0
\end{align}
for all $b \geq 0$.
By the mean value theorem we obtain
\begin{align}
\label{rmean01}
\widehat{r}_{\pm1 -1}(k,\pm k_0,k \mp k_0,b) &\,=\, i \big(j_1 \partial_k\omega(\theta_0(k,b)k,b) + \partial_k \omega(\pm k_0 -\theta_1(k,b)k,b)\big) k \,
\end{align}
with $\theta_0(k,b),\theta_1(k,b) \in [0,1]$, for all $k \in \R$. 
Hence,  there exist functions $\gamma_1 \in C^{0} (\R^{+}_0 \setminus \{ {b}_{\ast} \}, \R^{+})$ and $\widetilde{\delta}_1 \in C^{0} (\R^{+}_0 \setminus \{ {b}_{\ast} \} , (0,k_0/20))$
 such that 
\begin{equation} \label{0bound}
|\widehat{r}_{\pm 1 -1}(k,\pm k_0,k\mp k_0,b)| \geq \gamma_1(b)|k|
\end{equation} 
for $|k| \leq \widetilde{\delta}_1(b)$. Moreover, we have
\begin{align} \label{k0res}
& \widehat{r}_{-1 \pm1}(\pm k_0, \pm k_0, 0,b) = 0
\end{align}
for all $b \geq 0$. Using the mean value theorem again we obtain
\begin{align}
\label{rmean02}
& \widehat{r}_{-1 \pm1}(k,\pm k_0,k\mp k_0,b) 
\nonumber \\[1mm]
& \qquad  =\, -i \big(  \partial_k \omega(\pm k_0 +\theta_0(k,b)(k \mp k_0),b) + j_2 \partial_k \omega(\theta_1(k,b)(k \mp k_0),b) \big) (k \mp k_0)\,
\end{align}
with $\theta_0(k,b),\theta_1(k,b) \in [0,1]$, for all $k \in \R$. 
Hence,  there exist functions $\gamma_2 \in C^{0} (\R^{+}_0 \setminus \{ {b}_{\ast} \}, \R^{+})$ and $\widetilde{\delta}_2 \in C^{0} (\R^{+}_0 \setminus \{ {b}_{\ast} \} , (0,k_0/20))$
 such that 
\begin{equation} \label{k0bound}
|\widehat{r}_{-1 \pm1}(k,\pm k_0,k \mp k_0,b)| \geq \gamma_2(b) |k \mp k_0|
\end{equation} 
for $|k \mp k_0| \leq \widetilde{\delta}_2(b)$.

Since $k \mapsto \omega(k)$ is strictly monotonically increasing, $k \mapsto \widehat{r}_{j_1 j_2}(k,\pm k_0,k \mp k_0,b)$ has no other zeros if $j_1 = 1$ or $j_2=1$.  As discussed in the introduction the remaining zeros of $k \mapsto \widehat{r}_{-1-1}(k,\pm k_0,k \mp k_0,b)$
can be determined by analyzing the zeros of the function $\widehat{r}: \R \times \R^+_0 \to \R$
with
\begin{align}
\widehat{r}(k,b) = \omega(k,b) - \omega(k-k_0,b) - \omega(k_0,b) \,.
\end{align}
Because $k \mapsto \omega(k,b)$ is odd there holds  
\begin{align} 
\widehat{r}(k_0/2+k,b) = \widehat{r}(k_0/2-k,b) 
\end{align}
for all $k\in \R$ and all $b \in \R^+_0$ such that it is sufficient to analyze $\widehat{r}$ for $k \in  [k_0/2, \infty)$.
In the following, we present a quantitative description of the behavior of $\widehat{r}$ that is illustrated in Figure \ref{fig2}.

We have
\begin{align} \label{dkr}
& \partial_k\widehat{r}(k,b) = \partial_k \omega(k,b) - \partial_k \omega(k-k_0,b) =  \partial_k \omega(k,b) - \partial_k \omega(k_0-k,b) 
\end{align}
for all $k\in  [k_0/2, \infty)$ and all $b \in \R^+_0$. 

Using
\begin{align}
\label{astanh}
\mathrm{tanh}(k) &\,=\, \mathrm{sgn}(k) + \mathcal{O}(e^{-2|k|})\,,\\[2mm]
\frac{d}{dk}\, \mathrm{tanh}(k) &\,=\, \mathcal{O}(e^{-2|k|}) \label{asdktanh}
\end{align}
for $|k| \rightarrow \infty$, we deduce 
\begin{align} \label{asomega}
\omega(k,b) &\,=\, \mathrm{sgn}(k)\, |k|^{1/2} (1+bk^{2})^{1/2} + \mathcal{O}(e^{-|k|}) \,,\\[2mm] 
\label{asdkomega}
\partial_k \omega(k,b) &\,=\, \frac{1+3bk^{2}}{2 (1+bk^{2})} \, |k|^{-1/2} (1+bk^{2})^{1/2}  + \mathcal{O}(e^{-|k|})\,,
\\[2mm]
\label{asdk2omega}
\partial^{2}_k \omega(k,b) &\,=\, \mathcal{O}(|k|^{-3/2} (1+bk^{2})^{1/2}) 
\end{align}
for $|k| \rightarrow \infty$ uniformly with respect to $b \lesssim 1$. By Taylor's theorem we have
\begin{align} \label{tayloromega}
\omega(k,b)-\omega(k-k_0,b) = \partial_k \omega(k,b) k_0 - \frac12 \partial^{2}_k \omega(k+ \theta(k)(k-k_0),b) k_0^{2}\,
\end{align}
with $\theta(k) \in [0,1]$. Due to \eqref{asomega}--\eqref{tayloromega},
we obtain 
\begin{equation} \label{rinfty0}
\widehat{r}(k,0) \to - \omega(k_0,0) \;\;\, \text{for}\,\, k \to \infty
\end{equation}
and for all $b > 0\,$:
\begin{equation} \label{rinfty}
\widehat{r}(k,b) \to \infty \;\;\, \text{for}\,\, k \to \infty\,.
\end{equation}
If $b \geq 1/3$, then there holds $\partial^2_k \omega(k,b) >0$ for all $k >0$.  Because of $\omega(0,b) =0$ this implies $\omega(k_0/2,b) < \omega(k_0,b)/2$ for all $b \geq 1/3$ and
therefore $\widehat{r}(k_0/2,b) <0$ for all $b \geq 1/3$. Due to \eqref{dkr}, the positivity of $\partial^2_k \omega$ also yields
$\partial_k \widehat{r}(k,b) > 0$ for all $k > k_0/2$  and all $b \geq 1/3$.  Moreover, since 
$\partial^n_k\widehat{r}$ is continuous with respect to $k$ and $b$ for any $n \in \N_0$, there exist a constant $b_1 \in (0, 1/3)$ and functions $\gamma_j \in C^{0} ((b_1,\infty), \R^+)$, $j \in \{3,4,5\}$, and $\widetilde{\delta}_3 \in C^{0} ((b_1,\infty), (0,k_0/20))$  such that there holds
\begin{align}
 &\widehat{r}(k,b) \leq -\gamma_3(b) 
\end{align}
for all $b \in (b_1,\infty)$ and all $k\in [k_0/2, k_0-\widetilde{\delta}_3(b)]$,
\begin{align}
& \partial_k \widehat{r}(k,b) \geq \gamma_4(b) 
\end{align}
for all $b \in (b_1,\infty)$ and all $k\in [k_0-\widetilde{\delta}_3(b), k_0+\widetilde{\delta}_3(b)  ]$, as well as
\begin{align}
& \widehat{r}(k,b) \geq \gamma_5(b) 
\end{align}
for all $b \in (b_1,\infty)$ and all $k\in [k_0+\widetilde{\delta}_3(b), \infty)$; compare Figure \ref{fig2}, Panel (i)--(iii).

If $b =0$, then there holds $\partial^2_k \omega(k,b) < 0$ for all $k >0$, which implies $\omega(k_0/2,0) > \omega(k_0,0)/2$ and
therefore $\widehat{r}(k_0/2,0) > 0$. Due to \eqref{dkr}, the negativity of $\partial^2_k \omega$ also yields
$\partial_k \widehat{r}(k,0) < 0$ for all $k > k_0/2$. Hence,
there exist ${b}_2 \in (0, b_1)$,  $\gamma_j \in C^{0} ([0,{b}_2), \R^+)$ for $j \in \{6,7\}$, $\gamma_8 \in \R^+$ and $\widetilde{\delta}_4 \in C^{0} ([0,{b}_2), (0,k_0/20))$  such that there holds
\begin{align}
 & \widehat{r}(k,b) \geq \gamma_6(b) 
\end{align}
for all $b \in [0,{b}_2)$ and all $k\in [k_0/2, k_0-\widetilde{\delta}_4(b)]$,
\begin{align} \label{dkrminus}
& \partial_k \widehat{r}(k,b) \leq -\gamma_7(b)
\end{align}
for all $b \in [0,{b}_2)$ and all $k\in [k_0-\widetilde{\delta}_4(b), k_0+\widetilde{\delta}_4(b)  ]$, as well as
\begin{align}
& \widehat{r}(k,0) \leq -\gamma_8
\end{align}
for all  $k\in [k_0+\widetilde{\delta}_4(b), \infty)$; compare Figure \ref{fig2}, Panel (vii)--(ix).

Because of $\widehat{r}(k_0,b)=0$, \eqref{rinfty}, \eqref{dkrminus} and the intermediate value theorem there exist functions ${k}_1, {k}_2 \in C^{0} ((0,{b}_2), (k_0, \infty))$ with
\begin{align}
& \widehat{r}({k}_1(b),b) =  0\,, \\[1mm]
& \partial_k \widehat{r}({k}_2(b),b) =  0
\end{align}
for all $b \in (0, {b}_2)$.

Since there exists a strictly monotonically decreasing function ${k}_3 \in C^{0} ((0,1/3),(0, \infty))$ with 
${k}_3(b) \to 0$ for  $b \to 1/3$ and
${k}_3(b) \to \infty$ for $b \to 0$
such that for all $b \in (0,1/3)$ there holds
\begin{align} \label{omegaconc}
& \partial^{2}_k \omega(k,b) < 0 \;\;\, \text{if}\,\, 0 < k < {k}_3(b)\,,
\\[2mm]
& \partial^{2}_k \omega(k,b) = 0 \;\,\, \text{if}\,\, k = {k}_3(b)\,,\\[2mm] \label{omegaconv}
& \partial^{2}_k \omega(k,b) > 0 \;\;\, \text{if}\,\, k > {k}_3(b)
\,,
\end{align}
the mean value theorem yields 
${k}_2(b) > {k}_3(b)$ for all $b \in (0,{b}_2)$.
Moreover, because there exists a unique function ${k}_4 \in C^{0} ((0,1/3))$ with $k_4(b) > k_3(b)$ and
\begin{align}
\partial^{3}_k \omega(k_4(b),b) = 0 
\end{align}
for all $b \in (0,1/3)$ and since
\begin{align} 
& \partial^2_k\widehat{r}(k,b) = \partial^2_k \omega(k,b) - \partial^2_k \omega(k-k_0,b) 
\end{align}
for all $k\in  [k_0/2, \infty)$ and all $b \in \R^+_0$, the function $k \mapsto \partial^2_k\widehat{r}(k,b)$ can have at most one zero on $(k_2(b),\infty)$ for all $b \in (0,{b}_2)$. Because of $\widehat{r}(k_0,b)=0$, \eqref{rinfty},   \eqref{dkrminus} and the mean value theorem it follows that the function ${k}_2$ and therefore also the function ${k}_1$ is unique.

Let $\widetilde{k}(c,b)= 4c^{2}\mathrm{tanh}(k_0)/9k_0b$ for $c >0$. Using \eqref{asdkomega}--\eqref{tayloromega} yields
\begin{align}
\lim\limits_{b \to 0}   \partial_k \omega(\widetilde{k}(c,b),b)\, k_0 = c\, ( k_0\, \mathrm{tanh}(k_0))^{1/2} = c\, \omega(k_0,0)\,,
\end{align}
which implies
\begin{align} \label{k1}
& {k}_1(b) = \frac{4\, \mathrm{tanh}(k_0)}{9k_0b} (1+o(1)) \;\;\, \text{for}\,\,\, b  \to 0\,.
\end{align}
Moreover, due to \eqref{q11}--\eqref{frakm}, there exist  $b_0 \in (0, {b}_2]$ and functions $C_{TW\!I} \in C^{0} ([0,b_0), (1,\infty))$,
$\delta_1 \in C^{0} ([0,b_0), (0,1))$ with $ \delta_1(b) < 1-(20(k_1(b)-k_0))^{-1}k_0$ for all $b \in (0,b_0)$,
\begin{align}
& \lim\limits_{b \to 0} \delta_1(b) \geq \frac12 \,,\\[2mm]
& \lim\limits_{b \to 0} C_{TW\!I}(b) = 1 
\end{align}
such that
\begin{align} \label{CTWI}
& \frac{1}{C_{TW\!I}(b)} \leq \frac{-\widehat{q}_{j_1 j_1}(\pm {k},\pm k_0,\pm ({k}-k_0))}{\widehat{q}_{j_1 j_1}(\mp ({k}- k_0),\pm k_0,\mp {k})} \leq  C_{TW\!I}(b) 
\end{align}
for all $ j_1 \in \{-1,-2\}$, $b \in (0,b_0)$ and $k\in [{k}_1(b) - \delta_1(b)({k}_1(b)-k_0), {k}_1(b) + \delta_1(b)({k}_1(b)-k_0)]$. Furthermore, there exist functions $\gamma_j \in C^{0} ([0,b_0), \R^{+})$, $j \in \{9,10\}$, such that 
\begin{align} \label{k1est1}
& \widehat{r}(k,b) \leq -\gamma_{9}(b) \,.
\end{align} 
for all $b \in (0,b_0)$ and all $k \in [k_0 +\widetilde{\delta}_4(b), {k}_1(b) - \delta_1(b)({k}_1(b)-k_0)/2]$, as well as
\begin{align} \label{k1est2}
& \widehat{r}(k,b) \geq \gamma_{10}(b) \,. 
\end{align}
for all $b \in (0,b_0)$ and all $k \in [{k}_1(b) + \delta_1(b)({k}_1(b)-k_0)/2, \infty)$; compare Figure \ref{fig2}, Panel (vii)--(viii).

By using the method from \cite{SSZ15} to compute the values of the coefficients in the TWI systems belonging to the two-dimensional water wave problem with finite depth in the Eulerian formulation and in the arc length formulation, one can additionally show that one can choose $b_0=b_2$.
We remark that this property of $b_0$ is not needed for the proof of our error estimates.

Finally, we define 
\begin{align*}
\mathcal{B}:= \{b \in [0,b_0) \cup (b_1, \infty):\, k \mapsto \omega(k,b)  \,\,\mathrm{satisfies}\,\, \eqref{nrb1},\, \eqref{nrb3}\,\, \mathrm{and}\,\, \eqref{nrb4} \,\,\mathrm{with}\,\, M=6. \}\,.  
\end{align*}
and $\delta_0= \delta_0(b)$ by 
\begin{align*}
\delta_{0}(b)= 
\begin{cases}
\min \{ \widetilde{\delta}_j(b): j \in \{0,1,2,3\} \} & \quad {\rm if }\;b \in (b_1,\infty)   \,,\\[2mm]
\min \{ \widetilde{\delta}_j(b): j \in \{0,1,2,4\} \} & \quad{\rm if}\;  b \in [0,b_0)\,.
\end{cases}
\end{align*}

Now, we are prepared to prove
\begin{theorem}
\label{arcresult}
Let $k_0 >0$ and $s_A \geq {10}$.  Then for all ${\tau}_0, C_0 > 0$ there exists an  $\epsilon_0 > 0$ and a function $C \in C^0(\mathcal{B},\R^+)$
such that for all  $b \in \mathcal{B}$ and all solutions ${\tA}_{-11}^0 \in
C^0([0,{\tau}_0],H^{s_A}(\R,\C))$ of the NLS equation (\ref{nlsderive})
with 
$$ 
\sup_{\tau \in [0,{\tau}_0]} \| {\tA}_{-11}^0(\cdot,\tau) \|_{H^{s_A}(\R,\C)} \leq C_0
$$  
the following holds. For all  $\epsilon \in (0,\epsilon_0)$ there exists a solution 
$$
\mathcal{U} \in
C^0([0,{\tau}_0 \epsilon^{-2}], (L^{2}(\R,\R))^2 \times (H^{s}(\R,\R))^2)
$$
of \eqref{w-1}--\eqref{deltaneu}, where $s=s_A-2$, which satisfies
  $$\sup_{t\in[0,{ {\tau}_0} {\epsilon^{-2}}]} \Big\| \mathcal{U}(\cdot,t) -
    {\epsilon} \left( \ba{c}
{\tPsi}^0_{1}  \\[1mm] \partial^{2}_{\ualpha}  {\tPsi}^0_{1}    \ea \right)\!(\cdot,t)\Big\|_{(L^{2})^2 \times(H^{s})^2} 
\leq C(b) {\epsilon^{3/2}}\,,$$
where
$$
{\tPsi}^0_{1}(\alpha,t) \,=\, {\tA}_{-11}^0 (\eps
(\alpha-c_g t),\eps^2 t)
\,\EE^{1} \!\left( \ba{c} 
1 \\ 0  \ea \right) + \mathrm{c.c.}
\,. 
$$
\end{theorem} 

\subsection{The construction of the energy}

For the proof of Theorem \ref{arcresult} we introduce the energy 
\begin{equation} 
\mathcal{E}_{\check{s}} = E_{1,0} + \mathcal{E}_{2,\check{s}} \,
\end{equation}
with $0 \leq \check{s} \leq s := s_A-2$, 
\begin{equation} 
\mathcal{E}_{2,\check{s}} = \sum_{l=0}^{\check{s}} E_{2,l} \,,
\end{equation}
\begin{equation}
E_{j,0} = \sum_{j_1\in \{\pm j\}} \frac{1}{2}   \int_{\R}  \check{R}_{j_1}\, \rho_{j_1}^{0} \check{R}_{j_1}  \,d{\alpha} \,, \\[2mm]
\end{equation}
\begin{equation}
\label{ch R}
\check{R}_{j_1} = {R}_{j_1} + \eps \sum_{j_2\in \{\pm j_1\}}  \mathcal{N}_{j_1 j_2}(\psi_c,R_{j_2}) + \eps^{2}  \sum_{j_2\in \{\pm j_1\}} \mathcal{T}_{j_1 j_2}(\psi_c, \psi_c,R_{j_2}) 
\end{equation}
and 
\begin{equation} 
E_{2,l} = \sum_{j_1\in \{\pm2\}} \Big( \frac{1}{2}   \int_{\R} \partial_{\alpha}^{l} R_{j_1}\, \rho_{j_1}^{l} \partial_{\alpha}^{l} R_{j_1}   \,d{\alpha}  +     \veps \sum_{j_2 \in \{\pm j_1\}} \int_{\R}\partial_{\alpha}^{l} R_{j_1}\, \rho_{j_1}^{l}  \partial_{\alpha}^{l} N_{j_1 j_2}(\psi_c,R_{j_2})\,d{\alpha} \Big)
\end{equation}
for $1 \leq l \leq s$.

Here, $\rho_{j_1}^{l}$ is defined by its symbol
\begin{align}
\widehat{\rho}_{j_1}^{\,l}(k)= 
\begin{cases}
1+ \displaystyle \sum\limits_{\ell \in \{\pm1\}} \widehat{\rho}_{j_1 \ell}^{\,l}(k) & \quad{\rm if}\; b \in (0,b_0) \; \mathrm{and} \,\, \mathrm{sgn}(j_1)\,= -1   \,,\\[2mm]
1 & \quad {\rm otherwise }  \,,
\end{cases}
\end{align}
where
\begin{align}
& \widehat{\rho}_{j_1 \ell}^{\,l}(k) = \Big( \frac{-\widehat{q}_{j_1 j_1}(-k +\ell k_0,\ell k_0,-k)}{\widehat{q}_{j_1 j_1}( k,\ell k_0,k-\ell k_0)} \Big(\frac{-k+\ell k_0}{k}\Big)^{2l} -1 \Big)\,\widehat{\xi}_1 \Big( \frac{k+\ell (k_{1}- k_0)}{k_1- k_0}\Big),
\end{align}
$k_1= k_1(b)$ is as above
and
$\widehat{\xi}_{i} \in C^{\infty}_c(\R,\R)$, $i \in \{0,1\}$, satisfies 
\begin{align}
& \widehat{\xi}_{i}(k) \;=\; \left\{ \begin{array}{ll}
1 & \quad {\rm if }\; |k| \leq \delta_{i}/2 \,,
\\[2mm]
 0  & \quad {\rm if }\; |k| \geq \delta_{i}  \,, \\[2mm]
\widehat{\xi}_{i}(|k|) \in [0,1]  & \quad {\rm otherwise }\,
  \end{array} \right. 
\end{align}  
with $\delta_{i}=\delta_{i}(b)$ as above.
Because of \eqref{CTWI} there exist a constant $C_{\rho} \geq 1$ such that there holds 
\begin{equation}
\label{rhobound}
C_{\rho}^{-1} \,\leq\,  \widehat{\rho}_{j_1}^{\,l} \,\leq\, C_{\rho}
\end{equation}
for all $j_1 \in \{\pm 1,\pm2\}$ and all $0 \leq l \leq s$ uniformly on compact subsets of $\mathcal{B}$.

The functions $N_{j_1 j_2}$, $\mathcal{N}_{j_1 j_2}$ and $\mathcal{T}_{j_1 j_2}$ are defined as follows. We set 
\begin{align}
&  {N}_{j_1 j_2}(\varphi,f) = \left\{ \begin{array}{ll}
\!   {N}_{j_1 j_2}^1(\varphi,f)  
& \quad\! \mathrm{if}\;  |j_1| = 1,
\\[4mm]
\!    {N}_{j_1 j_2}^1(\varphi,f) +  {N}_{j_1 j_2}^2(\varphi,\partial_{\alpha}^{-1} f)   & \quad \! \mathrm{if}\;  |j_1| = 2\,,
\end{array} \right. 
\end{align}  
where
\begin{equation}
{N}_{j_1 j_2}^j(\varphi,g) \;=\; \sum_{\ell \in \{\pm1\}} {N}_{j_1 j_2 \ell}^j({\varphi}_{\ell},g)
\end{equation}
with
\begin{align}
& \widehat{N}_{j_1 j_2 \ell}^j({\varphi}_{\ell},g)(k) 
\;=\;  \int_{\R} \widehat{n}_{j_1 j_2 \ell}^{\, j}(k)\, \widehat{{\varphi}}_{\ell}(k-m)\, \widehat{g}(m)\,dm\,,\\[2mm]  
\label{nkdef}
& \widehat{n}_{j_1 j_2 \ell}^{\,j}(k) \;=\;  
\frac{ \widehat{\mathfrak{q}}_{j_1 j_2}^{\,j}(k,\ell k_0,k-\ell k_0)}{
\widehat{r}_{j_1 j_2}(k,\ell k_0,k-\ell k_0)} \; \frac{\widehat{\zeta}_{j_1 j_2 \ell}(k)\, (\widehat{\vartheta} - \eps \widehat{\xi}_0)(k- \ell k_0)}{\widehat{\vartheta}(k)}\,,\\[3mm]
& \widehat{\mathfrak{q}}_{j_1 j_2}^{1}(k,k-m,m)  = \left\{ \begin{array}{ll}
\!  \widehat{q}_{j_1 j_2}(k,k-m,m)
& \quad\! \mathrm{if}\;  |j_1| = 1,
\\[3mm]
\!   \widehat{q}_{j_1 j_2}(k,k-m,m)  - \widehat{q}_{j_1 j_2}^{2,4}(k,k-m,m)      & \quad \! \mathrm{if}\;  |j_1| = 2\,,
\end{array} \right. 
\\[3mm]
& \widehat{\mathfrak{q}}_{j_1 j_2}^{2}(k,k-m,m)  = im\,\widehat{q}_{j_1 j_2}^{2,4}(k,k-m,m) \,,
\\[3mm]
& \widehat{r}_{j_1 j_2}(k,k-m,m) \;=\; i\, (\, \mathrm{sgn}(j_1)\, \omega(k) + \omega(k-m) - \mathrm{sgn}(j_2)\, \omega(m))  
\end{align}
and
\begin{align}
& \widehat{\zeta}_{j_1 j_2 \ell}(k) = \left\{ \begin{array}{ll}
\! 1- \widehat{\xi}_1 \Big( \dfrac{k-\ell k_1}{k_1- k_0}\Big)\! - \widehat{\xi}_1 \Big( \dfrac{k+\ell (k_1-k_0)}{k_1- k_0}\Big)
& \quad\! \mathrm{if}\; b \in (0,b_0)\; \mathrm{and} 
\\ & \quad\! \mathrm{sgn}(j_1)=\mathrm{sgn}(j_2)= -1,
\\[3mm]
\! 1 & \quad \! {\rm otherwise },
\end{array} \right. 
\end{align}  
as well as
\begin{align}
&  \mathcal{N}_{j_1 j_2}(\varphi,f) = \left\{ \begin{array}{ll}
\!   \mathcal{N}_{j_1 j_2}^1(\varphi,f)  
& \quad\! \mathrm{if}\;  |j_1| = 1,
\\[4mm]
\!    \mathcal{N}_{j_1 j_2}^1(\varphi,f) +  \mathcal{N}_{j_1 j_2}^2(\varphi,\partial_{\alpha}^{-1} f)   & \quad \! \mathrm{if}\;  |j_1| = 2\,,
\end{array} \right. 
\end{align}  
where
\begin{equation}
{\mathcal{N}}_{j_1 j_2}^i(\varphi,g) \;=\; \sum_{j=0}^{1} \mathcal{N}_{j_1 j_2}^{i,j}({\varphi},g)
\end{equation}
with
\begin{align}
& \widehat{\mathcal{N}}_{j_1 j_2}^{i,0}({\varphi},g)(k) 
\;=\;  \int_{\R} \widehat{n}_{j_1 j_2}^{\,i,0}(k,k-m,m)\, \widehat{{\varphi}}(k-m)\, \widehat{g}(m)\,dm\,,\\[2mm]  
& \widehat{n}_{j_1 j_2}^{\,i,0}(k,k-m,m) \;=\;  
\widehat{P}_{0, \delta_{0}}\!(k)\; \frac{ \widehat{\mathfrak{q}}_{j_1 j_2}^{\,i}(k,k-m,m)}{
\widehat{r}_{j_1 j_2}(k,k-m,m)} \; \frac{\widehat{\vartheta}(m)}{\widehat{\vartheta}(k)}
\end{align}
and
\begin{align}
& \widehat{\mathcal{N}}_{j_1 j_2}^{i,1}({\varphi},g)(k) 
\;=\; \widehat{P}_{\delta_{0},\infty}(k)\, \widehat{N}_{j_1 j_2}^{i}({\varphi},g)(k) \,.
\end{align}  
Moreover, we set
\begin{equation}
{\mathcal{T}}_{j_1 j_2}(g,h,f) \;=\; \sum_{j=|j_1|}^{2} {\mathcal{T}}_{j_1 j_2}^{j}(g,h,f)\,,
\end{equation}
where 
\begin{align}
& {\mathcal{T}}^{j}_{j_1 j_2}(g,h,f) \;=\; \sum_{\ell \in \{\pm1\}} {\mathcal{T}}_{j_1 j_2 \ell}^{j}(g_{\ell},h_{\ell},f)\,,
\\[2mm]
& \widehat{\mathcal{T}}_{j_1 j_2 \ell}^{j}(g_{\ell},h_{\ell},f)(k) \;=\;  \int_{\R} \int_{\R} \widehat{\tau}^{j}_{j_1 j_2 \ell}(k)\, \widehat{g}_{\ell}(k-m)\,
\widehat{h}_{\ell}(m-n)\, \widehat{f}(n)\,dn dm\,,  
\\[2mm]
\label{taudef1}
& \widehat{\tau}^{1}_{j_1 j_2 \ell}(k) \;=\;\widehat{P}_{0,\delta_{0}}\!(k)\, \frac{ \widehat{c}_{j_1 \ell \ell j_2}(k,\ell k_0,\ell k_0,k -2\ell k_0)}{\widehat{\nu}_{j_1 j_2}(k, \ell k_0,\ell k_0,k -2\ell k_0)} \; \frac{ \widehat{\vartheta}(k- 2\ell k_0)}{\widehat{\vartheta}(k)}\,,
\\[2mm]
& \widehat{\tau}^{2}_{j_1 j_2 \ell}(k) \;=\; \sum_{j_3 \in \{\pm1\}}  \widehat{P}_{0,\delta_{0}}\!(k) \, \frac{ \widehat{q}_{j_1 j_2 j_3}(k,\ell k_0, k-\ell k_0,\ell k_0, k -2 \ell k_0)}{\widehat{r}_{j_1 j_2 j_3}(k,\ell k_0,k-\ell k_0,   \ell k_0,k -2\ell k_0)} \; \frac{ \widehat{\vartheta}(k- 2\ell k_0)}{\widehat{\vartheta}(k)}\,,
\\[2mm]
& \widehat{\nu}_{j_1 j_2}(k,k-m,m-n,n) \;=\; i\,(\mathrm{sgn}(j_1)\, \omega(k)+\omega(k-m)+ \omega(m-n) -\mathrm{sgn}(j_2)\,\omega(n)) 
\,,
\\[2mm]
& \widehat{r}_{j_1 j_2 j_3}(k,k-m,m,m-n,n) \;=\; \widehat{r}_{j_1 j_3}(k,k-m,m)\, \widehat{\nu}_{j_1 j_2}(k,k-m,m-n,n)\,,
\\[2mm]
& \widehat{q}_{j_1 j_2 j_3}(k,k-m,m,m-n,n) \;=\; \widehat{q}_{j_1 j_3}(k,k-m,m)\,
\widehat{q}_{j_3 j_2}(m,m-n,n)
\end{align}
and $\widehat{c}_{j_1 \ell \ell j_2}$ is given by
\begin{equation}
\label{taudef6}
\int_{\R} \int_{\R} \widehat{c}_{j_1 \ell \ell j_2}(k, k-m,m-n,n)\, \widehat{g}_{\ell}(k-m)\,
\widehat{h}_{\ell}(m-n)\, \widehat{f}(n)\,dn dm
\;=\; \widehat{\mathcal{C}}_{j_1 \ell \ell j_2}(g_{\ell},h_{\ell},f)(k)
\,. 
\end{equation}
 Finally, the functions $p_{\ell}$ with $p \in \{\varphi,g,h\}$ are defined by $\widehat{p}_{\ell} = \widehat{p}\, \chi_{\R^+_0}$ if $\ell=1$ and $\widehat{p}_{\ell} = \widehat{p}\, \chi_{\R^-_0}$ if $\ell=-1$.

As explained in the introduction, the construction of the energy $\mathcal{E}_s$ is inspired by the method of normal-form transforms, where the normal-form transform is incorporated directly in $\mathcal{E}_s$. The normal-form transform we use
is an extension of a normal-form transform of the form \eqref{inft} and \eqref{nkernel2} in order to handle the non-trivial resonances being present in the two-dimensional water wave problem with finite depth. 

The weight function $\widehat{\vartheta}$ and the correction function $\eps \widehat{\xi}_0$ are included to handle the non-trivial resonances at $\pm k_0$.
The trilinear mappings $\mathcal{T}_{j_1 j_2}$ are constructed in such a way that they generate terms in the evolution equation of $\mathcal{E}_s$ which cancel all the terms of order $\mathcal{O}(\eps)$ in the evolution equation which are caused by the fact that 
$\widehat{\vartheta}^{-1}$ is of order $\mathcal{O}(\eps^{-1})$ for $|k| \leq \delta_{0}$.
The weight functions $\widehat{\rho}_{j_1 \ell}^{\,l}$ and the correction functions           $\widehat{\zeta}_{j_1 j_2 \ell}$ are included to control the additional non-trivial resonances. Their form is motivated by the conserved quantity \eqref{staben}. The factor $1/(k_1-k_0)$ in the definition of $\widehat{\zeta}_{j_1 j_2 \ell}$ is chosen in such a way that we obtain error estimates which are uniform with respect to $b$ as $b$ and $\eps$ go to $0$.

The following lemma will allow us to show that it is sufficient for our goals that $\widehat{\rho}_{j_1 \ell}^{\,l}$, $\widehat{n}^{j}_{j_1 j_2 \ell}$ and $\widehat{\tau}^{j}_{j_1 j_2 \ell}$ in $\mathcal{E}_s$ depend only on $k$ and not on $k$ and $m$ like the kernel \eqref{nkernel}. 

\begin{lemma}
\label{int-kerne}
Let $p \in \R$, $g \in C^2(\R,\C)$ have a compactly supported Fourier transform, $\mathfrak{K}: \R^3 \to \C \cup \{\infty\}$ be uniformly bounded for all $(k,l,m) \in \R^3$ with $l-p,k-m-p \in \supp \widehat{g}$ and 
$f \in H^{s}(\R,\C)$ for $s \geq 0$.\\[2mm]
{\bf a)} If $l \mapsto \mathfrak{K}(k,l,m)$ is Lipschitz continuous in some neighborhood
of $p$ with a Lipschitz constant $L$ being independent of $k$ and $m$ if $k-m-p \in \supp \widehat{g}$, then there exist constants $\eps_0 >0$ and $C_2 >0$ with $C_2 \lesssim L  \Vert \widehat{g} \Vert_{L^{1}(s+1)}$ such that for all $\eps \in (0, \eps_0)$ there holds
\begin{equation} \label{IK1}
\Big \| \int \big(\mathfrak{K} (\cdot, \cdot-m,m) - \mathfrak{K} (\cdot, p ,m)\big)\, \veps^{-1} \widehat{g}\Big(\frac{ \cdot - m - p}{\varepsilon}\Big) \widehat{f}(m)\, dm\, \Big \|_{L^{2}(s)}
 \le C_2\, \varepsilon \| f \|_{H^{s}} \,.
\end{equation}
{\bf b)} If $m \mapsto \mathfrak{K}(k,k-m,m)$ is Lipschitz continuous for all $m \in \R$ for which  $k-m-p \in \supp \widehat{g}$ with a Lipschitz constant $L$ being independent of $k$ and $m$, then there exist a constant $C_3 >0$ with $C_3 \lesssim L  \Vert \widehat{g} \Vert_{L^{1}(s+1)}$ such that for all $\eps  > 0$ there holds
\begin{equation}  \label{IK2}
\Big \| \int \big(\mathfrak{K} (\cdot, \cdot-m,m)- \mathfrak{K}(\cdot, \cdot - m ,\cdot-p)\big)\, \veps^{-1} \widehat{g}\Big(\frac{ \cdot - m - p}{\varepsilon} \Big) \widehat{f}(m)\, dm\, \Big \|_{L^{2}(s)}
\le C_3\, \varepsilon \| f \|_{H^{s}} \,.
\end{equation}
\end{lemma}

\textbf{Proof.}
The lemma is proven analogously as Lemma 3.5 in \cite{DSW12}.
\qed  
\medskip

In order to verify Theorem \ref{arcresult} we would like to prove that the energy $\mathcal{E}_s$ controls the error $\mathcal{R}$ in $(L^{2}(\R,\R))^{2} \times (H^{s}(\R,\R))^{2}$ and remains bounded for all $t \in [0,\tau_0  \veps^{-2}]$. However, we will see that the right-hand side of the evolution equation for $\mathcal{E}_s$ contains terms which cannot be estimated by a multiple of $ \eps^{2}(\mathcal{E}_s+1)$ such that Gronwall's inequality cannot be applied to obtain the desired bound. Nevertheless, these terms can be identified as time derivatives of time dependent integrals. Hence, by adding these integrals to $\mathcal{E}_s$
we will obtain a new energy $\widetilde{\mathcal{E}_s}$ at the end of Subsection 4.4 and this energy will have the desired properties. We will prove

\begin{lemma}
\label{finenest}
For sufficiently small $\eps >0$, $\widetilde{\mathcal{E}_s}$ satisfies
\begin{align}
\label{tEest}
& \frac{d}{dt} \widetilde{\mathcal{E}_s}\, \lesssim\, \veps^{2} (\widetilde{\mathcal{E}_s}+1) \,,
\\[2mm] 
\label{tEnorm}
& \widetilde{\mathcal{E}_s} \lesssim \|\mathcal{R}_{1}\|_{(L^{2})^2}^{2} + \|\mathcal{R}_{2}\|_{(H^{{s}})^2}^{2} + \eps^{5} \,,\\[2mm] \label{normtE}
&  \|\mathcal{R}_{1}\|_{(L^{2})^2}^{2} + \|\mathcal{R}_{2}\|_{(H^{{s}})^2}^{2} \lesssim \widetilde{\mathcal{E}_s} + \eps^{5} \,,
\end{align}
as long as $\eps^{5/2} \Vert \mathcal{R}_{1} \Vert_{(L^{2})^2},\, \eps^{5/2} \Vert \mathcal{R}_{2} \Vert_{(H^{{s}})^2} \ll 1$,
uniformly on compact subsets of $\mathcal{B}$.
\end{lemma}

With the help of the energy $\widetilde{\mathcal{E}_s}$ it is also possible to give an alternative proof of local existence and uniqueness of solutions $\mathcal{R}$ to \eqref{Rmp1}--\eqref{R2R1p} in Sobolev spaces  without using the local existence and uniqueness results for the water wave problem \eqref{ytilde}--\eqref{ctilde}.  In this proof, the energy $\widetilde{\mathcal{E}_s}$ takes the role of the energy used in the proof of the local existence and uniqueness result for the water wave problem in the arc length formulation in \cite{AM05}. 

Because $\widetilde{\mathcal{E}_s}$ remains $\mathcal{O}(1)$-bounded for a timespan of order $\mathcal{O}(\eps^{-2})$, the solutions $\mathcal{R}$ to \eqref{Rmp1}--\eqref{R2R1p} even exist for this long timespan and belong to $(L^{2}(\R,\R))^2 \times (H^{s}(\R,\R))^2$. The Fourier truncation procedure from Section \ref{sec3} allows us to choose $s=s_A-2$ such that the 
error $\mathcal{R}$ has the same Sobolev regularity as the solution of the NLS equation.

Since the energy estimates \eqref{tEest}--\eqref{normtE} are only valid as long as $\eps^{5/2} \Vert \mathcal{R}_{1} \Vert_{(L^{2})^2} \ll 1$ and $\eps^{5/2} \Vert \mathcal{R}_{2} \Vert_{(H^{{s}})^2} \ll 1$, they are not sufficient to guarantee global existence of solutions to \eqref{Rmp1}--\eqref{R2R1p} and global existence of small oscillating wave packet-like solutions to the two-dimensional water wave problem with finite depth
by iterating the arguments of the proof of local existence in time. 

For the three-dimensional water wave problem with infinite depth, for the three-dimensional water wave problem with finite depth and either gravity or surface tension as well as for the two-dimensional water wave problem with infinite depth and either gravity or surface tension it is possible to
combine energy estimates of this type (so-called high-order energy estimates) with dispersive decay estimates to establish global existence of solutions for small data.  A detailed explanation of this procedure is given, for example, in \cite{KP83} and global existence results for the water wave problem obtained with the help of this procedure can be found, for example, in \cite{AD15, DIPP17, GMS09, IT14a, IT14b, IP15, IP18, Wa16, Wa19, Wu09}.

However, for the two-dimensional water wave problem with infinite depth, gravity and surface tension the optimal dispersive decay rate is not strong enough to ensure global existence of solutions and for the two-dimensional water wave problem with finite
depth as well as for the three-dimensional water wave problem with finite depth, gravity and surface tension such dispersive decay estimates cannot be expected because of the existence of solitary waves solutions, which are localized traveling waves of permanent form. In all these cases, the global existence of solutions is still an open problem. 
For a more detailed discussion of the existence issue for solutions to the water wave problem we refer, for example, to  
the review article \cite{D18} and the references therein.

\subsection{Fundamental properties of the energy}

Now, we show several fundamental properties of the operators $N_{j_1j_2}$,  $\mathcal{N}_{j_1j_2}$ and 
${\mathcal{T}}_{j_1 j_2}$, which will be mandatory for the proof of our energy estimates.

\begin{lemma}
\label{lem31}
The operators $N_{j_1j_2}^{i}$ have the following properties:\\[2mm]
{\bf a)} Fix $\varphi \in L^2(\R,\R)$ with $\mathrm{supp}\, \widehat{\varphi} = \mathrm{supp}\, \widehat{\psi}_c$. Then $f \mapsto N_{jj}^i(\varphi,f)$ defines a continuous linear map from $H^1(\R,\R)$ into $L^2(\R,\R)$ and $f \mapsto N_{j-j}^{i}(\varphi,f)$ a continuous linear map from $H^{(1-(|j|-1))/2}(\R,\R)$ into $L^2(\R,\R)$. Furthermore, 
 there exists a constant $C>0$ with $C \lesssim \Vert \widehat{\varphi} \Vert_{L^{1}}$ such that for all $f\in H^{(1-(|j|-1))/2}(\R,\R)$, all $g\in H^{1}(\R,\R)$ and all $h\in H^{1+(1-(|j|-1))/2}$ $(\R,\R)$  and    $p\in H^{2}(\R,\R)$ there holds
\begin{align}
\label{N0}
\Vert N_{jj}^i(\varphi, g) \Vert_{L^{2}} & \leq\, C \varepsilon^{-1} \Vert g \Vert_{H^1}\,, \\[2mm]
\label{N0b}
\Vert N_{j-j}^i(\varphi, f) \Vert_{L^{2}} & \leq \, C \varepsilon^{-1} \Vert f \Vert_{H^{(1-(|j|-1))/2}}\,,
\\[2mm]
\label{thetaN}
\Vert \mathcal{P} N_{jj}^i(\varphi, g) \Vert_{L^{2}} & \leq \, C \Vert g \Vert_{H^1}\,,\\[2mm]
\label{thetaNb}
\Vert \mathcal{P} N_{j-j}^i(\varphi, f) \Vert_{L^{2}} & \leq \, C   \Vert f \Vert_{H^{(1-(|j|-1))/2}}\,,
\\[2mm]
\label{d/dxN}
\Vert \partial_{\alpha} N_{jj}^i(\varphi, p) \Vert_{L^{2}} & \leq \, C  \Vert p \Vert_{H^2}\,,\\[2mm]
\label{d/dxNb}
\Vert \partial_{\alpha} N_{j-j}^i(\varphi, h) \Vert_{L^{2}}  & \leq \, C  \Vert h \Vert_{H^{{1+(1-(|j|-1))/2}}}
\end{align}
uniformly on compact subsets of $\mathcal{B}$, where $\mathcal{P}= P_{\delta_{0}, \infty}$ or $\mathcal{P}= \vartheta$.
\\[2mm]
{\bf b)}  Let $\varphi$ be as in a). Then for all $f \in L^2(\R,\R)$ there holds
\begin{align}
\label{P_0 N}
P_{0, \delta_{0}}N^{i}_{j_1j_2}(\varphi, P_{0, \delta_{0}}f)=0\,.
\end{align}
\end{lemma}
{\bf Proof.}
a) The key step of the proof is to discuss systematically the kernels
$\widehat{n}_{j_1 j_2 \ell}^{i}$ for all $i \in \{1,2\}$. We start by analyzing the behavior of $\widehat{n}_{j_1 j_2 \ell}^{i}$ in a neighborhood of the zeros of the factor $\widehat{r}_{j_1 j_2}$ in the denominator. As shown above, we have zeros at
\begin{itemize}
\item $k=0$ if $\mathrm{sgn}(j_2) = -1$.
\\
[1mm]
Because of \eqref{q11}--\eqref{frakm} we have
\begin{equation}
|\widehat{P}_{0,\delta_0}(k)\, \widehat{\mathfrak{q}}_{j_1 j_2}^{\,i}(k,\ell k_0,k-\ell k_0)| \lesssim |k|
\end{equation} 
uniformly with respect to $b \lesssim 1$. Hence, due to 
\eqref{0bound},
the singularity of $\widehat{n}_{j_1 j_2 \ell}^{i}$ at $k=0$ can be removed and then  $\widehat{P}_{0,\delta_0} \widehat{n}_{j_1 j_2 \ell}^{i}$ is bounded uniformly on compact subsets of $\mathcal{B}$. However, because of
\begin{equation}
\widehat{P}_{0,\delta_0}(k)\, \widehat{\vartheta}^{-1} (k) = \mathcal{O}(\veps^{-1})
\end{equation}
we have
\begin{align}
\label{nz0}
\widehat{P}_{0,\delta_0}(k)\, \widehat{n}_{j_1 j_2 \ell}^{i}(k) = \mathcal{O}(\veps^{-1})
\end{align}
uniformly on compact subsets of $\mathcal{B}$.
\item $k= \ell k_0$ if $\mathrm{sgn}(j_1) = -1$.
\\[1mm]
By construction of $(\widehat{\vartheta}- \eps \widehat{\xi}_0)(k - \ell k_0)$
we have
\begin{equation}
|\widehat{P}_{0,\delta_0}(k- \ell k_0)\, (\widehat{\vartheta} - \eps \widehat{\xi}_0)(k - \ell k_0)| \leq  (1+\eps)|k-\ell k_0|/\delta_{0} \,.
\end{equation}
Hence, because of \eqref{k0bound} and since $\sigma$ is differentiable with respect to $k$, where $\partial_k \sigma$ depends continuously on $b$,
the singularity of $\widehat{n}_{j_1 j_2 \ell}^{i}$ at $k= \ell k_0$ can be removed 
such that we obtain
\begin{align}
\label{nzk0}
\widehat{P}_{0,\delta_0}(k- \ell k_0)\, \widehat{n}_{j_1 j_2 \ell}^{i}(k) = \mathcal{O}({1})
\end{align}
uniformly on compact subsets of $\mathcal{B}$.
\item $k= \ell k_1$ and $k=-\ell(k_1-k_0)$  if $b \in (0,b_2)$ and $\mathrm{sgn}(j_1)=\mathrm{sgn}(j_2)=-1$.
\\[1mm]
The function $\widehat{\zeta}_{j_1j_2\ell}$ is constructed in such a way that
the singularities of $\widehat{n}_{j_1 j_2 \ell}^{i}$ at $k= \ell k_1$ and $k=-\ell(k_1-k_0)$ can be removed and then, due to \eqref{k1est1}--\eqref{k1est2} and the fact that $k \mapsto \omega(k,b)$ is odd, we obtain
\begin{align}
\label{nzadd1}
& \widehat{P}_{0,\delta_1(k_1-k_0)}(k- \ell k_1)\,\widehat{n}_{j_1 j_2 \ell}^{i}(k) = \mathcal{O}(|\widehat{\mathfrak{q}}_{j_1 j_2}^{\,i}(k,\ell k_0,k-\ell k_0)|)\,, \\[2mm]
\label{nzadd2}
& \widehat{P}_{0,\delta_1(k_1-k_0)}(k+ \ell (k_1-k_0))\,\widehat{n}_{j_1 j_2 \ell}^{i}(k) = \mathcal{O}(|\widehat{\mathfrak{q}}_{j_1 j_2}^{\,i}(k,\ell k_0,k-\ell k_0)|)
\end{align}
uniformly on compact subsets of $\mathcal{B}$.
\end{itemize}

Next, we analyze the asymptotic behavior of the kernels $\widehat{n}_{j_1 j_2 \ell}^{i}$ for $|k| \rightarrow \infty$. Because of \eqref{astanh}--\eqref{asdktanh} we have 
\begin{align}  \label{assigma}
\sigma(k,b) &\,=\, |k|^{1/2} (1+bk^{2})^{1/2} + \mathcal{O}(e^{-|k|}) 
\,,\\[2mm] 
\label{asdksigma}
\partial_k \sigma(k,b) &\,=\, \mathrm{sgn}(k)\, \frac{1+3bk^{2}}{2(1+bk^{2})}\, |k|^{-1/2}  (1+bk^{2})^{1/2} + \mathcal{O}(e^{-|k|})
\end{align}
 for $|k| \rightarrow \infty$ uniformly with respect to $b \lesssim 1$. 
Inserting \eqref{assigma} in \eqref{q2}--\eqref{q3} yields
\begin{align} \label{asq2}
(ik)^{-1}\, \widehat{q}_{j_1j_2}^{2,2}(k,\ell k_0,k-\ell k_0) &\,=\, \mathcal{O}(|k|^{-1/2}(1+bk^2)^{-1/2}) \,,\\[2mm]
\label{asq3}
(ik)^{-1}\, \widehat{q}_{j_1j_2}^{2,3}(k,\ell k_0,k-\ell k_0) &\,=\, \mathcal{O}(b^{1/2}|k|^{-1/2}) 
\end{align}
for $|k| \rightarrow \infty$ uniformly with respect to $b \lesssim 1$. 
Furthermore, with the help of \eqref{rootsigma} and the mean value theorem we derive
\begin{align}
& (ik)^{-1}\, \widehat{q}_{j_1j_2}^{2,4}(k,\ell k_0,k-\ell k_0) \nonumber \\[1mm] & \qquad =  \mathrm{sgn}(j_1)\, \sigma^{-1}( \ell k_0) (\ell k_0)^2\, \frac{ \sigma(k- \theta(k)(k-\ell k_0)) \, \partial_k \sigma(k- \theta(k)(k-\ell k_0))}{(\sigma(k)+\sigma(\ell k_0))\,(k-\ell k_0)} 
\end{align}
with $\theta(k) \in [0,1]$ for all $k \in \R \setminus \{\ell k_0\}$. Because of \eqref{assigma}--\eqref{asdksigma} we obtain
\begin{align} \label{asq4}
(ik)^{-1}\, \widehat{q}_{j_1j_2}^{2,4}(k,\ell k_0,k-\ell k_0) &\,=\,  
\mathcal{O}(|k|^{-3/2}(1+ bk^2)^{1/2})
\end{align}
for $|k| \rightarrow \infty$ uniformly with respect to $b \lesssim 1$. 

Due to \eqref{q11}--\eqref{q1}, \eqref{q5}, \eqref{asq2}--\eqref{asq3} and \eqref{asq4}, we conclude
\begin{align} \label{asqjj}
(ik)^{-1}\, \widehat{\mathfrak{q}}_{jj}^{1}(k,\ell k_0,k-\ell k_0) &\,=\, -1 + \mathcal{O}(|k|^{-1/2}) \,,\\[2mm] \label{asqj-j}
(ik)^{-1}\, \widehat{\mathfrak{q}}_{j-j}^{1}(k,\ell k_0,k-\ell k_0) &\,=\, \mathcal{O}(|k|^{-(|j|-1)/2}) \,, \\[2mm] \label{asfq2}
(ik)^{-1}\, \widehat{\mathfrak{q}}_{j_1j_2}^{2}(k,\ell k_0,k-\ell k_0) &\,=\,  
\mathcal{O}(|k|^{-1/2}(1+ bk^2)^{1/2})
\end{align}
for $|k| \rightarrow \infty$ uniformly with respect to $b \lesssim 1$. 

Moreover, by the mean value theorem we have
\begin{align}
\label{rmean1}
\widehat{r}_{j j}(k,\ell k_0,k-\ell k_0) &\,=\, i(\omega(\ell k_0)+ \mathrm{sgn}(j) \omega'(k-\theta(k)\ell k_0)\,\ell k_0)\,, \\[2mm]
\label{rmean2}
\widehat{r}_{j -j}(k,\ell k_0,k-\ell k_0)  &\,=\, i(\omega(\ell k_0)+ 2\, \mathrm{sgn}(j) \omega(k)-\mathrm{sgn}(j) \omega'(k-\theta(k)\ell k_0)\, \ell k_0)\, 
\end{align}
with $\theta(k) \in [0,1]$ for all $k \in \R$. Due to \eqref{asomega}--\eqref{asdkomega} and \eqref{k1est1}--\eqref{k1est2}, this implies 
\begin{align}
\label{asrjj}
\widehat{\zeta}_{jj\ell}(k)\,(\widehat{r}_{j j}(k,\ell k_0,k-\ell k_0))^{-1} &\,=\, \mathcal{O}((1+|k|^{-1/2}(1+ bk^2)^{1/2})^{-1} ) \,,\\[2mm] 
\label{asrj-j}
(\widehat{r}_{j -j}(k,\ell k_0,k-\ell k_0))^{-1}  &\,= \mathcal{O}(|k|^{-1/2}(1+ bk^2)^{-1/2})
\end{align}
for $|k| \rightarrow \infty$ uniformly on compact subsets of $\mathcal{B}$. 

Now, using \eqref{theta_m-k}, \eqref{asqjj}--\eqref{asqj-j}, \eqref{asfq2} and \eqref{asrjj}--\eqref{asrj-j}, we obtain
\begin{align}
\label{asnjj}
(ik)^{-1}\, \widehat{n}_{j j \ell}^{1}(k) & \,=\, \mathcal{O}((1+|k|^{-1/2}(1+ bk^2)^{1/2})^{-1})\,, \\[2mm]
\label{asnj-j}
\widehat{n}_{j -j \ell}^{1}(k) &\,=\, \mathcal{O}(|k|^{(1-(|j|-1))/2}(1+ bk^2)^{-1/2})\,,
\\[2mm]
\label{asn2jj}
(ik)^{-1}\, \widehat{n}_{j j \ell}^{2}(k) & \,=\, \mathcal{O}(1)\,, \\[2mm]
\label{asn2j-j}
\widehat{n}_{j -j \ell}^{2}(k) &\,=\, \mathcal{O}(1)
\end{align}
for $|k| \to \infty$ uniformly on compact subsets of $\mathcal{B}$.
Hence, combining \eqref{theta-1-eps-1}--\eqref{theta-1}, \eqref{nz0}, \eqref{nzk0},   \eqref{nzadd1}--\eqref{nzadd2}, \eqref{asnjj}--\eqref{asn2j-j} and Young's inequality for convolutions,
we arrive at \eqref{N0}--\eqref{d/dxNb}. 
Since
\begin{equation*}
\widehat{n}_{j_1 j_2 -\ell}^{i}(-k)= \overline{\widehat{n}_{j_1j_2 \ell}^{i}(k)} 
\end{equation*}
and $\varphi$ is real-valued, $f \mapsto N_{jj}^{i}(\varphi,f)$ is a continuous linear map from $H^1(\R,\R)$ into $L^2(\R,\R)$ and $f \mapsto N_{j-j}^{i}(\varphi,f)$ a continuous linear map from $H^{(1-(|j|-1))/2}(\R,\R)$ into $L^2(\R,\R)$, such that
we have proven all assertions of a).
\\[2mm]
b) is a direct consequence of
\begin{align*}
\widehat{P}_{0,\delta_{0}}\!(k)\, \widehat{P}_{0,\delta_{0}}\!(m)\, \chi_c(k-m)=0\,.
\end{align*}
$\hfill \Box$

\begin{lemma}
\label{lem31b}
Let $R_{-1}, R_1 \in L^2(\R,\R)$, $R_{- 2}, R_2 \in H^s(\R,\R)$ and $1 \leq l \leq s-1$. Then there holds
\begin{align} \label{nf}
& \int_{\R} \partial_{\alpha}^{l} R_{j_1} \, \rho_{j_1}^{l} \partial_{\alpha}^{l} \big(\mathrm{sgn}(j_1) i \omega N_{j_1 j_2}(\psi_c,R_{j_2})+ N_{j_1 j_2}(i \omega \psi_c,R_{j_2})
\nonumber \\[1mm] 
& \quad
- \mathrm{sgn}(j_2) N_{j_1 j_2}(\psi_c,i \omega R_{j_2})\big)\, d\alpha 
\nonumber \\[2mm] 
& \qquad\quad  =\, \int_{\R} \partial_{\alpha}^{l} R_{j_1}\, {\rho}_{j_1}^{l} \partial_{\alpha}^{l} Q_{j_1 j_2}(\psi_c) R_{j_2}\, d\alpha 
 + \eps \int_{\R} \partial_{\alpha}^{l} R_{j_1}\, {\rho}_{j_1}^{l} \partial_{\alpha}^{l} {Y}_{j_1 j_2}(\psi_c,R_{j_2})\, d\alpha \,
\end{align}
with
\begin{align} \label{nf-rest}
& \Big| \sum_{\genfrac{}{}{0pt}{}{j_1 \in \{\pm 1\},} {j_2 \in \{\pm j_1\}}}  \int_{\R} \partial_{\alpha}^{l} R_{j_1}\, {\rho}_{j_1}^{l}  \partial_{\alpha}^{l} {Y}_{j_1 j_2}(\psi_c, R_{j_2})\, d\alpha \Big| \,\lesssim\, \| \mathcal{R}_{1}\|_{(L^{2})^{2}}^2 + \| \mathcal{R}_{2}\|_{(H^{\max \{2,l\}})^{2}}^2 + \eps^{4} \,,
\end{align}
as long as $\eps^{5/2} \Vert \mathcal{R}_{1} \Vert_{(L^{2})^2},\, \eps^{5/2} \Vert \mathcal{R}_{2} \Vert_{(H^{\max \{2,l\}})^2} \ll 1$, and
\begin{align} \label{nf-rest2}
& \Big| \sum_{\genfrac{}{}{0pt}{}{j_1 \in \{\pm 2\},} {j_2 \in \{\pm j_1\}}}  \int_{\R} \partial_{\alpha}^{l} R_{j_1}\, {\rho}_{j_1}^{l}  \partial_{\alpha}^{l} {Y}_{j_1 j_2}(\psi_c, R_{j_2})\, d\alpha \Big| 
\nonumber \\[2mm] 
& \qquad\quad 
\lesssim\,  \| \mathcal{R}_{2} \|_{(H^{l})^2} \, (\| \mathcal{R}_{1}\|_{(L^{2})^{2}} + \| \mathcal{R}_{2}\|_{(H^{l+1})^{2}} + \eps^{2}) \,,
\end{align}
as long as $\eps^{5/2} \Vert \mathcal{R}_{1} \Vert_{(L^{2})^2},\, \eps^{5/2} \Vert \mathcal{R}_{2} \Vert_{(H^{l+1})^2} \ll 1$, uniformly on compact subsets of $\mathcal{B}$.
\end{lemma}
{\bf Proof.}
Let 
\begin{align*} 
\eps {Y}_{j_1 j_2}(\psi_c,R_{j_2}) :=
&\;\, \mathrm{sgn}(j_1) i \omega N_{j_1 j_2}(\psi_c,R_{j_2})+ N_{j_1 j_2}(i \omega \psi_c,R_{j_2})
\nonumber \\[1mm] 
&\;- \mathrm{sgn}(j_2) N_{j_1 j_2}(\psi_c,i \omega R_{j_2})  - Q_{j_1 j_2}(\psi_c)R_{j_2} \,. 
\end{align*}
Then we have
\begin{align}
\widehat{Y}_{j_1 j_2}(\psi_c,R_{j_2})(k)  & =\, \eps^{-1} \sum_{\mu=1}^{|j_1|} \sum_{\ell \in \{\pm1\}} \int_{\R} 
\widehat{\vartheta}^{-1}(k)\,K_{j_1 j_2 \ell}^{\mu} (k,k-m,m)\,
\widehat{\psi}_{\ell}(k-m)  \nonumber \\
& \hspace{3.3cm} \times (im)^{-(\mu -1)}\,\widehat{R}_{j_2}(m)\,dm\,
\end{align}
with
\begin{align}
& K_{j_1 j_2 \ell}^{\mu}(k,k-m,m)  \nonumber \\[1mm] & \qquad =\, \frac{\widehat{r}_{j_1 j_2}(k,k-m,m)\, \widehat{\zeta}_{j_1 j_2 \ell}(k)\, \widehat{\mathfrak{q}}_{j_1 j_2}^{\mu}(k, \ell k_0,k-\ell k_0)}{\widehat{r}_{j_1 j_2}(k,\ell k_0,k-\ell k_0)}
\, (\widehat{\vartheta}-\eps \widehat{\xi}_0)(k-\ell k_0) \nonumber
\\[2mm]
& \qquad \qquad \!\!\! - \widehat{\mathfrak{q}}_{j_1 j_2}^{\mu}(k,k-m,m)\,\widehat{\vartheta}(m) \,. \label{KLdef}
\end{align}
We split $K_{j_1 j_2 \ell}^{\mu}$ into
\begin{align} \label{KL1}
K_{j_1 j_2 \ell}^{\mu}(k,k-m,m) = &\,\, {K}_{j_1 j_2 \ell}^{\mu}(k,k-m,m)-{K}_{j_1 j_2 \ell}^{\mu}(k,\ell k_0,m) \nonumber \\[2mm]
&\,\, + {K}_{j_1 j_2 \ell}^{\mu}(k,\ell k_0,m)- {K}_{j_1 j_2 \ell}^{\mu}(k,\ell k_0,k-\ell k_0)\nonumber  \\[2mm]
&\,\, + {K}_{j_1 j_2 \ell}^{\mu}(k,\ell k_0,k-\ell k_0)\,. 
\end{align}
We have
\begin{align} \label{KL2}
{K}_{j_1 j_2 \ell}^{\mu}(k,\ell k_0,k-\ell k_0)= &\,\, -(1 - \widehat{\zeta}_{j_1 j_2 \ell}(k)+ \eps \widehat{\xi}_0(k-\ell k_0) )\, \widehat{\mathfrak{q}}_{j_1 j_2}^{\mu}(k, \ell k_0,k-\ell k_0) \,,
\end{align}
where $\widehat{\zeta}_{j_1 j_2 \ell}(k) =1$ if $b \in \mathcal{B} \setminus (0,b_0)$ or $\mathrm{sgn}(j_1) =1$ or $\mathrm{sgn}(j_2) =1$. 

Let $ b \in (0,b_0)$, $j \in \{-1,-2\}$, 
\begin{equation*}
\widehat{\Xi}_{j\ell}^{l \mu}(k,m) \,=\, -\widehat{\vartheta}^{-1}(k) \widehat{\rho}_{j}^{l}(k) k^{2l} (1 - \widehat{\zeta}_{j j \ell}(k))\, \widehat{\mathfrak{q}}_{jj}^{\mu}(k, \ell k_0,k-\ell k_0) (im)^{-(\mu -1)}
\end{equation*}
for all $k,m \in \R$ and the function $g_j^{l \mu}$ be defined by its Fourier transform
\begin{equation*}
\widehat{g}_j^{l\mu}(k) \,=\, \sum_{\ell \in \{\pm1\}} \int_{\R} 
(-ik)^{-l}\, \widehat{\Xi}_{j\ell}^{l\mu}(k,m)\,
\widehat{\psi}_{\ell}(k-m)\, \widehat{R}_{j}(m)\,dm 
\end{equation*}
for all $k \in \R$. 
Then, due to \eqref{dx-2R2}--\eqref{dx-1R2} and Young's inequality for convolutions, we have
$g_j^{l \mu} \in L^2(\R,\R)$
and with the help of Fubini's theorem we deduce
\begin{align*}
& \sum_{\ell \in \{\pm1\}} \int_{\R} \int_{\R} \overline{\widehat{R}_j(k)}\; \widehat{\Xi}_{j\ell}^{l \mu}(k,m)\,
\widehat{\psi}_{\ell}(k-m)\, \widehat{R}_{j}(m)\,dm dk \\[2mm] 
& \qquad =\, \frac12 \sum_{\ell \in \{\pm1\}} \int_{\R} \int_{\R} \overline{ \widehat{R}_j(k)}\; \widehat{\Xi}_{j\ell}^{l \mu}(k,m)\,
\widehat{\psi}_{\ell}(k-m)\, \widehat{R}_{j}(m)\,dm dk  \\[2mm] 
& \qquad \quad \,\,+ \frac12 \sum_{\ell \in \{\pm1\}} \int_{\R} \int_{\R}  {\widehat{R}_j(k)}\,\overline{\widehat{\Xi}_{j\ell}^{l \mu}(k,m)\,\widehat{\psi}_{\ell}(k-m)\,\widehat{R}_{j}(m)}\,dm dk   
\\[2mm]
& \qquad =\, \frac12 \sum_{\ell \in \{\pm1\}} \int_{\R} \int_{\R} \overline{ \widehat{R}_j(k)}\; \widehat{\Xi}_{j\ell}^{l \mu}(k,m)\,
\widehat{\psi}_{\ell}(k-m)\, \widehat{R}_{j}(m)\,dm dk  \\[2mm] 
& \qquad \quad \,\,+ \frac12 \sum_{\ell \in \{\pm1\}} \int_{\R} \int_{\R} \overline{\widehat{R}_j(k)\, \widehat{\Xi}_{j\ell}^{l \mu}(m,k)\,\widehat{\psi}_{\ell}(m-k)}\, {\widehat{R}_{j}(m)}\,dm dk  
\\[2mm]
& \qquad =\, \frac12 \sum_{\ell \in \{\pm1\}} \int_{\R} \int_{\R}  \overline{\widehat{R}_j(k)}\; \widehat{\Xi}_{j\ell}^{l \mu}(k,m)\,
\widehat{\psi}_{\ell}(k-m)\, \widehat{R}_{j}(m)\,dm dk  \\[2mm] 
& \qquad \quad \,\,+ \frac12 \sum_{\ell \in \{\pm1\}} \int_{\R} \int_{\R} \overline{\widehat{R}_j(k)}\; {\widehat{\Xi}_{j-\ell}^{l \mu}(-m,-k)}\;
{\widehat{\psi}_{-\ell}(k-m)}\,{\widehat{R}_{j}(m)}\,dm dk 
\\[1mm]
& \qquad =\, \frac12 \sum_{\ell \in \{\pm1\}} \int_{\R} \int_{\R} \overline{\widehat{R}_j(k)}\,\big( \widehat{\Xi}_{j\ell}^{l \mu}(k,k-\ell k_0) +  \widehat{\Xi}_{j\ell}^{l \mu}(-k+\ell k_0,-k)\big)\,    
\\
& \hspace{3.65cm} \times \widehat{\psi}_{\ell}(k-m)\, \widehat{R}_{j}(m)\,dm dk  
\end{align*}
\begin{align*} 
& \qquad \quad \,\,+ \frac12 \sum_{\ell \in \{\pm1\}} \int_{\R} \int_{\R} \overline{\widehat{R}_j(k)}\, \big({\widehat{\Xi}_{j\ell}^{l \mu}(k,m)} -  \widehat{\Xi}_{j\ell}^{l \mu}(k,k-\ell k_0) \big)\,   
\\
& \hspace{4.05cm} \times 
{\widehat{\psi}_{\ell}(k-m)}\, {\widehat{R}_{j}(m)}\,dm dk 
\\[2mm]
& \qquad \quad \,\,+ \frac12 \sum_{\ell \in \{\pm1\}} \int_{\R} \int_{\R} \overline{\widehat{R}_j(k)}\, \big({\widehat{\Xi}_{j\ell}^{l \mu}(-m,-k)} -  \widehat{\Xi}_{j\ell}^{l \mu}(-k+\ell k_0,-k) \big)\,   
\\
& \hspace{4.05cm} \times 
{\widehat{\psi}_{\ell}(k-m)}\, {\widehat{R}_{j}(m)}\,dm dk \,.
\end{align*}
Since 
$\widehat{\xi}_1(-k)= \widehat{\xi}_1(k)$ for all $k  \in \R$, we have by construction of $\widehat{\rho}_{j}^{l}$ and $\widehat{\Xi}_{j \ell}^{l \mu}$:
\begin{align*}
& - \sum_{\mu=1}^{|j|}\, \widehat{\Xi}_{j\ell}^{l \mu}(k,k-\ell k_0) + \widehat{\Xi}_{j\ell}^{l \mu}(-k+\ell k_0,-k) \\[-1mm]
& \qquad =\, \widehat{q}_{jj}(k, \ell k_0, k - \ell k_0)\, k^{2l}\, \widehat{\xi}_1\Big(\frac{k-\ell k_1}{k_1-k_0}\Big)  \\[1mm]
& \qquad \quad\, -\widehat{q}_{jj}(-k+ \ell k_0, \ell k_0, -k)\, (-k+\ell k_0)^{2l}\,  
\widehat{\xi}_1\Big(\frac{k+ \ell (k_1-k_0)}{k_1-k_0}\Big)
\\[1mm]
& \qquad \quad\, +  \widehat{q}_{jj}(-k +\ell k_0, \ell k_0, -k)\, (-k+\ell k_0)^{2l}\,  
\widehat{\xi}_1\Big(\frac{-k- \ell (k_1-k_0)}{k_1-k_0}\Big)
\\[1mm]
& \qquad \quad\, -  \widehat{q}_{jj}(k, \ell k_0, k- \ell k_0)\, k^{2l}\, \widehat{\xi}_1\Big(\frac{-k+\ell k_1}{k_1-k_0}\Big)  \\[1mm]
& \qquad =\, 0 \,,
\end{align*}
 which yields
\begin{align} \label{KL3}
& \sum_{\mu=1}^{|j|} \sum_{\ell \in \{\pm1\}} \int_{\R} \int_{\R} \overline{\widehat{R}_j(k)}\; \widehat{\Xi}_{j\ell}^ {l\mu}(k,m)\,
\widehat{\psi}_{\ell}(k-m)\, \widehat{R}_{j}(m)\,dm dk  \nonumber \\[1mm] 
& \qquad \quad =\, \frac12 \sum_{\mu=1}^{|j|} \sum_{\ell \in \{\pm1\}} \int_{\R} \int_{\R} \overline{\widehat{R}_j(k)}\, \big({\widehat{\Xi}_{j\ell}^{l \mu}(k,m)} -  \widehat{\Xi}_{j\ell}^{l \mu}(k,k-\ell k_0) \big)\,   
\nonumber \\
& \hspace{4.05cm} \times 
{\widehat{\psi}_{\ell}(k-m)}\, {\widehat{R}_{j}(m)}\,dm dk 
 \nonumber 
 \\[3mm]
&  \qquad \quad \quad \,\,+ \frac12 \sum_{\mu=1}^{|j|} \sum_{\ell \in \{\pm1\}} \int_{\R} \int_{\R} \overline{\widehat{R}_j(k)}\, \big({\widehat{\Xi}_{j\ell}^{l \mu}(-m,-k)} -  \widehat{\Xi}_{j\ell}^{l \mu}(-k+\ell k_0,-k) \big)\,   
\nonumber \\ 
& \hspace{4.4cm} \times 
{\widehat{\psi}_{\ell}(k-m)}\, {\widehat{R}_{j}(m)}\,dm dk \,.
\end{align}
Because of \eqref{pm1-ans}, \eqref{pm1-ans2}, \eqref{asdkomega}, \eqref{asqjj}--\eqref{asfq2} and \eqref{asnjj}--\eqref{asn2j-j} the functions 
\begin{align*}
& (k,k-m,m) \mapsto (1+k^{2l})^{-1}\,{\widehat{\Xi}_{j\ell}^{l \mu}}(k,m)\,, \\[1mm]
& (k,k-m,m) \mapsto (1+k^{2l})^{-1}\,{\widehat{\Xi}_{j\ell}^{l \mu}}(-m,-k)\,,\\[1mm]
& (k,k-m,m) \mapsto (1+|k|)^{-(\mu+1)/2}\, \widehat{\vartheta}^{-1}(k)\,{K}_{j_1 j_2 \ell}^{\mu}(k,k-m,m)   
\end{align*}
have the properties of $\mathfrak{K}$ in Lemma \ref{int-kerne} a) and b), respectively, such that we can use \eqref{RES3}, \eqref{dxR1R2}--\eqref{dx-2R2}, \eqref{IK1}--\eqref{IK2},   \eqref{KL1}--\eqref{KL3} and Young's inequality for convolutions 
to obtain \eqref{nf}--\eqref{nf-rest} uniformly on compact subsets of $\mathcal{B}$. 
$\hfill \Box$

\begin{lemma}
\label{lem44}
The operators $\mathcal{N}_{j_1j_2}^{i}$ have the following properties:\\[2mm]
{\bf a)} Fix $\varphi \in L^2(\R,\R)$ with $\mathrm{supp}\, \widehat{\varphi} = \mathrm{supp}\, \widehat{\psi}_c$. Then $f \mapsto \mathcal{N}_{jj}^i(\varphi,f)$ defines a continuous linear map from $H^1(\R,\R)$ into $L^2(\R,\R)$ and $f \mapsto \mathcal{N}_{j-j}^{i}(\varphi,f)$ a continuous linear map from $H^{(1-(|j|-1))/2}(\R,\R)$ into $L^2(\R,\R)$. Furthermore, 
 there exists a constant $C>0$ with $C \lesssim \Vert \widehat{\varphi} \Vert_{L^{1}}$ such that for all $f\in H^{(1-(|j|-1))/2}(\R,\R)$ and all $g\in H^{1}(\R,\R)$ there holds
\begin{align}
\label{cN0}
\Vert \mathcal{N}_{jj}^i(\varphi, g) \Vert_{L^{2}} & \leq\, C \varepsilon^{-1} \Vert g \Vert_{H^1}\,, \\[2mm]
\label{cN0b}
\Vert \mathcal{N}_{j-j}^i(\varphi, f) \Vert_{L^{2}} & \leq \, C \varepsilon^{-1} \Vert f \Vert_{H^{(1-(|j|-1))/2}}\,,
\\[2mm]
\label{PcN}
\Vert P_{\delta_{0}, \infty} \mathcal{N}_{jj}^i(\varphi, g) \Vert_{L^{2}} & \leq \, C \Vert g \Vert_{H^1}\,,
\\[2mm]
\label{PcNb}
\Vert P_{\delta_{0}, \infty} \mathcal{N}_{j-j}^i(\varphi, f) \Vert_{L^{2}} & \leq \, C  \Vert f \Vert_{H^{(1-(|j|-1))/2}}\,,
\end{align}
uniformly on compact subsets of $\mathcal{B}$.
\\[2mm]
{\bf b)} Let $\varphi$ be as in a). 
Then 
for all $f \in L^2(\R,\R)$ there holds
\begin{align}
\label{P_0 cN}
P_{0, \delta_{0}} \mathcal{N}^{i}(\varphi, P_{0, \delta_{0}}f)=0\,.
\end{align}
\end{lemma}
\textbf{Proof.}
a) The first step of the proof is to analyze the behavior of $\widehat{n}_{j_1 j_2}^{i,0}$ for all $i \in \{1,2\}$ in a neighborhood of the zeros of the factor $\widehat{r}_{j_1 j_2}$ in the denominator. Due to the localization of the supports of $\widehat{P}_{0,\delta_{0}}$ and  
$\widehat{\varphi}$, it is sufficient to consider only the zeros satisfying $|k| \leq \delta_{0}$ and $|k-m \mp k_0|  \leq \delta_{0}$. As shown above, the only zeros satisfying $|k| \leq \delta_{0}$ and $|k-m \mp k_0|  \leq \delta_{0}$ are $(0,\pm k_0, \mp k_0)$, which appear if 
$\mathrm{sgn}(j_2) = -1$. By the same arguments as those in the proof of Lemma \ref{lem31} it follows that the singularities of $\widehat{n}^{i,0}_{j_1 j_2}$ at $(0,\pm k_0, \mp k_0)$ can be removed and then there holds
\begin{align}
\label{cnz0}
\widehat{P}_{0, \delta_{0}}(k)\, \widehat{P}_{0, \delta_{0}}(k-m \mp k_0)\, \widehat{n}_{j_1 j_2}^{i,0}(k,k-m,m) = \mathcal{O}(\veps^{-1})
\end{align}
uniformly on compact subsets of $\mathcal{B}$. Hence, because of \eqref{thetaN}--\eqref{thetaNb} all assertions of a) are valid.
\\[2mm]
b) is again a direct consequence of
\begin{align*}
\widehat{P}_{0,\delta_{0}}\!(k)\, \widehat{P}_{0,\delta_{0}}\!(m)\, \chi_c(k-m)=0\,.
\end{align*}
$\hfill \Box$

\begin{lemma}
\label{lem44b}
Let $R_{-1}, R_1 \in L^2(\R,\R)$ and $R_{- 2}, R_2 \in H^2(\R,\R)$. Then there holds
\begin{align} \label{cnf}
& \int_{\R} R_{j_1}\,  \rho_{j_1}^{0}  \big(\mathrm{sgn}(j_1) i \omega \mathcal{N}_{j_1 j_2}(\psi_c,R_{j_2})+ \mathcal{N}_{j_1 j_2}(i \omega \psi_c,R_{j_2})
\nonumber \\[1mm] 
& \quad
- \mathrm{sgn}(j_2) \mathcal{N}_{j_1 j_2}(\psi_c,i \omega R_{j_2})\big) \, d\alpha 
\nonumber \\[2mm] 
& \qquad\quad  =\, \int_{\R} R_{j_1}\, {\rho}_{j_1}^{0} Q_{j_1 j_2}(\psi_c)R_{j_2}\, d\alpha 
 + \eps \int_{\R} R_{j_1}\, {\rho}_{j_1}^{0} \mathcal{Y}_{j_1 j_2}(\psi_c,R_{j_2})\, d\alpha \,
\end{align}
with
\begin{align} \label{cnf-rest}
& \Big| \sum_{\genfrac{}{}{0pt}{}{j_1 \in \{\pm 1, \pm 2\},} {j_2 \in \{\pm j_1\}}}                    \int_{\R} R_{j_1}\, {\rho}_{j_1}^{0} \mathcal{Y}_{j_1 j_2}(\psi_c,R_{j_2})\, d\alpha \Big| 
\,  \lesssim\,
\| \mathcal{R}_{1}\|_{(L^{2})^{2}}^2 + \| \mathcal{R}_{2}\|_{(H^{2})^{2}}^2 + \eps^{4} \,,
\end{align}
as long as $\eps^{5/2} \Vert \mathcal{R}_{1} \Vert_{(L^{2})^2},\, \eps^{5/2} \Vert \mathcal{R}_{2} \Vert_{(H^{2})^2} \ll 1$, uniformly on compact subsets of $\mathcal{B}$.
\end{lemma}
\textbf{Proof.}
Let
\begin{align*} 
\eps \mathcal{Y}_{j_1 j_2}(\psi_c,R_{j_2}) :=
&\;\, \mathrm{sgn}(j_1) i \omega \mathcal{N}_{j_1 j_2}(\psi_c,R_{j_2})+ \mathcal{N}_{j_1 j_2}(i \omega \psi_c,R_{j_2})- \mathrm{sgn}(j_2) \mathcal{N}_{j_1 j_2}(\psi_c,i \omega R_{j_2}) 
\nonumber \\[1mm] 
&\; - Q_{j_1 j_2}(\psi_c)R_{j_2} \,. 
\end{align*}
Then we have
\begin{align}
\widehat{\mathcal{Y}}_{j_1 j_2}(\psi_c,R_{j_2})(k) \,
& =\, \eps^{-1} \sum_{\mu=1}^{|j_1|} \sum_{\ell \in \{\pm1\}} \int_{\R} 
\widehat{P}_{\delta_{0},\infty}(k)\,K_{j_1 j_2 \ell}^{\mu} (k,k-m,m)\,
\widehat{\psi}_{\ell}(k-m)  \nonumber \\
& \hspace{3.3cm} \times (im)^{-(\mu -1)}\,\widehat{R}_{j_2}(m)\,dm\,,
\end{align}
where $K_{j_1 j_2 \ell}^{\mu}$ is defined by \eqref{KLdef}. 
Hence, the assertion of the lemma follows by the same arguments as those from the proof of Lemma \ref{lem31b}. 
$\hfill \Box$
\\[5mm]
Let  $i \in \{1,2\}$, $l \geq 1$, $\varphi \in L^2(\R,\R)$ with $\mathrm{supp}\, \widehat{\varphi} = \mathrm{supp}\, \widehat{\psi}_c$, $f:= (f_{-2}, f_2)^T \in (H^1(\R,\R))^2$,
\begin{equation*}
(\rho^l \mathfrak{N}^{i})(\varphi) f  = 
\left(
\begin{array}{cc} (\rho^l \mathfrak{N}^{i})_{-2-2}(\varphi) & (\rho^l \mathfrak{N}^{i})_{-22}(\varphi)  \\ (\rho^l \mathfrak{N}^{i})_{2-2}(\varphi) & (\rho^l \mathfrak{N}^{i})_{22}(\varphi)
\end{array}
\right)  \left(
\begin{array}{c} f_{-2} \\ f_{2}
\end{array}
\right)
\end{equation*}
with
\begin{equation*}
 (\rho^l \mathfrak{N}^{i})_{j_1j_2}(\varphi) f_{j_2} \,=\,  \rho_{j_1}^l \vartheta N_{j_1j_2}^{i}(\varphi) f_{j_2}\, =\,  \rho_{j_1}^l \vartheta N_{j_1j_2}^{i}(\varphi, f_{j_2})
\end{equation*}
for $j_1,j_2  \in \{ \pm 2 \}$. Then there holds
\begin{align*}
((\widehat{\rho^l \mathfrak{N}^{i}})_{j_1 j_2}(\varphi)f_{j_2})(k) &\,=\, \sum_{\ell \in \{\pm 1\}} \int_{\R} (\widehat{\rho^l \mathfrak{n}^{i}})_{j_1 j_2 \ell}(k)\,\widehat{\varphi}_{\ell}(k-m)\, \widehat{f}_{j_2}(m)\,dm 
\end{align*}
with
\begin{equation*}
(\widehat{\rho^l \mathfrak{n}^{i}})_{j_1 j_2 \ell}(k) \,=\, \widehat{\rho}_{j_1}^l(k)\,\widehat{\vartheta}(k) \,\widehat{n}^{1}_{j_1 j_2 \ell}(k) \,.
\end{equation*}
Moreover, let
\begin{equation*}
(\rho^l \mathfrak{N}^{i})^{\ast}(\varphi) = \left(
\begin{array}{cc} (\rho^l \mathfrak{N}^{i})^{\ast}_{-2-2}(\varphi) & (\rho^l \mathfrak{N}^{i})^{\ast}_{-22}(\varphi)  \\ (\rho^l \mathfrak{N}^{i})^{\ast}_{2-2}(\varphi) & (\rho^l \mathfrak{N}^{i})^{\ast}_{22}(\varphi)
\end{array}
\right) 
\end{equation*}
be the adjoint operator of $(\rho^l \mathfrak{N}^{i})(\varphi)$. That means, we have 
\begin{equation*}
\langle\, (\rho^l \mathfrak{N}^{i})(\varphi)f,\, g \,\rangle_{(L^2)^2} \,=\,
\langle\, f,\, (\rho^l \mathfrak{N}^{i})^{\ast}(\varphi)\, g\, \rangle_{(L^2)^2} 
\end{equation*}
for all $f,g \in (H^1(\R,\R))^2$,
as well as
\begin{align*}
(\rho^l \mathfrak{N}^{i})^{\ast}(\varphi) &=\, \left(
\begin{array}{cc} ( \rho_{-2}^l \vartheta N_{-2-2}^{i})^{\ast}(\varphi) & ( \rho_{-2}^l \vartheta N_{-22}^{i})^{\ast}(\varphi)  \\ ( \rho_{2}^l \vartheta N_{2-2}^{i})^{\ast}(\varphi) & ( \rho_{2}^l \vartheta N_{22}^{i})^{\ast}(\varphi)
\end{array}
\right)^T \\[2mm]
&=\,   \left(
\begin{array}{cc} ( \rho_{-2}^l \vartheta N_{-2-2}^{i})^{\ast}(\varphi) & ( \rho_{2}^l \vartheta N_{2-2}^{i})^{\ast}(\varphi)  \\ ( \rho_{-2}^l \vartheta N_{-22}^{i})^{\ast}(\varphi) & ( \rho_{2}^l \vartheta N_{22}^{i})^{\ast}(\varphi)
\end{array}
\right)
\,,
\end{align*}
where
\begin{equation*}
\langle\,  \rho_{j_1}^l \vartheta N_{j_1j_2}^{i}(\varphi)h,\, p\, \rangle_{L^2} \,=\,
\langle\, h,\, ( \rho_{j_1}^l \vartheta N_{j_1j_2}^{i})^{\ast}(\varphi)\,p\, \rangle_{L^2} 
\end{equation*}
for all $h,p \in H^1(\R,\R)$. Then
\begin{equation*}
(\rho^l \mathfrak{N}^{i})^{s}(\varphi) = \left(
\begin{array}{cc} (\rho^l \mathfrak{N}^{i})^{s}_{-2-2}(\varphi) & (\rho^l \mathfrak{N}^{i})^{s}_{-22}(\varphi)  \\ (\rho^l \mathfrak{N}^{i})^{s}_{2-2}(\varphi) & (\rho^l \mathfrak{N}^{i})^{s}_{22}(\varphi)
\end{array}
\right) 
 := \frac12 ((\rho^l \mathfrak{N}^{i})(\varphi) + (\rho^l \mathfrak{N}^{i})^{\ast}(\varphi))
\end{equation*}
denotes the symmetric part and
\begin{equation*}
(\rho^l \mathfrak{N}^{i})^{a}(\varphi) = \left(
\begin{array}{cc} (\rho^l \mathfrak{N}^{i})^{a}_{-2-2}(\varphi) & (\rho^l \mathfrak{N}^{i})^{a}_{-22}(\varphi)  \\ (\rho^l \mathfrak{N}^{i})^{a}_{2-2}(\varphi) & (\rho^l \mathfrak{N}^{i})^{a}_{22}(\varphi)
\end{array}
\right) 
:= \frac12 ((\rho^l \mathfrak{N}^{i})(\varphi) - (\rho^l \mathfrak{N}^{i})^{\ast}(\varphi))
\end{equation*}
the antisymmetric part of $(\rho^l \mathfrak{N}^{i})(\varphi)$. There holds
\begin{align*}
((\widehat{\rho^l \mathfrak{N}^{i}})^{\ast}_{j_1 j_2}(\varphi)f_{j_2})(k) &\,=\, \sum_{\ell \in \{\pm 1\}} \int_{\R} (\widehat{\rho^l \mathfrak{n}^{i}})^{\ast}_{j_1 j_2 \ell}(m)\,\widehat{\varphi}_{\ell}(k-m)\, \widehat{f}_{j_2}(m)\,dm
\end{align*}
with
\begin{equation*}
(\widehat{\rho^l \mathfrak{n}^{i}})^{\ast}_{j_1 j_2 \ell}(m) \,=\, (\widehat{\rho^l \mathfrak{n}^{i}})_{j_2 j_1 \ell}(-m) \,=\,
\widehat{\rho}_{j_2}^l(-m)\,\widehat{\vartheta}(-m) \,\widehat{n}^{i}_{j_2 j_1 \ell}(-m)
\end{equation*}
and therefore
\begin{align*}
((\widehat{\rho^l \mathfrak{N}^{i}})^{s}_{j_1 j_2 }(\varphi)f_{j_2})(k) &\,=\, \sum_{\ell \in \{\pm 1\}} \int_{\R} (\widehat{\rho^l \mathfrak{n}^{i}})^{s}_{j_1 j_2 \ell}(k,m)\,\widehat{\varphi}_{\ell}(k-m)\,\widehat{f}_{j_2}(m)\,dm\,
\end{align*}
with
\begin{align*}
(\widehat{\rho^l \mathfrak{n}^{i}})_{j_1 j_2 \ell}^{s}(k,m) &\,=\, \frac12 \big(  (\widehat{\rho^l \mathfrak{n}^{i}})_{j_1 j_2 \ell}(k) + (\widehat{\rho^l \mathfrak{n}^{i}})_{j_2 j_1 \ell}(-m)  \big) 
\end{align*}
as well as
\begin{align*}
((\widehat{\rho^l \mathfrak{N}^{i}})^{a}_{j_1 j_2 }(\varphi)f_{j_2})(k) &\,=\, \sum_{\ell \in \{\pm 1\}} \int_{\R} (\widehat{\rho^l \mathfrak{n}^{i}})^{a}_{j_1 j_2 \ell}(k,m)\,\widehat{\varphi}_{\ell}(k-m)\,\widehat{f}_{j_2}(m)\,dm\,
\end{align*}
with
\begin{align*}
(\widehat{\rho^l \mathfrak{n}^{i}})_{j_1 j_2 \ell}^{a}(k,m) &\,=\, \frac12 \big(  (\widehat{\rho^l \mathfrak{n}^{i}})_{j_1 j_2 \ell}(k) - (\widehat{\rho^l \mathfrak{n}^{i}})_{j_2 j_1 \ell}(-m)  \big) \,.
\end{align*}
More generally, for any densely defined linear operator $L$ on a Hilbert space we denote its adjoint operator by 
$L^{\ast}$, its symmetric part by $L^{s}$ and its antisymmetric part by $L^{a}$.

\begin{lemma}
\label{lem3sym}
The operators $(\rho^l \mathfrak{N}^{1})^{s}_{j_1j_2}$ with $j_1, j_2  \in \{\pm 2\}$ and $l \geq 1$ have the following properties:\\[2mm]
{\bf a)} Fix $\varphi \in L^2(\R,\R)$ with $\mathrm{supp}\, \widehat{\varphi} = \mathrm{supp}\, \widehat{\psi}_c$. Then $f \mapsto (\rho^l \mathfrak{N}^{1})_{j_1j_2}^s(\varphi)f$ can be extended to a continuous linear map from  $L^2(\R,\R)$ into $L^2(\R,\R)$  and
 there exists a constant $C_1>0$ with $C_1 \lesssim \Vert \widehat{\varphi} \Vert_{L^{1}}$ such that for all $f\in L^2(\R,\R)$ there holds
\begin{align}
\label{N0s}
\Vert (\rho^l \mathfrak{N}^{1})_{j_1j_2}^s(\varphi) f \Vert_{L^{2}} & \leq \, C_1\, \Vert f \Vert_{L^2}
\end{align}
uniformly on compact subsets of $\mathcal{B}$. Moreover, for all $\psi \in H^{2+p}(\R,\R)$ with $p >1/2$ there exists a constant ${C}_2>0$ with ${C}_2 \lesssim \Vert \widehat{\varphi} \Vert_{L^{1}} \Vert \partial_{\alpha} \psi \Vert_{H^{1+p}} $ such that for all $f\in L^2(\R,\R)$ there holds
\begin{align}
\label{N0skomm}
\Vert [\psi, (\rho^l \mathfrak{N}^{1})_{j_1j_2}^s(\varphi)] f \Vert_{H^{1}} & \leq \, {C}_2\, \Vert f \Vert_{L^2}
\end{align}
uniformly on compact subsets of $\mathcal{B}$. 
\\[2mm]
{\bf b)} 
For all $\varphi \in L^2(\R,\R)$ with $\mathrm{supp}\, \widehat{\varphi} = \mathrm{supp}\, \widehat{\psi}_c$ there exists a symmetric linear operator $G_j^l(\varphi): L^2(\R,\R) \to L^2(\R,\R)$ such that
\begin{equation}
\label{Gj}
(\rho^l \mathfrak{N}^{1})_{j-j}^s(\varphi)f = G_{j}^l(\varphi)f + M_{j}^l(\varphi)f   
\end{equation}
and there exists a constant $C_3>0$ with $C_3 \lesssim \Vert \widehat{\varphi} \Vert_{L^{1}}$ such that for all $f\in L^2(\R,\R)$ there holds
\begin{equation}
\label{Mj}
\|M_{j}^l(\varphi)f\|_{H^{1/2}}  \leq C_3\,  \|f\|_{L^2}
\end{equation}
uniformly on compact subsets of $\mathcal{B}$. Moreover, for all $\psi \in H^{2+p}(\R,\R)$ with $p >1/2$ there exists a constant ${C}_4>0$ with ${C}_4 \lesssim \Vert \widehat{\varphi} \Vert_{L^{1}} \Vert \partial_{\alpha} \psi \Vert_{H^{1+p}} $ such that for all $f\in L^2(\R,\R)$ there holds
\begin{align}
\label{Gjkomm}
\Vert [\psi, G_{j}^l(\varphi)] f \Vert_{H^{1}} & \leq \, {C}_4\, \Vert f \Vert_{L^2}
\end{align}
uniformly on compact subsets of $\mathcal{B}$. 
\end{lemma}
{\bf Proof.}
a) Let $j \in \{\pm 2\}$. We have
\begin{align} \label{p1qjj}
\partial_1\,\sum_{\mu=1}^{3} \widehat{q}_{jj}^{2,\mu}(k,\ell k_0,k-\ell k_0) \,&=\, -i + \mathcal{O}(|k|^{-1/2}) \,,
\\
\partial_1\,\sum_{\mu=1}^{3} \widehat{q}_{j-j}^{2,\mu}(k,\ell k_0,k-\ell k_0) \,&=\, \mathcal{O}(|k|^{-1/2}) \,,
\\
\partial_1^{2}\,\sum_{\mu=1}^{3} \widehat{q}_{j_1j_2}^{2,\mu}(k,\ell k_0,k-\ell k_0) \,&=\, 0 \,,\\ \label{p3qj1j2}
\partial_3\,\sum_{\mu=1}^{3} \widehat{q}_{j_1j_2}^{2,\mu}(k,\ell k_0,k-\ell k_0) \,&=\, \mathcal{O}(|k|^{-1/2}) \,,\\
\partial_3^{2}\,\sum_{\mu=1}^{3} \widehat{q}_{j_1j_2}^{2,\mu}(k,\ell k_0,k-\ell k_0) \,&=\, \mathcal{O}(|k|^{-3/2}) \,,\\ \label{p1p3q}
\partial_1 \partial_3\,\sum_{\mu=1}^{3} \widehat{q}_{j_1j_2}^{2,\mu}(k,\ell k_0,k-\ell k_0) \,&=\, \partial_3 \partial_1\, \,\sum_{\mu=1}^{3} \widehat{q}_{j_1j_2}^{2,\mu}(k,\ell k_0,k-\ell k_0) \,=\,\mathcal{O}(|k|^{-3/2}) 
\end{align}
for $|k| \rightarrow \infty$ uniformly with respect to $b \lesssim 1$. 

Because of \eqref{q5}--\eqref{frakm}, \eqref{oddq}--\eqref{odd-jq}, \eqref{k1}, \eqref{asqjj}--\eqref{asqj-j} and
\eqref{p1qjj}--\eqref{p1p3q} we obtain by
construction of $\widehat{\rho}_{j_1}^l$ and $\widehat{\zeta}_{j_1 j_2 \ell}$ for all $j \in \{\pm 2\}$, $\beta \in \{1,2\}$, $\gamma \in \{0,1\}$ and $n \in \{1,3\}$:
\begin{align} \label{dkweight}
 \Big (\frac{d^{\beta}}{dk^{\beta}}\,(\widehat{\rho}_j^l\, \widehat{\zeta}_{j j \ell}) \Big)(-k)\, \partial_n^{\gamma} \widehat{\mathfrak{q}}^{1}_{jj}(-k,\ell k_0,-k + \ell k_0)  \,&=\, \mathcal{O}(|k|^{-(\beta+\gamma-1)})
\,, \\[2mm]
 \Big (\frac{d^{\beta}}{dk^{\beta}}\,(\widehat{\rho}_{-j}^l\, \widehat{\zeta}_{-j j \ell})\Big)(-k)\, \partial_n^{\gamma} \widehat{\mathfrak{q}}^{1}_{-jj}(-k,\ell k_0,-k + \ell k_0)  \,&=\, \mathcal{O}(|k|^{-(\beta+\gamma-1/2)})
\end{align}
for $|k| \rightarrow \infty$ uniformly with respect to $b \lesssim 1$. 

Moreover, since $\omega$ is odd, we have
\begin{align} \label{rswitch}
\widehat{r}_{j_2 j_1}(-k,\ell k_0,-k -\ell k_0) 
&\,=\, \widehat{r}_{j_1 j_2}(k + \ell k_0 ,\ell k_0,k) \,.
\end{align}

Let $\widehat{r}_{j_1 j_2 \ell}: \R \to \C$ defined by $\widehat{r}_{j_1 j_2 \ell}(k) = \widehat{r}_{j_1 j_2}(k,\ell k_0,k - \ell k_0)$. Then there holds for all $\beta \in \{1,2\}$:
\begin{align} \label{dkrjjell}
\frac{d^{\beta}}{dk^{\beta}}\, \widehat{r}_{j j \ell}(k)\,&=\, \mathcal{O}(|k|^{-(\beta+1/2)}(1+bk^{2})^{1/2})\,,\\[2mm] \label{dkrj-jell}
\frac{d^{\beta}}{dk^{\beta}}\, \widehat{r}_{j -j \ell}(k)\,&=\, \mathcal{O}(|k|^{-(\beta-1/2)}(1+bk^{2})^{1/2})
\end{align}
for $|k| \rightarrow \infty$ uniformly with respect to $b$ in compact subsets of $\mathcal{B}$. 

With the help of Taylor's theorem as well as \eqref{q5}--\eqref{frakm}, \eqref{oddq}--\eqref{odd-jq}, \eqref{asqjj}--\eqref{asqj-j}, \eqref{asrjj}--\eqref{asrj-j} and \eqref{p1qjj}--\eqref{dkrj-jell} we derive
\begin{align} \label{nuspm}
&(\widehat{\rho^l \mathfrak{n}^{1}})_{j_1 j_2\ell}^s(k,m)\, \chi_{\ell}(k-m) 
\nonumber \\[2mm]
& \qquad  =\,
(\widehat{\rho^l \mathfrak{n}^{1}})_{j_1 j_2\ell}^{s,1}(k) + (\widehat{\rho^l \mathfrak{n}^{1}})_{j_1 j_2\ell}^{s,2}(k)\,(k-m) 
+ (\widehat{\rho^l \mathfrak{n}^{1}})_{j_1 j_2\ell}^{s,3}(k)\,(k-m) 
\nonumber \\[2mm]
& \qquad \quad\;\,
+  (\widehat{\rho^l \mathfrak{n}^{1}})_{j_1 j_2\ell}^{s,4}(k)\,
(k-m-2 \ell k_0) + (\widehat{\rho^l \mathfrak{n}^{1}})_{j_1 j_2\ell}^{s,5}(k)\,
(k-m-\ell k_0) 
\nonumber 
\\[2mm]
& 
\qquad \quad\;\,
+\delta_{j_1j_2}\, \mathcal{O}(|k|^{-1}(1+|k|^{-1/2}(1+bk^{2})^{1/2})^{-1})
\nonumber \\[2mm]
& 
\qquad \quad\;\,
+\delta_{j_1-j_2}\, \mathcal{O}(|k|^{-3/2}(1+bk^{2})^{-1/2})
\end{align}
for $|k| \to \infty$ uniformly with respect to $ m \in \R$ and $b$ in compact subsets of $\mathcal{B}$, where
\begin{align} 
(\widehat{\rho^l \mathfrak{n}^{1}})_{j_1 j_2\ell}^{s,1}(k) \,& =\,
\frac{\widehat{\rho}_{j_1}^l(k)\,\widehat{\zeta}_{j_1 j_2 \ell}(k)\,\widehat{\mathfrak{q}}^{1}_{j_1 j_2}(k,\ell k_0,k -\ell k_0)}{2 \widehat{r}_{j_1 j_2 \ell}(k)}
\nonumber \\[2mm]
& \quad \,\, + \frac{\widehat{\rho}_{j_2}^l(-k)\,\widehat{\zeta}_{j_2 j_1 \ell}(-k)\,\widehat{\mathfrak{q}}^{1}_{j_2 j_1}(-k,\ell k_0,-k +\ell k_0)}{2 \widehat{r}_{j_1 j_2 \ell}(k)} \,,
\\[2mm]
(\widehat{\rho^l \mathfrak{n}^{1}})_{j_1 j_2\ell}^{s,2}(k) \,& =\,
- \frac{\widehat{\rho}_{j_2}^l(-k)\,\widehat{\zeta}_{j_2 j_1 \ell}(-k)\, \partial_1 \widehat{\mathfrak{q}}^{1}_{j_2 j_1}(-k,\ell k_0,-k +\ell k_0)}{2 \widehat{r}_{j_1 j_2 \ell}(k)} \,,
\\[2mm]
(\widehat{\rho^l \mathfrak{n}^{1}})_{j_1 j_2\ell}^{s,3}(k) \,& =\,
- \frac{ \big(\textstyle \frac{d}{dk}\,(\widehat{\rho}_{j_2}^l\, \widehat{\zeta}_{j_2 j_1 \ell})\big)(-k)\, \widehat{\mathfrak{q}}^{1}_{j_2 j_1}(-k,\ell k_0,-k +\ell k_0)}{2 \widehat{r}_{j_1 j_2 \ell}(k)} \,,
\\[2mm]
(\widehat{\rho^l \mathfrak{n}^{1}})_{j_1 j_2\ell}^{s,4}(k) \,& =\,
- \frac{\widehat{\rho}_{j_2}^l(-k)\,\widehat{\zeta}_{j_2 j_1 \ell}(-k)\,\partial_3 \widehat{\mathfrak{q}}^{1}_{j_2 j_1}(-k,\ell k_0,-k +\ell k_0)}{2 \widehat{r}_{j_1 j_2 \ell}(k)} \,,
\\[2mm] \label{nuspm2}
(\widehat{\rho^l \mathfrak{n}^{1}})_{j_1 j_2\ell}^{s,5}(k) \,& =\,
\frac{\widehat{\rho}_{j_2}^l(-k)\,\widehat{\zeta}_{j_2 j_1 \ell}(-k)\,\widehat{\mathfrak{q}}^{1}_{j_2 j_1}(-k,\ell k_0,-k +\ell k_0)}{2 \widehat{r}_{j_1 j_2 \ell}(k)}\,
 \frac{\textstyle \frac{d}{dk}\, \widehat{r}_{j_1 j_2 \ell}(k)}{ \widehat{r}_{j_1 j_2 \ell}(k)} \,,
\end{align}
and consequently
\begin{align} 
\label{nus}
(\widehat{\rho^l \mathfrak{n}^{1}})_{j j\ell}^s(k,m)\, \chi_{\ell}(k-m) &\,=\, \mathcal{O}((1+|k|^{-1/2}(1+bk^2)^{1/2})^{-1})\,, \\[2mm]
\label{nus-}
(\widehat{\rho^l \mathfrak{n}^{1}})_{j -j \ell}^s(k,m)\, \chi_{\ell}(k-m) &\,=\, \mathcal{O}((1+bk^2)^{-1/2})
\end{align}
for $|k| \to \infty$ uniformly with respect to $ m \in \R$ and $b$ in compact subsets of $\mathcal{B}$, which implies \eqref{N0s}.

To prove the second assertion of a) we consider
\begin{align} \label{fcomm}
& \mathcal{F}([\psi, (\rho^l \mathfrak{N}^{1})_{j_1j_2}^s(\varphi)] f)(k) \nonumber \\[2mm]
& \qquad  =\, \sum_{\ell \in \{\pm 1\}} \Big( \int_{\R} \int_{\R} \widehat{\psi}(k-m)\, (\widehat{\rho^l \mathfrak{n}^{1}})_{j_1 j_2\ell}^s(m,n)\, \widehat{\varphi}_{\ell}(m-n)\, \widehat{f}(n)\,dndm
\nonumber \\[2mm]
& \qquad \qquad  \qquad\,\;\,\, - \int_{\R} \int_{\R} (\widehat{\rho^l \mathfrak{n}^{1}})_{j_1 j_2\ell}^s(k,m)\, \widehat{\varphi}_{\ell}(k-m)\, \widehat{\psi}(m-n)\, \widehat{f}(n)\,dndm \Big).
\end{align}
Expanding the kernels $(\widehat{\rho^l \mathfrak{n}^{1}})_{j_1 j_2\ell}^s(m,n)$ and $(\widehat{\rho^l \mathfrak{n}^{1}})_{j_1 j_2\ell}^s(k,m)$ by \eqref{nuspm}--\eqref{nuspm2}, 
rewriting the factors $(\widehat{\rho^l \mathfrak{n}^{1}})_{j_1 j_2\ell}^{s,i}(m)$ for $i= 1, \ldots, 5$  as 
\begin{align*} 
(\widehat{\rho^l \mathfrak{n}^{1}})_{j_1 j_2\ell}^{s,i}(m) =
(\widehat{\rho^l \mathfrak{n}^{1}})_{j_1 j_2\ell}^{s,i}(k) + \big((\widehat{\rho^l \mathfrak{n}^{1}})_{j_1 j_2\ell}^{s,i}(m)-(\widehat{\rho^l \mathfrak{n}^{1}})_{j_1 j_2\ell}^{s,i}(k)\big) 
\end{align*}
and using 
\eqref{q5}--\eqref{frakm}, \eqref{oddq}--\eqref{odd-jq}, \eqref{asqjj}--\eqref{asqj-j}, \eqref{asrjj}--\eqref{asrj-j}, \eqref{p1qjj}--\eqref{dkrj-jell}, the mean value theorem and the fact that the convolution is a commutative operation
yields 
\eqref{N0skomm}.
\\[2mm]
b) Let $G_j^l(\varphi)$ defined by 
\begin{equation}
{G}_{j}^l(\varphi)f \;=\; \frac12   \big(  (\rho^l \mathfrak{N}^{1})^s_{j-j}(\varphi)f +  (\rho^l \mathfrak{N}^{1})^s_{-jj}(\varphi)f \big)
\end{equation}
for all $f \in L^2(\R,\R)$ and ${M}_{j}^l(\varphi):= (\rho^l \mathfrak{N}^{1})^s_{j-j}(\varphi)-{G}^l_{j}(\varphi)$.
Then, 
$G_j^l(\varphi)$ is symmetric. Moreover, we have
\begin{equation} \label{Mjl}
{M}_{j}^l(\varphi)f \;=\; \frac12   \big(  (\rho^l \mathfrak{N}^{1})^s_{j-j}(\varphi)f -  (\rho^l \mathfrak{N}^{1})^s_{-jj}(\varphi)f \big) 
\end{equation}
for all $f \in L^2(\R,\R)$. There holds
\begin{align}
& (\widehat{\rho^l \mathfrak{n}^{1}})_{j -j \ell}^{s}(k,m) - (\widehat{\rho^l \mathfrak{n}^{1}})_{-j j \ell}^{s}(k,m)
\nonumber \\[2mm]
& \qquad =\,
\frac12 \big(\widehat{\rho}_{j}^l(k)\,\widehat{\vartheta}(k) \,\widehat{n}^{1}_{j -j \ell}(k) - \widehat{\rho}_{-j}^l(k)\,\widehat{\vartheta}(k) \, \widehat{n}^{1}_{-jj \ell}(k)\big) 
\nonumber \\[2mm]
& \qquad \quad \,\, +
\frac12 \big(\widehat{\rho}_{j}^l(-m)\,\widehat{\vartheta}(-m) \,\widehat{n}^{1}_{j -j \ell}(-m) - \widehat{\rho}_{-j}^l(-m)\,\widehat{\vartheta}(-m) \, \widehat{n}^{1}_{-j j \ell}(-m)\big) 
\end{align}
and
\begin{align}
& 
\widehat{\rho}_{j}^l(k)\,\widehat{\vartheta}(k) \,\widehat{n}^{1}_{j -j \ell}(k) - \widehat{\rho}_{-j}^l(k)\,\widehat{\vartheta}(k) \, \widehat{n}^{1}_{-jj \ell}(k)
\nonumber \\[2mm]
& \qquad =\,
\Big(\widehat{\rho}_{j}^l(k)\,\frac{\widehat{\mathfrak{q}}^{1}_{j -j}(k,\ell k_0,k -\ell k_0)}{\widehat{r}_{j -j}(k,\ell k_0,k -\ell k_0)} - \widehat{\rho}_{-j}^l(k)\,\frac{\widehat{\mathfrak{q}}^{1}_{-j j}(k,\ell k_0,k -\ell k_0)}{\widehat{r}_{-j j}(k,\ell k_0,k -\ell k_0)} \Big)
\nonumber \\[2mm]
& \qquad \quad \,\, \times
 (\widehat{\vartheta}- \eps \widehat{\xi}_0)(k-\ell k_0) \nonumber \\[2mm]
& \qquad =\,
\sum_{i=1}^{3}\, \widehat{\mathfrak{m}}_{j\ell}^{l,i}(k)\, (\widehat{\vartheta}- \eps \widehat{\xi}_0)(k-\ell k_0) \,,
\end{align}
where
\begin{align}
\widehat{\mathfrak{m}}_{j\ell}^{l,1}(k) \,& = \, \widehat{\rho}_{j}^l(k)\,\frac{\widehat{\mathfrak{q}}^{1}_{j -j}(k,\ell k_0,k -\ell k_0) + \widehat{\mathfrak{q}}^{1}_{-j j}(k,\ell k_0,k -\ell k_0) }{\widehat{r}_{j -j}(k,\ell k_0,k -\ell k_0)} \,,
\\[3mm]
\widehat{\mathfrak{m}}_{j\ell}^{l,2}(k) \,& = \, - \widehat{\rho}_{j}^l(k)\,\frac{\widehat{r}_{j -j}(k,\ell k_0,k -\ell k_0) + \widehat{r}_{-j j}(k,\ell k_0,k -\ell k_0)}{\widehat{r}_{j -j}(k,\ell k_0,k -\ell k_0) \,\widehat{r}_{-j j}(k,\ell k_0,k -\ell k_0)   } \; \widehat{\mathfrak{q}}^{1}_{-j j}(k,\ell k_0,k -\ell k_0)\,,
\\[3mm] \label{frakm3}
\widehat{\mathfrak{m}}_{j\ell}^{l,3}(k) \,& = \, \big(\widehat{\rho}_{j}^l(k)
- \widehat{\rho}_{-j}^l(k) \big) \, \frac{\widehat{\mathfrak{q}}^{1}_{-j j}(k,\ell k_0,k -\ell k_0)}{\widehat{r}_{-j j}(k,\ell k_0,k -\ell k_0)} \,.
\end{align}
Using \eqref{q5}, \eqref{odd-jq} and \eqref{asrj-j} we obtain
\begin{align} \label{frakm1est}
\widehat{\mathfrak{m}}_{j\ell}^{l,1}(k) = \mathcal{O}(|k|^{-1/2}(1+ bk^2)^{-1/2})
\end{align}
for $|k| \rightarrow \infty$ uniformly on compact subsets of $\mathcal{B}$. Because of
\begin{align}
\widehat{r}_{j -j}(k,\ell k_0,k -\ell k_0) + \widehat{r}_{-j j}(k,\ell k_0,k -\ell k_0) = 2 i \,\omega(\ell k_0) \,,
\end{align}
\eqref{asqj-j} and \eqref{asrj-j} we conclude
\begin{align} \label{frakm2est}
\widehat{\mathfrak{m}}_{j\ell}^{l,2}(k) = \mathcal{O}(|k|^{-1/2}(1+ bk^2)^{-1})
\end{align}
for $|k| \rightarrow \infty$ uniformly on compact subsets of $\mathcal{B}$. Furthermore, we have
\begin{align}
 \widehat{\rho}_{j}^l(k)
- \widehat{\rho}_{-j}^l(k) = \begin{cases}  - \mathrm{sgn}(j)
\displaystyle \sum\limits_{\nu \in \{\pm1\}} \widehat{\rho}_{-|j| \nu}^{l}(k) & \quad{\rm if}\; b \in (0,b_0)   \,,\\[3mm]
0 & \quad {\rm otherwise }  \,,
\end{cases}
\end{align}
and
\begin{align}
\widehat{\rho}_{-|j| \nu}^{l}(k) &= 
\left( \, -\frac{\widehat{q}_{-|j| -|j|}(-k+ \nu k_0,\nu k_0,-k) + \widehat{q}_{-|j| -|j|}(k,\nu k_0,k-\nu k_0)}{ \widehat{q}_{-|j| -|j|}(k,\nu k_0,k-\nu k_0)} \right.
\nonumber \\[2mm]
& \left. \qquad  +  \frac{-\widehat{q}_{-|j| -|j|}(-k+ \nu k_0,\nu k_0,-k) }{ \widehat{q}_{-|j| -|j|}(k,\nu k_0,k-\nu k_0)} \, \Big( \Big(1-\frac{\nu k_0}{k} \Big)^{2l} -1 \Big) \right) \nonumber \\[2mm]
& \quad \; \times \widehat{\xi}_1\Big(\frac{k+\nu(k_1-k_0)}{k_1-k_0}\Big)
\nonumber \\[2mm]
& =\, \mathcal{O}(|k|^{-1})
\end{align} 
for $|k| \rightarrow \infty$ uniformly on compact subsets of $\mathcal{B} \cap (0,b_0)$, where the last equality holds because of \eqref{q5}, \eqref{oddq}, \eqref{k1}, \eqref{CTWI}, \eqref{asq4}, \eqref{asqjj}, \eqref{p1qjj}, \eqref{p3qj1j2}, 
\begin{align} \label{p1o3q4}
\partial_{n}\, \widehat{q}_{j_1j_2}^{2,4}(k,\ell k_0,k-\ell k_0) &\,=\,  
\mathcal{O}(|k|^{-3/2}(1+ bk^2)^{1/2})
\end{align}
for $n \in \{1,3\}$ and $|k| \rightarrow \infty$ uniformly with respect to $b \lesssim 1$ and the mean value theorem.
Hence, due to \eqref{asqj-j} and \eqref{asrj-j}, we obtain
\begin{align} \label{frakm3est}
\widehat{\mathfrak{m}}_{j\ell}^{l,3}(k) = \mathcal{O}(|k|^{-1}(1+ bk^2)^{-1/2})
\end{align}
for $|k| \rightarrow \infty$ uniformly on compact subsets of $\mathcal{B}$.

Now, combining \eqref{Mjl}--\eqref{frakm3}, \eqref{frakm1est}, \eqref{frakm2est} and \eqref{frakm3est}, we arrive at \eqref{Mj}. 
 
Finally, \eqref{Gjkomm} is proven in an analogous manner as \eqref{N0skomm}. 
\qed

\begin{lemma}
\label{lem3anti}
The operators $(\rho^l \mathfrak{N}^{i})^{a}_{j_1j_2}$ with $i \in \{1,2\}$, $j_1,j_2  \in \{\pm 2\}$ and $l \geq 1$ have the following properties:\\[2mm] 
{\bf a)} Fix $\varphi \in L^2(\R,\R)$ with $\mathrm{supp}\, \widehat{\varphi} = \mathrm{supp}\, \widehat{\psi}_c$. Then $f \mapsto (\rho^l \mathfrak{N}^{1})^{a}_{jj}(\varphi)f$ defines a continuous linear map from $H^1(\R,\R)$ into $L^2(\R,\R)$ and $f \mapsto (\rho^l \mathfrak{N}^{1})^{a}_{j-j}(\varphi)f$ can be extended to a continuous linear map from $L^2(\R,\R)$ into $L^2(\R,\R)$. Furthermore, 
 there exists a constant $C_1>0$ with $C_1 \lesssim \Vert \widehat{\varphi} \Vert_{L^{1}}$ such that for all $f\in L^2(\R,\R)$ and all $g\in H^1(\R,\R)$ there holds
\begin{align}
\label{thetaNa}
\Vert (\rho^l \mathfrak{N}^{1})^{a}_{jj}(\varphi) g \Vert_{L^{2}} & \leq \, C_1\, \Vert g \Vert_{H^1}\,,\\[2mm]
\label{thetaNab}
\Vert (\rho^l \mathfrak{N}^{1})^{a}_{j-j}(\varphi) f \Vert_{L^{2}} & \leq C_1 \,  \Vert f \Vert_{L^2}
\end{align}
uniformly on compact subsets of $\mathcal{B}$.
\\[2mm]
{\bf b)} 
Let $\varphi$ be as in a). Then 
there exists a constant $C_2>0$ with $C_2 \lesssim \Vert \widehat{\varphi} \Vert_{L^{1}}$ such that for all $g\in H^1(\R,\R)$
there holds
\begin{equation}
\label{Gaj}
(\rho^l \mathfrak{N}^{1})_{jj}^a(\varphi) g = (\partial_{\alpha}^{-1} \rho^l \mathfrak{N}^{1})_{jj}^s(\varphi) \partial_{\alpha} g + \widetilde{M}_{j}^{l}(\varphi) g \, 
\end{equation}
with
\begin{equation}
\label{Maj}
\|\widetilde{M}_{j}^{l}(\varphi) g\|_{L^{2}}  \leq C_2\,  \| g\|_{L^2}
\end{equation}
uniformly on compact subsets of $\mathcal{B}$. Moreover, for all $\psi \in H^{2+p}(\R,\R)$ with $p >1/2$ there exists a constant ${C}_3>0$ with ${C}_3 \lesssim \Vert \widehat{\varphi} \Vert_{L^{1}} \Vert \partial_{\alpha} \psi \Vert_{H^{1+p}} $ such that for all $f \in L^2(\R,\R)$ there holds
\begin{align}
\label{Gajkomm}
\Vert [\psi, (\partial_{\alpha}^{-1}\rho^l \mathfrak{N}^{1})_{jj}^s(\varphi)] f \Vert_{H^{1}} & \leq \, {C}_3\, \Vert f \Vert_{L^2}
\end{align}
uniformly on compact subsets of $\mathcal{B}$. 
\\[2mm]
{\bf c)} 
Let $\varphi$ be as in a). Then $f \mapsto (\rho^l \mathfrak{N}^{2})^{a}_{jj}(\varphi)f$ defines a continuous linear map from $L^2(\R,\R)$ into $L^{2}(\R,\R)$ and 
 there exists a constant $C_4>0$ with $C_4 \lesssim \Vert \widehat{\varphi} \Vert_{L^{1}}$ such that for all $f\in L^2(\R,\R)$ there holds
\begin{align}
\label{thetaN2a}
\Vert (\rho^l \mathfrak{N}^{2})^{a}_{jj}(\varphi) f \Vert_{L^{2}} & \leq \, C_4\, \Vert f \Vert_{L^2}
\end{align}
uniformly on compact subsets of $\mathcal{B}$.
\end{lemma}
{\bf Proof.}
a) follows directly from Lemma \ref{lem31} a) and Lemma \ref{lem3sym} a).
\\[2mm]
b) There holds
\begin{align} \label{na-1}
& (\widehat{\rho^l \mathfrak{n}^{1}})^{a}_{jj\ell}(k,m)  \nonumber
\\[2mm]
& \qquad  =\, \frac12  \big( (ik)^{-1} \, \widehat{\rho}_{j}^{l}(k)\, \widehat{\vartheta}(k)\,\widehat{n}^{1}_{jj\ell}(k)  +  (i(-m))^{-1} \, \widehat{\rho}_{j}^{l}(-m)\, \widehat{\vartheta}(-m)\,\widehat{n}^{1}_{jj\ell}(-m) \big)\, im   \nonumber
\\[2mm]
& \qquad \quad \,\, + \frac12\, \widehat{\rho}_{j}^{l}(k)\, \widehat{\vartheta}(k)\, (ik)^{-1}\, \widehat{n}^{1}_{jj\ell}(k)\, i(k-m)  \,,
\end{align}
which, due to \eqref{asnjj}, implies \eqref{Gaj}--\eqref{Maj}.

\eqref{Gajkomm} is proven in an analogous manner as \eqref{N0skomm}.
\\[2mm]
c) By construction of $\widehat{\rho}_{j}^l \widehat{\zeta}_{j j \ell}$ and by 
using Taylor's theorem as well as \eqref{oddq}, \eqref{k1}, \eqref{asfq2}, \eqref{asrjj}, \eqref{rswitch}--\eqref{dkrjjell}, \eqref{p1o3q4} and
\begin{align} 
\partial_{n_1} \partial_{n_2} \, \widehat{q}_{j_1j_2}^{2,4}(k,\ell k_0,k-\ell k_0)
&\,=\, \mathcal{O}(|k|^{-5/2}(1+ bk^2)^{1/2})  
\end{align}
for $n_1, n_2 \in \{1,3\}$ and $|k| \rightarrow \infty$ uniformly with respect to $b \lesssim 1$
we obtain
\begin{align} 
&(\widehat{\rho^l \mathfrak{n}^{2}})_{jj \ell}^a(k,m)\, \chi_{\ell}(k-m) 
\nonumber \\[2mm]
& \qquad  =\, \frac{\widehat{\rho}_{j}^l(-k)\,\widehat{\zeta}_{j j \ell}(-k)\, \partial_1 \widehat{\mathfrak{q}}^{2}_{j j}(-k,\ell k_0,-k +\ell k_0)}{2 \widehat{r}_{j j \ell}(k)} \,(k-m) 
\nonumber \\[2mm]
& \qquad \quad\;\, + \frac{ \big(\textstyle \frac{d}{dk}\,(\widehat{\rho}_{j}^l\, \widehat{\zeta}_{j j \ell})\big)(-k)\, \widehat{\mathfrak{q}}^{2}_{j j}(-k,\ell k_0,-k +\ell k_0)}{2 \widehat{r}_{j j \ell}(k)} \,(k-m) 
\nonumber 
\end{align}
\begin{align}
& 
\qquad \quad\;\,
+  \frac{\widehat{\rho}_{j}^l(-k)\,\widehat{\zeta}_{j j \ell}(-k)\,\partial_3 \widehat{\mathfrak{q}}^{2}_{j j}(-k,\ell k_0,-k +\ell k_0)}{2 \widehat{r}_{j j \ell}(k)} \,
(k-m-2 \ell k_0) 
\nonumber \\[2mm]
& 
\qquad \quad\;\,
- \frac{\widehat{\rho}_{j}^l(-k)\,\widehat{\zeta}_{j j \ell}(-k)\,\widehat{\mathfrak{q}}^{2}_{j j}(-k,\ell k_0,-k +\ell k_0)}{2 \widehat{r}_{j j \ell}(k)}\,
 \frac{\textstyle \frac{d}{dk}\, \widehat{r}_{j j \ell}(k)}{ \widehat{r}_{j j \ell}(k)} \,
(k-m-\ell k_0) 
\nonumber \\[2mm]
& 
\qquad \quad\;\,
+\mathcal{O}(|k|^{-1})
\nonumber \\[3mm] \label{nua}
& 
\qquad  
=\, \mathcal{O}(1)
\end{align}
for $|k| \to \infty$ uniformly with respect to $ m \in \R$ and $b$ in compact subsets of $\mathcal{B}$, which implies \eqref{thetaN2a}.
\qed

\begin{lemma}
\label{Tlem}
The operators $\mathcal{T}_{j_1j_2}$ have the following properties:\\[2mm]
{\bf a)}
Fix functions $g,h$ with $\widehat{g}, \widehat{h} \in L^1(\mathbb{R}, \C)$ and 
$\mathrm{supp} \, \widehat{g} = \mathrm{supp} \, \widehat{h} = \mathrm{supp} \, \widehat{\psi}_c$. Then $f \mapsto \mathcal{T}_{j_1j_2}(g,h,f)$ defines a continuous linear map from $L^2(\R,\C)$ into $L^2(\R,\C)$, and for all $f \in L^2(\R,\C)$ we have
\begin{align}
\label{T-abschaetzung}
\Vert \mathcal{T}_{j_1j_2}(g,h,f) \Vert_{L^2} \lesssim {\varepsilon}^{-1} \Vert \widehat{g}\Vert_{L^1} \Vert \widehat{h}\Vert_{L^1} \Vert  f\Vert_{L^2}
\end{align}
uniformly on compact subsets of $\mathcal{B}$.
If $f,g$ and $h$ are real-valued, then $\mathcal{T}_{j_1j_2}(g,h,f)$ is also real-valued.
\\[2mm]
{\bf b)}
For sufficiently small $\eps >0$ we have 
\begin{align}
\label{T-wahl-1}
&\mathrm{sgn}(j_1) i \omega \mathcal{T}_{j_1j_2}^1(\psi_c,\psi_c,R_{j_2})+
\mathcal{T}_{j_1j_2}^1(i \omega \psi_c,\psi_c,R_{j_2})\\
& +\mathcal{T}_{j_1j_2}^1(\psi_c,i \omega \psi_c,R_{j_2}) 
- \mathrm{sgn}(j_2) \mathcal{T}_{j_1j_2}^1(\psi_c,\psi_c,i \omega R_{j_2}) \nonumber \\
&\qquad\quad = {\varepsilon}^{-1} \mathcal{C}_{j_1j_2}(\psi_c, \psi_c) R_{j_2} + Y_{j_1j_2}^1(\psi_c,\psi_c, R_{j_2})\, \nonumber
\end{align}
with
\begin{equation}
\label{T-rest-1}
\|Y_{j_1j_2}^1(\psi_c,\psi_c, R_{j_2})\|_{L^2} = \mathcal{O}( \Vert R_{j_2} \Vert_{L^2} + \Vert R_{2j_2} \Vert_{H^{3/2}})
\end{equation}
uniformly on compact subsets of $\mathcal{B}$ if $j_1,j_2 \in \{\pm 1\}$, and
\begin{align}
\label{T-wahl}
&\mathrm{sgn}(j_1) i \omega \mathcal{T}_{j_1j_2}^2(\psi_c,\psi_c,R_{j_2})+
\mathcal{T}_{j_1j_2}^2(i \omega \psi_c,\psi_c,R_{j_2})\\
& +\mathcal{T}_{j_1j_2}^2(\psi_c,i \omega \psi_c,R_{j_2}) 
- \mathrm{sgn}(j_2) \mathcal{T}_{j_1j_2}^2(\psi_c,\psi_c,i \omega R_{j_2}) \nonumber \\
&\qquad\quad = \sum_{j_3 \in \{\pm |j_1|\}} \mathcal{N}_{j_1j_3}(\psi_c, Q_{j_3j_2}(\psi_c) R_{j_2} ) + Y_{j_1j_2}^2(\psi_c,\psi_c, R_{j_2})\, \nonumber
\end{align}
with
\begin{equation}
\label{T-rest}
\|Y_{j_1j_2}^2(\psi_c,\psi_c, R_{j_2})\|_{L^2} = \mathcal{O}( \Vert R_{j_2} \Vert_{H^2})
\end{equation}
uniformly on compact subsets of $\mathcal{B}$.
\\[2mm]
{\bf c)}
For all $f \in L^2(\R,\C)$ we have
\begin{align}
\label{P_0 T}
P_{\delta_{0},\infty}\mathcal{T}_{j_1j_2}(\psi_c, \psi_c,f)=0\,.
\end{align}
\end{lemma}

\textbf{Proof.}
a) Because of \eqref{nrb4}, \eqref{theta-1-eps-1}, \eqref{taudef1}--\eqref{taudef6} and \eqref{nz0} we have
\begin{align*}
\|\widehat{\tau}_{j_1j_2 \ell}^{j}\|_{L^{\infty}} = \mathcal{O}(\varepsilon^{-1})
\end{align*}
uniformly on compact subsets of $\mathcal{B}$ for $j_1 \in \{\pm1,\pm2\}$, $j_2 \in \{\pm j_1\}$, $|j_1| \leq j \leq 2$  and $\ell \in \{ \pm 1\}$. With the help of 
Young's inequality for convolutions we obtain
\begin{align*}
\Vert \mathcal{T}_{j_1j_2}(g,h,f) \Vert_{L^2} &\lesssim \sum_{\genfrac{}{}{0pt}{}{|j_1| \leq j \leq 2,} {\ell \in \{ \pm 1\}} } \|\widehat{\tau}_{j_1j_2 \ell}^{j}\|_{L^{\infty}}  \Vert \widehat{g}_{\ell}\Vert_{L^1} \Vert  \widehat{h}_{\ell} \Vert_{L^1} \Vert  f\Vert_{L^2} \lesssim  \varepsilon^{-1} \Vert \widehat{g} \Vert_{L^1} \Vert  \widehat{h} \Vert_{L^1} \Vert  f\Vert_{L^2} 
\end{align*}
uniformly on compact subsets of $\mathcal{B}$. Furthermore, since
\begin{equation*}
\widehat{\tau}_{j_1 j_2 -\ell}^{j}(-k)= \overline{\widehat{\tau}_{j_1j_2 \ell}^{j}(k)}\,, 
\end{equation*}
we conclude that $\mathcal{T}_{j_1j_2}(g,h,f)$ is real-valued if $f,g$ and $h$ are real-valued.
\\[2mm]
b) To prove the first assertion of b), we first show that there holds 
\begin{align}
\label{Ndom1}
{\varepsilon}^{-1} \mathcal{C}_{j_1j_2}(\psi_c, \psi_c) R_{j_2}
= \sum_{\ell \in \{\pm 1\}} {\varepsilon}^{-1} P_{0, \delta_{0}} 
\mathcal{C}_{j_1j_2}(\psi_{\ell}, \psi_{\ell}, R_{j_2})
+ \mathcal{O}(\Vert R_{j_2} \Vert_{L^2} + \Vert R_{2j_2} \Vert_{H^{3/2}})
\end{align}
uniformly on compact subsets of $\mathcal{B}$ such that it is sufficient to prove that the $L^{2}$-norm of 
\begin{align*}
\widetilde{Y}^1_{j_1j_2} :=\sum_{\ell \in \{\pm 1\}} \Big(  
&\mathrm{sgn}(j_1) i \omega \mathcal{T}_{j_1j_2 \ell}^{1}(\psi_{\ell},\psi_{\ell},R_{j_2})+
\mathcal{T}_{j_1j_2 \ell}^{1}(i \omega \psi_{\ell},\psi_{\ell},R_{j_2})\\
& +\mathcal{T}_{j_1j_2 \ell}^{1}(\psi_{\ell},i \omega \psi_{\ell},R_{j_2}) 
- \mathrm{sgn}(j_2) \mathcal{T}_{j_1j_2 \ell}^{1}(\psi_{\ell},\psi_{\ell},i \omega R_{j_2})\\[2mm]
&- {\varepsilon}^{-1} P_{0, \delta_{0}} \mathcal{C}_{j_1j_2}(\psi_{\ell}, \psi_{\ell}) R_{j_2}
\Big)
\end{align*}
is of order $\mathcal{O}(\Vert R_{j_2} \Vert_{L^2})$ uniformly on compact subsets of $\mathcal{B}$, which we will obtain
by construction of $\mathcal{T}_{j_1j_2 \ell}^{1}$ and because of Lemma \ref{int-kerne}.

To verify \eqref{Ndom1}, we split ${\varepsilon}^{-1} \mathcal{C}_{j_1j_2}(\psi_c, \psi_c) R_{j_2}$ into
\begin{align*}
{\varepsilon}^{-1} \mathcal{C}_{j_1j_2}(\psi_c, \psi_c) R_{j_2}
 =& \sum_{\ell \in \{\pm 1\}} {\varepsilon}^{-1} P_{0, \delta_{0}} \mathcal{C}_{j_1j_2}(\psi_{\ell}, \psi_{\ell}) R_{j_2}
 + \sum_{\ell \in \{\pm 1\}} {\varepsilon}^{-1} P_{0, \delta_{0}} 
\mathcal{C}_{j_1j_2}(\psi_{\ell}, \psi_{-\ell}) R_{j_2} \\[1mm]
&+ {\varepsilon}^{-1} P_{\delta_{0}, \infty}\, \mathcal{C}_{j_1j_2}(\psi_c, \psi_c) R_{j_2}
\,.
\end{align*}
Due to \eqref{RES3}, \eqref{theta-1-eps-0}, \eqref{c1} and \eqref{c2} we have
\begin{align*}
\Vert {\varepsilon}^{-1} P_{\delta_{0}, \infty}\, \mathcal{C}_{j_1j_2}(\psi_c, \psi_c) R_{j_2}\Vert_{L^2} = \mathcal{O}(\Vert R_{j_2} \Vert_{L^2} + \Vert R_{2j_2} \Vert_{H^{3/2}})
\end{align*}
uniformly on compact subsets of $\mathcal{B}$.

It follows from \eqref{ytilde}--\eqref{ctilde} that each summand of 
$\mathcal{C}_{j_1j_2}(\psi_c, \psi_c) R_{j_2}$ contains at least one $\alpha$-derivative. Using this fact as well as $\psi_c \partial_{\alpha} \psi_c = (\partial_{\alpha} (\psi_c)^2)/2$ and  the inequality $|n| \leq |k| + |k-n|$, we obtain
\begin{align*}
&\;{\varepsilon}^{-1} \mathcal{F}(P_{0, \delta_{0}} \mathcal{C}_{j_1j_2}(\psi_{\ell}, \psi_{-\ell}) R_{j_2})(k) \\[2mm]
=&\;
\widehat{P}_{0, \delta_{0}}\!(k) \! \int_{\mathbb{R}}\int_{\mathbb{R}}  \widehat{\vartheta}^{-1}
(k)\, \widehat{c}_{j_1 \ell -\ell j_2}(k,k-m,m-n,n)\,
\widehat{\psi}_{\ell}(k-m)\, \widehat{\psi}_{-\ell}(m-n)\, \widehat{\vartheta}(n) \widehat{R}_{j_2}(n)\, dn dm\,
\end{align*}
with
\begin{align}
\label{cderiv}
|\widehat{c}_{j_1 \ell -\ell j_2}(k,k-m,m-n,n)|\;
\lesssim\; |k|+|k-n|
\end{align}
uniformly on compact subsets of $\mathcal{B}$.
\eqref{theta-1-eps-1}, \eqref{theta-1} and \eqref{cderiv} as well as Fubini's theorem yield
\begin{align*}
&\;
\big|\widehat{P}_{0, \delta_{0}}(k) \int_{\mathbb{R}}\int_{\mathbb{R}} \widehat{\vartheta}^{-1}
(k)\, \widehat{c}_{j_1 \ell -\ell j_2}(k,k-m,m-n,n)\,
\widehat{\psi}_{\ell}(k-m) \, \widehat{\psi}_{-\ell}(m-n) \,\widehat{\vartheta}(n)\, \widehat{R}_{j_2}(n)\, dn dm \big| \\[2mm]
&\; \qquad \lesssim  \int_{\mathbb{R}}\int_{\mathbb{R}} \widehat{P}_{0, \delta_{0}}(k) \, \big(1+ \big|\textstyle \frac{k-n}{\veps}\big| \big)\,
|\widehat{\psi}_{\ell}(k-m)| \,|\widehat{\psi}_{-\ell}(m-n)|\, |\widehat{R}_{j_2}(n)|\, dn dm 
\\[2mm]
&\; \qquad = \int_{\mathbb{R}} \widehat{P}_{0, \delta_{0}}(k) \,
\big(1+ \big|\textstyle \frac{k-n}{\veps}\big| \big)\,
(|\widehat{\psi}_{\ell}| \ast |\widehat{\psi}_{-\ell}|)(k-n) \,|\widehat{R}_{j_2}(n)|\, dn
\end{align*}
uniformly on compact subsets of $\mathcal{B}$. Because of \eqref{pm1-ans} and \eqref{pm1-ans2}, the function $|\widehat{\psi}_{\ell}| \ast |\widehat{\psi}_{-\ell}|$ is strongly concentrated near $0$, more precisely, it
has a compact support being independent of $\veps$ and is 
of the form
\begin{equation*}
(|\widehat{\psi}_{\ell}| \ast |\widehat{\psi}_{-\ell}|)(k) = {\varepsilon}^{-1}  
\widehat{g}({\varepsilon}^{-1} k)
\end{equation*}
for all $k \in \R$, with $\widehat{g} \in L^{1}(1)$. 
Hence, by using \eqref{RES3} and Young's inequality for convolutions, we conclude
\begin{align*}
\| {\varepsilon}^{-1} P_{0, \delta_{0}} \mathcal{C}_{j_1j_2}(\psi_{\ell}, \psi_{-\ell}) R_{j_2}
\|_{L^{2}} =&\,\mathcal{O}(\Vert R_{j_2} \Vert_{L^2})
\end{align*}
uniformly on compact subsets of $\mathcal{B}$.
Therefore, we have verified
\eqref{Ndom1}.

To estimate $ \Vert \widetilde{Y}^{1}_{j_1j_2} \Vert_{L^2}$, we use
\begin{equation*}
\widehat{\widetilde{Y}^1_{j_1 j_2}}(k) \,= \sum_{\ell \in \{\pm1\}} \int_{\R} \int_{\R} 
K_{j_1 j_2}^{1}(k,k-m,m-n,n, \ell k_0, k-2\ell k_0)\,
\widehat{\psi}_{\ell}(k-m)\,  \widehat{\psi}_{\ell}(m-n)\,  \widehat{R}_{j_2}(n)\,dndm\,
\end{equation*}
with
\begin{align*}
& K_{j_1 j_2 \ell}^{1}(k,k-m,m-n,n) \\[1mm] & \qquad =\, \widehat{P}_{0, \delta_{0}}(k)\, \frac{\widehat{\nu}_{j_1 j_2}(k,k-m,m-n,n)\, \widehat{c}_{j_1 \ell \ell j_2}(k, \ell k_0,\ell k_0, k-2\ell k_0)\, \widehat{\vartheta}(k-2\ell k_0)}{\widehat{\nu}_{j_1 j_2}(k,\ell k_0,\ell k_0, k-2\ell k_0)\,\widehat{\vartheta}(k)}\\[2mm]
& \qquad \qquad \!\!\! - \widehat{P}_{0, \delta_{0}}(k)\, \widehat{\vartheta}^{-1}(k)\,\widehat{c}_{j_1 \ell \ell j_2}(k,k-m,m-n,n)\,\widehat{\vartheta}(n) \,.
\end{align*}
We split $K_{j_1 j_2\ell}^{1}$ into
\begin{align*}
&K_{j_1 j_2 \ell}^{1}(k,k-m,m-n,n) 
\\[1mm]
&\qquad = \, {K}^{1}_{j_1 j_2 \ell}(k,k-m,m-n,n)- {K}^{1}_{j_1 j_2 \ell}(k,\ell k_0,m-n,n) \\[1mm]
&\qquad\quad \, + {K}^{1}_{j_1 j_2 \ell}(k,\ell k_0,m-n,n)- {K}^{1}_{j_1 j_2 \ell}(k,\ell k_0,\ell k_0,n) 
\\[1mm]
&\qquad\quad\, + {K}^{1}_{j_1 j_2 \ell}(k,\ell k_0,\ell k_0,n)- {K}^{1}_{j_1 j_2 \ell}(k,\ell k_0,\ell k_0, k-2\ell k_0)\\[1mm]
&\qquad\quad\, +  {K}^{1}_{j_1 j_2 \ell}(k,\ell k_0,\ell k_0, k-2\ell k_0) \,.
\end{align*}
Since 
\begin{align*}
{K}^{1}_{j_1 j_2 \ell}(k,\ell k_0,\ell k_0, k-2\ell k_0) = 0\,,
\end{align*}
we deduce by applying \eqref{RES3} and Lemma \ref{int-kerne} that there holds
\begin{align*}
\Vert \widetilde{Y}^{1}_{j_1j_2}  \Vert_{L^2} =  \mathcal{O}(\Vert {R}_{j_2} \Vert_{L^2})
\end{align*}
uniformly on compact subsets of $\mathcal{B}$.
Hence, we have proven the first assertion of b).

To prove the second assertion of b), we first show that there holds
\begin{align}
\label{Ndom}
\mathcal{N}_{j_1j_3}(\psi_c, Q_{j_3j_2}(\psi_c) R_{j_2})
= \sum_{\ell \in \{\pm 1\}} P_{0, \delta_{0}} \mathcal{N}_{j_1j_3}(\psi_{\ell}, Q_{j_3j_2}(\psi_{\ell}) R_{j_2}) + \mathcal{O}(\Vert R_{j_2} \Vert_{H^2})
\end{align}
uniformly on compact subsets of $\mathcal{B}$ such that it is sufficient to prove that the $L^{2}$-norm of 
\begin{align*}
\widetilde{Y}^2_{j_1j_2} :=\sum_{\ell \in \{\pm 1\}} \Big(  
&\mathrm{sgn}(j_1) i \omega \mathcal{T}_{j_1j_2 \ell}^{2}(\psi_{\ell},\psi_{\ell},R_{j_2})+
\mathcal{T}_{j_1j_2 \ell}^{2}(i \omega \psi_{\ell},\psi_{\ell},R_{j_2})\\
& +\mathcal{T}_{j_1j_2 \ell}^{2}(\psi_{\ell},i \omega \psi_{\ell},R_{j_2}) 
- \mathrm{sgn}(j_2) \mathcal{T}_{j_1j_2 \ell}^{2}(\psi_{\ell},\psi_{\ell},i \omega R_{j_2})\\[2mm]
&- \sum_{j_3 \in \{\pm |j_1|\}} P_{0, \delta_{0}} \mathcal{N}_{j_1j_3}(\psi_{\ell}, Q_{j_3j_2}(\psi_{\ell}) R_{j_2})
\Big)
\end{align*}
is of order $\mathcal{O}(\Vert R_{j_2} \Vert_{L^2})$ uniformly on compact subsets of $\mathcal{B}$, which we will obtain
by construction of $\mathcal{T}_{j_1j_2}^2$ and because of Lemma \ref{int-kerne}.

To verify \eqref{Ndom}, we split $\mathcal{N}_{j_1j_3}(\psi_c, Q_{j_3j_2}(\psi_c) R_{j_2})$ into
\begin{align*}
\mathcal{N}_{j_1j_3}(\psi_c, Q_{j_3j_2}(\psi_c) R_{j_2})
 =& \sum_{\ell \in \{\pm 1\}} P_{0, \delta_{0}} \mathcal{N}_{j_1j_3}(\psi_{\ell}, Q_{j_3j_2}(\psi_{\ell}) R_{j_2}) \\
 &+ \sum_{\ell \in \{\pm 1\}} P_{0, \delta_{0}} \mathcal{N}_{j_1j_3}(\psi_{\ell}, Q_{j_3j_2}(\psi_{-\ell}) R_{j_2}) \\
&+ P_{\delta_{0}, \infty} \mathcal{N}_{j_1j_3}(\psi_c, Q_{j_3j_2}(\psi_c) R_{j_2})
\,.
\end{align*}
Because of \eqref{RES3}, \eqref{q}--\eqref{frakm}, \eqref{thetaN}--\eqref{thetaNb} and Young's inequality for convolutions we conclude that the $L^{2}$-norm of the last summand is of order $\mathcal{O}(\Vert R_{j_2} \Vert_{H^2})$ uniformly on compact subsets of $\mathcal{B}$.
Moreover, due to Lemma \ref{int-kerne} and \eqref{P_0 cN}, there holds
\begin{align*}
&\;\mathcal{F}\,P_{0, \delta_{0}}\,\mathcal{N}_{j_1j_3}(\psi_{\ell}, Q_{j_3j_2}(\psi_{j}) R_{j_2})
(k) \\[2mm]
&\; \qquad =\,
\int_{\mathbb{R}}\int_{\mathbb{R}} {K}_{j_1j_2j_3}(k,\ell k_0,k-\ell k_0,jk_0, k-(\ell+j)k_0)
\widehat{\psi}_{\ell}(k-m) \widehat{\psi}_{j}(m-n) \widehat{R}_{j_2}(n)\, dn dm \\[2mm]
&\; \qquad \quad  + \mathcal{O}(\Vert R_{j_2} \Vert_{L^2})
\end{align*}
uniformly on compact subsets of $\mathcal{B}$, where
\begin{align*}
&\;{K}_{j_1j_2j_3}(k,\ell k_0,k-\ell k_0, jk_0, k-(\ell+j)k_0) \\[2mm]
&\; \qquad =\, \widehat{P}_{0, \delta_{0}}(k)\, \frac{\widehat{q}_{j_1 j_2 j_3}(k,\ell k_0,k-\ell k_0,j k_0,k-(\ell+j)k_0)\,  \widehat{\vartheta}(k-(\ell+j) k_0)}{
\widehat{r}_{j_1j_3}(k,\ell k_0,k-\ell k_0)\,\widehat{\vartheta}(k)}\,.
\end{align*}

If $j=-\ell$, then we have
\begin{align*}
&\;{K}_{j_1j_2j_3}(k,\ell k_0,k-\ell k_0, jk_0, k-(\ell+j)k_0) \\[2mm]
&\; \qquad =\,{K}_{j_1j_2j_3}(k,\ell k_0,k-\ell k_0,-\ell k_0, k) \\[2mm]
&\; \qquad =\, \widehat{P}_{0, \delta_{0}}(k)\, \frac{\widehat{q}_{j_1 j_2 j_3}(k,\ell k_0,k-\ell k_0,-\ell k_0,k)\,  \widehat{\vartheta}(k)}{
\widehat{r}_{j_1j_3}(k,\ell k_0,k-\ell k_0)\,\widehat{\vartheta}(k)}
\end{align*}
and the factor $\widehat{\vartheta}(k)$ in the denominator is canceled by the same factor in the numerator, ${K}_{j_1j_2j_3}(k,\ell k_0,k-\ell k_0,-\ell k_0,k)$ contains no factors which are of order $\mathcal{O}(\eps^{-1})$ such that there holds
\begin{align*}
\sup_{k \in \R}\, |{K}_{j_1j_2j_3}(k,\ell k_0,k-\ell k_0,-\ell k_0,k)| = \mathcal{O}(1)
\end{align*} 
uniformly on compact subsets of $\mathcal{B}$.
Hence, we obtain 
\begin{align*}
\|\mathcal{F}\, \mathcal{N}_{j_1j_3}(\psi_{\ell}, Q_{j_3j_2}(\psi_{-\ell}) R_{j_2})\|_{L^{2}} =&\,
 \mathcal{O}(\Vert R_{j_2} \Vert_{L^2})
\end{align*}
uniformly on compact subsets of $\mathcal{B}$ such that we have verified
\eqref{Ndom}.

We have
\begin{align*}
\widehat{\widetilde{Y}^{2}_{j_1j_2}}(k)
\,& = \,
 \sum_{\ell \in \{\pm 1\}} \int_{\mathbb{R}}\int_{\mathbb{R}} 
K^{2}_{j_1j_2}(k,k-m,m-n,n,\ell k_0,k-\ell k_0, k-2\ell k_0)\,
\widehat{\psi}_{\ell}(k-m)\, \times \\[1mm]
&\; \qquad \qquad \qquad \; \times\, 
\widehat{\psi}_{\ell}(m-n)\,\widehat{R}_{j_2}(n)\, dn dm \\[2mm]
&\; \quad  + \mathcal{O}(\Vert R_{j_2} \Vert_{L^2})
\end{align*}
uniformly on compact subsets of $\mathcal{B}$, where
\begin{align*}
&\;{K}^{2}_{j_1j_2}(k,k-m,m-n,n,\ell k_0,k-\ell k_0, k-2\ell k_0) \\[2mm]
&\; \qquad = \,\sum_{j_3 \in \{\pm 1\}} \Big(  \widehat{P}_{0, \delta_{0}}(k)\,  \frac{    \widehat{q}_{j_1j_2j_3}(k,\ell k_0, k-\ell k_0, \ell k_0, k-2\ell k_0)\,\widehat{\vartheta}(k-2\ell k_0)}{ 
\widehat{r}_{j_1j_3}(k,\ell k_0,k-\ell k_0)\, \widehat{\nu}_{j_1j_2}(k,\ell k_0, \ell k_0,k-2\ell k_0)\, \widehat{\vartheta}(k)}\, \times \\[2mm]
&\; \qquad \quad \,\qquad \qquad \times \big(\widehat{\nu}_{j_1j_2}(k,k-m,m-n,n)-  \widehat{\nu}_{j_1j_2}(k,\ell k_0, \ell k_0,k-2\ell k_0)\big) \Big).
\end{align*}

Now, we can apply again Fubini's theorem, Young's inequality for convolutions and Lemma \ref{int-kerne} 
to obtain
\begin{align*}
\Vert \widetilde{Y}^{2}_{j_1j_2}  \Vert_{L^2} =  \mathcal{O}(\Vert {R}_{j_2} \Vert_{L^2})
\end{align*}
uniformly on compact subsets of $\mathcal{B}$.
Hence, we have proven the second assertion of b).
\\[2mm]
c) follows directly by the definition of $\mathcal{T}_{j_1j_2}$.
\qed
\\[2mm]
Now, we are able to compare our energy with Sobolev norms of the error. We obtain

\begin{lemma}
\label{checkRlem}
For sufficiently small $\eps > 0$,
we have
\begin{equation}
\label{checkR gegen R}
\sum_{j \in \{\pm 1, \pm 2 \}} \Vert \check{R}_{j} \Vert_{L^2}^2 
\lesssim \sum_{j \in \{\pm 1, \pm 2 \}}  \Vert R_{j} \Vert_{L^{2}}^2  +  \Vert \partial_{\alpha} \mathcal{R}_{2} \Vert_{(H^{1})^2}^2  + \eps^4    \,,
\end{equation}
\begin{equation}
\label{R gegen checkR}
\sum_{j \in \{\pm 1, \pm 2\}} \Vert R_{j} \Vert_{L^2}^2 \lesssim   \sum_{j \in \{\pm 1, \pm 2\}} \Vert \check{R}_{j} \Vert_{L^2}^2 + \veps \,  \Vert \partial_{\alpha} \mathcal{R}_{2} \Vert_{(H^{1})^2}^2  + \eps^5 \,,
\end{equation}
as long as $\eps^{5/2} \Vert \mathcal{R}_{1} \Vert_{(L^{2})^2},\, \eps^{5/2} \Vert \mathcal{R}_{2} \Vert_{(H^{2})^2} \ll 1$, uniformly on compact subsets of $\mathcal{B}$.
\end{lemma}

\textbf{Proof.}
Estimate \eqref{checkR gegen R} follows from the estimates \eqref{dxR1R2}--\eqref{dx-2R2}, \eqref{cN0}--\eqref{cN0b} and \eqref{T-abschaetzung}.
\\
To prove \eqref{R gegen checkR} we introduce
$R_{j_1}^0:=P_{0, \delta_{0}} R_{j_1}$, $\check{R}_{j_1}^0:=P_{0, \delta_{0}}\check{R}_{j_1}$, $R_{j_1}^1:=P_{ \delta_{0}, \infty}R_{j_1}$, $\check{R}_{j_1}^1:=P_{ \delta_{0}, \infty}\check{R}_{j_1}$ 
and split $R_{j_1}$, $\check{R}_{j_1}$ into $R_{j_1}=R_{j_1}^0+R_{j_1}^1$ and $\check{R}_{j_1}=\check{R}_{j_1}^0+\check{R}_{j_1}^1$.
Because of \eqref{P_0 cN} and \eqref{P_0 T}, $R_{j_1}^0$ satisfies 
\begin{align}
& R_{j_1}^0 + \varepsilon^{2} \sum_{{j_2} \in \{ \pm j_1 \} } \mathcal{T}_{j_1 j_2}(\psi_c, \psi_c, R_{j_2}^0) \nonumber \\[2mm]
& \qquad  = \check{R}_{j_1}^0 - \varepsilon P_{0,\delta_{0}} \sum_{j_2 \in \{ \pm j_1 \} }  \mathcal{N}_{j_1 j_2} (\psi_c, R_{j_2}^1)- \varepsilon^{2} \sum_{j_2 \in \{ \pm j_1\} } \mathcal{T}_{j_1 j_2}(\psi_c, \psi_c, R_{j_2}^1)\,.
\end{align}
Multiplying this equation with $R_{j_1}^0$, integrating, summing over $j_1 \in  \{ \pm 1, \pm 2\}$ and using \eqref{T-abschaetzung} yields
\begin{align}
\label{R0est}
\sum_{j_1 \in \{\pm 1, \pm 2\}} \Vert R_{j_1}^0 \Vert_{L^2}
& \lesssim \sum_{j_1 \in \{\pm 1, \pm 2\}} \Vert \check{R}_{j_1}^0 \Vert_{L^2} + \Vert R_{j_1}^1 \Vert_{L^2}
\end{align}
uniformly on compact subsets of $\mathcal{B}$ for sufficiently small $\eps > 0$. Moreover, $R_{j_1}^1$ satisfies 
\begin{equation}
R_{j_1}^1 + \varepsilon P_{\delta_{0},\infty} \sum_{j_2 \in \{ \pm j_1 \} } \mathcal{N}_{j_1 j_2}(\psi_c, R_{j_2}^1)
= \check{R}_{j_1}^1 - \varepsilon P_{\delta_{0},\infty} \sum_{j_2 \in \{ \pm j_1 \} }  \mathcal{N}_{j_1 j_2}(\psi_c, R_{j_2}^0)\,.
\end{equation}
Multiplying this equation with $R_{j_1}^1$, integrating, summing over $j_1 \in  \{ \pm 1, \pm 2\}$ and using \eqref{PcN}--\eqref{PcNb} yields
\begin{align*}
\sum_{j_1 \in \{\pm 1, \pm 2\}} \Vert R_{j_1}^1 \Vert_{L^2}^2 
& \lesssim  \sum_{j_1 \in \{\pm 1, \pm 2\}} \Vert \check{R}_{j_1}^1 \Vert_{L^2} \Vert R_{j_1}^1 \Vert_{L^2} \\[2mm]
& \quad
+  \varepsilon \! \sum_{\genfrac{}{}{0pt}{}{j_1 \in \{\pm 1, \pm 2\},} {j_2 \in \{\pm j_1\}}} \! (\Vert R_{j_2}^0 \Vert_{L^2} + \Vert R_{j_2}^1 \Vert_{H^{1}}) \Vert R_{j_1}^1 \Vert_{L^2} \\[2mm]
& \quad
+  \varepsilon \! \sum_{j_1 \in \{\pm 2\}} \! (\Vert \partial_{\alpha}^{-1} R_{-2}  \Vert_{L^2}  +   \Vert \partial_{\alpha}^{-1} R_2 \Vert_{L^2}) \Vert R_{j_1}^1 \Vert_{L^2}
\end{align*}
uniformly on compact subsets of $\mathcal{B}$. With the help of \eqref{dxR1R2}--\eqref{dx-2R2}, \eqref{R0est} we deduce
\begin{align}
\label{R1est}
\sum_{j_1 \in \{\pm 1, \pm 2\}} \Vert R^1_{j_1} \Vert_{L^2}^2
& \, \lesssim \sum_{j_1 \in \{\pm 1, \pm 2 \}} \Vert \check{R}_{j_1}^1 \Vert_{L^2}^2 + \varepsilon  \, \Vert R_{j_1}^0 \Vert_{L^2}^2 +  \eps \, \Vert R_{j_1}^1 \Vert_{H^{1}}^2 
\nonumber  \\[2mm]
&\, \quad +  \eps\,  \Vert \partial_{\alpha}^{-1} R_{-2}  \Vert_{L^2}^2  + \eps \, \Vert \partial_{\alpha}^{-1} R_2 \Vert_{L^2}^2
\nonumber  \\[2mm]
& \, \lesssim \sum_{j_1 \in \{\pm 1, \pm 2 \}} \Vert \check{R}_{j_1} \Vert_{L^2}^2 + \veps \, \Vert \partial_{\alpha} \mathcal{R}_{2} \Vert_{(H^1)^2}^2 + \eps^5
\end{align}
and
\begin{align}
\label{R0est2}
\sum_{j_1 \in \{\pm 1, \pm 2\}} \Vert R_{j_1}^0 \Vert_{L^2}^2
& \, \lesssim  \sum_{j_1 \in \{\pm 1, \pm 2 \}}   \Vert \check{R}_{j_1} \Vert_{L^2}^2 + \veps \, \Vert \partial_{\alpha} \mathcal{R}_{2} \Vert_{(H^1)^2}^2 + \eps^5
\end{align} 
uniformly on compact subsets of $\mathcal{B}$ for sufficiently small $\eps > 0$. Combining \eqref{R1est} and \eqref{R0est2} yields \eqref{R gegen checkR}.
\qed
\\[3mm]
The analysis of $\mathcal{E}_{2,s}$ will be simplified by

\begin{lemma}
\label{int-helfer}
Let $f \in H^{l}(\R,\R)$ and $g \in H^{m}(\R,\R)$ with $l,m \ge 0$. Then we have 
\begin{align}
\label{int und theta}
\int_{\R}\partial_{\alpha}^{l} f\, \partial_{\alpha}^{m} \vartheta g\, d{\alpha} = \int_{\R}\partial_{\alpha}^{l} f\,  \partial_{\alpha}^{m} g\,d{\alpha}+ \mathcal{O}(\Vert f\Vert_{L^2} \Vert g\Vert_{L^2})\,,
\end{align}
\begin{align}
\label{int und theta-1}
\int_{\R}\partial_{\alpha}^{l} f\, \partial_{\alpha}^{m+1} \vartheta^{-1} g\, d{\alpha} = \int_{\R}\partial_{\alpha}^{l} f\,  \partial_{\alpha}^{m+1} g\, d{\alpha}+\mathcal{O}(\Vert f\Vert_{L^2} \Vert g\Vert_{L^2})\,.
\end{align}
\end{lemma}

\textbf{Proof.}
The proof is analogous to the proof of Lemma 4.4 in \cite{DH18}.
\qed
\\[3mm]
We obtain

\begin{lemma} \label{cor32}
For sufficiently small $\eps > 0$ and $2 \leq \check{s} \leq s$, we have
\begin{align}  \label{Enorm}
& \mathcal{E}_{\check{s} } \lesssim \|\mathcal{R}_{1}\|_{(L^{2})^2}^{2} + \|\mathcal{R}_{2}\|_{(H^{\check{s}})^2}^{2} + \eps^{5} \,,\\[2mm] \label{normE}
&  \|\mathcal{R}_{1}\|_{(L^{2})^2}^{2} + \|\mathcal{R}_{2}\|_{(H^{\check{s}})^2}^{2} \lesssim \mathcal{E}_{{\check{s}}} + \eps^{5} \,,
\end{align}
as long as $\eps^{5/2} \Vert \mathcal{R}_{1} \Vert_{(L^{2})^2},\, \eps^{5/2} \Vert \mathcal{R}_{2} \Vert_{(H^{\check{s}})^2} \ll 1$, uniformly on compact subsets of $\mathcal{B}$.
\end{lemma}

\textbf{Proof.}
Because of \eqref{rhobound}, \eqref{dx-2R2}--\eqref{dx-1R2}, Leibniz's rule, Lemma \ref{lem31} and Lemma \ref{int-helfer} there holds
\begin{align}
E_{2,l} = & \frac12 \sum_{j_1\in \{\pm2\}} \int_{\R} \partial_{\alpha}^{l} R_{j_1}\, \rho_{j_1}^{l} \partial_{\alpha}^{l} R_{j_1}   \,d{\alpha}  + \eps \langle \partial_{\alpha}^{l} \mathcal{R}_2, (\rho^{l} \mathfrak{N}^{1})(\psi_c) 
\partial_{\alpha}^{l} \mathcal{R}_2 \rangle_{(L^2)^2} 
\nonumber
\\[1mm] 
& + \eps\, \mathcal{O}( \|\mathcal{R}_{1}\|_{(L^{2})^2}^{2} + \|\mathcal{R}_2\|_{(H^{\max \{2,l\}})^2}^2 + \eps^4)
\end{align}
uniformly on compact subsets of $\mathcal{B}$ for all $1 \leq l \leq s$. Moreover, we have
\begin{equation}
 \langle \partial_{\alpha}^{l} \mathcal{R}_2, (\rho^{l} \mathfrak{N}^{1})(\psi_c) \partial_{\alpha}^{l} \mathcal{R}_2 \rangle_{(L^2)^2} = \langle \partial_{\alpha}^{l} \mathcal{R}_2, (\rho^{l} \mathfrak{N}^{1})^s(\psi_c)
\partial_{\alpha}^{l} \mathcal{R}_2 \rangle_{(L^2)^2} \,.
\end{equation}
Hence, due to Lemma \ref{lem3sym}, Lemma \ref{checkRlem}, Lemma \ref{int-helfer},  \eqref{dxR1R2} and \eqref{rhobound},  we obtain \eqref{Enorm}--\eqref{normE} uniformly on compact subsets of $\mathcal{B}$.
\qed

\subsection{The energy estimates}

Now, we are prepared to estimate $\frac{d}{dt} \mathcal{E}_{s}$. First, we show
\medskip

\begin{lemma}
For sufficiently small $\eps > 0$, we have
\begin{equation}
\label{dt E_0}
\frac{d}{dt} \mathcal{E}_{0}  \lesssim
 \varepsilon^2 ( \mathcal{E}_{2} +1) \,,
\end{equation}
as long as $\eps^{5/2} \Vert \mathcal{R}_{1} \Vert_{(L^{2})^2},\, \eps^{5/2} \Vert \mathcal{R}_{2} \Vert_{(H^{\check{s}})^2} \ll 1$, uniformly on compact subsets of $\mathcal{B}$.
\end{lemma}

\textbf{Proof.}
Because of \eqref{Rmps} and \eqref{ch R}  and since $\rho_{j_1}^0$ is symmetric, we have
\begin{align*}
\frac{d}{dt} \mathcal{E}_{0} = \sum_{j_1 \in \{\pm1,\pm2\}} \int_{\mathbb{R}} {\check{R}_{j_1}}\, \rho_{j_1}^0\, \partial_t\check{R}_{j_1}\, d\alpha \,
\end{align*}
with
\begin{align*}
\partial_t\check{R}_{j_1} =\,  \mathrm{sgn}(j_1) i \omega \check{R}_{j_1} + \eps^2 \sum_{k=1}^{12} F_{j_1}^k(\Psi, \mathcal{R})\,,
\end{align*}
where
\begin{align*}
F_{j_1}^1(\Psi, \mathcal{R}) \,=\,&\,\, \eps^{-1} \sum_{j_2 \in \{\pm j_1\}} \Big(  Q_{j_1 j_2}(\psi_c){R}_{j_2} - 
\mathrm{sgn}(j_1) i \omega \mathcal{N}_{j_1 j_2}(\psi_c,{R}_{j_2})\\
&\qquad\qquad\qquad\!\! - \mathcal{N}_{j_1 j_2}(i \omega \psi_c,{R}_{j_2})
+ \mathrm{sgn}(j_2) \mathcal{N}_{j_1 j_2}(\psi_c,i \omega {R}_{j_2}) \Big)
\,,
\\[1mm]
F_{j_1}^2(\Psi, \mathcal{R}) \,=\,&\,\, \varepsilon^{-9/2} \vartheta^{-1}   \mathrm{res}_{\mathrm{sgn}(j_1)}(\eps \Psi)\,,
\\[1mm]
F_{j_1}^3(\Psi, \mathcal{R}) \,=\,&\,\, \eps^{-1} \sum_{j_2 \in \{\pm j_1\}} \mathcal{N}_{j_1 j_2}( \partial_t \psi_c + i \omega \psi_c, R_{j_2})\,,
\\[1mm]
F_{j_1}^4(\Psi, \mathcal{R}) \,=\,&\,\, \eps^{-7/2} \sum_{j_2 \in \{\pm j_1\}} 
\mathcal{N}_{j_1 j_2}
(\psi_c,\vartheta^{-1} \mathrm{res}_{\mathrm{sgn}(j_2)}(\eps \Psi))\,, 
\\[1mm]
F_{j_1}^5(\Psi, \mathcal{R}) \,=\,&\,\, \sum_{j_2 \in \{\pm j_1\}} \Big(
\eps^{-1}  \mathcal{C}_{j_1 j_2}
(\psi_c,\psi_c){R}_{j_2} +  \sum_{j_3 \in \{\pm1\}} \mathcal{N}_{j_1 j_3}
(\psi_c,Q_{j_3 j_2}(\psi_c){R}_{j_2})
\\[1mm]
&\qquad\qquad\, - \mathrm{sgn}(j_1) i \omega \mathcal{T}_{j_1j_2}(\psi_c,\psi_c,R_{j_2})-
\mathcal{T}_{j_1j_2}(i \omega \psi_c,\psi_c,R_{j_2}) 
\\[1mm]
&\qquad\qquad\, - 
\mathcal{T}_{j_1j_2}(\psi_c,i \omega \psi_c,R_{j_2}) 
+ \mathrm{sgn}(j_2) \mathcal{T}_{j_1j_2}(\psi_c,\psi_c,i \omega R_{j_2})
\Big)\,,
\\[1mm]
F_{j_1}^6(\Psi, \mathcal{R}) \,=\,&\,\,  \sum_{j_2 \in \{\pm j_1\}}
\mathcal{T}_{j_1 j_2}(\partial_t \psi_c + i \omega \psi_c, \psi_c,R_{j_2}) +\mathcal{T}_{j_1 j_2}(\psi_c,\partial_t \psi_c + i\omega \psi_c,R_{j_2}) \,,
\\[2mm]
F_{j_1}^7(\Psi, \mathcal{R}) \,=\,&\,\,   \mathcal{M}_{j_1}(\Psi, \mathcal{R})\,,
\\[2mm]
F_{j_1}^8(\Psi, \mathcal{R}) \,=\,&\,\,  \sum_{j_2, j_3 \in \{\pm j_1\}}  \mathcal{N}_{j_1 j_3}
(\psi_c, \mathcal{C}_{j_3 j_2}(\psi_c,\psi_c){R}_{j_2})\,,
\\[1mm]
F_{j_1}^{9}(\Psi, \mathcal{R}) \,=\, &\,\, \eps \sum_{j_2 \in \{\pm j_1\}} 
\mathcal{N}_{j_1 j_2}
(\psi_c, \mathcal{M}_{j_2}(\Psi, \mathcal{R})) \,,
\\[1mm]
F_{j_1}^{10}(\Psi, \mathcal{R}) \,=\, &\,\,\eps \sum_{j_2,j_3 \in \{\pm j_1\}}  \mathcal{T}_{j_1 j_3}
(\psi_c,\psi_c,  Q_{j_3 j_2}(\psi_c){R}_{j_2} +  \mathcal{C}_{j_3 j_2}(\psi_c,\psi_c){R}_{j_2})\,,
\\[1mm]
F_{j_1}^{11}(\Psi, \mathcal{R}) \,=\, &\,\,\eps^{2} \sum_{j_2 \in \{\pm j_1\}}
\mathcal{T}_{j_1 j_2}(\psi_c,\psi_c,  \mathcal{M}_{j_2} (\Psi, \mathcal{R})) \,,
\\[1mm]
F_{j_1}^{12}(\Psi, \mathcal{R}) \,=\, &\,\,\eps^{-5/2} \sum_{j_2 \in \{\pm j_1\}}
\mathcal{T}_{j_1 j_2}(\psi_c,\psi_c,  \vartheta^{-1} \mathrm{res}_{\mathrm{sgn}(j_2)}(\eps \Psi))
\,.
\end{align*}
Due to the skew symmetry of $i\omega$  we obtain
\begin{align*}
\frac{d}{dt} \mathcal{E}_{0} = \eps^{2} \, \sum_{k=1}^{12} \,\sum_{j_1 \in \{\pm1,\pm2\}}\, \int_{\mathbb{R}} {\check{R}_{j_1}}\, \rho_{j_1}^0\,F_{j_1}^{k}(\Psi, \mathcal{R})  \, d\alpha  \,.
\end{align*}
Because of \eqref{dxR1R2}--\eqref{dx-2R2}, \eqref{rhobound}, \eqref{cN0}--\eqref{cN0b}, \eqref{cnf}--\eqref{cnf-rest} and \eqref{T-abschaetzung} we deduce
\begin{align*}
\Big| \sum_{j_1 \in \{\pm 1, \pm2\}} \int_{\mathbb{R}} {\check{R}_{j_1}}\, \rho_{j_1}^0\,F_{j_1}^{1}(\Psi, \mathcal{R})  \, d\alpha  \Big | \lesssim  
\| \mathcal{R}_{1}\|_{(L^{2})^{2}}^2 + \| \mathcal{R}_{2}\|_{(H^{2})^{2}}^2 + \eps^{4}   
\end{align*}
and with the help of
the Cauchy-Schwarz inequality and \eqref{rhobound} we conclude
\begin{align*}
\frac{d}{dt} \mathcal{E}_{0}
\,&\lesssim\,   \eps^2 \Big ( \| \mathcal{R}_{1}\|_{(L^{2})^{2}}^2 + \| \mathcal{R}_{2}\|_{(H^{2})^{2}}^2 + \eps^{4} +  \mathcal{E}_{0}  + \sum_{k=2}^{12} \,\sum_{j_1 \in \{\pm1,\pm2\}}\, \Vert F_{j_1}^{k}(\Psi, \mathcal{R}) \Vert_{L^2}^2  \Big)\,.
\end{align*}
Because of \eqref{RES1}, \eqref{theta-1-eps-1}--\eqref{theta-1-eps-0}, \eqref{cN0}--\eqref{PcNb},
\eqref{P_0 cN} and \eqref{T-abschaetzung} we have
\begin{align*}
& \Vert F_{j_1}^{2}(\Psi, \mathcal{R}) \Vert_{L^2}^2   \lesssim 1\,, \\[2mm]
& \Vert F_{j_1}^{4}(\Psi, \mathcal{R}) \Vert_{L^2}^2 + \Vert F_{j_1}^{12}(\Psi, \mathcal{R}) \Vert_{L^2}^2  \lesssim \eps^2\,.
\end{align*}
Due to \eqref{dxR1R2}--\eqref{dx-2R2} and \eqref{T-wahl-1}--\eqref{T-rest} we deduce
\begin{equation*}
\Vert F_{j_1}^{5}(\Psi, \mathcal{R}) \Vert_{L^2}^2   \lesssim \| \mathcal{R}_{1}\|_{(L^{2})^{2}}^2 + \| \mathcal{R}_{2}\|_{(H^{2})^{2}}^2 + \eps^{4}   
  \,.
\end{equation*}
Furthermore, \eqref{dtpsi}, \eqref{dxR1R2}--\eqref{dx-2R2}, \eqref{cN0}--\eqref{cN0b} and \eqref{T-abschaetzung} yield
\begin{equation*}
\Vert F_{j_1}^{3}(\Psi, \mathcal{R}) \Vert_{L^2}^2 +  \Vert F_{j_1}^{6}(\Psi, \mathcal{R}) \Vert_{L^2}^2   \lesssim  \| \mathcal{R}_{1}\|_{(L^{2})^{2}}^2 + \| \mathcal{R}_{2}\|_{(H^{2})^{2}}^2 + \eps^{4}   
  \,.
\end{equation*}
Using \eqref{calM1}, \eqref{dxR1R2}--\eqref{dx-2R2}, \eqref{cN0}--\eqref{cN0b} and \eqref{T-abschaetzung}, we obtain
\begin{equation*}
\Vert F_{j_1}^{7}(\Psi, \mathcal{R}) \Vert_{L^2}^2  + \Vert F_{j_1}^{9}(\Psi, \mathcal{R}) \Vert_{L^2}^2 + \Vert F_{j_1}^{11}(\Psi, \mathcal{R}) \Vert_{L^2}^2 \lesssim 
 \| \mathcal{R}_{1}\|_{(L^{2})^{2}}^2 + \| \mathcal{R}_{2}\|_{(H^{2})^{2}}^2 + \eps^{4}    
  \,. 
\end{equation*}
Finally, \eqref{theta-1-eps-1}--\eqref{theta-1-eps-0}, \eqref{c1}, \eqref{dxR1R2}--\eqref{dx-2R2}, \eqref{q}--\eqref{q13}, \eqref{c2},
\eqref{cN0}--\eqref{PcNb},
\eqref{P_0 cN} and \eqref{T-abschaetzung} imply
\begin{equation*}
\Vert F_{j_1}^{8}(\Psi, \mathcal{R}) \Vert_{L^2}^2  + \Vert F_{j_1}^{10}(\Psi, \mathcal{R}) \Vert_{L^2}^2 \lesssim \| \mathcal{R}_{1}\|_{(L^{2})^{2}}^2 + \| \mathcal{R}_{2}\|_{(H^{2})^{2}}^2 + \eps^{4}   
  \,.
\end{equation*}
All bounds are uniform on compact subsets of $\mathcal{B}$.

Hence, because of \eqref{normE} we arrive at
\begin{equation*}
\frac{d}{dt} \mathcal{E}_{0}  \lesssim
 \varepsilon^2( {\mathcal{E}}_{2}
  +1)
\end{equation*}
 as long as $\eps^{5/2} \Vert \mathcal{R}_{1} \Vert_{(L^{2})^2},\, \eps^{5/2} \Vert \mathcal{R}_{2} \Vert_{(H^{\check{s}})^2} \ll 1$, uniformly on compact subsets of $\mathcal{B}$.
\qed
\\[5mm]
For the estimates of $\frac{d}{dt} {E}_{2,l}$ with $1 \leq l \leq s$ we use 
\smallskip

\begin{lemma} \label{lem33c}
Let $s \geq 0$ and $A: H^s(\R,\R) \to  L^2(\R,\R)$ be an antisymmetric linear operator. 
\\[2mm]     
{\bf a)} Let $S_j: L^2(\R,\R) \to  L^2(\R,\R)$, $j \in \{1,2\}$, be symmetric linear operators with $S_jH^s(\R,\R)\subseteq H^s(\R,\R)$. Then for all $f \in H^s(\R,\R)$ there holds
\begin{align}
\label{dx1} 
\int_{\R} S_1S_2f\, A f\, d{\alpha}
=&\, -\frac{1}{2}\, \int_{\R} f\, [A, S_1S_2] f\, d{\alpha} -\frac{1}{2}\, \int_{\R} f\, A [S_1, S_2] f\, d{\alpha} .
\end{align}
{\bf b)} Let $S, S_{j}: L^2(\R,\R) \to  L^2(\R,\R)$, $j \in \{\pm 1\}$, be  symmetric linear operators with $SH^s(\R,\R),$ $S_{j}H^s(\R,\R) \subseteq H^s(\R,\R)$. Then for all $f \in H^s(\R,\R)$ there holds
\begin{align}
\label{dx2}
& \sum_{j \in \{\pm1\}} \int_{\R} S S_{j} f_{j}\, Af_{-j}\,d{\alpha} \nonumber\\
& \qquad  =\, \frac{1}{2}\, \int_{\R} (f_{-1}-f_1)\, S(S_{-1}-S_{1})A(f_{-1}+f_1)\, \,d\alpha \nonumber \\[2mm]  
& \qquad \quad\,\, - \frac{1}{2}\, \sum_{j \in \{\pm1\}} \Big(\int_{\R} f_{j}\, [A, S_{j}S] f_{-j}\, d{\alpha} + \int_{\R} f_{-j}\, A[S, S_{j}] f_{j}\,  d{\alpha} \Big)
\nonumber \\[1mm]
& \qquad \quad\,\, -  \sum_{j \in \{\pm1\}} \frac{j}{4} \Big(\int_{\R} f_{j}\, [A,(S_{-1}-S_{1})S] f_{j}\, d{\alpha} + \int_{\R} f_{j}\, A[S_{-1}-S_1,S] f_{j}\, d{\alpha}  \Big).
\end{align}
{\bf c)} Let $L_j: L^2(\R,\R) \to  L^2(\R,\R)$, $j \in \{\pm1\}$, be linear operators  
with $L_{j}H^s(\R,\R) \subseteq H^s(\R,\R)$ and $L_1^{\ast} = L_{-1}$. Moreover, let $S: L^2(\R,\R) \to  L^2(\R,\R)$ be a symmetric linear operator with $SH^s(\R,\R)\subseteq H^s(\R,\R)$. Then for all $f \in H^s(\R,\R)$ there holds
\begin{align}
\label{dx3} 
\sum_{j \in \{\pm1\}}\, \int_{\R} SL_{j} f_{j}\, Af_{-j}\,d{\alpha}
=&\, -\int_{\R} f_1\, [A, SL_{-1}] f_{-1}\, d{\alpha} - \int_{\R} f_{-1}\, A[S, L_{1}]f_{1} \, d{\alpha} . 
\end{align}
\end{lemma}

{\bf Proof.}
We have
\begin{align*}
  \int_{\R} S_1S_2f\, A f\, d\ualpha\; 
  =\; & \int_{\R} S_2S_1f\, Af\, d\ualpha + \int_{\R} [S_1, S_2] f\, Af\, d{\alpha} 
\\[2mm]
=\; & \int_{\R} f\, S_1S_2Af\, d\ualpha + \int_{\R} [S_1, S_2] f\, Af\, d{\alpha} 
\\[2mm]
 =\; & \int_{\R} f\, AS_1S_2 f\, d\ualpha - \int_{\R} f\, [A, S_1S_2] f\, d{\alpha} + \int_{\R} [S_1, S_2] f\, Af\, d{\alpha} 
\\[2mm]
=\; & -\int_{\R} S_1S_2f\, Af\, d\ualpha - \int_{\R} f\, [A, S_1S_2] f\, d{\alpha} 
- \int_{\R} f\, A [S_1, S_2] f\, d{\alpha} \,,
\end{align*}
which implies \eqref{dx1}, and
\begin{align*}
&\sum_{j \in \{\pm1\}}\, \int_{\R} SS_{j} f_{j}\, Af_{-j}\,d{\alpha}\;\\
& \qquad =\, \frac{1}{2}\, \sum_{j \in \{\pm1\}} \Big(\int_{\R} S_{j}S f_{j}\, Af_{-j}\,d{\alpha} + \int_{\R} SS_{j} f_{j}\, Af_{-j}\,d{\alpha} +
\int_{\R} [S, S_{j}] f_{j}\, Af_{-j}\, d{\alpha} \Big) 
\\[1mm]
& \qquad  =\, \frac{1}{2}\, \sum_{j \in \{\pm1\}} \Big(\int_{\R} S_{j}S f_{j}\, Af_{-j}\,d{\alpha} + \int_{\R} f_{j}\, AS_{j}S f_{-j}\, d\ualpha \Big) 
\\[1mm]
& \qquad \quad\,\,  - \frac{1}{2}\, \sum_{j \in \{\pm1\}} \Big(\int_{\R} f_{j}\, [A, S_{j}S] f_{-j}\, d{\alpha} + \int_{\R} f_{-j}\, A[S, S_{j}] f_{j}\,  d{\alpha} \Big)
\\[1mm] 
& \qquad =\, \frac{1}{2}\, \sum_{j \in \{\pm1\}} \Big(\int_{\R} S_{j}S f_{j}\, Af_{-j}\,d{\alpha} - \int_{\R}S_{j}S f_{-j}\, A f_{j}\, d\ualpha \Big) 
\\[1mm]
& \qquad \quad\,\, - \frac{1}{2}\, \sum_{j \in \{\pm1\}} \Big(\int_{\R} f_{j}\, [A, S_{j}S] f_{-j}\, d{\alpha} + \int_{\R} f_{-j}\, A[S, S_{j}] f_{j}\,  d{\alpha} \Big)
\\[2mm]
& \qquad  =\, \frac{1}{2}\, \Big( \int_{\R} (S_{-1}-S_{1})S\, f_{-1}\, Af_{1}\,d\ualpha - \int_{\R} (S_{-1}-S_{1})Sf_1\, A f_{-1}\,d\ualpha \Big) \\[2mm]  
& \qquad \quad\,\, - \frac{1}{2}\, \sum_{j \in \{\pm1\}} \Big(\int_{\R} f_{j}\, [A, S_{j}S] f_{-j}\, d{\alpha} + \int_{\R} f_{-j}\, A[S, S_{j}] f_{j}\,  d{\alpha} \Big)
\\[2mm]
& \qquad  =\, \frac{1}{2}\, \Big( \int_{\R} (S_{-1}-S_{1})S\, f_{-1}\, A(f_{-1}+f_{1})\,d\alpha 
 - \int_{\R} (S_{-1}-S_{1})S\,  f_1\, A (f_{-1}+f_1)\,dx \Big)
\\[2mm]  
& \qquad \quad\,\, - \frac{1}{2}\, \sum_{j \in \{\pm1\}} \Big(\int_{\R} f_{j}\, [A, S_{j}S] f_{-j}\, d{\alpha} + \int_{\R} f_{-j}\, A[S, S_{j}] f_{j}\,  d{\alpha} \Big)
\end{align*}
\begin{align*}
& \qquad \quad\,\, -  \sum_{j \in \{\pm1\}}  \frac{j}{4}  \Big(\int_{\R} f_{j}\, [A,(S_{-1}-S_{1})S] f_{j}\, d{\alpha} + \int_{\R} f_{j}\, A[S_{-1}-S_{1},S] f_{j}\, d{\alpha}  \Big)        
\\[2mm]
& \qquad  =\, \frac{1}{2}\, \int_{\R} (f_{-1}-f_1)\, S(S_{-1}-S_{1})A(f_{-1}+f_1) \,d\alpha 
\\[2mm]
& \qquad \quad\,\, - \frac{1}{2}\, \sum_{j \in \{\pm1\}} \Big(\int_{\R} f_{j}\, [A, S_{j}S] f_{-j}\, d{\alpha} + \int_{\R} f_{-j}\, A[S, S_{j}] f_{j}\,  d{\alpha} \Big)
\\[1mm]
& \qquad \quad\,\, - \frac{1}{2}\, \sum_{j \in \{\pm1\}} \Big(\int_{\R} f_{j}\, [A, S_{j}S] f_{-j}\, d{\alpha} + \int_{\R} f_{-j}\, A[S, S_{j}] f_{j}\,  d{\alpha} \Big)
\\[1mm]
& \qquad \quad\,\, -  \sum_{j \in \{\pm1\}} \frac{j}{4} \Big(\int_{\R} f_{j}\, [A,(S_{-1}-S_{1})S] f_{j}\, d{\alpha} + \int_{\R} f_{j}\, A[S_{-1}-S_1,S] f_{j}\, d{\alpha}  \Big).
\end{align*}
Moreover, we deduce
\begin{align*}
\sum_{j \in \{\pm1\}}\, \int_{\R} SL_{j} f_{j}\, Af_{-j}\,d{\alpha}\;
=\;&  \int_{\R} SL_{-1} f_{-1}\, Af_{1}\,d{\alpha} + \int_{\R} L_1Sf_{1}\, A f_{-1}\, d\ualpha 
\\
& + \int_{\R} [S, L_{1}] f_{1} \, A f_{-1} \, d{\alpha}
\\[2mm]
=\;&  \int_{\R} SL_{-1} f_{-1}\, Af_{1}\,d{\alpha} + \int_{\R} f_{1}\, ASL_{-1} f_{-1}\, d\ualpha \\[2mm]
& - \int_{\R} f_{1}\, [A, SL_{-1}] f_{-1}\, d{\alpha}
- \int_{\R} f_{-1}\, A[S, L_{1}]f_{1} \, d{\alpha} 
 \\[2mm]
 =\;& - \int_{\R} f_{1}\, [A, SL_{-1}] f_{-1}\, d{\alpha} - \int_{\R} f_{-1}\, A[S, L_{1}]f_{1} \, d{\alpha} \,. 
\\[-7mm]
\end{align*}
\qed
\\[2mm]
From now on, let $1 \leq l \leq s$. We compute
\begin{align*}
\frac{d}{dt} E_{2,l} =& \sum_{j_1 \in \{\pm2\}} \Big( \int_{\R} \partial^{l}_{\alpha} R_{j_1}\, \rho_{j_1}^{l} \partial_t\partial^{l}_{\alpha} R_{j_1} \,d{\alpha}\,
+   \veps \sum_{j_2 \in \{\pm2\}} \Big( \int_{\R} \partial_t\partial^{l}_{\alpha} R_{j_1}\, \rho_{j_1}^{l} \partial^{l}_{\alpha} N_{j_1j_2}(\psi_c)R_{j_2}\,d{\alpha}
\\[2mm] & 
+ \int_{\R} \partial^{l}_{\alpha} R_{j_1}\, \rho_{j_1}^{l} \partial^{l}_{\alpha} N_{j_1j_2}(\psi_c)\partial_t R_{j_2}\,d{\alpha}\,
+ \int_{\R} \partial^{l}_{\alpha} R_{j_1}\, \rho_{j_1}^{l} \partial^{l}_{\alpha} N_{j_1j_2}(\partial_t \psi_c) R_{j_2} \,d{\alpha} \Big)\Big)\,.
\end{align*}
Using \eqref{Rmps} we obtain
\begin{align*}
\frac{d}{dt} E_{2,l} 
=\;& \sum_{j_1 \in \{\pm2\}} \Big(\,\mathrm{sgn}(j_1) \int_{\R} \partial^{l}_{\alpha} R_{j_1}\,\rho_{j_1}^{l} i\omega \partial^{l}_{\alpha} R_{j_1} \,d{\alpha}\\[1mm] &
\qquad\qquad\! + \int_{\R} \partial^{l}_{\alpha} R_{j_1}\,\rho_{j_1}^{l} \veps^{-5/2} \partial^{l}_{\alpha} \vartheta^{-1} \mathrm{res}_{{j_1}}\!(\veps\Psi)\,d{\alpha}\, \Big)
\\[2mm] 
& 
+\; \veps \sum_{j_1, j_2 \in \{\pm2\}}\! \Big(\, \int_{\R} \partial^{l}_{\alpha} R_{j_1}\,\rho_{j_1}^{l} \partial^{l}_{\alpha} Q_{j_1j_2}(\psi_c)R_{j_2} \,d{\alpha} \\[1mm]
& 
\qquad\qquad\qquad\,\,\, + \int_{\R} \mathrm{sgn}(j_1) i \omega \partial^{l}_{\alpha} R_{j_1}\,\rho_{j_1}^{l} \partial^{l}_{\alpha} N_{j_1j_2}(\psi_c) R_{j_2}\,d{\alpha}
\end{align*}
\begin{align*}
&\qquad\qquad\qquad\,\,\, + \int_{\R} \mathrm{sgn}(j_2) \partial^{l}_{\alpha} R_{j_1}\,\rho_{j_1}^{l} \partial^{l}_{\alpha} N_{j_1j_2}(\psi_c) i \omega R_{j_2}\,d{\alpha}
\\[2mm]
&\qquad\qquad\qquad\,\,\, - \int_{\R} \partial^{l}_{\alpha} R_{j_1}\,\rho_{j_1}^{l} \partial^{l}_{\alpha} N_{j_1j_2}(i \omega \psi_c) R_{j_2}\,d{\alpha}
\\[2mm]
&\qquad\qquad\qquad\,\,\,
+ \int_{\R} \partial^{l}_{\alpha} R_{j_1}\,\rho_{j_1}^{l} \partial^{l}_{\alpha} N_{j_1j_2}(\partial_t \psi_c +i \omega \psi_c) R_{j_2} \,d{\alpha}
\\[2mm] 
& \qquad\qquad\qquad\,\,\,
- \int_{\R} \veps^{-5/2} \partial^{l+1}_{\alpha} \vartheta^{-1} \mathrm{res}_{{j_1}}\!(\veps\Psi)\,\rho_{j_1}^{l} \partial^{l-1}_{\alpha} N_{j_1j_2}(\psi_c)R_{j_2}\,d{\alpha}
\\[2mm]& \qquad\qquad\qquad\,\,\,
+ \int_{\R} \partial^{l}_{\alpha} R_{j_1}\,\rho_{j_1}^{l} \partial^{l}_{\alpha} N_{j_1j_2}(\psi_c)\veps^{-5/2} \vartheta^{-1} \mathrm{res}_{{j_2}}\!(\veps\Psi)\,d{\alpha}\, \Big)
\\[2mm]
&
+\; \veps^{2}  \sum_{j_1,j_2,j_3 \in \{\pm2\}}\!  \Big(\, \int_{\R}  \partial^{l}_{\alpha} Q_{j_1j_3}({\psi_c}) R_{j_3}\,\rho_{j_1}^{l} \partial^{l}_{\alpha} N_{j_1j_2}(\psi_c) R_{j_2}\,d{\alpha}\\[1mm]
&\qquad\qquad\qquad\qquad\! + \int_{\R} \partial^{l}_{\alpha} R_{j_1}\,\rho_{j_1}^{l} \partial^{l}_{\alpha} N_{j_1j_2}(\psi_c) Q_{j_2j_3}(\psi_c) R_{j_3}\,d{\alpha} \,\Big) 
\\[2mm]
&
+\;\veps^{2} \sum_{j_1, j_2 \in \{\pm2\}} 
\! \Big(\, \int_{\R} \partial^{l}_{\alpha} R_{j_1}\,\rho_{j_1}^{l} \partial^{l}_{\alpha} W_{j_1j_2}(\Psi, \mathcal{R})R_{j_2}\,d{\alpha}
\\[1mm]
& 
\qquad\qquad\qquad\,\,\,\,\,\, + \int_{\R} \partial^{l}_{\alpha} R_{j_1}\, \rho_{j_1}^{l} \partial_{\alpha}^{l} \mathcal{M}_{j_1}(\Psi, \mathcal{R})\,d{\alpha}\,\Big)
\\[2mm]
&
+\; \veps^{3}  \sum_{j_1,j_2,j_3 \in \{\pm2\}}\!  \Big(\, \int_{\R}  \partial^{l}_{\alpha} W_{j_1j_3}(\Psi, \mathcal{R}) R_{j_3}\, \rho_{j_1}^{l} \partial^{l}_{\alpha} N_{j_1j_2}(\psi_c) R_{j_2}\,d{\alpha}\\[1mm]
&\qquad\qquad\qquad\qquad\! + \int_{\R} \partial^{l}_{\alpha} R_{j_1}\, \rho_{j_1}^{l} \partial^{l}_{\alpha} N_{j_1j_2}(\psi_c) W_{j_2j_3}(\Psi, \mathcal{R}) R_{j_3}\,d{\alpha} \,\Big)
\\[2mm]
&
+\; \veps^{3}  \sum_{j_1,j_2 \in \{\pm2\}}\!  \Big(\,
\int_{\R}  \partial^{l}_{\alpha} \mathcal{M}_{j_1}(\Psi, \mathcal{R})\, \rho_{j_1}^{l}\partial^{l}_{\alpha} N_{j_1j_2}(\psi_c) R_{j_2}\,d{\alpha}\\[1mm]
&\qquad\qquad\qquad\,\,\,\,\,\, + \int_{\R} \partial^{l}_{\alpha} R_{j_1}\,\rho_{j_1}^{l} \partial^{l}_{\alpha} N_{j_1j_2}(\psi_c) \mathcal{M}_{j_2}(\Psi, \mathcal{R})\,d{\alpha} \,\Big)\,. 
\end{align*}

Due to the skew symmetry of $i\omega$  and the symmetry of $\rho_{j_1}^{l}$ the first integral equals zero. 
Moreover, because of 
\eqref{RES1}, \eqref{dtpsi}, \eqref{theta-1-eps-1}--\eqref{theta-1}, \eqref{rhobound}, \eqref{N0}--\eqref{P_0 N}, \eqref{N0s}, \eqref{normE} 
and \eqref{dx1}
the sum of the second, the seventh, the eighth and the ninth line can be bounded by $C\veps^{3}(\mathcal{E}_{s}+1)$ for a constant $C>0$, as long as $\eps^{5/2} \Vert \mathcal{R}_{1} \Vert_{(L^{2})^2},\, \eps^{5/2} \Vert \mathcal{R}_{2} \Vert_{(H^{{s}})^2} \ll 1$, uniformly on compact subsets of $\mathcal{B}$. 
Hence, using \eqref{rhobound}, \eqref{nf}--\eqref{nf-rest}, \eqref{int und theta}--\eqref{int und theta-1} and \eqref{normE}, 
we obtain
\begin{align*}
\frac{d}{dt} E_{2,l} =\;& \sum_{j=1}^{7} I_{j} + \veps^{2}\,\mathcal{O}(\mathcal{E}_{s}+1) \,,
\end{align*}
as long as $\eps^{5/2} \Vert \mathcal{R}_{1} \Vert_{(L^{2})^2},\, \eps^{5/2} \Vert \mathcal{R}_{2} \Vert_{(H^{{s}})^2} \ll 1$, uniformly on compact subsets of $\mathcal{B}$, where
\begin{align*}
I_1 =\;& -\veps^{2} \sum_{j_1,j_2  \in \{\pm2\}}  \int_{\R} \partial^{l}_{\alpha} R_{j_1}\, \rho_{j_1}^{l} \partial^{l}_{\alpha} \vartheta Y_{j_1j_2}({\psi}_c, R_{j_2})\,d{\alpha}\,, 
\end{align*}
\begin{align*}
I_2 =\;& \veps^{2} \sum_{j_1,j_2  \in \{\pm2\}} \Big(\, \int_{\R}  \partial^{l}_{\alpha} Q_{j_1j_2}({\psi}_c) R_{j_2}\, \rho_{j_1}^{l} \partial_{\alpha}^{l} \vartheta N_{j_1j_1}(\psi_c) R_{j_1}\,d{\alpha}\\[1mm]
&\qquad\qquad\qquad\!\! +  \int_{\R} \partial^{l}_{\alpha} R_{j_1}\, \rho_{j_1}^{l} \partial^{l}_{\alpha} \vartheta N_{j_1j_1}(\psi_c) Q_{j_1j_2}({\psi}_c) R_{j_2}\,d{\alpha}\, \Big), 
\\[3mm]
I_3 =\; 
& \veps^{2} \sum_{j_1, j_2 \in \{\pm2\}} 
\! \Big(\, \int_{\R} \partial^{l}_{\alpha} Q_{j_1j_2}({\psi}_c) R_{j_2} \, \rho_{j_1}^{l} \partial_{\alpha}^{l} \vartheta N_{j_1-j_1}(\psi_c) R_{-j_1}\,d{\alpha}\; \\[2mm]
&\qquad\qquad\qquad\!\! + \int_{\R} \partial^{l}_{\alpha} R_{j_1}\, \rho_{j_1}^{l} \partial^{l}_{\alpha} \vartheta N_{j_1-j_1}(\psi_c) Q_{-j_1j_2}({\psi}_c) R_{j_2}\,d{\alpha}\, \Big),
\\[3mm]
I_4 =\; 
& \veps^{2} \sum_{j_1, j_2 \in \{\pm2\}} 
\! \Big(\, \int_{\R} \partial^{l}_{\alpha} R_{j_1}\,\rho_{j_1}^{l} \partial^{l}_{\alpha} W_{j_1j_2}(\Psi, \mathcal{R})R_{j_2}\,d{\alpha}
\\[1mm]
& 
\qquad\qquad\qquad\!\!\! + \int_{\R} \partial^{l}_{\alpha} R_{j_1}\, \rho_{j_1}^{l}  \partial_{\alpha}^{l} \mathcal{M}_{j_1}(\Psi, \mathcal{R})\,d{\alpha}\,\Big),
\\[3mm] 
I_5 =\;  &  \veps^{3} \sum_{j_1,j_2  \in \{\pm2\}} \Big(\, \int_{\R}  \partial^{l}_{\alpha} W_{j_1j_2}({\Psi}, \mathcal{R}) R_{j_2}\, \rho_{j_1}^{l} \partial_{\alpha}^{l} \vartheta N_{j_1j_1}(\psi_c) R_{j_1}\,d{\alpha}\\[1mm]
&\qquad\qquad\qquad\!\!\! + \int_{\R} \partial^{l}_{\alpha} R_{j_1}\, \rho_{j_1}^{l} \partial^{l}_{\alpha} \vartheta N_{j_1j_1}(\psi_c) W_{j_1j_2}({\Psi}, \mathcal{R}) R_{j_2}\,d{\alpha} \,\Big),
\\[3mm] 
I_6 =\;  &  \veps^{3} \sum_{j_1,j_2  \in \{\pm2\}} \Big(\,
 \int_{\R} \partial^{l}_{\alpha} W_{j_1j_2}({\Psi}, \mathcal{R}) R_{j_2} \, \rho_{j_1}^{l} \partial_{\alpha}^{l} \vartheta N_{j_1-j_1}(\psi_c) R_{-j_1}\,d{\alpha}\; \\[2mm]
&\qquad\qquad\qquad\!\!\! + \int_{\R} \partial^{l}_{\alpha} R_{j_1}\, \rho_{j_1}^{l} \partial^{l}_{\alpha} \vartheta N_{j_1-j_1}(\psi_c) W_{-j_1j_2}({\Psi}, \mathcal{R}) R_{j_2}\,d{\alpha}\, \Big), 
\\[3mm]
I_7 =\;  &  \veps^{3}  \sum_{j_1,j_2 \in \{\pm2\}}\!  \Big(\,
\int_{\R}  \partial^{l}_{\alpha} \mathcal{M}_{j_1}(\Psi, \mathcal{R})\, \rho_{j_1}^{l} \partial^{l}_{\alpha} \vartheta N_{j_1j_2}(\psi_c) R_{j_2}\,d{\alpha}\\[1mm]
&\qquad\qquad\qquad\!\!\! + \int_{\R} \partial^{l}_{\alpha} R_{j_1}\, \rho_{j_1}^{l} \partial^{l}_{\alpha} \vartheta N_{j_1j_2}(\psi_c) \mathcal{M}_{j_2}(\Psi, \mathcal{R})\,d{\alpha} \,\Big). 
\end{align*}
First, we analyze $I_2$. To extract all terms with more than $l$ spatial derivatives falling on $R_{2}$ or $R_{-2}$ we use Leibniz's rule, integration by parts, \eqref{q1}--\eqref{frakm}, \eqref{thetaN}, \eqref{N0s}, \eqref{thetaNa}, \eqref{Gaj}--\eqref{Maj} and \eqref{thetaN2a} to obtain
\begin{align*}
I_2
=\;& \veps^{2} \sum_{j_1, j_2 \in \{\pm2\}} \Big(\,  \int_{\R}  \partial^{l}_{\alpha} Q_{j_1j_2}({\psi}_c) R_{j_2}\, \rho_{j_1}^{l} \vartheta N_{j_1j_1}(\psi_c) \partial^{l}_{\alpha} R_{j_1}\,d{\alpha}\\[1mm]
&\qquad\qquad\qquad+ l  \int_{\R} \partial^{l}_{\alpha} Q_{j_1j_2}({\psi}_c) R_{j_2}\, \rho_{j_1}^{l} \vartheta N_{j_1j_1}^{1}(\partial_{\alpha} \psi_c) \partial^{l-1}_{\alpha} R_{j_1}\,d{\alpha}\\[3mm]
&\qquad\qquad\qquad+ \int_{\R} \partial^{l}_{\alpha} R_{j_1} \, \rho_{j_1}^{l} \vartheta N_{j_1j_1}(\psi_c)  \partial^{l}_{\alpha} Q_{j_1j_2}({\psi}_c) R_{j_2}\,d{\alpha} \\[3mm]
&\qquad\qquad\qquad+ l  \int_{\R} \partial^{l}_{\alpha} R_{j_1}\, \rho_{j_1}^{l} \vartheta N_{j_1j_1}^{1}(\partial_{\alpha} \psi_c) \partial^{l-1}_{\alpha} Q_{j_1j_2}({\psi}_c) R_{j_2}\,d{\alpha}
\,\Big)
\end{align*}
\begin{align*}
&+\, \veps^{2}\,\mathcal{O}(\mathcal{E}_{s}+1)
\\[2mm]
\;& + 2 \veps^{2} \sum_{j_1, j_2 \in \{\pm2\}} \Big(\,  \int_{\R}  \partial^{l}_{\alpha} Q_{j_1j_2}({\psi}_c) R_{j_2}\, (\rho^{l} \mathfrak{N}^{1})_{j_1j_1}^s (\psi_c) \partial^{l}_{\alpha} R_{j_1}\,d{\alpha}
\\[1mm]
&\qquad\qquad\qquad\quad\;\; - \int_{\R}  \partial^{l-1}_{\alpha} Q_{j_1j_2}({\psi}_c) R_{j_2}\, (\rho^{l} \mathfrak{N}^{2})_{j_1j_1}^a (\psi_c) \partial^{l}_{\alpha} R_{j_1}\,d{\alpha}
\\[3mm]
&\qquad\qquad\qquad\quad\;\; + l  \int_{\R} \partial^{l}_{\alpha} Q_{j_1j_2}({\psi}_c) R_{j_2}\, (\rho^{l} \mathfrak{N}^{1})_{j_1j_1}^a (\partial_{\alpha} \psi_c) \partial^{l-1}_{\alpha} R_{j_1}\,d{\alpha}
\,\Big)
\\[2mm]
&+\, \veps^{2}\,\mathcal{O}(\mathcal{E}_{s}+1)
\\[2mm]
 =\;& \sum_{i=1}^{4}\,  2\veps^{2} \sum_{j_1, j_2 \in \{\pm2\} } \int_{\R}  Q^i_{j_1j_2} (\psi_c)  \partial^{l}_{\alpha} R_{j_2}\, S_{j_1}^l (\psi_c) \partial^{l}_{\alpha} R_{j_1}\,d{\alpha} 
\,+\, \veps^{2}\,\mathcal{O}(\mathcal{E}_{s}+1)
\\[1mm]
=:\; & \sum_{i=1}^4 \, I_2^i \,+\, \veps^{2}\,\mathcal{O}(\mathcal{E}_{s}+1) \,,
\end{align*}
as long as $\eps^{5/2} \Vert \mathcal{R}_{1} \Vert_{(L^{2})^2}, \eps^{5/2} \Vert \mathcal{R}_{2} \Vert_{(H^{{s}})^2} \ll 1$, uniformly on compact subsets of $\mathcal{B}$, where 
\begin{equation*}
S_{j_1}^l (\psi_c) = (\rho^{l} \mathfrak{N}^{1})_{j_1j_1}^s (\psi_c) + l (\partial^{-1}_{\alpha} \rho^{l} \mathfrak{N}^{1})_{j_1j_1}^s(\partial_{\alpha}\psi_c)\,.
\end{equation*}
Due to \eqref{int und theta}--\eqref{int und theta-1}, we have
\begin{align*}
 I_2^1
=\;& -2\veps^{2}\sum_{j_1 \in \{\pm2\}} \int_{\R}  \partial_{\alpha} ( \psi_c\, \partial^{l}_{\alpha} R_{j_1})\, S_{j_1}^l (\psi_c) \partial^{l}_{\alpha} R_{j_1}\,d{\alpha} \;+\; \veps^{2} \mathcal{O}(\mathcal{E}_{s}+1)
\\[2mm]
=\;& -2 \veps^{2}\sum_{j_1 \in \{\pm2\}} \int_{\R} \partial^{l+1}_{\alpha} R_{j_1}\, \psi_c\, S_{j_1}^l (\psi_c) \partial^{l}_{\alpha} R_{j_1}\,d{\alpha} \;+\; \veps^{2} \mathcal{O}(\mathcal{E}_{s}+1)\,.
\end{align*}
Since $\psi_c$ is real-valued and $S_{j_1}^l  (\psi_c)$ symmetric, using \eqref{N0skomm}, \eqref{Gajkomm} and \eqref{dx1} yields
\begin{align*}
 I_2^1
=\;& \veps^{2} \sum_{j_1 \in \{\pm2\}} \int_{\R} \partial^{l}_{\alpha} R_{j_1}\, [\partial_{\alpha},\,  \psi_c\, S_{j_1}^l (\psi_c)] \partial^{l}_{\alpha} R_{j_1}\,d{\alpha} \;+\; \veps^{2} \mathcal{O}(\mathcal{E}_{s}+1)
\\[2mm]
=\;& \veps^{2} \sum_{j_1 \in \{\pm2\}} \Big( \int_{\R} \partial^{l}_{\alpha} R_{j_1}\, \partial_{\alpha}  \psi_c\, S_{j_1}^l  (\psi_c) \partial^{l}_{\alpha} R_{j_1}\,d{\alpha}  \\[1mm]
&\qquad\qquad\quad  + \int_{\R} \partial^{l}_{\alpha} R_{j_1}\, \psi_c\, S_{j_1}^l  (\partial_{\alpha} \psi_c) \partial^{l}_{\alpha} R_{j_1}\,d{\alpha} \Big)
+ \veps^{2} \mathcal{O}(\mathcal{E}_{s}+1)\\[2mm]
=\;& \veps^{2} \mathcal{O}(\mathcal{E}_{s}+1) \,,
\end{align*}
as long as $\eps^{5/2} \Vert \mathcal{R}_{1} \Vert_{(L^{2})^2}, \eps^{5/2} \Vert \mathcal{R}_{2} \Vert_{(H^{{s}})^2} \ll 1$, uniformly on compact subsets of $\mathcal{B}$.

Because of \eqref{N0skomm}, \eqref{Gajkomm}, \eqref{int und theta}--\eqref{int und theta-1} and
\eqref{dx1}--\eqref{dx2} we have
\begin{align*}
 I_2^2
=\;& \veps^{2} \sum_{j_1 \in \{\pm2\}}  \int_{\R}  (K_0 \sigma^{-1} \partial_{\alpha} \psi_c)  \sigma^{-1} \partial_{\alpha}^{l+1} (R_{-2}-R_{2}) S_{j_1}^l  (\psi_c) \partial^{l}_{\alpha} R_{j_1}\,d{\alpha} + \veps^{2} \mathcal{O}(\mathcal{E}_{s}+1)            
\end{align*}
\begin{align*}
=\;& \veps^{2} \sum_{j_1 \in \{\pm2\}}  \int_{\R}  \mathrm{sgn}(j_1) (K_0 \sigma^{-1} \partial_{\alpha} \psi_c)  \sigma^{-1} \partial_{\alpha}^{l+1} R_{-j_1} S_{j_1}^l  (\psi_c) \partial^{l}_{\alpha} R_{j_1}\,d{\alpha} + \veps^{2} \mathcal{O}(\mathcal{E}_{s}+1)            
\\[2mm]
=\;& - \frac{\veps^{2}}{2} \int_{\R} \partial_{\alpha}^{l} (R_{-2}-R_{2}) (K_0 \sigma^{-1} \partial_{\alpha} \psi_c) (S_{-2}^l +S_2^l ) (\psi_c) \sigma^{-1} \partial^{l+1}_{\alpha} ( R_{-2}+  R_{2})\,d{\alpha} \\[2mm]
\;& +\; \veps^{2} \mathcal{O}(\mathcal{E}_{s}+1) \,,
\end{align*}
as long as $\eps^{5/2} \Vert \mathcal{R}_{1} \Vert_{(L^{2})^2},\, \eps^{5/2} \Vert \mathcal{R}_{2} \Vert_{(H^{{s}})^2} \ll 1$. There holds
\begin{align}
\label{domega} 
\sigma^{-1} \partial_{\alpha}  
=\;& (1-b \partial_{\alpha}^2)^{-1} i\omega \,. 
\end{align} 
Moreover, \eqref{Rmp2}--\eqref{dx-1R2} yield 
\begin{align} 
i\omega ( R_{-2}+  R_{2})
\,=&\,  
- \partial_t (R_{-2}- R_{2}) 
\nonumber
\\[2mm] 
&\, 
- \veps \vartheta^{-1} \partial_{\alpha}(( \psi_c +\eps  \partial^{-2}_{\alpha} g_{+}(\Psi_2^h, \mathcal{R}_2))\, \vartheta (R_{-2}-R_{2})) 
\nonumber
\\[2mm] 
&\,  - \frac{\veps}{2} \vartheta^{-1} \partial_{\alpha} ([\sigma, \partial^{-2}_{\alpha} \vartheta(R_{-2}+ R_{2})] \sigma^{-1} (\partial_{\alpha}^2 \psi_c + \eps g_{-}(\Psi_2^h, \mathcal{R}_2)))
\nonumber
\\[2mm] 
&\,  
+ \veps f(\Psi,\mathcal{R}_{2}) \,,
\label{dxdt}
\end{align}
with a function $f$ satisfying
\begin{align} 
\label{dxdt-est}
\|f(\Psi,\mathcal{R}_{2})\|_{H^{s}} \,=\, \mathcal{O}((\mathcal{E}_s +1)^{1/2})\,,
\end{align}
as long as $\eps^{5/2} \Vert \mathcal{R}_{1} \Vert_{(L^{2})^2}, \eps^{5/2} \Vert \mathcal{R}_{2} \Vert_{(H^{{s}})^2} \ll 1$. 
Hence, with the help of \eqref{sigma-id}, \eqref{sigmak1}, \eqref{dx-2R2},  \eqref{asnjj}, \eqref{N0skomm}, \eqref{nus}, \eqref{Gajkomm}, \eqref{int und theta}--\eqref{int und theta-1} and \eqref{dx1} we obtain
\begin{align*}
 I_2^2
=\;&  \frac{\veps^{2}}{2}  \int_{\R}  \partial_{\alpha}^{l} (R_{-2}-R_{2}) (K_0  \sigma^{-1} \partial_{\alpha} \psi_c) (S_{-2}^l +S_2^l ) (\psi_c) (1-b \partial_{\alpha}^2)^{-1} \partial_t \partial^{l}_{\alpha} ( R_{-2}-  R_{2})\,d{\alpha} 
\\[2mm]
\;& +\;
\frac{\veps^{3}}{2}  \int_{\R}  \partial_{\alpha}^{l} (R_{-2}-R_{2}) (K_0  \sigma^{-1} \partial_{\alpha} \psi_c) (S_{-2}^l +S_2^l ) (\psi_c) (1-b \partial_{\alpha}^2)^{-1} \,       \\[2mm]
\;& \hspace{1.5cm} \times 
\big(( \psi_c +\eps  \partial^{-2}_{\alpha} g_{+}(\Psi_2^h, \mathcal{R}_2))\,   \partial^{l+1}_{\alpha} ( R_{-2}- R_{2})\big)\,d{\alpha} 
\\[2mm]
\;& +\;
\frac{\veps^{3}}{4}  \int_{\R}  \partial_{\alpha}^{l} (R_{-2}-R_{2}) (K_0  \sigma^{-1} \partial_{\alpha} \psi_c) (S_{-2}^l +S_2^l ) (\psi_c) (1-b \partial_{\alpha}^2)^{-1} \,        \\[2mm]
\;& \hspace{1.5cm} \times 
\big(\sigma^{-1} (\partial_{\alpha}^2 \psi_c + \eps g_{-}(\Psi_2^h, \mathcal{R}_2))\,\sigma \partial^{l-1}_{\alpha} ( R_{-2}+ R_{2})\big)\,d{\alpha}
\\[2mm]
\;& +\; \veps^{2} \mathcal{O} (\mathcal{E}_s +1 )
\\[2mm]
=\;&  \frac{\veps^{2}}{4}  \frac{d}{dt} \int_{\R} \partial_{\alpha}^{l} (R_{-2}-R_{2}) (K_0  \sigma^{-1} \partial_{\alpha}\psi_c) (S_{-2}^l +S_2^l ) (\psi_c) (1-b \partial_{\alpha}^2)^{-1}\partial^{l}_{\alpha} ( R_{-2}-  R_{2}) \,d{\alpha}  \\[2mm]
\;& + \; \veps^{2} \mathcal{O}(\mathcal{E}_s +1 )\,,
\end{align*}
as long as $\eps^{5/2} \Vert \mathcal{R}_{1} \Vert_{(L^{2})^2}, \eps^{5/2} \Vert \mathcal{R}_{2} \Vert_{(H^{{s}})^2} \ll 1$,
uniformly on compact subsets of $\mathcal{B}$.

For $I_2^3$ and $I_2^4$ we can also use  \eqref{sigma-id}, \eqref{sigmak1}, \eqref{dx-2R2}, \eqref{asnjj}, \eqref{N0skomm}, \eqref{nus}, \eqref{Gajkomm}, \eqref{int und theta}--\eqref{int und theta-1}  and \eqref{dx1}--\eqref{dx2} to deduce
\begin{align*}
 I_2^3
=\;& - \veps^{2} \sum_{j_1 \in \{\pm2\}} \int_{\R} b (\sigma^{-1} \partial^{2}_{\alpha} \psi_c) K_0  \sigma^{-1} \partial^{l+2}_{\alpha} (R_{-2}-R_2)\, S_{j_1}^l (\psi_c) \partial^{l}_{\alpha} R_{j_1}\,d{\alpha} + \veps^{2} \mathcal{O}(\mathcal{E}_{2,s}) 
\\[2mm]
=\;& \frac{\veps^{2}}{2} \int_{\R}  \partial_{\alpha}^{l} (R_{-2}-R_{2}) (\sigma^{-1} \partial^{2}_{\alpha} \psi_c) (S_{-2}^l +S_2^l ) (\psi_c) b K_0  \sigma^{-1} \partial^{l+2}_{\alpha}
( R_{-2}+  R_{2})\,d{\alpha}\\[3mm]
\;& + \veps^{2} \mathcal{O}(\mathcal{E}_{s}+1)  
\\[3mm]
=\;& \veps^{2} \mathcal{O}(\mathcal{E}_{s}+1)  
\end{align*}
and
\begin{align*}
 I_2^4
=\;& \veps^{2} \sum_{j_1\in \{\pm2\}} \int_{\R}  \mathrm{sgn}(j_1) ([\sigma, \partial^{l-1}_{\alpha} (R_{-2}+R_2)] \sigma^{-1} \partial^{2}_{\alpha} \psi_c)
\, S_{j_1}^l (\psi_c) \partial^{l}_{\alpha} R_{j_1}\,d{\alpha} \\[2mm]
&
+ \veps^{2} \mathcal{O}(\mathcal{E}_{s}+1) 
\\[2mm]
=\;& \veps^{2} \sum_{j_1\in \{\pm2\}} \int_{\R}  \mathrm{sgn}(j_1)  (\sigma^{-1} \partial^{2}_{\alpha} \psi_c) \sigma \partial^{l-1}_{\alpha} (R_{-2}+R_2)\, S_{j_1}^l (\psi_c) \partial^{l}_{\alpha} R_{j_1}\,d{\alpha} \\[2mm]
&
+ \veps^{2} \mathcal{O}(\mathcal{E}_{s}+1) 
\\[2mm]
=\;& - \frac{\veps^{2}}{2} \int_{\R}  \partial_{\alpha}^{l} (R_{-2}-R_{2})  (\sigma^{-1} \partial^{2}_{\alpha} \psi_c) (S_{-2}^l +S_2^l ) (\psi_c) \sigma \partial^{l-1}_{\alpha} (R_{-2}+R_2)\,d{\alpha} \\[3mm]
& + \veps^{2} \mathcal{O}(\mathcal{E}_{s}+1) 
\\[3mm]
=\;& \veps^{2} \mathcal{O}(\mathcal{E}_{s}+1)  \,,
\end{align*}
as long as $\eps^{5/2} \Vert \mathcal{R}_{1} \Vert_{(L^{2})^2},\, \eps^{5/2} \Vert \mathcal{R}_{2} \Vert_{(H^{{s}})^2} \ll 1$,
uniformly on compact subsets of $\mathcal{B}$.

Next, we examine $I_3$. Using Leibniz's rule, integration by parts, \eqref{q1}--\eqref{frakm}, \eqref{thetaNb}, \eqref{N0s},  
\eqref{Gj}--\eqref{Mj} and \eqref{int und theta}--\eqref{int und theta-1}we conclude
\begin{align*}
I_3
=\;& \veps^{2} \sum_{j_1, j_2 \in \{\pm2\}} \Big(\,   \int_{\R} \partial^{l}_{\alpha} Q_{j_1j_2}({\psi}_c) R_{j_2} \, \rho_{j_1}^{l} \vartheta 
N_{j_1-j_1}^{1} (\psi_c) \partial_{\alpha}^{l} R_{-j_1}\,d{\alpha}\; 
\\[3mm]
&\qquad\qquad\qquad+ \int_{\R} \partial^{l}_{\alpha} R_{-j_1}\, \rho_{-j_1}^{l} \vartheta N_{-j_1 j_1}^{1} (\psi_c) \partial_{\alpha}^{l} Q_{j_1j_2}({\psi}_c) R_{j_2}\,d{\alpha}
\,\Big)
+ \veps^{2}\,\mathcal{O}(\mathcal{E}_{s}+1)
\\[2mm]
 =\;& \sum_{i=1}^{4}\,  2\veps^{2} \sum_{j_1, j_2 \in \{\pm2\} } \int_{\R} Q^i_{j_1j_2} (\psi_c)  \partial^{l}_{\alpha} R_{j_2}\, (\rho^{l} \mathfrak{N}^{1})_{j_1-j_1}^s (\psi_c) \partial^{l}_{\alpha} R_{-j_1}\,d{\alpha}
+ \veps^{2}\,\mathcal{O}(\mathcal{E}_{s}+1)
\\[2mm]
=\;&  2\veps^{2} \sum_{j_1, j_2 \in \{\pm2\} } \int_{\R} Q^1_{j_1j_2} (\psi_c)  \partial^{l}_{\alpha} R_{j_2}\, (\rho^{l} \mathfrak{N}^{1})_{j_1-j_1}^s (\psi_c) \partial^{l}_{\alpha} R_{-j_1}\,d{\alpha} 
\\[2mm] 
&+
\sum_{i=2}^{4}\,  2\veps^{2} \sum_{j_1, j_2 \in \{\pm2\} } \int_{\R} Q^i_{j_1j_2} (\psi_c)  \partial^{l}_{\alpha} R_{j_2}\, G_{j_1}^l  (\psi_c) \partial^{l}_{\alpha} R_{-j_1}\,d{\alpha}
+ \veps^{2}\,\mathcal{O}(\mathcal{E}_{s}+1)
\end{align*}
\begin{align*} 
=:\; & \sum_{i=1}^4 \, I_3^i \,+\, \veps^{2}\,\mathcal{O}(\mathcal{E}_{s}+1)\,,
\end{align*}
as long as $\eps^{5/2} \Vert \mathcal{R}_{1} \Vert_{(L^{2})^2}, \eps^{5/2} \Vert \mathcal{R}_{2} \Vert_{(H^{{s}})^2} \ll 1$, uniformly on compact subsets of $\mathcal{B}$. 

Because of \eqref{N0skomm}, \eqref{int und theta}--\eqref{int und theta-1} and \eqref{dx3} we obtain
\begin{align*}
 I_3^1
=\;&  -2\veps^{2}\sum_{j_1 \in \{\pm2\}} \int_{\R} \partial^{l+1}_{\alpha} R_{-j_1}\, \psi_c\, (\rho^{l} \mathfrak{N}^{1})_{-j_1j_1}^s(\psi_c) \partial^{l}_{\alpha} R_{j_1}\,d{\alpha} + \veps^{2} \mathcal{O}(\mathcal{E}_{s}+1) \\[2mm]
=\;& 2\veps^{2}  \int_{\R} \partial^{l}_{\alpha} R_{2}\, [\partial_{\alpha}, \psi_c\, (\rho^{l} \mathfrak{N}^{1})_{2-2}^s(\psi_c)] \partial^{l}_{\alpha} R_{-2}\,d{\alpha} + \veps^{2} \mathcal{O}(\mathcal{E}_{s}+1) 
\\[2mm]
=\;& 2\veps^{2}  \int_{\R} \partial^{l}_{\alpha} R_{2}\, \partial_{\alpha} \psi_c (\rho^{l} \mathfrak{N}^{1})_{2-2}^s(\psi_c) \partial^{l}_{\alpha} R_{-2}\,d{\alpha}  
\\[2mm]
& + 2\veps^{2}  \int_{\R} \partial^{l}_{\alpha} R_{2}\, \psi_c\, (\rho^{l} \mathfrak{N}^{1})_{2-2}^s( \partial_{\alpha} \psi_c) \partial^{l}_{\alpha} R_{-2}\,d{\alpha} 
+ \veps^{2} \mathcal{O}(\mathcal{E}_{s}+1) \\[2mm]
=\;& \veps^{2} \mathcal{O}(\mathcal{E}_{s}+1) \,,
\end{align*}
as long as $\eps^{5/2} \Vert \mathcal{R}_{1} \Vert_{(L^{2})^2},\, \eps^{5/2}\Vert \mathcal{R}_{2} \Vert_{(H^{{s}})^2} \ll 1$, uniformly on compact subsets of $\mathcal{B}$.

Due to \eqref{Gjkomm}, \eqref{int und theta}--\eqref{int und theta-1},
\eqref{dx1}--\eqref{dx2} and \eqref{domega}--\eqref{dxdt-est}, we have
\begin{align*}
 I_3^2
=\;& \veps^{2}  \sum_{j_1 \in \{\pm2\}}  \int_{\R}  (K_0 \sigma^{-1} \partial_{\alpha} \psi_c)  \sigma^{-1} \partial_{\alpha}^{l+1} (R_{-2}-R_{2}) G_{j_1}^l  (\psi_c) \partial^{l}_{\alpha} R_{-j_1}\,d{\alpha} + \veps^{2} \mathcal{O}(\mathcal{E}_{2,s})            \\[2mm]
=\;& - \frac{\veps^{2}}{2} \int_{\R} \partial_{\alpha}^{l} (R_{-2}-R_{2}) (K_0 \sigma^{-1} \partial_{\alpha} \psi_c) (G_{-2}^l + G_{2}^l ) (\psi_c) \sigma^{-1} \partial^{l+1}_{\alpha} ( R_{-2}+  R_{2})\,d{\alpha} \\[3mm]
\;& +\; \veps^{2} \mathcal{O}(\mathcal{E}_{s}+1) 
\end{align*}
\begin{align*}
=\;&  \frac{\veps^{2}}{4}  \frac{d}{dt} \int_{\R} \partial_{\alpha}^{l} (R_{-2}-R_{2}) (K_0  \sigma^{-1} \partial_{\alpha}\psi_c) (G_{-2}^l + G_{2}^l ) (\psi_c) (1-b \partial_{\alpha}^2)^{-1} \partial^{l}_{\alpha} ( R_{-2}-  R_{2}) \,d{\alpha}\\[3mm]
\;& +  \veps^{2} \mathcal{O}(\mathcal{E}_s +1 )\,,
\end{align*}
as long as $\eps^{5/2} \Vert \mathcal{R}_{1} \Vert_{(L^{2})^2},\, \eps^{5/2} \Vert \mathcal{R}_{2} \Vert_{(H^{{s}})^2} \ll 1$, uniformly on compact subsets of $\mathcal{B}$.

For $I_3^3$ and $I_3^4$ we  can use \eqref{sigma-id}, \eqref{sigmak1}, \eqref{dx-2R2}, \eqref{asomega}, \eqref{asqj-j}, \eqref{Gjkomm}, \eqref{nus-}, \eqref{int und theta}--\eqref{int und theta-1} and \eqref{dx1}--\eqref{dx2} again to deduce
\begin{align*}
 I_3^3
=\;& - \veps^{2} \sum_{j_1 \in \{\pm2\}} \int_{\R} b (\sigma^{-1} \partial^{2}_{\alpha} \psi_c) K_0  \sigma^{-1} \partial^{l+2}_{\alpha} (R_{-2}-R_2)\, G_{j_1}^l (\psi_c) \partial^{l}_{\alpha} R_{-j_1}\,d{\alpha}\\[2mm]
\;& + \veps^{2} \mathcal{O}(\mathcal{E}_{s}+1) 
\\[2mm]
=\;& \frac{\veps^{2}}{2} \int_{\R} \partial_{\alpha}^{l} (R_{-2}-R_{2}) (\sigma^{-1} \partial^{2}_{\alpha} \psi_c) (G_{-2}^l + G_{2}^l )(\psi_c)\, b K_0  \sigma^{-1} \partial^{l+2}_{\alpha}
( R_{-2}+  R_{2})\,d{\alpha}\\[2mm]
\;& + \veps^{2} \mathcal{O}(\mathcal{E}_{s}+1)  
\\[2mm]
=\;& \veps^{2} \mathcal{O}(\mathcal{E}_{s}+1)  
\end{align*}
and
\begin{align*}
 I_3^4
=\;& \veps^{2} \sum_{j_1\in \{\pm2\}} \int_{\R}  \mathrm{sgn}(j_1) ([\sigma, \partial^{l-1}_{\alpha} (R_{-2}+R_2)] \sigma^{-1} \partial^{2}_{\alpha} \psi_c)  \, G_{j_1}^l (\psi_c) \partial^{l}_{\alpha} R_{-j_1}\,d{\alpha}\\[2mm]
\;&  + \veps^{2} \mathcal{O}(\mathcal{E}_{s}+1) 
\\[2mm]
=\;&- \frac{\veps^{2}}{2} \int_{\R}  \partial_{\alpha}^{l} (R_{-2}-R_{2})  (\sigma^{-1} \partial^{2}_{\alpha} \psi_c) (G_{-2}^l + G_{2}^l )(\psi_c) \sigma \partial^{l-1}_{\alpha} (R_{-2}+R_2)\,d{\alpha} \\[2mm]
\;& + \veps^{2} \mathcal{O}(\mathcal{E}_{s}+1) 
 \\[2mm]
=\;&  \veps^{2} \mathcal{O}(\mathcal{E}_{s}+1) \,,
\end{align*}
as long as $\eps^{5/2} \Vert \mathcal{R}_{1} \Vert_{(L^{2})^2},\, \eps^{5/2} \Vert \mathcal{R}_{2} \Vert_{(H^{{s}})^2} \ll 1$, uniformly on compact subsets of $\mathcal{B}$.

Now, we investigate $I_4$. Because of \eqref{estcalMmp2} we have
\begin{align*}
I_4
 =\;& \sum_{\genfrac{}{}{0pt}{}{j_1, j_2 \in \{\pm2\},} {i \in \{1, \ldots, 4\}}} \veps^{2}\int_{\R} \rho_{j_1}^{l} \partial^{l}_{\alpha} W^i_{j_1j_2} (\Psi, \mathcal{R}) R_{j_2}\, \partial^{l}_{\alpha} R_{j_1}\,d{\alpha} + \veps^{2} \mathcal{O}(\mathcal{E}_{s}+1)  
\\[2mm]
=:\;& \sum_{i=1}^4 \, I_4^i + \veps^{2} \mathcal{O}(\mathcal{E}_{s}+1) 
\,,
\end{align*}
as long as $\eps^{5/2} \Vert \mathcal{R}_{1} \Vert_{(L^{2})^2},\, \eps^{5/2} \Vert \mathcal{R}_{2} \Vert_{(H^{{s}})^2} \ll 1$, uniformly on compact subsets of $\mathcal{B}$.
Analogously to the case of $I_2$ we conclude
\begin{align*}
  I_4^1
=\;& -\veps^{2}\sum_{j_1 \in \{\pm2\}} \int_{\R} \partial_{\ua}^{-2} g_{+}(\Psi_{2}^{h}, \mathcal{R}_{2})\, \partial^{l+1}_{\alpha} R_{j_1}\, \partial^{l}_{\alpha} R_{j_1}\,d{\alpha}
\\[2mm]\;& 
\,+\,\veps^{2} \mathcal{O}(\mathcal{E}_s +1)
\end{align*}
\begin{align*}
=\;& \frac{\veps^{2}}{2} \sum_{j_1 \in \{\pm2\}} \int_{\R} \partial_{\ua}^{-1} g_{+}(\Psi_{2}^{h}, \mathcal{R}_{2})\, (\partial^{l}_{\alpha} R_{j_1})^{2}\,d{\alpha} 
\,+\,\veps^{2} \mathcal{O}(\mathcal{E}_s +1 ) \\[2mm]
=\;& \; \veps^{2} \mathcal{O}(\mathcal{E}_s +1)
\end{align*}
and
\begin{align*}
  I_4^2
=\;& -\frac{\veps^{2}}{2} \int_{\R} \partial_{\alpha}^{l}(R_{-2}-R_{2}) (K_0 \sigma^{-1} \partial_{\ua}^{-1} g_{-}(\Psi_{2}^{h}, \mathcal{R}_{2})+c(\Psi, \mathcal{R}))
\,       \\[2mm]
\;& \hspace{1.35cm} \times 
\sigma^{-1} \partial^{l+1}_{\alpha} ( R_{-2}+ R_{2})\,d{\alpha} 
+ \; \veps^{2} \mathcal{O}(\mathcal{E}_s +1)\\[2mm]
=\;&  \frac{\veps^{2}}{4}  \frac{d}{dt} \int_{\R} \partial_{\alpha}^{l} (R_{-2}-R_{2}) (K_0  \sigma^{-1} \partial_{\ua}^{-1} g_{-}(\Psi_{2}^{h}, \mathcal{R}_{2})+c(\Psi, \mathcal{R})) 
 \\[2mm]
\;& \hspace{1.35cm} \times
(1-b \partial_{\alpha}^2)^{-1} \partial^{l}_{\alpha} ( R_{-2}-  R_{2}) \,d{\alpha}   
+ \; \veps^{2} \mathcal{O}(\mathcal{E}_s +1)\,,
\end{align*}
as long as $\eps^{5/2} \Vert \mathcal{R}_{1} \Vert_{(L^{2})^2},\, \eps^{5/2} \Vert \mathcal{R}_{2} \Vert_{(H^{{s}})^2} \ll 1$, uniformly on compact subsets of $\mathcal{B}$.
Furthermore, we deduce
\begin{align*}
  I_4^3
=\;& \frac{\veps^{2}}{2} \int_{\R}  \partial_{\alpha}^{l}(R_{-2}-R_{2}) (\sigma^{-1} g_{-}(\Psi_{2}^{h}, \mathcal{R}_{2})) 
\,bK_0 \sigma^{-1} \partial_{\alpha}^{l+2} (R_{-2}+R_{2})\,d{\alpha} \\[2mm]  
\;&
+ \; \veps^{2} \mathcal{O}(\mathcal{E}_s +1)\\[2mm]
=\;& - \frac{\veps^{2}}{4}  \frac{d}{dt} \int_{\R} \partial_{\alpha}^{l} (R_{-2}-R_{2})(\sigma^{-1}  g_{-}(\Psi_{2}^{h}, \mathcal{R}_{2}))\, b K_0 \partial_{\alpha}  \\[2mm]
\;& \hspace{1.8cm} \times
(1-b \partial_{\alpha}^2)^{-1} \partial^{l}_{\alpha} ( R_{-2}-  R_{2}) \,d{\alpha} 
+ \; \veps^{2} \mathcal{O}(\mathcal{E}_s +1)
\end{align*}
as long as $\eps^{5/2} \Vert \mathcal{R}_{1} \Vert_{(L^{2})^2},\, \eps^{5/2} \Vert \mathcal{R}_{2} \Vert_{(H^{{s}})^2} \ll 1$, uniformly on compact subsets of $\mathcal{B}$. 
Moreover, due to \eqref{sigma-id}, \eqref{sigmak1}, \eqref{dx-2R2}, we have
\begin{align*}
  I_4^4
=\;& -\frac{\veps^{2}}{2} \int_{\R} \sigma\partial_{\ua}^{l-1}(R_{ -2}+R_{2}) (\sigma^{-1} g_{-}(\Psi_{2}^{h}, \mathcal{R}_{2}))\,
\partial_{\alpha}^{l} (R_{-2}-R_{2}) \,d{\alpha}        
\\[2mm]
&+ \; \veps^{2} \mathcal{O}(\mathcal{E}_s +1)
\end{align*}
and because of 
\begin{align}
\sigma \partial_{\alpha}^{-1}
=\;& (-K_0  \partial_{\alpha})^{-1}  i\omega  
\end{align}
and \eqref{dxdt}--\eqref{dxdt-est} we obtain
\begin{align*}
I_4^4
=\;& \frac{\veps^{2}}{4}  \frac{d}{dt} \int_{\R} \partial_{\alpha}^{l} (R_{-2}-R_{2})\,  \sigma^{-1}  g_{-}(\Psi_{2}^{h}, \mathcal{R}_{2}) 
(-K_0 \partial_{\alpha})^{-1} \partial^{l}_{\alpha} ( R_{-2}-  R_{2}) \,d{\alpha}   \\[2mm]
& + \; \veps^{2} \mathcal{O}(\mathcal{E}_s +1)\,,
\end{align*}
as long as $\eps^{5/2} \Vert \mathcal{R}_{1} \Vert_{(L^{2})^2}, \eps^{5/2} \Vert \mathcal{R}_{2} \Vert_{(H^{{s}})^2} \ll 1$,
uniformly on compact subsets of $\mathcal{B}$.

Next, we consider $I_5$ and $I_6$. Analogously to the cases of $I_2$ and $I_3$ we obtain 
\begin{align*}
I_5 
&\,=\, \sum_{\genfrac{}{}{0pt}{}{j_1, j_2 \in \{\pm2\},} {i \in \{1, \ldots, 4\}}} 2\veps^{3} \int_{\R}   W^i_{j_1j_2} (\Psi, \mathcal{R})  \partial^{l}_{\alpha} R_{j_2}\, S_{j_1}^l  (\psi_c) \partial^{l}_{\alpha} R_{j_1}\,d{\alpha} \\[2mm]
&\,\quad\, + \veps^{3}\,\mathcal{O}(\mathcal{E}_s +1)
\\[2mm]
&\;=:\, \sum_{i=1}^4 \, I_5^i +\, \veps^{3}\,\mathcal{O}(\mathcal{E}_s +1) \,,
\\[2mm]
I_6
&\,=\, \sum_{\genfrac{}{}{0pt}{}{j_1, j_2 \in \{\pm2\},} {i \in \{1, \ldots, 4\}}} 2\veps^{3} \int_{\R}   W^i_{j_1j_2} (\Psi, \mathcal{R})  \partial^{l}_{\alpha} R_{j_2}\, (\rho^{l} \mathfrak{N}^{1})_{j_1-j_1}^s (\psi_c) \partial^{l}_{\alpha} R_{-j_1}\,d{\alpha} \\[2mm]
&\,\quad\,+ \veps^{3}\,\mathcal{O}(\mathcal{E}_s +1)
\\[2mm]
&\,=\,  \sum_{i=1}^4 \, I_6^i +\, \veps^{3}\,\mathcal{O}(\mathcal{E}_s +1)
\end{align*}
as long as $\eps^{5/2} \Vert \mathcal{R}_{1} \Vert_{(L^{2})^2},\, \eps^{5/2}\Vert \mathcal{R}_{2} \Vert_{(H^{{s}})^2} \ll 1$,
uniformly on compact subsets of $\mathcal{B}$
and therefore
\begin{align*}
  I_5^1
=\;& \veps^{3}\sum_{j_1 \in \{\pm2\}} \Big( \int_{\R} \partial^{l}_{\alpha} R_{j_1}\, \partial_{\ua}^{-1} g_{+}(\Psi_{2}^{h}, \mathcal{R}_{2})\, S_{j_1}^l  (\psi_c) \partial^{l}_{\alpha} R_{j_1}\,d{\alpha}  \\[1mm]
&\qquad\qquad\quad  + \int_{\R} \partial^{l}_{\alpha} R_{j_1}\, \partial_{\ua}^{-2} g_{+}(\Psi_{2}^{h}, \mathcal{R}_{2})\, S_{j_1}^{l} (\partial_{\alpha} \psi_c) \partial^{l}_{\alpha} R_{j_1}\,d{\alpha} \Big)\\[2mm]
\;&  + \veps^{3} \mathcal{O}(\mathcal{E}_s +1)\\[2mm]
=\;& \veps^{3} \mathcal{O}(\mathcal{E}_s +1)\,, 
\\[3mm]
  I_6^1=
\;& 2\veps^{3}  \int_{\R} \partial^{l}_{\alpha} R_{2}\, \partial_{\ua}^{-1} g_{+}(\Psi_{2}^{h}, \mathcal{R}_{2}) (\rho^{l} \mathfrak{N}^{1})_{2-2}^s(\psi_c) \partial^{l}_{\alpha} R_{-2}\,d{\alpha}  \\[2mm]
& + 2\veps^{3}  \int_{\R} \partial^{l}_{\alpha} R_{2}\, \partial_{\ua}^{-2} g_{+}(\Psi_{2}^{h}, \mathcal{R}_{2}) (\rho^{l} \mathfrak{N}^{1})_{2-2}^s( \partial_{\alpha} \psi_c) \partial^{l}_{\alpha} R_{-2}\,d{\alpha} \\[2mm]
& +  \veps^{3} \mathcal{O}(\mathcal{E}_s +1) \\[2mm]
=\;& \veps^{3} \mathcal{O}(\mathcal{E}_s +1) 
\,,
\\[3mm]
 I_5^2
=\;&  \frac{\veps^{3}}{4} \frac{d}{dt} \int_{\R} \partial_{\alpha}^{l} (R_{-2}-R_{2})\, (K_0  \sigma^{-1} \partial_{\ua}^{-1} g_{-}(\Psi_{2}^{h}, \mathcal{R}_{2})+c(\Psi, \mathcal{R}))\,  (S_{-2}^l + S_{2}^l ) (\psi_c)   \\[2mm]
\;& \hspace{1.35cm} \times (1-b \partial_{\alpha}^2)^{-1}\, \partial^{l}_{\alpha} ( R_{-2}-  R_{2}) \,d{\alpha}+  \veps^{3} \mathcal{O}(\mathcal{E}_s +1)
\,,
\\[3mm]
 I_6^2
=\;&  \frac{\veps^{3}}{4} \frac{d}{dt} \int_{\R} \partial_{\alpha}^{l} (R_{-2}-R_{2})\, (K_0  \sigma^{-1} \partial_{\ua}^{-1} g_{-}(\Psi_{2}^{h}, \mathcal{R}_{2})+c(\Psi, \mathcal{R}))\, (G_{-2}^l + G_{2}^l )(\psi_c)   \\[2mm]
\;& \hspace{1.35cm} \times (1-b \partial_{\alpha}^2)^{-1}\, \partial^{l}_{\alpha} ( R_{-2}-  R_{2}) \,d{\alpha}+  \veps^{3} \mathcal{O}(\mathcal{E}_s +1)
\,,
\\[3mm]
 I_5^3
=\;& \frac{\veps^{3}}{2} \int_{\R}  
\partial_{\alpha}^{l} (R_{-2}-R_{2}) (\sigma^{-1} g_{-}(\Psi_{2}^{h}, \mathcal{R}_{2}))(S_{-2}^l +S_2^l ) (\psi_c) b K_0  \sigma^{-1} \partial^{l+2}_{\alpha}
( R_{-2}+  R_{2})\,d{\alpha}\\[2mm]
\;& + \veps^{3} \mathcal{O}(\mathcal{E}_s +1)\\[2mm]
=\;& \veps^{3} \mathcal{O}(\mathcal{E}_s +1 )
\,,
\\[3mm]
 I_6^3
=\;& 
\frac{\veps^{3}}{2} \int_{\R}  
\partial_{\alpha}^{l} (R_{-2}-R_{2}) (\sigma^{-1} g_{-}(\Psi_{2}^{h}, \mathcal{R}_{2}))(G_{-2}^l +G_2^l ) (\psi_c) b K_0  \sigma^{-1} \partial^{l+2}_{\alpha}
( R_{-2}+  R_{2})\,d{\alpha}\\[2mm]
\;& 
+ \veps^{3} \mathcal{O}(\mathcal{E}_s +1)\\[2mm]
=\;& \veps^{3} \mathcal{O}(\mathcal{E}_s +1 )
\,,
\\[3mm]
 I_5^4
=\;& - \frac{\veps^{3}}{2}  \int_{\R} 
\partial_{\alpha}^{l} (R_{-2}-R_{2})  (\sigma^{-1} g_{-}(\Psi_{2}^{h}, \mathcal{R}_{2})) (S_{-2}^l +S_2^l ) (\psi_c) \sigma \partial^{l-1}_{\alpha} (R_{-2}+R_2)\,d{\alpha} 
\\[2mm]
\;&     
+  \veps^{3} \mathcal{O}(\mathcal{E}_s  +1)\\[2mm]
=\;& 
\veps^{3} \mathcal{O}(\mathcal{E}_s +1)
\,,
\end{align*}
\begin{align*}
 I_6^4
=\;& - \frac{\veps^{3}}{2} \int_{\R} 
\partial_{\alpha}^{l} (R_{-2}-R_{2})  (\sigma^{-1} g_{-}(\Psi_{2}^{h}, \mathcal{R}_{2})) (G_{-2}^l +G_2^l ) (\psi_c) \sigma \partial^{l-1}_{\alpha} (R_{-2}+R_2)\,d{\alpha} 
\\[2mm]
\;&     
+  \veps^{3} \mathcal{O}(\mathcal{E}_s +1)\\[2mm]
=\;& 
\veps^{3} \mathcal{O}(\mathcal{E}_s +1)
\end{align*}
as long as $\eps^{5/2} \Vert \mathcal{R}_{1} \Vert_{(L^{2})^2},\, \eps^{5/2} \Vert \mathcal{R}_{2} \Vert_{(H^{{s}})^2} \ll 1$,
uniformly on compact subsets of $\mathcal{B}$.

Next, we estimate $I_7$. Due to \eqref{estcalMmp2}, \eqref{thetaNb} and \eqref{N0s}, we deduce
\begin{align*}
I_7 =\;  &  \veps^{3} \sum_{j_1  \in \{\pm2\}} \Big(\, \int_{\R}  \partial^{l}_{\alpha} \mathcal{M}_{j_1}({\Psi}, \mathcal{R})\, \partial_{\alpha}^{l} \rho_{j_1}^{l} \vartheta N_{j_1j_1}(\psi_c) R_{j_1}\,d{\alpha}
\; \\[2mm]
&\qquad\qquad\qquad\!\!\! + \int_{\R} \partial^{l}_{\alpha} R_{j_1}\, \partial^{l}_{\alpha} \rho_{j_1}^{l} \vartheta N_{j_1j_1}(\psi_c) \mathcal{M}_{j_1}({\Psi}, \mathcal{R})\,d{\alpha}\, \Big) \\[2mm]
& + \veps^{3} \mathcal{O}(\mathcal{E}_s +1)
\\[2mm]
=\; & \veps^{3} \sum_{j_1  \in \{\pm2\}}  \int_{\R}  \partial^{l}_{\alpha} \mathcal{M}_{j_1}({\Psi}, \mathcal{R})\, (\rho^{l} \mathfrak{N}^{1})_{j_1j_1}^s (\psi_c) \partial_{\alpha}^{l} R_{j_1}\,d{\alpha} 
+ \veps^{3} \mathcal{O}(\mathcal{E}_s +1)\\[2mm]
=\;& 
\veps^{3} \mathcal{O}(\mathcal{E}_s +1) \,,
\end{align*}
as long as $\eps^{5/2} \Vert \mathcal{R}_{1} \Vert_{(L^{2})^2},\, \eps^{5/2} \Vert \mathcal{R}_{2} \Vert_{(H^{{s}})^2} \ll 1$,
uniformly on compact subsets of $\mathcal{B}$.

Finally, we bound $I_1$ by improving the estimate \eqref{nf-rest2}. 
It follows from the proof of this estimate that
\begin{align*}
 I_1
=\;& \sum_{\nu=1}^{2} I_{1,1}^{\nu}   + \sum_{\nu=1}^{2} I_{1,2}^{\nu}  + \sum_{\nu=2}^{4} I_{1,3}^{\nu} + I_{1,4}+  \veps^{2} \mathcal{O}(\mathcal{E}_s+1)\,
\end{align*}
with
\begin{align*}
I_{1,1}^{\nu}
=\;& -\eps \sum_{\genfrac{}{}{0pt}{}{j_1,j_2 \in \{\pm 2\},} {\ell \in \{ \pm 1\}} }
\int_{\R} \int_{\R} \overline{(ik)^{l}\, \widehat{R}_{j_1}(k)}\, \big(\widehat{r}_{j_1j_2}(k,k-m,m)-\widehat{r}_{j_1j_2}(k,\ell k_0,m)\big)\,
\\[-3mm]
& \hspace{3.25cm} \times (ik)^{l}\, (\widehat{\rho^l \mathfrak{n}^{\nu}})_{j_1j_2 \ell}(k)\, 
\widehat{\psi}_{\ell}(k-m)\,(im)^{-(\nu-1)}\, \widehat{R}_{j_2}(m)\, dmdk\,,\\[3mm]
I_{1,2}^{\nu}
=\;& -\eps \sum_{\genfrac{}{}{0pt}{}{j_1,j_2 \in \{\pm 2\},} {\ell \in \{ \pm 1\}} }
\int_{\R} \int_{\R} \overline{(ik)^{l}\, \widehat{R}_{j_1}(k)}\, \big(\widehat{r}_{j_1j_2}(k,\ell k_0,m)-\widehat{r}_{j_1j_2}(k,\ell k_0,k-\ell k_0)\big)\,  
\\[-3mm]
& \hspace{3.25cm} \times (ik)^{l}\,  (\widehat{\rho^l \mathfrak{n}^{\nu}})_{j_1j_2 \ell}(k)\, 
\widehat{\psi}_{\ell}(k-m)\,(im)^{-(\nu-1)}\, \widehat{R}_{j_2}(m)\, dmdk\,,
\\[3mm]
I_{1,3}^{\nu}
=\;& \eps \sum_{\genfrac{}{}{0pt}{}{j_1,j_2 \in \{\pm 2\},} {\ell \in \{ \pm 1\}} }
\int_{\R} \int_{\R} \overline{(ik)^{l}\, \widehat{R}_{j_1}(k)}\, \big(\widehat{q}^{2,\nu}_{j_1j_2}(k,k-m,m)-\widehat{q}^{2,\nu}_{j_1j_2}(k,\ell k_0,m)\big)\, \
\\[-3mm]
& \hspace{2.8cm} \times (ik)^{l}\,
\widehat{\psi}_{\ell}(k-m)\, \widehat{\vartheta}(m)\,\widehat{R}_{j_2}(m)\, dmdk\,,
\\[3mm]
I_{1,4}
=\;& \eps \sum_{\genfrac{}{}{0pt}{}{j_1,j_2 \in \{\pm 2\},} {\ell \in \{ \pm 1\}} }
\int_{\R} \int_{\R} \overline{(ik)^{l}\, \widehat{R}_{j_1}(k)}\, \big(\widehat{\mathfrak{q}}^1_{j_1j_2}(k,\ell k_0,m)-\widehat{\mathfrak{q}}^1_{j_1j_2}(k,\ell k_0,k-\ell k_0)\big)\, 
\\[-3mm]
& \hspace{2.8cm} \times (ik)^{l}\,
\widehat{\psi}_{\ell}(k-m)\, \widehat{\vartheta}(m)\,\widehat{R}_{j_2}(m)\, dmdk\,.
\end{align*}
Because of 
\begin{align} \label{omega2l}
\widehat{r}_{j_1j_2}(k,k-m,m)-\widehat{r}_{j_1j_2}(k,\ell k_0,m) =  i \big(\omega(k-m) - \omega(\ell k_0)\big) 
\end{align}
as well as \eqref{asnjj}--\eqref{asnj-j} and Lemma \ref{int-kerne} there holds
\begin{align*}
I_{1,1}^1
=\;& -\eps \sum_{\genfrac{}{}{0pt}{}{j_1 \in \{\pm 2\},} {\ell \in \{ \pm 1\}} }
\int_{\R} \int_{\R} \overline{(ik)^{l}\, \widehat{R}_{j_1}(k)}\,  i \big(\omega(k-m) - \omega(\ell k_0)\big) 
\\[-2mm]
& \hspace{2.85cm} \times  
 (\widehat{\rho^l \mathfrak{n}^{1}})_{j_1j_1 \ell}(k)\,
\widehat{\psi}_{\ell}(k-m)\,(im)^{l}\,\widehat{R}_{j_1}(m)\, dmdk 
\\[2mm]
& + \veps^{2} \mathcal{O}(\mathcal{E}_s +1) \,.
\end{align*}
For symmetry reasons it follows
\begin{align*}
I_{1,1}^1
=\;& -\eps \sum_{\genfrac{}{}{0pt}{}{j_1 \in \{\pm 2\},} {\ell \in \{ \pm 1\}} }
\int_{\R} \int_{\R} \overline{(ik)^{l}\, \widehat{R}_{j_1}(k)}\,  i \big(\omega(k-m) - \omega(\ell k_0)\big) 
\\[-2mm]
& \hspace{2.85cm} \times  
 (\widehat{\rho^l \mathfrak{n}^{1}})_{j_1j_1 \ell}^s(k)\,
\widehat{\psi}_{\ell}(k-m)\,(im)^{l}\,\widehat{R}_{j_1}(m)\, dmdk 
\\[2mm]
& + \veps^{2} \mathcal{O}(\mathcal{E}_s +1) \,,
\end{align*}
such that, due to \eqref{nus} and Lemma \ref{int-kerne}, we conclude
\begin{align*}
I_{1,1}^1
=\;&  \veps^{2} \mathcal{O}(\mathcal{E}_s+1) \,,
\end{align*}
as long as $\eps^{5/2} \Vert \mathcal{R}_{1} \Vert_{(L^{2})^2},\, \eps^{5/2} \Vert \mathcal{R}_{2} \Vert_{(H^{{s}})^2} \ll 1$, uniformly on compact subsets of $\mathcal{B}$.

Furthermore, by using \eqref{dx-2R2}, \eqref{asn2jj}--\eqref{asn2j-j}, \eqref{omega2l} and Lemma \ref{int-kerne} we obtain
\begin{align*}
I_{1,1}^2
=\;&  \veps^{2} \mathcal{O}(\mathcal{E}_s+1) \,,
\end{align*}
as long as $\eps^{5/2} \Vert \mathcal{R}_{1} \Vert_{(L^{2})^2},\, \eps^{5/2} \Vert \mathcal{R}_{2} \Vert_{(H^{{s}})^2} \ll 1$, uniformly on compact subsets of $\mathcal{B}$.

Also because of \eqref{asnjj}--\eqref{asnj-j} and Lemma \ref{int-kerne} we deduce
\begin{align*}
I_{1,2}^1
=\;& -\eps \sum_{\genfrac{}{}{0pt}{}{j_1 \in \{\pm 2\},} {\ell \in \{ \pm 1\}} }
\int_{\R} \int_{\R} \overline{(ik)^{l}\, \widehat{R}_{j_1}(k)}\, \big(\widehat{r}_{j_1j_1}(k,\ell k_0,m)-\widehat{r}_{j_1j_1}(k,\ell k_0,k-\ell k_0)\big)\, 
\\[-2mm]
& \hspace{2.9cm} \times (\widehat{\rho^l \mathfrak{n}^{1}})_{j_1j_1 \ell}(k)\, 
\widehat{\psi}_{\ell}(k-m)\,(im)^{l}\,\widehat{R}_{j_1}(m)\, dmdk 
\\[2mm]
& + \veps^{2} \mathcal{O}(\mathcal{E}_s+1)\,.
\end{align*}
We split the integral kernel into
\begin{align*}
& \big(\widehat{{r}}_{j_1j_1}(k,\ell k_0,m)-\widehat{r}_{j_1j_1}(k,\ell k_0,k-\ell k_0)\big)\, 
(\widehat{\rho^l \mathfrak{n}^{1}})_{j_1j_1 \ell}(k)
\\[2mm]
& \qquad 
= \big(\widehat{\mathfrak{r}}_{j_1\ell}^{s}(k,m) + \widehat{\mathfrak{r}}_{j_1\ell}^{a}(k,m)\big) \,
(\widehat{\rho^l \mathfrak{n}^{1}})_{j_1j_1 \ell}(k)
\end{align*}
with
\begin{align*}
\widehat{\mathfrak{r}}_{j_1\ell}^{s}(k,m) =\;& -\frac 12 \,\mathrm{sgn}(j_1)\, i\, \Big( \big( \omega(m) - \omega(k-\ell k_0)\big) + \big(\omega(-k) - \omega(-m-\ell k_0)\big) \Big),\\[2mm]
\widehat{\mathfrak{r}}_{j_1\ell}^{a}(k,m) =\;& -\frac 12 \,\mathrm{sgn}(j_1)\, i\, \Big( \big( \omega(m) - \omega(k-\ell k_0)\big) - \big(\omega(-k) - \omega(-m-\ell k_0)\big) \Big).
\end{align*}
By the mean value theorem we have
\begin{align}
& \big( \omega(m) - \omega(k-\ell k_0)\big) \pm \big(\omega(-k) - \omega(-m-\ell k_0)\big) \nonumber \\[2mm]
& \qquad = \big( \omega(m) - \omega(k-\ell k_0)\big) \mp \big(\omega(k) - \omega(m+\ell k_0)\big) \nonumber \\[2mm]
& \qquad = - \omega'(m +\theta_0(k,m,\ell)(k-m-\ell k_0))\, (k-m-\ell k_0) \nonumber \\[2mm]
& \qquad \quad \mp \omega'(k -\theta_1(k,m,\ell)(k-m-\ell k_0))\, (k-m-\ell k_0)\,
\end{align}
with $\theta_0(k,m,\ell), \theta_1(k,m,\ell) \in [0,1]$. 
Hence, we obtain
\begin{align} \label{frakrs}
\widehat{\mathfrak{r}}_{j_1\ell}^{s}(k,m)\, \chi_{\ell}(k-m)  = \mathcal{O}(|k|^{-1/2}(1+bk^{2})^{1/2})\, (k-m-\ell k_0)
\end{align}
for $|k| \to \infty$ uniformly with respect to $m \in \R$ and $b \lesssim 1$ and,
by using the mean value theorem once more, 
\begin{align} \label{frakra}
\widehat{\mathfrak{r}}_{j_1\ell}^{a}(k,m)\, \chi_{\ell}(k-m)  = \mathcal{O}(|k|^{-3/2}(1+bk^{2})^{1/2})\, (k-m-\ell k_0)
\end{align}
for $|k| \to \infty$ uniformly with respect to $m \in \R$ and $b \lesssim 1$.
Consequently, due to \eqref{asnjj}, \eqref{nus} and Lemma \ref{int-kerne}, we conclude
\begin{align*}
I_{1,2}^1
=\;& -\frac{\veps}{2}  \sum_{\genfrac{}{}{0pt}{}{j_1 \in \{\pm 2\},} {\ell \in \{ \pm 1\}} }
\int_{\R} \int_{\R} \overline{(ik)^{l} \widehat{R}_{j_1}(k)}\,\widehat{\mathfrak{r}}_{j_1\ell}^{s}(k,m) 
\, (\widehat{\rho^l \mathfrak{n}^{1}})^s_{j_1j_1 \ell}(k)\, \widehat{\psi}_{\ell}(k-m)
\\[-2mm]
& \hspace{2.95cm} \times
(im)^{l}\widehat{R}_{j_1}(m)\, dmdk 
\\[2mm]
& -\frac{\veps}{2}  \sum_{\genfrac{}{}{0pt}{}{j_1 \in \{\pm 2\},} {\ell \in \{ \pm 1\}} }
\int_{\R} \int_{\R} \overline{(ik)^{l} \widehat{R}_{j_1}(k)}\,\widehat{\mathfrak{r}}_{j_1\ell}^{a}(k,m) 
\, (\widehat{\rho^l \mathfrak{n}^{1}})^a_{j_1j_1 \ell}(k)\, \widehat{\psi}_{\ell}(k-m)
\\[-2mm]
& \hspace{2.95cm} \times
(im)^{l}\widehat{R}_{j_1}(m)\, dmdk 
\\[2mm]
& + \veps^{2} \mathcal{O}(\mathcal{E}_s+1)
\\[2mm]
=\;& \veps^{2} \mathcal{O}(\mathcal{E}_s+1) \,,
\end{align*}
as long as $\eps^{5/2} \Vert \mathcal{R}_{1} \Vert_{(L^{2})^2},\, \eps^{5/2} \Vert \mathcal{R}_{2} \Vert_{(H^{{s}})^2} \ll 1$, uniformly on compact subsets of $\mathcal{B}$.

Furthermore, because of \eqref{dx-2R2}, \eqref{thetaN2a}, \eqref{frakrs}--\eqref{frakra} and
Lemma \ref{int-kerne}, we have
\begin{align*}
I_{1,2}^2
=\;& -\frac{\veps}{2}  \sum_{\genfrac{}{}{0pt}{}{j_1 \in \{\pm 2\},} {\ell \in \{ \pm 1\}} }
\int_{\R} \int_{\R} \overline{(ik)^{l} \widehat{R}_{j_1}(k)}\,\widehat{\mathfrak{r}}_{j_1\ell}^{s}(k,m) 
\, (\widehat{\rho^l \mathfrak{n}^{2}})^s_{j_1j_1 \ell}(k)\, \widehat{\psi}_{\ell}(k-m)
\\[-2mm]
& \hspace{2.95cm} \times
(im)^{l-1}\widehat{R}_{j_1}(m)\, dmdk 
\\[1mm]
=\;& \veps^{2} \mathcal{O}(\mathcal{E}_s+1) \,.
\end{align*}
Hence, using 
\begin{align*}
(\widehat{\rho^l \mathfrak{n}^2})^{s}_{j_1j_1\ell}(k,m) = \;& \frac12 \big( (ik)^{-1} \,
(\widehat{\rho^l \mathfrak{n}^{2}})_{j_1j_1 \ell}(k) - (i(-m))^{-1} \, (\widehat{\rho^l \mathfrak{n}^{2}})_{j_1j_1 \ell}(-m) \big) im 
\\[2mm]
& + \frac12 (ik)^{-1} \,
(\widehat{\rho^l \mathfrak{n}^{2}})_{j_1j_1 \ell}(k) \, i(k-m) 
\end{align*}
as well as \eqref{pm1-ans}, \eqref{pm1-ans2}, \eqref{asn2jj}, \eqref{frakrs} and Young's inequality for convolutions, we obtain
\begin{align*}
I_{1,2}^2 =\;& \veps^{2} \mathcal{O}(\mathcal{E}_s+1) \,,
\end{align*}
as long as $\eps^{5/2} \Vert \mathcal{R}_{1} \Vert_{(L^{2})^2},\, \eps^{5/2} \Vert \mathcal{R}_{2} \Vert_{(H^{{s}})^2} \ll 1$, uniformly on compact subsets of $\mathcal{B}$.

Analogously to the case of $I_4^{2}$ and $I_4^{3}$, we deduce
\begin{align*}
I_{1,3}^{2}
=\;& \frac{\eps^{2}}{4} \frac{d}{dt} \int_{\R}  \partial_{\alpha}^{l} (R_{-2}-R_2) L_1(\psi_c) (1-b \partial_{\alpha}^2)^{-1} \partial_{\alpha}^{l} (R_{-2}-R_2)\,d\alpha + \veps^{2} \mathcal{O}(\mathcal{E}_s+1)\,,
\\[3mm]
I_{1,3}^{3}
=\;& - \frac{\eps^{2}}{4} \frac{d}{dt} \int_{\R}  \partial_{\alpha}^{l} (R_{-2}-R_2) L_2(\psi_c) b K_0 \partial_{\alpha}(1-b \partial_{\alpha}^2)^{-1}   \partial_{\alpha}^{l} (R_{-2}-R_2)\,d\alpha 
\\[2mm]
\;& 
+ \veps^{2} \mathcal{O}(\mathcal{E}_s+1)\,,
\end{align*}
as long as $\eps^{5/2} \Vert \mathcal{R}_{1} \Vert_{(L^{2})^2},\, \eps^{5/2} \Vert \mathcal{R}_{2} \Vert_{(H^{{s}})^2} \ll 1$, uniformly on compact subsets of $\mathcal{B}$, where
\begin{align*}
\widehat{L}_1(\psi_c)(k) =\;&  \sum_{\ell \in \{\pm 1\}} \veps^{-1}(\widehat{K}_0(k)  \sigma^{-1}(k) ik - \widehat{K}_0(\ell k_0)  \sigma^{-1}(\ell k_0) i\ell k_0)\,\widehat{\psi}_{\ell}(k)\,,\\[2mm]
\widehat{L}_2(\psi_c)(k) =\;&  -\sum_{\ell \in \{\pm 1\}} \veps^{-1}( \sigma^{-1}(k)k^{2} - \sigma^{-1}(\ell k_0)(\ell k_0)^{2})\,\widehat{\psi}_{\ell}(k)
\,.
\end{align*}
To bound $I_{1,3}^{4}$ we split the integral kernel into
\begin{align*}
& \widehat{q}^{2,\nu}_{j_1j_2}(k,k-m,m)-\widehat{q}^{2,\nu}_{j_1j_2}(k,\ell k_0,m)
\\[2mm] 
& \qquad 
= \frac12\,  \mathrm{sgn}(j_1)\, ik\, i\, \frac{\sigma(k)-\sigma(k-m)}{k-(k-m)} \, \big(\sigma^{-1}(k-m) (k-m)^2 - \sigma^{-1}(\ell k_0)(\ell k_0)^{2} \big)\, (im)^{-1}
\\[2mm] 
& \qquad \quad \,
+ \frac12\,  \mathrm{sgn}(j_1)\, ik\, i\, \Big(\frac{\sigma(k)-\sigma(k-m)}{k-(k-m)} - \frac{\sigma(k)-\sigma(\ell k_0)}{k-\ell k_0}\Big)\, \sigma^{-1}(\ell k_0)(\ell k_0)^{2}\, (im)^{-1} \,.
\end{align*}
With the help of Taylor's theorem we obtain
\begin{align*}
\frac{\sigma(k)-\sigma(k-m)}{k-(k-m)} - \frac{\sigma(k)-\sigma(\ell k_0)}{k-\ell k_0} & \leq\, \frac12 \, \|\partial_k^{2} \sigma \|_{L^{\infty}} \, (k-m-\ell k_0)
\\[2mm]
& = \, \mathcal{O}(1)\, (k-m-\ell k_0)
\end{align*}
uniformly with respect to $k,m \in \R$ and $b \lesssim 1$. 
Because of Lemma \ref{int-kerne} we conclude
\begin{align*}
I_{1,3}^{4}
=\;& -\frac{\eps^{2}}{2} \int_{\R}   ( [\sigma,\partial_{\ua}^{l-1}(R_{ -2}+R_{2})] L_2(\psi_c)) \partial_{\alpha}^{l} (R_{-2}-R_2)\,d\alpha + \veps^{2} \mathcal{O}(\mathcal{E}_s+1) \,.
\end{align*}
Hence, we can proceed analogously to the case of $I_4^{4}$ to deduce
\begin{align*}
I_{1,3}^{4}
=\;& \frac{\eps^{2}}{4} \frac{d}{dt} \int_{\R}  \partial_{\alpha}^{l} (R_{-2}-R_2) L_2(\psi_c) (-K_0 \partial_{\alpha})^{-1} \partial_{\alpha}^{l} (R_{-2}-R_2)\,d\alpha + \veps^{2} \mathcal{O}(\mathcal{E}_s+1) \,,
\end{align*}
as long as $\eps^{5/2} \Vert \mathcal{R}_{1} \Vert_{(L^{2})^2},\, \eps^{5/2} \Vert \mathcal{R}_{2} \Vert_{(H^{{s}})^2} \ll 1$, uniformly on compact subsets of $\mathcal{B}$

Finally, due to the mean value theorem, we have
\begin{align*}
(\widehat{\mathfrak{q}}^1_{j_1j_2}(k,\ell k_0,m) - \widehat{\mathfrak{q}}^1_{j_1j_2}(k,\ell k_0,k - \ell k_0))\, \chi_{\ell}(k-m) \,= \, \mathcal{O}(|k|^{-1/2})\, (k-m-\ell k_0)
\end{align*}
for $|k| \to \infty$ uniformly with respect to $ m \in \R$ and $b  \lesssim 1$, such that with the help of Lemma \ref{int-kerne} we obtain
\begin{align*}
I_{1,4}
=\;& \veps^{2}  \mathcal{O}(\mathcal{E}_s+1) \,,
\end{align*}
as long as $\eps^{5/2} \Vert \mathcal{R}_{1} \Vert_{(L^{2})^2},\, \eps^{5/2} \Vert \mathcal{R}_{2} \Vert_{(H^{{s}})^2} \ll 1$, uniformly on compact subsets of $\mathcal{B}$.
\medskip

Now, we define our final energy $\widetilde{\mathcal{E}_s}$ by 
\begin{equation}
\widetilde{\mathcal{E}_s} := \mathcal{E}_s + \frac{1}{4} \, \veps^{2}\, \sum_{l=1}^s \sum_{i=1}^7 h_{l}^i\,
\end{equation}
with
\begin{align*}
h_{l}^1 \;=\;& -\int_{\R} \partial_{\alpha}^{l} (R_{-2}-R_{2})L_1(\psi_c)(1-b \partial_{\alpha}^2)^{-1} \partial^{l}_{\alpha} ( R_{-2}-  R_{2}) \,d{\alpha} \,, \\[2mm]
h_{l}^2 \;=\;&  \int_{\R}  \partial_{\alpha}^{l} (R_{-2}-R_{2}) L_2(\psi_c)b  K_0 \partial_{\alpha} (1-b \partial_{\alpha}^2)^{-1}      \partial^{l}_{\alpha}( R_{-2}-  R_{2})\,d{\alpha} \,, 
\\[2mm]
h_{l}^3 \;=\;& -  \int_{\R}  \partial_{\alpha}^{l} (R_{-2}-R_{2})  L_2(\psi_c) (-K_0 \partial_{\alpha})^{-1} \partial^{l}_{\alpha} (R_{-2}-R_2)\,d{\alpha} \,, 
\\[2mm]
h_{l}^4 \;=\;& - \int_{\R} \partial_{\alpha}^{l} (R_{-2}-R_{2}) (K_0  \sigma^{-1} \partial_{\alpha}\psi_c) (S_{-2}^l +S_2^l +G_{-2}^l +G_2^l ) (\psi_c) 
 \\[1mm]
\;& \hspace{0.9cm} \times (1-b \partial_{\alpha}^2)^{-1}   \partial^{l}_{\alpha} ( R_{-2}-  R_{2}) \,d{\alpha} \,, 
\\[2mm]
h_{l}^5 \;=\;& - \int_{\R} \partial_{\alpha}^{l} (R_{-2}-R_{2}) (K_0  \sigma^{-1} \partial_{\ua}^{-1} g_{-}(\Psi_{2}^{h}, \mathcal{R}_{2})+c(\Psi, \mathcal{R})) 
 \\[1mm]
\;& \hspace{0.9cm} \times (1+ {\eps} (S_{-2}^l +S_2^l +G_{-2}^l +G_2^l ) (\psi_c))
(1-b \partial_{\alpha}^2)^{-1} \partial^{l}_{\alpha} ( R_{-2}-  R_{2}) \,d{\alpha} \,, 
\\[2mm]
h_{l}^6 \;=\;& \int_{\R} \partial_{\alpha}^{l} (R_{-2}-R_{2})(\sigma^{-1}  g_{-}(\Psi_{2}^{h}, \mathcal{R}_{2}))\, b K_0 \partial_{\alpha} 
(1-b \partial_{\alpha}^2)^{-1} \partial^{l}_{\alpha} ( R_{-2}-  R_{2}) \,d{\alpha} \,, 
\\[2mm]
h_{l}^7 \;=\;&  - \int_{\R} \partial_{\alpha}^{l} (R_{-2}-R_{2})\,  \sigma^{-1}  g_{-}(\Psi_{2}^{h}, \mathcal{R}_{2}) 
(-K_0 \partial_{\alpha})^{-1} \partial^{l}_{\alpha} ( R_{-2}-  R_{2}) \,d{\alpha}  
\,.
\end{align*}

Then, for sufficiently small $\eps >0$, the energy $\widetilde{\mathcal{E}_s}$ satisfies
the estimates \eqref{tEest}--\eqref{normtE},
as long as $\eps^{5/2} \Vert \mathcal{R}_{1} \Vert_{(L^{2})^2}, \eps^{5/2} \Vert \mathcal{R}_{2} \Vert_{(H^{{s}})^2} \ll 1$,
uniformly on compact subsets of $\mathcal{B}$. Hence, we have proven Lemma \ref{finenest}.
\qed
\medskip

{\bf Proofs of Theorems \ref{arcresult} and \ref{mainresult}.}
If $\Vert \mathcal{R}_{1}|_{t=0} \Vert_{(L^{2})^2}, \Vert \mathcal{R}_{2}|_{t=0}  \Vert_{(H^{{s}})^2} \lesssim 1$, then Lemma \ref{finenest} allows us to use
Gronwall's inequality to obtain for sufficiently small $\veps>0$ the $\mathcal{O}(1)$-bounded\-ness of $\widetilde{\mathcal{E}_s}$ for all $t\in[0,\tau_0/\veps^{2}]$ uniformly on compact subsets of $\mathcal{B}$. Due to
\eqref{RES2} and \eqref{normtE}, Theorem \ref{arcresult} follows. Transferring the assertions of Theorem \ref{arcresult} into Eulerian coordinates finally yields Theorem \ref{mainresult}.
\qed
\\[3mm]
\textbf{Acknowledgment.} This work was supported by the Deutsche Forschungsgemeinschaft DFG under the grant DU 1198/2. The author thanks Max He{\ss} for discussions and the referee for useful comments.

\end{document}